\documentclass[11pt,a4paper,reqno]{amsart}

\pdfoutput=1

\usepackage[utf8]{inputenc}
\usepackage{epsfig,amsmath,latexsym,amssymb}
\usepackage{amsthm, tikz}
\usetikzlibrary{calc}
\usepackage{bm}
\usepackage[numbers,sort&compress,longnamesfirst]{natbib}
\usepackage{graphicx}
\usepackage{url}
\usepackage[margin=1.1in]{geometry}
\usepackage{color}
\usepackage[colorlinks,citecolor=blue,urlcolor=blue,filecolor=blue]{hyperref}
\usepackage{multirow}
\usepackage{enumitem}
\usepackage{lscape}
\usepackage{cases}
\usepackage{dsfont}
\usepackage[ruled,vlined]{algorithm2e}
\usepackage{algorithmic}
\usepackage{placeins}
\usepackage{stmaryrd}
\usepackage{mathrsfs}
\usepackage{caption}

\usepackage{amsfonts,bbm}

\newcommand{\Td}{\mathrm{d}}

\newcommand{\Bb}{\mathbf{b}}
\newcommand{\Bc}{\mathbf{c}}

\newcommand{\Bp}{\mathbf{p}}
\newcommand{\Bq}{\mathbf{q}}

\newcommand{\BA}{\mathbf{A}}

\newcommand{\BF}{\mathbf{F}}

\newcommand{\BI}{\mathbf{I}}

\newcommand{\BL}{\mathbf{L}}

\newcommand{\BP}{\mathbf{P}}
\newcommand{\BQ}{\mathbf{Q}}

\newcommand{\BU}{\mathbf{U}}

\newcommand{\BIe}{{\boldsymbol{e}}}
\newcommand{\BIf}{{\boldsymbol{f}}}
\newcommand{\BIg}{{\boldsymbol{g}}}

\newcommand{\BIx}{{\boldsymbol{x}}}
\newcommand{\BIy}{{\boldsymbol{y}}}
\newcommand{\BIz}{{\boldsymbol{z}}}

\newcommand{\BIY}{{\boldsymbol{Y}}}

\newcommand{\CB}{\mathcal{B}}
\newcommand{\CC}{\mathcal{C}}

\newcommand{\CF}{\mathcal{F}}

\newcommand{\CJ}{\mathcal{J}}

\newcommand{\CL}{\mathcal{L}}

\newcommand{\CP}{\mathcal{P}}

\newcommand{\CS}{\mathcal{S}}
\newcommand{\CT}{\mathcal{T}}

\newcommand{\CV}{\mathcal{V}}

\newcommand{\CY}{\mathcal{Y}}
\newcommand{\CZ}{\mathcal{Z}}

\newcommand{\TTC}{\mathtt{C}}

\newcommand{\TTE}{\mathtt{E}}

\newcommand{\TTP}{\mathtt{P}}

\newcommand{\TTS}{\mathtt{S}}

\newcommand{\TTV}{\mathtt{V}}

\newcommand{\TTc}{\mathtt{c}}

\newcommand{\TTp}{\mathtt{p}}

\newcommand{\TTs}{\mathtt{s}}

\newcommand{\ILambda}{\mathit{\Lambda}}

\newcommand{\BPi}{{\boldsymbol{\Pi}}}

\newcommand{\BOmega}{{\boldsymbol{\Omega}}}

\newcommand{\BILambda}{{\boldsymbol{\ILambda}}}

\newcommand{\Bzeta}{{\boldsymbol{\zeta}}}

\newcommand{\Btheta}{{\boldsymbol{\theta}}}

\newcommand{\Blambda}{{\boldsymbol{\lambda}}}

\newcommand{\Bnu}{{\boldsymbol{\nu}}}
\newcommand{\Bxi}{{\boldsymbol{\xi}}}

\newcommand{\Brho}{{\boldsymbol{\rho}}}

\newcommand{\Bomega}{{\boldsymbol{\omega}}}

\newcommand{\R}{\mathbb{R}} 
\newcommand{\N}{\mathbb{N}} 
\newcommand{\vecone}{\mathbf{1}} 
\newcommand{\veczero}{\mathbf{0}} 
\newcommand{\INDI}{\mathbbm{1}} 
\newcommand{\TRANSP}{\mathsf{T}} 
\newcommand{\DIFF}{\Td} 
\newcommand{\DIFFX}[1]{\,\Td{#1}} 
\newcommand{\DIFFM}[2]{\,{#1}({#2})}
\newcommand{\PROB}{\mathbb{P}} 
\newcommand{\EXP}{\mathbb{E}} 
\newcommand{\clos}{\mathrm{cl}} 
\newcommand{\support}{\mathrm{supp}} 
\newcommand{\conv}{\mathrm{conv}} 

\DeclareMathOperator*{\argmin}{arg\,min} 
\DeclareMathOperator*{\argmax}{arg\,max} 
\DeclareMathOperator*{\minimize}{\mathrm{minimize}} 
\DeclareMathOperator*{\maximize}{\mathrm{maximize}} 

\DeclareFontFamily{U}{mathx}{\hyphenchar\font45}
\DeclareFontShape{U}{mathx}{m}{n}{
      <5> <6> <7> <8> <9> <10>
      <10.95> <12> <14.4> <17.28> <20.74> <24.88>
      mathx10
      }{}
\DeclareSymbolFont{mathx}{U}{mathx}{m}{n}
\DeclareMathSymbol{\bigtimes}{1}{mathx}{"91}

\newtheorem{theorem}{Theorem}[section]
\newtheorem{corollary}[theorem]{Corollary}
\newtheorem{lemma}[theorem]{Lemma}

\newtheorem{remark}[theorem]{Remark}
\newtheorem{proposition}[theorem]{Proposition}
\newtheorem{example}[theorem]{Example}

\newtheorem{setting}[theorem]{Setting}

\usepackage{times}

\DontPrintSemicolon

\newcommand{\union}{\hspace{1pt}\cup\hspace{1pt}}
\newcommand{\intersection}{\hspace{1pt}\cap\hspace{1pt}}

\newcommand{\extremept}{\mathrm{extr}}
\newcommand{\recession}{\mathrm{recc}}
\newcommand{\lineality}{\mathrm{line}}
\newcommand{\ineq}{\mathrm{in}}
\newcommand{\eq}{\mathrm{eq}}
\newcommand{\DRO}{\mathrm{DRO}}
\newcommand{\OT}{\mathrm{OT}}

\newcommand{\LSIPprimal}[1]{{\normalfont\big(\hyperref[eqn:dro-dual-LSIP]{$\textsf{LSIP}\normalfont{\text{-}}{#1}$}\big)}}
\newcommand{\LPprimal}[1]{{\normalfont\big(\hyperref[eqn:dro-dual-LP-relax]{$\textsf{LP}\normalfont{\text{-}}\CV\normalfont{\text{-}}{#1}$}\big)}}
\newcommand{\LPdual}[1]{{\normalfont\big(\hyperref[eqn:dro-dual-LP-relax-dual]{$\textsf{LP}^*\normalfont{\text{-}}\CV\normalfont{\text{-}}{#1}$}\big)}}

\numberwithin{equation}{section}

\numberwithin{figure}{section}

\numberwithin{table}{section}

\begin{document}

\title[Numerical method for two-stage DRO with marginal constraints]{Numerical method for approximately optimal solutions \\ of two-stage distributionally robust optimization \\ with marginal constraints}

\author[A. Neufeld]{Ariel Neufeld}
\author[Q. Xiang]{Qikun Xiang}

\address{Division of Mathematical Sciences, Nanyang Technological University, 21 Nanyang Link, 637371 Singapore}
\email{ariel.neufeld@ntu.edu.sg}

\address{Division of Mathematical Sciences, Nanyang Technological University, 21 Nanyang Link, 637371 Singapore}
\email{qikun.xiang@ntu.edu.sg}

\date{}

\thanks{AN and QX
gratefully acknowledge the financial support by the MOE AcRF Tier 2 Grant \textit{MOE-T2EP20222-0013}.} 


\begin{abstract}
    We consider a general class of two-stage distributionally robust optimization (DRO) problems where
    the ambiguity set is constrained by fixed marginal probability laws that are not necessarily discrete.
    We derive primal and dual formulations of this class of problems and subsequently
    develop a numerical algorithm for computing approximate optimizers as well as approximate worst-case probability measures.
    Moreover, our algorithm computes both an upper bound and a lower bound for the optimal value of the problem, where the difference between the computed bounds provides a direct sub-optimality estimate of the computed solution. 
    Most importantly, the sub-optimality can be controlled to be arbitrarily close to~0 by appropriately choosing the inputs of the algorithm. 
    To demonstrate the effectiveness of the proposed algorithm, 
    we apply it to three prominent instances of two-stage DRO problems in task scheduling, multi-product assembly, and supply chain network design with edge failure.
    The ambiguity sets in these problem instances involve a large number of continuous or discrete marginals.
    The numerical results showcase that the proposed algorithm computes 
    high-quality robust decisions along with non-conservative sub-optimality estimates.
\end{abstract}

\maketitle

\section{Introduction}
\label{sec:intro}

Decision problems in uncertain environments are naturally present in many important application areas.
Examples of such problems include portfolio selection \citep{delage2010distributionally, gao2017data, mohajerin2018data}, inventory management \citep*{wang2016likelihood}, scheduling \citep*{chen2018distributionally, kong2020appointment, mak2015appointment}, resource allocation \citep*{wiesemann2012robust}, and transportation \citep*{bertsimas2010models, hu2020production, wang2021target, goh2025strategic}. 
For these decision problems, the two-stage stochastic programming model (see, e.g., \citep*{shapiro2009lectures}), in which a random event occurs between the two stages, is widely adopted. 
In this model, the decision maker makes the so-called \textit{here-and-now} decision in the first stage. 
Subsequently, after the outcome of the random event is observed, the decision maker makes the so-called \textit{wait-and-see} decision in the second stage which depends on the random outcome. 
Let us denote the first-stage decision by $\BIx$, denote the cost incurred in the first stage by $c_1(\BIx)$, denote the outcome of the random event by $\Bomega\in\BOmega$ where $\BOmega$ contains all possible outcomes, and denote the optimal cost in the second-stage decision problem by $C_2(\BIx,\Bomega)$. 
Moreover, the decision maker has access to a probability law $\mu$ of the random outcome $\Bomega$.
Hence, the first-stage decision is made by minimizing the sum of the first-stage cost and the expected second-stage cost, i.e., 
$\min_{\BIx}\big\{c_1(\BIx)+\int_{\BOmega}C_2(\BIx,\Bomega)\DIFFM{\mu}{\DIFF\Bomega}\big\}$. 
Since stochastic programming can be highly sensitive to the choice of the probability law $\mu$, robust optimization has been proposed as a conservative alternative, in which
the decision maker minimizes the worst-case cost where the random outcome $\Bomega$ can be any element of the set $\BOmega$, i.e., 
$\min_{\BIx}\big\{c_1(\BIx)+\sup_{\Bomega\in\BOmega}\{C_2(\BIx,\Bomega)\}\big\}$; 
see, e.g., \citep*{bental2009robust, simchilevi2019constraint, zeng2013solving}.

Compared to stochastic programming and robust optimization, distributionally robust optimization (DRO) \citep{bertsimas2010models, delage2010distributionally, goh2010distributionally,kuhn2025distributionally} achieves a balance between performance and robustness.  
In DRO, the decision maker specifies a collection of probability measures $\CP_{\BOmega}$ on $\BOmega$, termed the ambiguity set, which contains all plausible candidates of the probability law of the random outcome $\Bomega$, and subsequently minimizes the worst-case expected cost where the probability law $\mu$ of the random outcome can be any element of $\CP_{\BOmega}$, i.e., 
$\min_{\BIx}\Big\{c_1(\BIx)+\sup_{\mu\in\CP_{\BOmega}}\big\{\int_{\BOmega}C_2(\BIx,\Bomega)\DIFFM{\mu}{\DIFF\Bomega}\big\}\Big\}$. 
Consequently, DRO provides more robustness than stochastic programming while being less conservative than robust optimization.
The choice of the ambiguity set is central to the performance of DRO. 
A good ambiguity set should encode the prior beliefs of the decision maker, be rich enough to contain a good approximation of the true underlying probability law, 
and also allow tractable computation of the resulting DRO problem. 
A variety of ambiguity sets have been considered in the literature, including but not limited to those based on
moments \citep*{delage2010distributionally, goh2010distributionally, wiesemann2014distributionally, goh2025strategic, natarajan2023distributionally},
moments and dispersion \citep*{chen2024robust},
statistical divergences (e.g., the Kullback--Leibler divergence) \citep*{calafiore2007ambiguous, huang2021multi},
likelihood \citep*{wang2016likelihood},
Wasserstein distance \citep*{chen2020robust, gao2016distributionally, hanasusanto2018conic, mohajerin2018data, wozabal2012framework, wozabal2014robustifying, zhao2018data},
parametric families \citep*{iyengar2022hedging},
and 
fixed marginals \citep*{chen2018distributionally, gao2017data};
see also the comprehensive survey of \citet*{kuhn2025distributionally} about the types of ambiguity sets in the literature.

In this paper, we study a general class of two-stage DRO problems in which the first-stage cost $c_1(\cdot)=\langle\Bc_1,\cdot\,\rangle$ is linear 
and the second-stage decision problem is a linear programming (LP) problem where the right-hand side of the constraints has a jointly affine dependence on the first-stage decision $\BIx$ and the random outcome $\Bomega$ (see \ref{setts:dro-stage2} in Setting~\ref{sett:dro}). 
Moreover, we consider ambiguity sets based on fixed marginals, i.e., they contain all couplings of fixed one-dimensional probability measures (see \ref{setts:dro-uncertainty} in Setting~\ref{sett:dro}).
This class of two-stage DRO models contains prominent decision problems in operations research including but not limited to: scheduling problems with uncertain task durations \citep[Section~5.1]{chen2018distributionally}, multi-product assembly problems with uncertain demands (also known as the assemble-to-order system), and supply chain network design problems with uncertain demands and edge failure; see Examples~\ref{exp:scheduling}--\ref{exp:supplychain} in Section~\ref{ssec:examples}. 
The use of ambiguity set with fixed marginals is motivated by the experience that one typically has much less ambiguity about the marginal distributions of a multivariate uncertain quantity than about its dependence structure, as discussed by \citet*{eckstein2018robust, ren2021data, muller2025multivariate}. 
For example, when managing multiple types of risks, the probability laws of individual risk types can be modeled and estimated from historical data, for which there exist a plethora of parametric and non-parametric methods.
On the other hand, estimating the interdependence structure among risk types would require access to time-synchronized historical data, which is typically much more challenging to obtain and is often unavailable. 
Compared to moment-based constraints and Wasserstein distance based constraints with respect to a discrete measure, another advantage of ambiguity sets with marginal constraints is that they rule out highly unrealistic probability measures such as those supported on a small number of points; see, for example, \citep*[Proposition~1]{long2021supermodularity} and \citep[Theorem~3.3]{wozabal2012framework}. 
The goal of this paper is to develop a numerical algorithm for efficiently and accurately computing approximate optimizers as well as the associated approximate worst-case probability measures for this class of problems.

\subsection*{Literature review}
A widely adopted approach for solving two-stage optimization problems is called adaptive optimization (also known as adjustable optimization); see, e.g., \citep*{bental2004adjustable, bertsimas2015on, bertsimas2016duality, bertsimas2012on, bertsimas2018scalable, bertsimas2019adaptive, elhousni2021on, goh2010distributionally, xu2018copositive}. 
In adaptive optimization, rather than letting the second-stage decision be optimal given the first-stage decision and the uncertain quantities, one restricts the second-stage decision to depend on the uncertain quantities via pre-specified parametric decision rules, typically affine or piece-wise affine. 
\citet*{chen2020robust} proposed a framework for adaptive DRO with the so-called event-wise ambiguity set. Under the event-wise affine decision rule, the problem is computationally tractable, and conservative solutions can be computed using state-of-the-art commercial solvers. 
It has been empirically shown by \citet*{saif2021data} in a distributionally robust capacitated facility location problem that conservative solutions produced by adaptive DRO with an affine decision rule are comparable to exact solutions. 
Despite this, to the best of our knowledge, there is no theoretical bound for the sub-optimality of conservative solutions. 
In fact, the use of adaptive decision rules may even lead to infeasibility, as shown by \citet*[Equation~(13)]{bertsimas2019adaptive}. 

DRO problems with ambiguity sets constrained by marginals have been studied by \citet*{gao2017data}, \citet*{ren2021data}, and \citet*{chen2018distributionally}. 
\citet{gao2017data} have developed duality results for the inner maximization problem in DRO when the ambiguity set is subject to marginal constraints as well as a Wasserstein distance based constraint on its dependence structure. 
However, the computational tractability of the resulting dual formulation only holds under the restrictive assumptions that the (second-stage) cost function is the maximum of finitely many affine functions and that the marginals all have finite support.
\citet{ren2021data} considered two-stage DRO problems with marginal constraints where the marginals are discrete and subject to further ambiguities in the probabilities.
However, their LP formulation involves exponentially many constraints unless one restricts the support of the underlying probability measure to a small number of points.
\citet{chen2018distributionally} studied a special case of two-stage DRO problems 
and focused on deriving sufficient conditions for the polynomial time solvability of the inner maximization problem.
However, the sufficient conditions need to be established on a case-by-case basis, and it is difficult to check them for more sophisticated problems.
Moreover, they did not provide a concrete numerical algorithm for approximately solving the problem when the marginals are non-discrete.
Compared to \citep{chen2018distributionally}, the numerical algorithm that we develop is applicable to a broader class of two-stage DRO problems presented in Setting~\ref{sett:dro}, which subsumes the problems considered by \citet{chen2018distributionally} as a special case.

We would also like to highlight that the inner maximization problem in our model (i.e., the maximization over probability measures with respect to a fixed first-stage decision) is a special case of multi-marginal optimal transport problem \citep{benamou2021optimal, pass2015multi},
where the marginals are one-dimensional, and 
the cost function is the optimal value of an LP problem. 
\citet*{chen2020distributionally} have studied the distributionally robust max flows problem, which belongs to this problem class.
Moreover, instances from this class 
have also been applied to discrete choice modeling and assortment optimization; see, e.g., the studies of \citet*{natarajan2009persistency} and \citet*{ruan2022nonparametric}.
While most numerical methods for multi-marginal optimal transport and related problems with non-discrete marginals rely on discretization \citep*{carlier2015numerical, eckstein2019robust, guo2019computational, neufeld2021deep} and/or regularization techniques \citep*{aquino2019bounds, aquino2020minmax, eckstein2019computation, eckstein2018robust, henry2019martingale, neufeld2021model}, 
\citet{neufeld2022numerical} recently developed a numerical method that computes approximate optimizers as well as sub-optimality estimates of high-dimensional multi-marginal optimal transport problems without discretization or regularization. 

\subsection*{Contributions and organization of this paper}

The main contributions of this paper are summarized as follows. 
\begin{enumerate}[label=(\arabic*),beginpenalty=10000]
    \item
    We develop primal and dual formulations of a general class of two-stage DRO problems with marginal constraints,
    establish the respective strong duality,
    and prove the existence of primal and dual optimizers.
    Moreover, we derive the saddle point conditions for the primal and dual optimization problems which lead to sufficient and necessary conditions for their optimizers.
    These results are presented in Theorem~\ref{thm:primaldual}.

    \item
    We develop a scheme to construct both a sub-optimizer of the two-stage DRO problem and a sub-optimizer of the inner maximization problem from a sub-optimizer of the dual formulation,
    where the sub-optimality is given by the optimal value of a global maximization problem.
    These results are presented in Theorem~\ref{thm:dual-suboptimality} and Corollary~\ref{cor:dual-suboptimality-primal-coupling}.
    This gives rise to a numerical algorithm for computing an approximate optimizer of the two-stage DRO problem and an approximate worst-case probability measure.
    Our algorithm is guaranteed to converge and computes
    provable lower and upper bounds for the optimal value of the two-stage DRO problem as well as sub-optimality estimates for the computed approximate optimizers.
    In particular, the sub-optimality can be controlled to be arbitrarily small by suitably choosing the inputs of the algorithm.
    We present these properties of our algorithm in Theorem~\ref{thm:iterative}.

    \item
    We perform numerical experiments to demonstrate the proposed algorithm in three prominent two-stage DRO problems: 
    task scheduling, multi-product assembly, and supply chain network design with edge failure.
    The ambiguity sets in these problems involve a large number of marginals, 
    which include both discrete and continuous probability measures. 
    The numerical results show that the computed approximately optimizers are close to being optimal.
\end{enumerate}

The rest of this paper is organized as follows. 
Section~\ref{sec:dro} presents our two-stage DRO model with marginal constraints. 
Three prominent instances
of two-stage DRO problems subsumed by our model are discussed in Section~\ref{ssec:examples}. 
In Section~\ref{sec:primaldual},
we develop the primal and dual formulations and analyze their theoretical properties.
The numerical algorithm used for approximately solving the two-stage DRO problem as well as its theoretical properties are presented in Section~\ref{sec:algorithm}. 
Finally, in Section~\ref{sec:experiments}, we showcase
the performance of the proposed numerical algorithm in the three two-stage DRO problems discussed in
Section~\ref{ssec:examples}.
The appendix contains the proofs of the theoretical results in this paper. 

\section{Model for two-stage DRO problems with marginal constraints}
\label{sec:dro}

\subsection{Settings}
\label{ssec:settings}
Let us first introduce the notions and notations that are used throughout this paper. 
We denote by $\R$, $\R_{+}$, $\R_{-}$
the sets of real numbers, 
non-negative real numbers, and 
non-positive real numbers, respectively.
We assume all vectors to be column vectors and denote vectors and vector-valued functions by boldface symbols.
In particular, for $n\in\N$, we denote by $\veczero_n$ and $\vecone_n$ the vector in $\R^n$ with all entries equal to~0
and the vector in $\R^n$ with all entries equal to~1,
i.e., $\veczero_n:=(\underbrace{0,\ldots,0}_{n\mathrm{~times}})^\TRANSP$,
$\vecone_n:=(\underbrace{1,\ldots,1}_{n\mathrm{~times}})^\TRANSP$.
For $n\in\N$ and two vectors $\BIx=(x_1,\ldots,x_n)^\TRANSP\in\R^n$, $\BIy=(y_1,\ldots,y_n)^\TRANSP\in\R^n$, 
we let $\BIx\le\BIy$ and $\BIx\ge\BIy$ denote the corresponding component-wise inequalities, i.e., 
$\BIx\le\BIy \;\Leftrightarrow$ $x_1\le y_1,\;\ldots,\;x_n\le y_n$,
$\BIx\ge\BIy \;\Leftrightarrow$ $x_1\ge y_1,\;\ldots,\;x_n\ge y_n$.
Moreover, we denote by $\langle\,\cdot\,,\cdot\,\rangle$ the Euclidean dot product, i.e., $\langle\BIx,\BIy\rangle:=\BIx^\TRANSP\BIy$.
We call a subset of a Euclidean space a polyhedral convex set if it is the intersection of finitely many closed half-spaces. 
For a subset $A$ of a Euclidean space, we use $\clos(A)$ and $\conv(A)$ to denote the closure and convex hull of $A$, respectively. 
For a polyhedral convex set $P$, we let 
$\extremept(P)$ denote the set of extreme points of~$P$, 
let $\recession(P)$ denote the recession cone of $P$,
and let $\lineality(P):=\recession(P)\intersection {-\recession(P)}$ denote the lineality space of $P$; see, e.g., \citep[p.61 \& p.65]{rockafellar1970convex}. 
Furthermore, for $a,b\in\R$, we let $a\wedge b$ denote $\min\{a,b\}$, let $a\vee b$ denote $\max\{a,b\}$, and let $(a)^+$ denote $\max\{a,0\}$. 
For any closed set $\CY\subseteq\R^n$ with $n\in\N$, we let $\CB(\CY)$ denote the Borel subsets of $\CY$, 
let $\CP(\CY)$ denote the set of Borel probability measures on $\CY$,
and let $\CP_1(\CY)$ denote the set of Borel probability measures on $\CY$ with finite first moments.
Let $\support(\mu)$ denote the support of any $\mu\in\CP(\CY)$, and let 
$\delta_{\BIy}$ denote the Dirac measure at $\BIy\in\CY$.
For any $\mu\in\nobreak\CP(\CY)$ and any Borel measurable function $g:\CY\to\R$,
we denote by $g\sharp\mu$ the pushforward of~$\mu$ by~$g$,
i.e., $g\sharp\mu(B):=\mu\big(g^{-1}(B)\big)$ for all $B\in\CB(\R)$.
Moreover, for $n\in\N$ and for any lower semi-continuous and proper convex function $f:\R^d\to\R\union\{\infty\}$, 
we let $\partial f(\BIz):=\big\{\BIg\in\R^n: f(\BIx) \ge f(\BIz) + \langle\BIg,\BIx-\BIz\rangle\; \forall \BIx\in\R^n\big\}$ denote the subdifferential of~$f$ at~$\BIz$.
We let $\CC(\R)$ denote the set of continuous functions 
on $\R$, 
and let $\CC_{\mathrm{cvx}}(\R)$ denote 
the set of continuous and convex functions on $\R$. 
Lastly, for any $f\in\CC(\R)$ and any $z\in\R$,
we let $\partial_{-}f(z)$ and $\partial_{+} f$ denote the left- and right-derivatives of~$f$ at~$z$,
i.e., $\partial_{-}f(z):=\lim_{h\to0^{-}}\frac{f(z+h)-f(z)}{h}$,
$\partial_{+}f(z):=\lim_{h\to0^{+}}\frac{f(z+h)-f(z)}{h}$.
In particular, for any $f\in\CC_{\mathrm{cvx}}(\R)$,
it holds that $\partial f(z)=\big[\partial_{-}f(z), \partial_{+}f(z)\big]\ne\emptyset$ for all $z\in\R$.

We consider a two-stage distributionally robust optimization (DRO) problem with marginal constraints in which both the first-stage and the second-stage optimization problems have linear objectives and constraints. 
The concrete details of our settings are presented below in Setting~\ref{sett:dro}.

\begin{setting}
    \label{sett:dro}
    In the two-stage DRO problem with marginal constraints, we introduce the following assumptions and notations.
    \begin{enumerate}[label=\normalfont{(A\arabic*)},leftmargin=25pt]
        \item\label{setts:dro-uncertainty}
        The outcome of the random event is given by the space $\BOmega:={\Omega_1\times\cdots\times\Omega_d}$, where $d\in\nobreak\N\intersection [2,\infty)$ and, for $i=1,\ldots,d$, $\Omega_i:=[\underline{\omega}_i, \overline{\omega}_i]$ for some $-\infty<\underline{\omega}_i\le \overline{\omega}_i<\infty$.
        For $i=1,\ldots,d$, the marginal probability law for the $i$-th component of the random outcome is given by $\mu_i\in\CP(\Omega_i)$.
        Let $\Gamma(\mu_1,\ldots,\mu_d)$ denote the set of couplings of $\mu_1,\ldots,\mu_d$, defined as 
        \begin{align*}
        \hspace{33pt}\Gamma(\mu_1,\ldots,\mu_d):=\hspace{-0.15em}\Big\{\gamma\in\CP(\Omega_1\times\cdots\times \Omega_d):\text{the marginal of }\gamma\text{ on }\Omega_i\text{ is }\mu_i\; \forall 1\le i\le d\Big\}.
        \end{align*}

        \item\label{setts:dro-stage1}
        The first-stage decision variable is denoted by $\BIx\in\R_{+}^{k_1}$ where $k_1\in\N$.
        The first-stage cost function is given by $\R_{+}^{k_1}\ni \BIx\mapsto \langle\Bc_1,\BIx\rangle\in\R$ for some $\Bc_1\in\R^{k_1}$.
        The feasible set of the first-stage decision problem is given by the non-empty polyhedral convex set $\CS_1:=\big\{\BIx\in\nobreak\R_{+}^{k_1}:{\BL_{\ineq}\BIx\le\Bq_{\ineq}},\;\BL_{\eq}\BIx=\Bq_{\eq}\big\}$, where $n_{\ineq}\in\N_0$, $n_{\eq}\in\N_0$, $\BL_{\ineq}\in\R^{n_{\ineq}\times k_1}$, $\Bq_{\ineq}\in\nobreak\R^{n_{\ineq}}$, $\BL_{\eq}\in\R^{n_{\eq}\times k_1}$, $\Bq_{\eq}\in\R^{n_{\eq}}$.

        \item\label{setts:dro-stage2}
        The second-stage decision variable is denoted by $\BIz\in\R_{+}^{k_2}$ where $k_2\in\N$.
        The second-stage cost function is given by $\R_{+}^{k_2}\ni\BIz\mapsto\langle\Bc_2,\BIz\rangle\in\R$ for some $\Bc_2\in\R^{k_2}$.
        The feasible set of the second-stage decision problem is a non-empty polyhedral convex set $\CS_2(\BIx,\Bomega):=\big\{\BIz\in\R_{+}^{k_2}:\BA_{\ineq}\BIz\le {\BQ_{\ineq}\BIx+\BP_{\ineq}\Bomega+\Bb_{\ineq}},\; \BA_{\eq}\BIz=\nobreak \BQ_{\eq}\BIx+\nobreak\BP_{\eq}\Bomega+\nobreak\Bb_{\eq}\big\}$ that depends on the first-stage decision $\BIx\in\CS_1$ and the random outcome~$\Bomega\in\BOmega$,
        where $m_{\ineq}\in\nobreak\N_0$, 
        $m_{\eq}\in\N_0$,
        $m_{\ineq}+m_{\eq}>0$,
        $\BA_{\ineq}\in\R^{m_{\ineq}\times k_2}$, 
        $\BQ_{\ineq}\in\R^{m_{\ineq}\times k_1}$, 
        $\BP_{\ineq}\in\R^{m_{\ineq}\times d}$, 
        $\Bb_{\ineq}\in\R^{m_{\ineq}}$,
        ${\BA_{\eq}\in\R^{m_{\eq}\times k_2}}$, 
        $\BQ_{\eq}\in\R^{m_{\eq}\times k_1}$, 
        $\BP_{\eq}\in\R^{m_{\eq}\times d}$, 
        $\Bb_{\eq}\in\R^{m_{\eq}}$.
        The optimal value of the second-stage decision problem, also referred to as the optimal second-stage cost function, is denoted by $C_2(\BIx,\Bomega):=\inf_{\BIz\in \CS_2(\BIx,\Bomega)}\big\{\langle\Bc_2,\BIz\rangle\big\}$ for all $\BIx\in\CS_1$ and all $\Bomega\in\BOmega$.
        Moreover, we assume that $\alpha:=\sup_{\BIx\in\CS_1,\,\Bomega\in\BOmega}\big\{C_2(\BIx,\Bomega)\big\}<\infty$.
    \end{enumerate}
    We aim to approximately solve
    the following two-stage DRO problem:
    \begin{align}
    \begin{split}
    C_{\DRO}&:=\inf_{\BIx\in \CS_1}\bigg\{\langle\Bc_1,\BIx\rangle+\sup_{\mu\in\Gamma(\mu_1,\ldots,\mu_d)}\bigg\{\int_{\BOmega}C_2(\BIx,\Bomega)\DIFFM{\mu}{\DIFF\Bomega}\bigg\}\bigg\},
    \end{split}
    \tag{\textsf{DRO}}
    \label{eqn:dro}
    \end{align}
    where $C_{\DRO}$ denotes its optimal value.
\end{setting}

In the problem defined in Setting~\ref{sett:dro}, $\langle\Bc_1,\BIx\rangle$ corresponds to first-stage cost when the first-stage decision is~$\BIx$,
the term $\int_{\BOmega}C_2(\BIx,\Bomega)\DIFFM{\mu}{\DIFF\Bomega}$ corresponds to the expected optimal second-stage cost under the probability measure $\mu\in\CP(\BOmega)$, and \sloppy{$\sup_{\mu\in\Gamma(\mu_1,\ldots,\mu_d)}\big\{\int_{\BOmega}C_2(\BIx,\Bomega)\DIFFM{\mu}{\DIFF\Bomega}\big\}$} corresponds to the worst-case expected optimal second-stage cost when the probability measure can be any coupling of the given marginals $\mu_1,\ldots,\mu_d$. 
The problem \eqref{eqn:dro} minimizes the worst-case expected cost in both decision stages.

Under Setting~\ref{sett:dro}, the second-stage decision problem $C_2(\BIx,\Bomega)$ is a linear minimization problem in which the right-hand side of the equality and inequality constraints have a jointly affine dependence on the first-stage decision~$\BIx$ and the random outcome~$\Bomega$. 
This type of problem structure has been widely studied in the literature in the context of robust optimization (see, e.g., \citep{bertsimas2010models, bertsimas2012on, bertsimas2015on, bertsimas2016duality, bertsimas2018scalable, elhousni2021on, xu2018copositive}) and DRO (see, e.g., \citep{long2021supermodularity}).

In order to make the subsequent analyses tractable, let us first replace $C_2(\BIx,\Bomega)$ in \eqref{eqn:dro} with its dual LP problem. 

\begin{lemma}
    \label{lem:stage2-dual}
    Under Setting~\ref{sett:dro}, the following statements hold.
    \begin{enumerate}[label=(\roman*), leftmargin=20pt]
    \item\label{lems:stage2-dual-duality}
    The following duality holds:
    \begin{align*}
    C_2(\BIx,\Bomega)&=\sup_{(\Blambda_{\ineq}^\TRANSP,\Blambda_{\eq}^\TRANSP)^\TRANSP\in \CS_2^*}\big\{\langle\BQ_{\ineq}\BIx+\BP_{\ineq}\Bomega+\Bb_{\ineq},\Blambda_{\ineq}\rangle+\langle\BQ_{\eq}\BIx+\BP_{\eq}\Bomega+\Bb_{\eq},\Blambda_{\eq}\rangle\big\}\\
    &\hspace{23em}\forall \BIx\in \CS_1,\;\forall\Bomega\in\BOmega,
    \end{align*}
    where
    \begin{align}
    \CS_2^*:=\big\{(\Blambda_{\ineq}^\TRANSP,\Blambda_{\eq}^\TRANSP)^\TRANSP:\Blambda_{\ineq}\in\R^{m_{\ineq}}_{-},\;\Blambda_{\eq}\in\R^{m_{\eq}},\;\BA_{\ineq}^\TRANSP\Blambda_{\ineq}+\BA_{\eq}^\TRANSP\Blambda_{\eq}\le \Bc_2\big\}
    \label{eqn:stage2-dual-feasible-set}
    \end{align}
    is a polyhedral convex set in $\R^{m_{\ineq}+m_{\eq}}$.

    \item\label{lems:stage2-dual-finiteset}
    There exists a finite set $\CV\subset\CS^*_2$ that satisfies
    \begin{align*}
    C_2(\BIx,\Bomega)&=\max_{(\Blambda_{\ineq}^\TRANSP,\Blambda_{\eq}^\TRANSP)^\TRANSP\in \CV}\big\{\langle\BQ_{\ineq}\BIx+\BP_{\ineq}\Bomega+\Bb_{\ineq},\Blambda_{\ineq}\rangle+\langle\BQ_{\eq}\BIx+\BP_{\eq}\Bomega+\Bb_{\eq},\Blambda_{\eq}\rangle\big\} \\
    &\hspace{23em}\forall \BIx\in \CS_1,\;\forall\Bomega\in\BOmega.
    \end{align*}
    In particular, $\CV$ can be chosen to be equal to the finite set $\CV_{\extremept}:=\extremept\big(\CS^*_2\intersection \lineality(\CS^*_2)^{\perp}\big)$ or any finite superset of $\CV_{\extremept}$,
    where $\lineality(\CS^*_2)^{\perp}$ denotes the linear subspace of $\R^{m_{\ineq}+m_{\eq}}$ orthogonal to $\lineality(\CS^*_2)$.
    \end{enumerate}
\end{lemma}

Due to Lemma~\ref{lem:stage2-dual}, we present the following alternative setting that is equivalent to Setting~\ref{sett:dro} in order to simplify the subsequent analyses. 

\begin{setting}
In addition to \ref{setts:dro-uncertainty} and \ref{setts:dro-stage1} in Setting~\ref{sett:dro}, we introduce the following assumptions and notations.
\begin{enumerate}[label=\normalfont{(A\arabic**)},leftmargin=32pt]
\setcounter{enumi}{2}
\item \label{setts:dro-stage2-alt}
The optimal second-stage cost function is given by $C_2(\BIx,\Bomega):=\max_{\Blambda\in \CS_2^*}\big\{\langle\BQ\BIx+\BP\Bomega+\Bb,\Blambda\rangle\big\}$ for all $\BIx\in\CS_1$ and all $\Bomega\in\BOmega$, 
where $k^*_2\in\N$, $\BQ\in\R^{k^*_2\times k_1}$, $\BP\in\R^{k^*_2\times d}$, $\Bb\in\nobreak\R^{k^*_2}$, and $\CS^*_2\subset\R^{k^*_2}$ is a non-empty polyhedral convex set which is defined by (\ref{eqn:stage2-dual-feasible-set}) with 
$m_{\ineq}\in\nobreak\N_0$, 
$m_{\eq}\in\N_0$,
$m_{\ineq}+m_{\eq}=k^*_2$,
$k_2\in\N$,
$\BA_{\ineq}\in\R^{m_{\ineq}\times k_2}$, 
${\BA_{\eq}\in\R^{m_{\eq}\times k_2}}$,
$\Bc_2\in\R^{k_2}$. 
Moreover, we assume that $\alpha:=\sup_{\BIx\in\CS_1,\,\Bomega\in\BOmega}\big\{C_2(\BIx,\Bomega)\big\}<\infty$.
Furthermore, we denote $\CV_{\extremept}:=\extremept\big(\CS^*_2\intersection \lineality(\CS^*_2)^{\perp}\big)$, which is a finite subset of $\CS^*_2$.
\end{enumerate}
\label{sett:dro-alt}
\end{setting}

Under Setting~\ref{sett:dro-alt}, Lemma~\ref{lem:stage2-dual}\ref{lems:stage2-dual-finiteset} implies that
whenever $\CV$ is a finite set where $\CV_{\extremept}\subseteq\CV \subset\CS^*_2$,
it holds that $C_2(\BIx,\Bomega)={\max_{\Blambda\in \CV}\big\{\langle\BQ\BIx+\BP\Bomega+\Bb,\Blambda\rangle\big\}}$ 
for all $\BIx\in\nobreak\CS_1$ and all $\Bomega\in\BOmega$.

\begin{remark}\label{rmk:Chen-marginal-DRO}
    The marginal distribution model considered by \citet{chen2018distributionally} can be seen as a special case of our model in
    Setting~\ref{sett:dro-alt},
    where \citet{chen2018distributionally} only study the inner maximization problem in \eqref{eqn:dro}, i.e.,
    $\sup_{\mu\in\Gamma(\mu_1,\ldots,\mu_d)}\big\{\int_{\BOmega}C_2(\BIx_0,\Bomega)\DIFFM{\mu}{\DIFF\Bomega}\big\}$ with respect to a fixed first-stage decision $\BIx_0$,
    rather than the joint min-max problem in \eqref{eqn:dro}.
    Specifically, one can observe that Setting~\ref{sett:dro-alt} subsumes 
    the model of \citet{chen2018distributionally}
    when 
    $k_1=k^*_2=d$, 
    $\BQ=-\BI_d$,
    $\BP=\BI_d$ (where $\BI_d$ denotes the $d$-by-$d$ identity matrix), $\Bb=\veczero_d$, 
    and $\CS_1=\{\BIx_0\}$ 
    (i.e., the first-stage decision problem is trivial).
\end{remark}

\subsection{Examples}
\label{ssec:examples}

The model for two-stage DRO problems with marginal constraints introduced in Setting~\ref{sett:dro} covers a wide range of prominent decision problems in operations research,
including but not limited to:
task scheduling \citep{chen2018distributionally}, 
multi-product assembly \citep{atan2017assemble}, 
supply chain network design \citep{atamturk2007two, cheng2018two, matthews2019designing},
the newsvendor problem \citep{shapiro2009lectures, wang2016likelihood, wiesemann2014distributionally},
lot sizing on a network \citep{bertsimas2016duality, xu2018copositive}, 
and resource allocation in a temporal network~\citep{wiesemann2012robust}.
In the following, we discuss in detail three examples of the model: task scheduling, multi-product assembly, and supply chain network design with edge failure.

\begin{example}[Task scheduling \normalfont{\citep[Section~5.1]{chen2018distributionally}}]
\label{exp:scheduling}
In this problem, there are $d\in\N$ tasks arranged in a fixed order and need to be scheduled within a fixed time window $[0,T]$ for $T>0$. 
For $i=1,\ldots,d$, let $x_i\in\R_+$ denote the scheduled duration of the $i$-th task. Hence, the $i$-th task will be scheduled to begin at time $\sum_{j=1}^{i-1}x_j$.  
The actual duration of the $i$-th task, which is denoted by $\omega_i$, is a random variable with law $\mu_i\in\CP\big([\underline{\omega}_i,\overline{\omega}_i]\big)$ for $i=1,\ldots,d$, where $\underline{\omega}_i\ge 0$ and $\overline{\omega}_i>\underline{\omega}_i$ are lower and upper bounds for the duration of the $i$-th task. 
There is no information about the interdependence among the durations of tasks. 
It is assumed that the $(i+1)$-th task can only begin after the $i$-th task has been completed. Since the actual duration may be longer than the scheduled duration, the $i$-th task may incur a delay, denoted by $z_i$, which is defined recursively as follows: 
$z_1:=(\omega_1-x_1)^+\in\R_+$, 
$z_i:=(z_{i-1}+\omega_i-x_i)^+\in\R_+$ for $i=2,\ldots,d$, 
that is, 
$z_i$ is the difference between the actual and the scheduled completion time of the $i$-th task if it is positive, and otherwise $z_i:=0$. 
The objective of the task scheduling problem is to minimize a weighted total delay, i.e., $\sum_{i=1}^dc_{2,i}z_i$, where $c_{2,1},\ldots,c_{2,d}\in\R_+$ denote the weights.

Formulating this problem into our two-stage DRO model in Setting~\ref{sett:dro}, the first-stage decision is $\BIx:=(x_1,\ldots,x_d)^\TRANSP\in\R^d_{+}$, and we have $k_1:=d$, $\CS_1:=\big\{(x_1,\ldots,x_d)^\TRANSP\in\R_{+}^d:\sum_{i=1}^dx_i\le T\big\}$. 
Since no cost is incurred in the first stage, we have $\Bc_1:=\veczero_d$. 
The optimal second-stage cost function $C_2(\,\cdot\,,\cdot\,)$ is given by
\begin{align*}
C_2\big((x_1,\ldots,x_d)^\TRANSP,(\omega_1,\ldots,\omega_d)^\TRANSP\big)&:=\min\left\{\sum_{i=1}^dc_{2,i}z_i:
\begin{tabular}{l}
    $z_i\in\R_+\;\forall 1\le i\le d,$\\
    $z_1\ge \omega_1-x_1,$\\
    $z_i\ge z_{i-1}+\omega_i-x_i\;\forall 2\le i\le d$
\end{tabular}%
\!\right\}\\
& \hspace{70pt} \forall (x_1,\ldots,x_d)^\TRANSP\in\CS_1,\; \forall (\omega_1,\ldots,\omega_d)^\TRANSP\in\BOmega,
\end{align*}
which has the required form in \ref{setts:dro-stage2} with $k_2:=d$. 
Moreover, 
for any $(x_1,\ldots,x_d)^\TRANSP\in\CS_1$ and 
any $(\omega_1,\ldots,\omega_d)^\TRANSP\in\BOmega$,
observe that defining
$\tilde{z}_1:=(\omega_1-x_1)^+\in\R_+$, 
$\tilde{z}_i:=(\tilde{z}_{i-1}+\omega_i-x_i)^+\in\R_+$ for $i=2,\ldots,d$ 
yields an optimizer $(\tilde{z}_1,\ldots,\tilde{z}_d)^\TRANSP$ of $C_2\big((x_1,\ldots,x_d)^\TRANSP,(\omega_1,\ldots,\omega_d)^\TRANSP\big)$,
where $\tilde{z}_i\le \sum_{j=1}^{i}\omega_j\le \sum_{j=1}^{i}\overline{\omega}_j$ for $i=1,\ldots,d$.
It thus holds that $C_2\big((x_1,\ldots,x_d)^\TRANSP,(\omega_1,\ldots,\omega_d)^\TRANSP\big)\le \sum_{i=1}^{d}c_{2,i}\sum_{j=1}^{i}\overline{\omega}_j<\infty$ for all $(x_1,\ldots,x_d)^\TRANSP\in\CS_1$ and all $(\omega_1,\ldots,\omega_d)^\TRANSP\in\BOmega$,
which shows that the assumption~\ref{setts:dro-stage2} is satisfied with respect to some $\alpha\le\sum_{i=1}^{d}c_{2,i}\sum_{j=1}^{i}\overline{\omega}_j$. 
\end{example}

\begin{example}[Multi-product assembly]
\label{exp:multiproduct}
This example is the distributionally robust version of the multi-product assembly problem adapted from Chapter~1.3.1 of \citep{shapiro2009lectures}. This is also known as the assemble-to-order system; see, e.g., the review about assemble-to-order systems by \citet*{atan2017assemble} and the references therein.
In this problem, let us consider a manufacturer that produces $d\in\N$ products. The production of these products requires $k_1\in\N$ parts that need to be ordered from suppliers with per-unit prices $c_{1,1}>0,\ldots,c_{1,k_1}>0$. 
Producing each unit of product $i$ requires $u_{i,j}\ge0$ units of part $j$ for $i=1,\ldots,d$, $j=1,\ldots,k_1$. 
The demand for product $i$, denoted by $\omega_i$, is a random variable with law $\mu_i\in\CP\big([\underline{\omega}_i,\overline{\omega}_i]\big)$ for $i=1,\ldots,d$, where $\underline{\omega}_i\ge 0$ and $\overline{\omega}_i>\underline{\omega}_i$ are lower and upper bounds for the demand for product~$i$. 
There is no information about the interdependence among the demands for products. 
Once the demands for the products are known, the manufacturer decides the quantity of each product to produce. The quantity $z_i$ of product~$i$ produced shall not exceed its demand $\omega_i$. 
The production of each unit of product $i$ earns the manufacturer a revenue of $r_i>0$. 
After production, the unused parts $h_1,\ldots,h_{k_1}$ have salvage values $s_1,\ldots,s_{k_1}$ such that $0\le s_j<c_{1,j}$ for $j=1,\ldots,k_1$. 

Formulating this problem into our two-stage DRO model in Setting~\ref{sett:dro}, the first-stage decision is $\BIx:=(x_1,\ldots,x_{k_1})^\TRANSP\in\R_{+}^{k_1}$, which corresponds to the quantities of parts to be ordered from the suppliers. We have $\CS_1:=\R^{k_1}_+$ and $\Bc_1:=(c_{1,1},\ldots,c_{1,k_1})^\TRANSP\in\nobreak\R_{+}^{k_1}$. 
Let 
$\Bc_{2,\mathrm{reve}}:=-{(r_1,\ldots,r_d)^\TRANSP}\in\R_{-}^{d}$, 
$\Bc_{2,\mathrm{salv}}:=-{(s_1,\ldots,s_{k_1})^\TRANSP}\in\nobreak\R_{-}^{k_1}$, 
$\BIz_{\mathrm{prod}}:=(z_1,\ldots,z_d)^\TRANSP\in\R_{+}^d$, 
$\BIz_{\mathrm{salv}}:=(h_1,\ldots,h_{k_1})^\TRANSP\in\R_{+}^{k_1}$, 
and let 
$\BU\in\R^{d\times k_1}$ be a matrix with entries $[\BU]_{i,j}:=u_{i,j}$ for $i=1,\ldots,d$, $j=1,\ldots,k_1$. 
Then, the optimal second-stage cost function $C_2(\BIx,\Bomega)$ is given by
\begin{align*}
C_2(\BIx,\Bomega)&:=\min\left\{\langle\Bc_{2,\mathrm{reve}},\BIz_{\mathrm{prod}}\rangle+\langle\Bc_{2,\mathrm{salv}},\BIz_{\mathrm{salv}}\rangle:
\begin{tabular}{l}
    $\BIz_{\mathrm{prod}}\in\R_{+}^{d},\;\BIz_{\mathrm{salv}}\in\R_{+}^{k_1},$ \\
    $\BU^\TRANSP\BIz_{\mathrm{prod}}+\BIz_{\mathrm{salv}}=\BIx,\;\BIz_{\mathrm{prod}}\le\Bomega$
\end{tabular}\!\!
\right\}\\
&\hspace{286pt} \forall \BIx\in\CS_1,\; \forall \Bomega\in\BOmega,
\end{align*}
which has the required form in \ref{setts:dro-stage2} with $k_2:=k_1+d$.
Moreover, $C_2(\BIx,\Bomega)\le 0$ holds for all $\BIx\in \CS_1$ and all $\Bomega\in\BOmega$ and thus the assumption \ref{setts:dro-stage2} is satisfied with respect to some $\alpha\le0$.  
\end{example}

\begin{example}[Supply chain network design with edge failure]
\label{exp:supplychain}
This example is a distributionally robust supply chain network design problem inspired by the non-robust version in Chapter~1.5 of \citep{shapiro2009lectures}, as well as the studies of \citet*{atamturk2007two, cheng2018two}, and \citet*{matthews2019designing}.
In this problem, we consider a supply chain network $(\TTV,\TTE)$ of a type of goods in which the vertices $\TTV:=\TTS\union\TTP\union\TTC$ consist of the suppliers~$\TTS$, the processing facilities~$\TTP$, and the customers~$\TTC$. The edges $\TTE:=\TTE_{\TTS,\TTP}\union\TTE_{\TTP,\TTC}$ consist of edges $\TTE_{\TTS,\TTP}\subseteq\TTS\times\TTP$ from the suppliers to the processing facilities and edges $\TTE_{\TTP,\TTC}\subseteq\TTP\times\TTC$ from the processing facilities to the customers. 
Each supplier $\TTs\in\TTS$ can supply $b_{\TTs}>0$ amount of goods which is known prior to the first decision stage.
Each processing facility $\TTp\in\TTP$ has a maximum processing capability of~$\overline{x}_{\TTp}>0$ and is associated with an investment cost of $c_{1,\TTp}>0$ for each unit of processing capability. 
Each customer $\TTc\in\TTC$ has demand $\omega_{\TTc}$ which is a random variable with law $\mu_{\TTc}\in\CP\big([\underline{\omega}_{\TTc},\overline{\omega}_{\TTc}]\big)$ where $\underline{\omega}_{\TTc}\ge0$ and $\overline{\omega}_{\TTc}>\underline{\omega}_{\TTc}$ are the minimum and maximum demands of the customer $\TTc$. 
Each edge $(\TTs,\TTp)\in\TTE_{\TTS,\TTP}$ is associated with a cost $c_{2,\TTs,\TTp}>0$ which is the sum of the per-unit transportation cost along this edge and the per-unit processing cost at the processing facility~$\TTp$. 
Similarly, each edge $(\TTp,\TTc)\in\TTE_{\TTP,\TTC}$ is associated with a per-unit transportation cost $c_{2,\TTp,\TTc}>0$.
Moreover, each edge $(\TTs,\TTp)\in\TTE_{\TTS,\TTP}$ has a maximum transportation capacity of $b_{\TTs,\TTp}>0$ and each edge $(\TTp,\TTc)\in\TTE_{\TTP,\TTC}$ has a maximum transportation capacity of $b_{\TTp,\TTc}>0$. 
Furthermore, there are subsets of edges $\widetilde{\TTE}_{\TTS,\TTP}\subset\TTE_{\TTS,\TTP}$ and $\widetilde{\TTE}_{\TTP,\TTC}\subset\TTE_{\TTP,\TTC}$ that are susceptible to failure. 
For each edge $(\TTs,\TTp)\in\widetilde{\TTE}_{\TTS,\TTP}$, let $\omega_{\TTs,\TTp}$ be a Bernoulli random variable indicating the failure of $(\TTs,\TTp)$, i.e., $\omega_{\TTs,\TTp}=1$ indicates the failure of $(\TTs,\TTp)$ while $\omega_{\TTs,\TTp}=0$ indicates no failure.
Similarly, for each edge $(\TTp,\TTc)\in\widetilde{\TTE}_{\TTP,\TTC}$, let $\omega_{\TTp,\TTc}$ be a Bernoulli random variable indicating 
the failure of $(\TTp,\TTc)$. 
There is no information about the interdependence among the demands and the failure of the edges. 

In the first decision stage, the decision maker determines the amount of investment for the processing capability $x_{\TTp}\ge0$ of each processing facility $\TTp\in\TTP$. 
In the second decision stage, with the processing capabilities $(x_{\TTp})_{\TTp\in\TTP}$, the demands $(\omega_{\TTc})_{\TTc\in\TTC}$, and the edge failures $(\omega_{\TTs,\TTp})_{(\TTs,\TTp)\in\widetilde{\TTE}_{\TTS,\TTP}}$, $(\omega_{\TTp,\TTc})_{(\TTp,\TTc)\in\widetilde{\TTE}_{\TTP,\TTC}}$ known, the decision maker minimizes the total operational cost including the transportation costs and the processing costs. 
Let $\TTE_{\TTs,\cdot}:=\big\{\TTp'\in\TTP:(\TTs,\TTp')\in\TTE_{\TTS,\TTP}\big\}$, 
$\TTE_{\hspace{1pt}\cdot,\TTp}:=\big\{\TTs'\in\TTS:(\TTs',\TTp)\in\TTE_{\TTS,\TTP}\big\}$, 
$\TTE_{\TTp,\cdot}:=\big\{\TTc'\in\TTC:(\TTp,\TTc')\in\TTE_{\TTP,\TTC}\big\}$, 
$\TTE_{\hspace{1pt}\cdot,\TTc}:=\big\{\TTp'\in\TTP:(\TTp',\TTc)\in\TTE_{\TTP,\TTC}\big\}$ for all $\TTs\in\TTS$, $\TTp\in\TTP$, $\TTc\in\TTC$. 
The second-stage decision problem can be formulated as follows:
\begin{align*}
\minimize_{(z_{\TTs,\TTp}),\,(z_{\TTp,\TTc})}\quad & \textstyle\left(\sum_{(\TTs,\TTp)\in\TTE_{\TTS,\TTP}}c_{2,\TTs,\TTp}z_{\TTs,\TTp}\right)+\left(\sum_{(\TTp,\TTc)\in\TTE_{\TTP,\TTC}}c_{2,\TTp,\TTc}z_{\TTp,\TTc}\right)\\
\mathrm{subject~to}\quad 
& \textstyle\left(\sum_{\TTs\in\TTE_{\hspace{1pt}\cdot,\TTp}}z_{\TTs,\TTp}\right)-\left(\sum_{\TTc\in\TTE_{\TTp,\cdot}}z_{\TTp,\TTc}\right)=0 
    \hspace{126.8pt} \forall \TTp\in\TTP,\allowdisplaybreaks\\
& \textstyle\hspace{93.6pt}\sum_{\TTp\in\TTE_{\TTs,\cdot}}z_{\TTs,\TTp}\le b_{\TTs}
    \hspace{122.6pt}\forall \TTs\in\TTS,\allowdisplaybreaks\\
& \textstyle\hspace{92.6pt}\sum_{\TTp\in\TTE_{\hspace{1pt}\cdot,\TTc}}z_{\TTp,\TTc}\ge \omega_{\TTc} 
    \hspace{120.7pt}\forall \TTc\in\TTC,\allowdisplaybreaks\\
& \textstyle\hspace{92.6pt}\sum_{\TTs\in\TTE_{\hspace{1pt}\cdot,\TTp}}z_{\TTs,\TTp}\le x_{\TTp}
    \hspace{121.2pt}\forall \TTp\in\TTP,\allowdisplaybreaks\\
& \hspace{131.8pt} z_{\TTs,\TTp}\le b_{\TTs,\TTp} 
    \hspace{58.1pt} \forall (\TTs,\TTp)\in\TTE_{\TTS,\TTP}\setminus\widetilde{\TTE}_{\TTS,\TTP},\allowdisplaybreaks\\
& \hspace{131.8pt} z_{\TTp,\TTc}\le b_{\TTp,\TTc} 
    \hspace{58.1pt} \forall (\TTp,\TTc)\in\TTE_{\TTP,\TTC}\setminus\widetilde{\TTE}_{\TTP,\TTC},\allowdisplaybreaks\\
& \hspace{131.8pt} z_{\TTs,\TTp}\le b_{\TTs,\TTp}(1-\omega_{\TTs,\TTp}) 
    \hspace{40pt} \forall (\TTs,\TTp)\in\widetilde{\TTE}_{\TTS,\TTP},\allowdisplaybreaks\\
& \hspace{131.8pt} z_{\TTp,\TTc}\le b_{\TTp,\TTc}(1-\omega_{\TTp,\TTc}) 
    \hspace{40pt} \forall (\TTp,\TTc)\in\widetilde{\TTE}_{\TTP,\TTC},\\
& \hspace{64.8pt} z_{\TTs,\TTp}\in\R_+ \quad \forall (\TTs,\TTp)\in\TTE_{\TTS,\TTP}, \qquad 
z_{\TTp,\TTc}\in\R_+ \quad \forall (\TTp,\TTc)\in\TTE_{\TTP,\TTC}.
\end{align*}
In the above problem, the decision variables $(z_{\TTs,\TTp})_{(\TTs,\TTp)\in\TTE_{\TTS,\TTP}}\subset\R_+$ and $(z_{\TTp,\TTc})_{(\TTp,\TTc)\in\TTE_{\TTP,\TTC}}\subset\R_+$ represent the amount of goods flowing through the edges in the supply chain network. 
The objective $\big(\sum_{(\TTs,\TTp)\in\TTE_{\TTS,\TTP}}c_{2,\TTs,\TTp}z_{\TTs,\TTp}\big)+\big(\sum_{(\TTp,\TTc)\in\TTE_{\TTP,\TTC}}c_{2,\TTp,\TTc}z_{\TTp,\TTc}\big)$ corresponds to the total transportation and processing cost when transporting the goods from the suppliers to the processing facilities and then to the customers after processing. 
The constraints $\big(\sum_{\TTs\in\TTE_{\hspace{1pt}\cdot,\TTp}}z_{\TTs,\TTp}\big)-\big(\sum_{\TTc\in\TTE_{\TTp,\cdot}}z_{\TTp,\TTc}\big)=0$ $\forall \TTp\in\TTP$ correspond to the flow conservation condition requiring the amount of goods that flow into each processing facility $\TTp\in\TTP$ to be equal the amount of goods that flow out of $\TTp$. 
The constraints $\sum_{\TTp\in\TTE_{\TTs,\cdot}}z_{\TTs,\TTp}\le b_{\TTs}$ $\forall \TTs\in\TTS$ require the amount of goods that flow out of each supplier $\TTs\in\TTS$ to not exceed its total supply $b_{\TTs}$. 
The constraints $\sum_{\TTp\in\TTE_{\hspace{1pt}\cdot,\TTc}}z_{\TTp,\TTc}\ge \omega_{\TTc}$ $\forall \TTc\in\TTC$ require the amount of goods that flow into each customer $\TTc\in\TTC$ to meet its demand $\omega_{\TTc}$. 
The constraints $\sum_{\TTs\in\TTE_{\hspace{1pt}\cdot,\TTp}}z_{\TTs,\TTp}\le x_{\TTp}$ $\forall \TTp\in\TTP$ require the amount of goods that flow into each processing facility $\TTp\in\TTP$ to not exceed its processing capability $x_{\TTp}$. 
The constraints $z_{\TTs,\TTp}\le b_{\TTs,\TTp}$ $\forall (\TTs,\TTp)\in\TTE_{\TTS,\TTP}\setminus\widetilde{\TTE}_{\TTS,\TTP}$ and $z_{\TTp,\TTc}\le b_{\TTp,\TTc}$ $\forall (\TTp,\TTc)\in\TTE_{\TTP,\TTC}\setminus\widetilde{\TTE}_{\TTP,\TTC}$ require 
the amount of goods flowing through
each non-susceptible edge in the supply chain network
to be within its maximum capacity.
The constraints $z_{\TTs,\TTp}\le b_{\TTs,\TTp}(1-\omega_{\TTs,\TTp})$ $\forall (\TTs,\TTp)\in\widetilde{\TTE}_{\TTS,\TTP}$ and $0\le z_{\TTp,\TTc}\le b_{\TTp,\TTc}(1-\omega_{\TTp,\TTc})$ $\forall (\TTp,\TTc)\in\widetilde{\TTE}_{\TTP,\TTC}$ require
the amount of goods flowing through
each susceptible edge in the supply chain network
to be within its maximum capacity when it does not fail,
and to be zero in the event that it fails.
The overall objective of the decision maker is to minimize the total investment in the first stage, i.e., $\sum_{\TTp\in\TTP}c_{1,\TTp}x_{\TTp}$, and the total transportation and processing costs in the second stage. 
With appropriate vectorizations of the decision variables and the parameters, the second-stage decision problem can be represented as:
\begin{align}
\begin{split}
\minimize_{\BIz_{\TTS,\TTP},\,\BIz_{\TTP,\TTC}}\quad 
& \langle\Bc_{2,\TTS,\TTP},\BIz_{\TTS,\TTP}\rangle + \langle\Bc_{2,\TTP,\TTC},\BIz_{\TTP,\TTC}\rangle\\
\mathrm{subject~to}\quad 
& \BA_{\mathrm{ba},\TTS,\TTP}\BIz_{\TTS,\TTP} - \BA_{\mathrm{ba},\TTP,\TTC}\BIz_{\TTP,\TTC}=\veczero_{|\TTP|},\\
& \BA_{\mathrm{su}}\BIz_{\TTS,\TTP}\le\Bb_{\TTS}, \quad 
\BA_{\mathrm{de}}\BIz_{\TTP,\TTC}\ge\Bomega_{\TTC}, \quad
\BA_{\mathrm{pr}}\BIz_{\TTS,\TTP}\le\BIx_{\TTP},\\
& \BIz_{\TTS,\TTP}\le\Bb_{\TTS,\TTP}-\BF_{\TTS,\TTP}\Bomega_{\TTS,\TTP}, \quad
\BIz_{\TTP,\TTC}\le\Bb_{\TTP,\TTC}-\BF_{\TTP,\TTC}\Bomega_{\TTP,\TTC},\\
& \BIz_{\TTS,\TTP}\in\R^{|\TTE_{\TTS,\TTP}|}_{+}, \quad \BIz_{\TTP,\TTC}\in\R^{|\TTE_{\TTP,\TTC}|}_{+},
\end{split}%
\label{eqn:supplychain-stage2obj}%
\end{align}%
for suitably specified 
$\Bc_{2,\TTS,\TTP}\in\nobreak\R_{+}^{|\TTE_{\TTS,\TTP}|}$, 
$\Bc_{2,\TTP,\TTC}\in\nobreak\R_{+}^{|\TTE_{\TTP,\TTC}|}$, 
$\BA_{\mathrm{ba},\TTS,\TTP}\in\nobreak\R_{+}^{|\TTP|\times|\TTE_{\TTS,\TTP}|}$, 
$\BA_{\mathrm{ba},\TTP,\TTC}\in\nobreak\R_{+}^{|\TTP|\times|\TTE_{\TTP,\TTC}|}$, 
$\BA_{\mathrm{su}}\in\nobreak\R_{+}^{|\TTS|\times|\TTE_{\TTS,\TTP}|}$, 
$\Bb_{\TTS}\in\nobreak\R_{+}^{|\TTS|}$, 
$\BA_{\mathrm{de}}\in\nobreak\R_{+}^{|\TTC|\times|\TTE_{\TTS,\TTP}|}$, 
$\Bomega_{\TTC}\in\nobreak\R_{+}^{|\TTC|}$, 
$\BA_{\mathrm{pr}}\in\nobreak\R_{+}^{|\TTP|\times|\TTE_{\TTS,\TTP}|}$, 
$\BIx_{\TTP}\in\nobreak\R^{|\TTP|}_+$, 
$\Bb_{\TTS,\TTP}\in\nobreak\R_{+}^{|\TTE_{\TTS,\TTP}|}$, 
$\Bb_{\TTP,\TTC}\in\nobreak\R_{+}^{|\TTE_{\TTP,\TTC}|}$, 
$\BF_{\TTS,\TTP}\in\nobreak\R_{+}^{|\TTE_{\TTS,\TTP}|\times|\widetilde{\TTE}_{\TTS,\TTP}|}$, 
$\BF_{\TTP,\TTC}\in\R_{+}^{|\TTE_{\TTP,\TTC}|\times|\widetilde{\TTE}_{\TTP,\TTC}|}$, 
$\Bomega_{\TTS,\TTP}\in\{0,1\}^{|\widetilde{\TTE}_{\TTS,\TTP}|}$, 
$\Bomega_{\TTP,\TTC}\in\{0,1\}^{|\widetilde{\TTE}_{\TTP,\TTC}|}$, 
where $\Bc_{2,\TTS,\TTP}$, $\Bc_{2,\TTP,\TTC}$ represent the transportation and processing costs of the edges, 
$\Bb_{\TTS}$ represents the supplies, 
$\Bomega_{\TTC}$ represents the demands, 
$\BIx_{\TTP}$ represents the processing capabilities, 
$\Bb_{\TTS,\TTP}$, $\Bb_{\TTP,\TTC}$ represent the maximum transportation capacities of the edges, 
and $\Bomega_{\TTS,\TTP}$, $\Bomega_{\TTP,\TTC}$ represent the failure of the edges. 
To guarantee that (\ref{eqn:supplychain-stage2obj}) is feasible for all feasible first-stage decisions, we introduce the auxiliary variables $\BIx_{0,\TTS,\TTP}\in\R_{+}^{|\TTE_{\TTS,\TTP}|}$ and $\BIx_{0,\TTP,\TTC}\in\R_{+}^{|\TTE_{\TTP,\TTC}|}$ to the first-stage decision variables, and define 
\begin{align*}
\CS_1:=\left\{\!(\BIx_{\TTP}^\TRANSP,\BIx_{0,\TTS,\TTP}^\TRANSP,\BIx_{0,\TTP,\TTC}^\TRANSP)^{\TRANSP} \in \R_{+}^{|\TTP|+|\TTE|}:\!
\begin{tabular}{l}
    $\BIx_{\TTP}\le\overline{\BIx}_{\TTP},\;\BA_{\mathrm{ba},\TTS,\TTP}\BIx_{0,\TTS,\TTP}-\BA_{\mathrm{ba},\TTP,\TTC}\BIx_{0,\TTP,\TTC}=\veczero_{|\TTP|},$ \\
    $\BA_{\mathrm{su}}\BIx_{0,\TTS,\TTP}\le \Bb_{\TTS},\; 
    \BA_{\mathrm{de}}\BIx_{0,\TTP,\TTC}\ge \overline{\Bomega}_{\TTC},\;
    \BA_{\mathrm{pr}}\BIx_{0,\TTS,\TTP}\le \BIx_{\TTP},$ \\
    $\BIx_{0,\TTS,\TTP}\le\Bb_{\TTS,\TTP}-\BF_{\TTS,\TTP}\vecone_{|\widetilde{\TTE}_{\TTS,\TTP}|},\; 
    \BIx_{0,\TTP,\TTC}\le\Bb_{\TTP,\TTC}-\BF_{\TTP,\TTC}\vecone_{|\widetilde{\TTE}_{\TTP,\TTC}|}$
\end{tabular}\!\!\!
\right\},
\end{align*}
where $\overline{\BIx}_{\TTP}\in\R_{+}^{|\TTP|}$ and $\overline{\Bomega}_{\TTC}\in\R_{+}^{|\TTC|}$ are the vectorized version of $(\overline{x}_{\TTp})_{\TTp\in\TTP}$ and $(\overline{\omega}_{\TTc})_{\TTc\in\TTC}$.
Let us assume in addition that $\CS_1$ is non-empty. 

Formulating this problem into our two-stage DRO model in Setting~\ref{sett:dro}, we have $d:=|\TTC|+|\widetilde{\TTE}|$,  
$k_1:=|\TTP|+|\TTE|$, $k_2:=|\TTE|$, 
$(\Omega_i)_{i=1:|\TTC|}:=\big([\underline{\omega}_{\TTc},\overline{\omega}_{\TTc}]\big)_{\TTc\in\TTC}$, and $(\Omega_{|\TTC|+i})_{i=1:|\widetilde{\TTE}|}:=\{0,1\}^{|\widetilde{\TTE}|}$. 
Subsequently, let $C_2\big((\BIx_{\TTP}^\TRANSP,\BIx_{0,\TTS,\TTP}^\TRANSP,\BIx_{0,\TTP,\TTC}^\TRANSP)^\TRANSP,\allowbreak(\Bomega_{\TTC}^\TRANSP,\Bomega_{\TTS,\TTP}^\TRANSP,\Bomega_{\TTP,\TTC}^\TRANSP)^\TRANSP\big)$ denote the optimal value of (\ref{eqn:supplychain-stage2obj}), which corresponds to the optimal second-stage cost function. 
One can check that $C_2\big((\BIx_{\TTP}^\TRANSP,\BIx_{0,\TTS,\TTP}^\TRANSP,\BIx_{0,\TTP,\TTC}^\TRANSP)^\TRANSP,\allowbreak(\Bomega_{\TTC}^\TRANSP,\Bomega_{\TTS,\TTP}^\TRANSP,\Bomega_{\TTP,\TTC}^\TRANSP)^\TRANSP\big)$ has the required form in \ref{setts:dro-stage2}. 
Moreover, it holds that
\begin{align*}
    &C_2\big((\BIx_{\TTP}^\TRANSP,\BIx_{0,\TTS,\TTP}^\TRANSP,\BIx_{0,\TTP,\TTC}^\TRANSP)^\TRANSP,(\Bomega_{\TTC}^\TRANSP,\Bomega_{\TTS,\TTP}^\TRANSP,\Bomega_{\TTP,\TTC}^\TRANSP)^\TRANSP\big)\le \langle\Bc_{2,\TTS,\TTP},\Bb_{\TTS,\TTP}\rangle+\langle\Bc_{2,\TTP,\TTC},\Bb_{\TTP,\TTC}\rangle\\
    & \hspace{100pt} \forall (\BIx_{\TTP}^\TRANSP,\BIx_{0,\TTS,\TTP}^\TRANSP,\BIx_{0,\TTP,\TTC}^\TRANSP)^\TRANSP \in\CS_1, \; \forall (\Bomega_{\TTC}^\TRANSP,\Bomega_{\TTS,\TTP}^\TRANSP,\Bomega_{\TTP,\TTC}^\TRANSP)^\TRANSP\in\BOmega,
\end{align*}
and thus \ref{setts:dro-stage2} is satisfied with respect to some
$\alpha\le\langle\Bc_{2,\TTS,\TTP},\Bb_{\TTS,\TTP}\rangle+\langle\Bc_{2,\TTP,\TTC},\Bb_{\TTP,\TTC}\rangle<\infty$. 
\end{example}

\section{Primal and dual formulations}
\label{sec:primaldual}

In this section, we develop an iterative approximation scheme for \eqref{eqn:dro} with the goal of constructing an approximate optimizer of \eqref{eqn:dro} as well as an approximate optimizer of the inner maximization problem in \eqref{eqn:dro}.
Specifically, we reformulate \eqref{eqn:dro} into a minimization problem denoted by \eqref{eqn:dro-primal} and a maximization problem denoted by \eqref{eqn:dro-dual},
and we examine their optimality conditions.
Our iterative approximation scheme is inspired by the optimality conditions for \eqref{eqn:dro-dual}.
The results we develop in this section will serve as the theoretical basis of the concrete numerical algorithm in Section~\ref{sec:algorithm}.

Under Setting~\ref{sett:dro-alt}, let us first introduce the following notations which are used throughout the rest of the paper.
For $i=1,\ldots,d$,
let $\Bp_i\in\R^{k^*_2}$ denote the $i$-th column of the matrix $\BP\in\R^{k^*_2\times d}$
and let $\OT_{\mu_i}(\nu)\in\R$ denote the minimum optimal transportation cost between $\mu_i$ and any $\nu\in\CP_1(\R)$ with respect to the cost function $\Omega_i\times \R \ni (\omega,z)\mapsto -z\omega\in\R$, i.e.,
\begin{align}
    \mathrm{OT}_{\mu_i}(\nu):=\inf_{\gamma\in\Gamma(\mu_i,\nu)}\bigg\{\int_{\Omega_i\times\R}-z\omega\DIFFM{\gamma}{\DIFF \omega,\DIFF z}\bigg\} \qquad \forall \nu\in\CP_1(\R),\; \forall 1\le i\le d,
    \label{eqn:notation-rho-def}
\end{align}
where $\Gamma(\mu_i,\nu):=\big\{\gamma\in\CP(\Omega_i\times\R):$ the marginals of $\gamma$ on $\Omega_i$ and $\R$ are $\mu_i$ and $\nu$ respectively$\big\}$ denote the set of couplings of $\mu_i$ and $\nu$.

We now present the primal and dual formulations of \eqref{eqn:dro}.
\begin{theorem}[Primal and dual formulations of \eqref{eqn:dro}]
\label{thm:primaldual}
Under Setting~\ref{sett:dro-alt}, recall that $C_{\DRO}$ denotes the optimal value of \eqref{eqn:dro}.
Let \eqref{eqn:dro-primal} denote the following minimization problem:
\begin{align}
    \begin{split}
        \minimize_{\BIx,\,f_{1},\ldots,f_d} \quad & \langle\Bc_1,\BIx\rangle+\sum_{i=1}^d\int_{\Omega_i}\sup_{z\in\R}\big\{z\omega-f_i(z)\big\}\DIFFM{\mu_i}{\DIFF\omega}  \\
        \mathrm{subject~to} \quad & \BL_{\ineq}\BIx\le \Bq_{\ineq},\qquad \BL_{\eq}\BIx=\Bq_{\eq}, \\
        & \langle\BQ\BIx+\Bb,\Blambda\rangle + \sum_{i=1}^d f_{i}\big(\langle\Bp_i,\Blambda\rangle\big) \le 0 \quad \forall \Blambda\in\CS_2^*,  \\
        & \BIx\in\R^{k_1}_{+}, \qquad f_{i}\in\CC_{\mathrm{cvx}}(\R) \quad \forall 1\le i\le d,
    \end{split}
    \tag{\textsf{Primal}}
    \label{eqn:dro-primal}
\end{align}
and let \eqref{eqn:dro-dual} denote the following maximization problem:
\begin{align}
    \begin{split}
        \maximize_{\tau,\,\Bxi_{\ineq},\,\Bxi_{\eq}} \quad & \int_{\CS^*_2}\langle\Bb,\Blambda\rangle\DIFFM{\tau}{\DIFF\Blambda} 
        - \Bigg(\sum_{i=1}^d\OT_{\mu_i}\big(\langle\Bp_i,\cdot\,\rangle\sharp\tau\big)\Bigg)
        + \langle\Bq_{\ineq},\Bxi_{\ineq}\rangle+\langle\Bq_{\eq},\Bxi_{\eq}\rangle  \\
        \mathrm{subject~to} \quad & \Bc_1+\int_{\CS_2^*}\BQ^\TRANSP\Blambda \DIFFM{\tau}{\DIFF\Blambda} \ge \BL_{\ineq}^\TRANSP\Bxi_{\ineq} + \BL_{\eq}^\TRANSP\Bxi_{\eq}, \\
        & \tau\in\CP_1(\CS^*_2), \qquad \Bxi_{\ineq} \in \R^{n_{\ineq}}_{-},\qquad \Bxi_{\eq}\in\R^{n_{\eq}}.
    \end{split}
    \tag{\textsf{Dual}}
    \label{eqn:dro-dual}
\end{align}
Then, the following statements hold.\widowpenalty-1000
\begin{enumerate}[label=(\roman*),beginpenalty=10000]
    \item\label{thms:primaldual-strongduality}
    $C_{\mathrm{DRO}}<\infty$ and the optimal values of
    \eqref{eqn:dro-primal} and \eqref{eqn:dro-dual} are both equal to~$C_{\mathrm{DRO}}$. 
    
    \item\label{thms:primaldual-boundedness}
    $C_{\mathrm{DRO}}>-\infty$ if and only if there exist $\Blambda\in\CS^*_2$, $\Bxi_{\ineq}\in\R^{n_{\ineq}}_{-}$, $\Bxi_{\eq}\in\R^{n_{\eq}}$ which satisfy $\Bc_1+\BQ^\TRANSP\Blambda\ge\BL_{\ineq}^\TRANSP\Bxi_{\ineq}+\BL_{\eq}^\TRANSP\Bxi_{\eq}$.
    In particular, $(\delta_{\Blambda},\Bxi_{\ineq},\Bxi_{\eq})$ is feasible for \eqref{eqn:dro-dual} in this case.
\end{enumerate}
Now, suppose in addition that $C_{\mathrm{DRO}}>-\infty$. 
Then, the following statements hold.\widowpenalty-1000
\begin{enumerate}[label=(\roman*),beginpenalty=10000]
    \setcounter{enumi}{2}

    \item\label{thms:primaldual-primal-optimizer-existence}
    \eqref{eqn:dro-primal} has an optimizer.

    \item\label{thms:primaldual-primal-optimizer-properties}
    Let $(\bar{\BIx},\bar{f}_1,\ldots,\bar{f}_d)$ be an arbitrary optimizer of \eqref{eqn:dro-primal}
    and let $\bar{f}_i^*(\omega):=\sup_{z\in\R}\big\{z\omega-\bar{f}_i(z)\big\}$ $\forall\omega\in\Omega_i$ denote the convex conjugate of $\bar{f}_i$ for $i=1,\ldots,d$. Then, $\bar{\BIx}$ is an optimizer of \eqref{eqn:dro}.
    Moreover,
    for any $\mu\in\Gamma(\mu_1,\ldots,\mu_d)$,
    $\mu$ is optimal for the inner maximization problem in \eqref{eqn:dro} with respect to $\bar{\BIx}$, i.e., $\mu\in\argmax_{\mu'\in\Gamma(\mu_1,\ldots,\mu_d)}\big\{\!\int_{\BOmega}C_2(\bar{\BIx},\Bomega)\DIFFM{\mu'}{\DIFF\Bomega}\big\}$,
    if and only if 
    $\sum_{i=1}^{d}\bar{f}_i^*(\omega_i)=C_2(\bar{\BIx},\Bomega)$ for 
    $\mu$-almost every $\Bomega=(\omega_1,\ldots,\omega_d)^\TRANSP\in\BOmega$.

    \item\label{thms:primaldual-primal-approxoptimizer}
    For any $\varepsilon>0$ and any \mbox{$\varepsilon$-optimizer} $(\BIx_{\varepsilon},f_{1,\varepsilon},\ldots,f_{d,\varepsilon})$ of \eqref{eqn:dro-primal}, $\BIx_{\varepsilon}$~is an $\varepsilon$-optimizer of~\eqref{eqn:dro}. 

    \item\label{thms:primaldual-dual-optimizer-existence}
    \eqref{eqn:dro-dual} has an optimizer $(\bar{\tau},\bar{\Bxi}_{\ineq},\bar{\Bxi}_{\eq})$ which satisfies $\support(\bar{\tau})\subseteq\CV_{\extremept}$.

    \item\label{thms:primaldual-dual-optimizer-properties}
    Let $(\bar{\tau},\bar{\Bxi}_{\ineq},\bar{\Bxi}_{\eq})$ be an arbitrary optimizer of \eqref{eqn:dro-dual}.
    Then, there exists $\bar{\pi}\in\CP(\BOmega\times\CS^*_2)$ 
    with marginal $\mu_i$ on $\Omega_i$ for $i=1,\ldots,d$ and marginal $\bar{\tau}$ on $\CS^*_2$, where for any optimizer $\bar{\BIx}$ of \eqref{eqn:dro}, 
    it holds that $\Blambda\!\in\!\argmax_{\Blambda'}\!\big\{\hspace{-1pt}\langle\BQ\bar{\BIx}+\BP\Bomega+\Bb,\Blambda'\rangle\hspace{-1pt}\big\}$ for $\bar{\pi}$-almost every $(\Bomega,\Blambda)\hspace{-1pt}\in\hspace{-1pt}\BOmega\times\nobreak\CS^*_2$, and
    the marginal $\bar{\mu}$ of $\bar{\pi}$ on $\BOmega$ satisfies $\bar{\mu}\in\argmax_{\mu\in\Gamma(\mu_1,\ldots,\mu_d)}\!\big\{\!\int_{\BOmega}C_2(\bar{\BIx},\Bomega)\DIFFM{\mu}{\DIFF\Bomega}\big\}$.

    \item\label{thms:primaldual-saddle}
    Let $(\BIx,f_1,\ldots,f_d)$ be feasible for \eqref{eqn:dro-primal} and let $(\tau,\Bxi_{\ineq},\Bxi_{\eq})$ be feasible for \eqref{eqn:dro-dual}.
    Then, $(\BIx,f_1,\ldots,f_d)$ and $(\tau,\Bxi_{\ineq},\Bxi_{\eq})$ are optimizers of \eqref{eqn:dro-primal} and \eqref{eqn:dro-dual} 
    if and only if the following conditions are satisfied:
    \begin{enumerate}[label=(\alph*),leftmargin=17pt]
        \item $\langle\BL_{\ineq}\BIx-\Bq_{\ineq}, \Bxi_{\ineq}\rangle=0$,
        
        \item $\big\langle\Bc_1+{\textstyle\int_{\CS^*_2}\BQ^\TRANSP\Blambda}\DIFFM{\tau}{\DIFF\Blambda}-\BL_{\ineq}^\TRANSP\Bxi_{\ineq}-\BL_{\eq}^{\TRANSP}\Bxi_{\eq},\BIx\big\rangle=0$,
        
        \item $\langle\BQ\BIx+\Bb, \Blambda\rangle+\sum_{i=1}^d f_i\big(\langle\Bp_i,\Blambda\rangle\big)= 0$ for $\tau$-almost every $\Blambda$,
        
        \item ${-\!\int_{\Omega_i}}\!\sup_{z\in\R}\!\big\{z\omega-f_i(z)\big\}\DIFFM{\mu_i}{\DIFF\omega}-\int_{\R}f_i\DIFFX{\langle\Bp_i,\cdot\,\rangle\sharp\tau}= \OT_{\mu_i}\!\big(\langle\Bp_i,\cdot\,\rangle\sharp\tau\big)$ for $i=1,\ldots,d$.
    \end{enumerate}

\end{enumerate}
\end{theorem}

In the above theorem, \eqref{eqn:dro-primal} and \eqref{eqn:dro-dual} are minimization and maximization problems that are easier to solve than the min-max problem \eqref{eqn:dro}, and
statement~\ref{thms:primaldual-strongduality} shows that we are able to bound the optimal value of \eqref{eqn:dro} from above and below once we obtain feasible solutions of \eqref{eqn:dro-primal} and \eqref{eqn:dro-dual}.
Statement~\ref{thms:primaldual-boundedness} implies that 
in order to guarantee that \eqref{eqn:dro} is bounded from below, it suffices to find one point $\Blambda\in\CS^*_2$ such that 
$\Bc_1+\BQ^\TRANSP\Blambda\ge\BL_{\ineq}^\TRANSP\Bxi_{\ineq}+\BL_{\eq}^\TRANSP\Bxi_{\eq}$ holds for some $\Bxi_{\ineq}\in\R^{n_{\ineq}}_{-}$, $\Bxi_{\eq}\in\R^{n_{\eq}}$.
In practice, the existence of such $\Blambda\in\CS^*_2$ can be checked by computing a feasible solution of an LP problem.
This technique will be used in Section~\ref{sec:algorithm} to initialize our numerical algorithm.
Statements~\ref{thms:primaldual-primal-optimizer-properties} and \ref{thms:primaldual-primal-approxoptimizer}
demonstrate that an exact or approximate optimizer of \eqref{eqn:dro-primal} directly yields an exact or approximate optimizer of \eqref{eqn:dro}.
Moreover, statements~\ref{thms:primaldual-primal-optimizer-properties} and \ref{thms:primaldual-dual-optimizer-properties}
provide interpretations of the optimizers of \eqref{eqn:dro-primal} and \eqref{eqn:dro-dual}.
Specifically, statement~\ref{thms:primaldual-primal-optimizer-properties} shows that the functions $\bar{f}_1,\ldots,\bar{f}_d\in\CC_{\mathrm{cvx}}(\R)$ in an optimizer $(\bar{\BIx},\bar{f}_1,\ldots,\bar{f}_d)$ of \eqref{eqn:dro-primal} characterize the worst-case probability measures with respect to the first-stage decision $\bar{\BIx}$,
and statement~\ref{thms:primaldual-dual-optimizer-properties} shows that the probability measure $\bar{\tau}\in\CP_1(\CS^*_2)$ in an optimizer $(\bar{\tau},\bar{\Bxi}_{\ineq},\bar{\Bxi}_{\eq})$ of \eqref{eqn:dro-dual} can be interpreted as a ``best response'' to a worst-case probability measure $\bar{\mu}\in\Gamma(\mu_1,\ldots,\mu_d)$ with respect to optimal first-stage decisions.
Lastly, statement~\ref{thms:primaldual-saddle} provides the saddle point conditions for the optimizers of \eqref{eqn:dro-primal} and \eqref{eqn:dro-dual}.
In particular, combining statement~\ref{thms:primaldual-saddle} with statements~\ref{thms:primaldual-primal-optimizer-existence} and \ref{thms:primaldual-dual-optimizer-existence}
yields sufficient and necessary optimality conditions for \eqref{eqn:dro-primal} and \eqref{eqn:dro-dual}.

Despite that \eqref{eqn:dro-dual} does not admit a unique optimizer in general, we are able to derive the following property about the optimizers of \eqref{eqn:dro-dual} under the additional assumption that $\mu_1,\ldots,\mu_d$ are absolutely continuous with respect to the Lebesgue measure.

\begin{corollary}[Property of dual optimizers]
    \label{cor:dual-optimizers}
    Let Setting~\ref{sett:dro-alt} hold, assume that $C_{\DRO}>-\infty$, and assume that $\mu_1\in\CP(\Omega_1),\ldots,\mu_d\in\CP(\Omega_d)$ are absolutely continuous with respect to the Lebesgue measure.
    Let $(\bar{\BIx},\bar{f}_1,\ldots,\bar{f}_d)$ be an arbitrary optimizer of \eqref{eqn:dro-primal} and let $(\bar{\tau},\bar{\Bxi}_{\ineq},\bar{\Bxi}_{\eq})$ be an arbitrary optimizer of \eqref{eqn:dro-dual}.
    For $i=1,\ldots,d$, let $\bar{f}_i^*(\omega):=\sup_{z\in\R}\big\{z\omega-\bar{f}_i(z)\big\}$ $\forall\omega\in\Omega_i$ denote the convex conjugate of $\bar{f}_i$.
    Then, it holds for $\mu_i$-almost every $\omega\in\Omega_i$ that
    $\partial \bar{f}_i^*(\omega)$ is a singleton.
    Let $\nabla \bar{f}_i^*:\Omega_i\to\R$ be a Borel measurable function such that $\nabla \bar{f}_i^*(\omega)\in\partial \bar{f}_i^*(\omega)$ for all $\omega\in\Omega_i$,
    i.e., $\nabla \bar{f}_i^*$ is a selector of $\partial \bar{f}_i^*$, which is uniquely defined $\mu_i$-almost everywhere.
    Subsequently, it holds that $\langle\Bp_i,\cdot\,\rangle\sharp\bar{\tau}=\nabla \bar{f}_i^*\sharp\mu_i$ for $i=1,\ldots,d$.
    In particular, for any two optimizers $(\bar{\tau},\bar{\Bxi}_{\ineq},\bar{\Bxi}_{\eq}),(\bar{\tau}',\bar{\Bxi}_{\ineq}',\bar{\Bxi}_{\eq}')$ of \eqref{eqn:dro-dual}, it holds that
    $\langle\Bp_i,\cdot\,\rangle\sharp\bar{\tau} = \langle\Bp_i,\cdot\,\rangle\sharp\bar{\tau}'$ for $i=1,\ldots,d$.
\end{corollary}

\begin{remark}\label{rmk:dual-marginal-uniqueness}%
    Corollary~\ref{cor:dual-optimizers} is a generalization of the so-called ``marginal uniqueness'' property in \citep[Corollary~3]{chen2018distributionally}.
    As shown by Theorem~\ref{thm:primaldual}\ref{thms:primaldual-dual-optimizer-properties}, the probability measure $\bar{\tau}\in\CP_1(\CS^*_2)$ in an optimizer $(\bar{\tau},\bar{\Bxi}_{\ineq},\bar{\Bxi}_{\eq})$ of \eqref{eqn:dro-dual}
    can be interpreted as the ``best response'' to a worst-case probability measure $\bar{\mu}\in\Gamma(\mu_1,\ldots,\mu_d)$ with respect to optimal first-stage decisions.
    Hence, for $i=1,\ldots,d$,
    $\langle\Bp_i,\cdot\,\rangle\sharp\bar{\tau}$ can be seen as the ``best response'' to the uncertain outcome $\omega_i$ with law $\mu_i$.
    Corollary~\ref{cor:dual-optimizers} states that $\BP\sharp\bar{\tau}$ has unique marginals for any optimizer $(\bar{\tau},\bar{\Bxi}_{\ineq},\bar{\Bxi}_{\eq})$ of \eqref{eqn:dro-dual}.
    However, this interpretation is more involved than the ``marginal uniqueness'' property in \citep[Corollary~3]{chen2018distributionally} due to the presence of the first-stage optimization over $\BIx\in\CS_1$,
    and due to the assumption that $\BP$ is not restricted to the identity matrix in our model.
\end{remark}

In the following, we shift our attention to developing a numerical method for constructing feasible and approximately optimal solutions of \eqref{eqn:dro-primal} and \eqref{eqn:dro-dual}.
One may observe from Theorem~\ref{thm:primaldual}\ref{thms:primaldual-dual-optimizer-existence} that 
the decision variable $\tau\in\CP_1(\CS^*_2)$ in
\eqref{eqn:dro-dual} can be restricted to those satisfying $\support(\tau)\subseteq\CV_{\extremept}$ without affecting the optimal value of \eqref{eqn:dro-dual}.
Consequently, the parametrization $\tau\leftarrow\sum_{\Blambda\in\CV_{\extremept}}\nu_{\Blambda}\delta_{\Blambda}$ via the probabilities $(\nu_{\Blambda})_{\Blambda\in\CV_{\extremept}}\subset\R_+$ where $\sum_{\Blambda\in\CV_{\extremept}}\nu_{\Blambda}=1$ 
results in a maximization problem in finite dimensions.
However, this approach is impractical because $|\CV_{\extremept}|$ can be prohibitively large and there is no computationally tractable way of enumerating the elements of $\CV_{\extremept}$, especially when $k^*_2$ is large. 
Despite this, 
it is plausible that there exists an exact or approximate optimizer $(\bar{\tau},\bar{\Bxi}_{\ineq},\bar{\Bxi}_{\eq})$ of \eqref{eqn:dro-dual} 
where $\support(\bar{\tau})$ is a finite set with manageable size.
This is indeed the case as we show later in the numerical experiments in Section~\ref{sec:experiments}.
Therefore, let us derive the sub-optimality bound for \eqref{eqn:dro-dual} in the following theorem, which is inspired by the optimality conditions in Theorem~\ref{thm:primaldual}\ref{thms:primaldual-saddle}.

\begin{theorem}[Sub-optimality bound for \eqref{eqn:dro-dual}]\label{thm:dual-suboptimality}
    Under Setting~\ref{sett:dro-alt}, assume that \mbox{$C_{\DRO}>-\infty$}.
    Let $(\tau,\Bxi_{\ineq},\Bxi_{\eq})$ be feasible for \eqref{eqn:dro-dual}
    where $\support(\tau)$ is a finite set,
    and let $\CV$ be an arbitrary finite set satisfying $\support(\tau)\subseteq\CV\subset\CS^*_2$.
    Consider the following constraints for 
    $\BIx\in\CS_1$, $f_1,\ldots,f_d\in\CC_{\mathrm{cvx}}(\R)$ with respect to a tolerance parameter $\varepsilon_{\OT}> 0$:
    \begin{subnumcases}{}
        \textstyle\langle\BQ\BIx+\Bb, \Blambda\rangle+\sum_{i=1}^d f_i\big(\langle\Bp_i,\Blambda\rangle\big)= 0 
        \hspace{161.8pt} \forall \Blambda\in\support(\tau),\label{eqn:dual-suboptimality-primal-feasibility-potential-eq}\\
        \textstyle\langle\BQ\BIx+\Bb, \Blambda\rangle+\sum_{i=1}^d f_i\big(\langle\Bp_i,\Blambda\rangle\big)\le 0 
        \hspace{143.8pt} \forall \Blambda\in\CV\setminus\support(\tau),\label{eqn:dual-suboptimality-primal-feasibility-potential-le}\\
        \langle\BL_{\ineq}\BIx-\Bq_{\ineq}, \Bxi_{\ineq}\rangle=0,\label{eqn:dual-suboptimality-primal-complementarity} \\
        \big\langle\Bc_1+{\textstyle\int_{\CS^*_2}\BQ^\TRANSP\Blambda}\DIFFM{\tau}{\DIFF\Blambda}-\BL_{\ineq}^\TRANSP\Bxi_{\ineq}-\BL_{\eq}^{\TRANSP}\Bxi_{\eq},\BIx\big\rangle=0,\label{eqn:dual-suboptimality-dual-complementarity}\\
        \mathrm{for~}i=1,\ldots,d: \nonumber\\
        \textstyle\quad {-\int_{\Omega_i}\sup_{z\in\R}\big\{z\omega-f_i(z)\big\}\DIFFM{\mu_i}{\DIFF\omega}}-\int_{\R}f_i\DIFFX{\langle\Bp_i,\cdot\,\rangle\sharp\tau}\ge \OT_{\mu_i}\big(\langle\Bp_i,\cdot\,\rangle\sharp\tau\big)-\varepsilon_{\OT},\label{eqn:dual-suboptimality-potential-optimizer}\\
        \quad \partial f_i(z)\subseteq \Omega_i
        \hspace{295.5pt} \forall z\in\R. \label{eqn:dual-suboptimality-potential-subgradient}
    \end{subnumcases}
    Moreover, for any $\BIx\in\CS_1$ and any $f_1,\ldots,f_d\in\CC_{\mathrm{cvx}}(\R)$,
    let $\mathrm{vio}(\BIx,f_1,\ldots,f_d)$ be defined as follows:
    \begin{align}
        \mathrm{vio}(\BIx,f_1,\ldots,f_d):=\sup_{\Blambda\in\CS^*_2}\Bigg\{\langle\BQ\BIx+\Bb, \Blambda\rangle+\sum_{i=1}^d f_i\big(\langle\Bp_i,\Blambda\rangle\big)\Bigg\}.
        \label{eqn:dual-suboptimality-vio-def}
    \end{align}
    Then, the following statements hold.
    \begin{enumerate}[label=(\roman*),beginpenalty=10000]
        \item\label{thms:dual-suboptimality-vio-positivity}
        For any $\BIx\in\CS_1$, $f_1,\ldots,f_d\in\CC_{\mathrm{cvx}}(\R)$ that satisfy (\ref{eqn:dual-suboptimality-primal-feasibility-potential-eq})--(\ref{eqn:dual-suboptimality-potential-subgradient}), it holds
        that $0\le \mathrm{vio}(\BIx,f_1,\ldots,f_d)<\infty$ and a maximizer of (\ref{eqn:dual-suboptimality-vio-def}) is attained in $\CV_{\extremept}$.

        \item\label{thms:dual-suboptimality-suboptimality}
        Assume that $\BIx\in\CS_1$, $f_1,\ldots,f_d\in\CC_{\mathrm{cvx}}(\R)$ satisfy (\ref{eqn:dual-suboptimality-primal-feasibility-potential-eq})--(\ref{eqn:dual-suboptimality-potential-subgradient}) and define
        $\varepsilon_{\mathrm{sub}}:=\mathrm{vio}(\BIx,f_1,\ldots,f_d)+d\varepsilon_{\OT}$.
        Then, $(\tau,\Bxi_{\ineq},\Bxi_{\eq})$ is an $\varepsilon_{\mathrm{sub}}$-optimizer of \eqref{eqn:dro-dual}.
        Moreover, 
        if $\mathrm{vio}(\BIx,f_1,\ldots,f_d)>0$, then
        the set 
        $\big\{\Blambda\in\CS^*_2:\langle\BQ\BIx+\Bb, \Blambda\rangle+\sum_{i=1}^d f_i\big(\langle\Bp_i,\Blambda\rangle\big)>\nobreak0\big\}$ 
        is disjoint from~$\CV$.
        
    \end{enumerate}
\end{theorem}

Theorem~\ref{thm:dual-suboptimality}\ref{thms:dual-suboptimality-suboptimality} allows one to quantify the sub-optimality of any $(\tau,\Bxi_{\ineq},\Bxi_{\eq})$ feasible for \eqref{eqn:dro-dual}
by first finding $\BIx\in\CS_1$, $f_1,\ldots,f_d\in\CC_{\mathrm{cvx}}(\R)$ that satisfy (\ref{eqn:dual-suboptimality-primal-feasibility-potential-eq})--(\ref{eqn:dual-suboptimality-potential-subgradient}) with respect to some $\varepsilon_{\OT}>0$
and then computing $\mathrm{vio}(\BIx,f_1,\ldots,f_d)+d\varepsilon_{\OT}$.
In particular, observe from Theorem~\ref{thm:primaldual}\ref{thms:primaldual-saddle} that if $(\tau,\Bxi_{\ineq},\Bxi_{\eq})$ is an optimizer of \eqref{eqn:dro-dual}, then for every optimizer $(\BIx,f_1,\ldots,f_d)$ of \eqref{eqn:dro-primal}, $\BIx,f_1,\ldots,f_d$ satisfy (\ref{eqn:dual-suboptimality-primal-feasibility-potential-eq})--(\ref{eqn:dual-suboptimality-potential-subgradient}) with respect to any $\varepsilon_{\OT}>0$.
The settings of Theorem~\ref{thm:dual-suboptimality} also yield 
approximate optimizers of \eqref{eqn:dro-primal}, \eqref{eqn:dro}, and the inner maximization problem in \eqref{eqn:dro}, as well as their sub-optimality bounds.
This is presented in the following corollary.

\begin{corollary}[Sub-optimality bounds for \eqref{eqn:dro-primal}, \eqref{eqn:dro}, and the inner maximization problem]\label{cor:dual-suboptimality-primal-coupling}
    Under the settings of Theorem~\ref{thm:dual-suboptimality},
    assume that $\BIx\in\CS_1$, $f_1,\ldots,f_d\in\CC_{\mathrm{cvx}}(\R)$ satisfy (\ref{eqn:dual-suboptimality-primal-feasibility-potential-eq})--(\ref{eqn:dual-suboptimality-potential-subgradient})
    with respect to some $\varepsilon_{\OT}>0$,
    and let $\varepsilon_{\mathrm{vio}}:=\mathrm{vio}(\BIx,f_1,\ldots,f_d)$,
    $\varepsilon_{\mathrm{sub}}:=\varepsilon_{\mathrm{vio}}+d\varepsilon_{\OT}$.
    \begin{enumerate}[label=(\roman*),beginpenalty=10000]
        \item\label{cors:dual-suboptimality-primal-coupling-primal}%
        Then, $(\BIx,f_1-\nobreak\varepsilon_{\mathrm{vio}},\allowbreak f_2,\ldots,f_d)$ is an $\varepsilon_{\mathrm{sub}}$-optimizer of \eqref{eqn:dro-primal}, and
        $\BIx$ is an $\varepsilon_{\mathrm{sub}}$-optimizer of \eqref{eqn:dro}.
    \end{enumerate}
    In the following, let $K:=\nobreak|\CV|$, enumerate $\CV=\{\Blambda_1,\ldots,\Blambda_K\}$, and for $i=1,\ldots,d$,
    let $\sigma_i:\{1,\ldots,K\}\to\{1,\ldots,K\}$ be a permutation
    satisfying $\big\langle\Bp_i,\Blambda_{\sigma_i(1)}\big\rangle\le \big\langle\Bp_i,\Blambda_{\sigma_i(2)}\big\rangle\le\nobreak \cdots \le \big\langle\Bp_i,\Blambda_{\sigma_i(K)}\big\rangle$.
    Let $(\mathscr{X},\mathscr{F},\PROB)$ be a probability space,
    and let $J:\mathscr{X}\to \{1,\ldots,K\}$ be a categorical random variable satisfying $\PROB[J=\nobreak j]=\tau\big(\{\Blambda_j\}\big)$ for $j=1,\ldots,K$.
    Moreover, let $(U_1,\ldots,U_d)^\TRANSP:\mathscr{X}\to[0,1]^d$ be a random variable that is independent of $J$, such that each of ${U_1,\ldots,U_d}$ is uniformly distributed on $[0,1]$.
    In particular, the distribution function of $(U_1,\ldots,U_d)^\TRANSP$ is a copula.
    Furthermore, for $i=1,\ldots,d$,
    let $G^{-}_{\mu_i}:[0,1]\to\Omega_i$ be defined as follows:
    \begin{align}
        G^{-}_{\mu_i}(u)&:=\inf\big\{t\in\Omega_i: \mu_i\big([\underline{\omega}_i,t]\big)\ge u\big\} && \forall u\in[0,1],
        \label{eqn:invcdf-def}
    \end{align}
    and let 
    $Y_i:\mathscr{X}\to\Omega_i$ be a random variable defined as follows:
    \begin{align*}
        Y_i:=G^{-}_{\mu_i}\!\left(\textstyle U_i\tau\big(\{\Blambda_{J}\}\big) + \sum_{j=1}^{\sigma_i^{-1}(J)-1}\tau\big(\{\Blambda_{\sigma_i(j)}\}\big)\!\right)\!.
    \end{align*}
    Lastly, let $\mu_{\tau}\in\CP(\BOmega)$ denote the law of the random variable $(Y_1,\ldots,Y_d)^\TRANSP:$ $\mathscr{X}\to\BOmega$.
    \begin{enumerate}[label=(\roman*),beginpenalty=10000]
        \setcounter{enumi}{1}
        \item\label{cors:dual-suboptimality-primal-coupling-meas}%
        Then, $\mu_{\tau}$ is an $\varepsilon_{\mathrm{sub}}$-optimizer of the inner maximization problem in \eqref{eqn:dro} with respect to $\BIx$,
        i.e., $\mu_{\tau}$ is an $\varepsilon_{\mathrm{sub}}$-approximate worst-case probability measure with respect to the first-stage decision $\BIx$.
    \end{enumerate}
\end{corollary}

In Corollary~\ref{cor:dual-suboptimality-primal-coupling},
$\mu_{\tau}$ is constructed to be the law of a random variable $(Y_1,\ldots,Y_d)^{\TRANSP}$ on a probability space $(\mathscr{X},\mathscr{F},\PROB)$.
Thus, independent samples from $\mu_{\tau}$ can be efficiently generated and used for the numerical approximation of $\mu_{\tau}$.
This is presented and discussed in Algorithm~\ref{alg:iterative} and Setting~\ref{sett:iterative}.

\section{Numerical algorithm}
\label{sec:algorithm}

Combining our results in Theorem~\ref{thm:primaldual}, Theorem~\ref{thm:dual-suboptimality}, and Corollary~\ref{cor:dual-suboptimality-primal-coupling} leads to an iterative approximation scheme for \eqref{eqn:dro},
where one begins with an arbitrary feasible solution $(\tau,\Bxi_{\ineq},\Bxi_{\eq})$ of \eqref{eqn:dro-dual} where $\tau$ has finite support, and iteratively adds points to the support of $\tau$ to obtain increasingly better approximate optimizers of \eqref{eqn:dro-dual}.
This process also produces increasingly better approximate optimizers of \eqref{eqn:dro-primal}, and hence of \eqref{eqn:dro}, as a consequence of Corollary~\ref{cor:dual-suboptimality-primal-coupling}.

Concretely, let us first consider a non-empty finite set $\CV\subset\CS^*_2$ and 
the restriction of \eqref{eqn:dro-dual} obtained by requiring $\support(\tau)\subseteq\CV$.
This allows us to parametrize $\tau$ by $\nu_{\Blambda}\leftarrow\tau\big(\{\Blambda\}\big)\in\nobreak\R_+$ $\forall\Blambda\in\nobreak\CV$
and reduce the decision space of \eqref{eqn:dro-dual} to a finite-dimensional one.
Observe that $\OT_{\mu_i}\big(\langle\Bp_i,\cdot\,\rangle\sharp\tau\big)=\max_{f\in\CC_{\mathrm{cvx}}(\R)}\Big\{{-\int_{\Omega_i}}\sup_{z\in\R}\big\{z\omega-f(z)\big\}\DIFFM{\mu_i}{\DIFF\omega}-\int_{\R}f\DIFFX{\langle\Bp_i,\cdot\,\rangle\sharp\tau}\Big\}$ for ${i=1,\ldots,d}$ and for any $\tau\in\CP_1(\CS^*_2)$ with $\support(\tau)\subseteq\CV$; see, e.g., Lemma~\ref{lem:OT-1d}\ref{lems:OT-1d-duality} with $\mu\leftarrow\mu_i$, $\nu\leftarrow \langle\Bp_i,\cdot\,\rangle\sharp\tau$.
This max-representation allows us to
linearize each of the nonlinear terms $\OT_{\mu_i}\big(\langle\Bp_i,\cdot\,\rangle\sharp\tau\big)$ 
in the objective of \eqref{eqn:dro-dual} and formulate the aforementioned restriction of \eqref{eqn:dro-dual} into the following linear semi-infinite programming (LSIP) problem:
\begin{align}
    \maximize_{\substack{(\nu_{\Blambda}),\,(\psi_i),\\ \Bxi_{\ineq},\,\Bxi_{\eq}}} \quad & 
    \Bigg(\sum_{\Blambda\in\CV} \langle\Bb,\Blambda\rangle\nu_{\Blambda}\Bigg) 
    - \Bigg(\sum_{i=1}^d \psi_i\Bigg)
    + \langle\Bq_{\ineq},\Bxi_{\ineq}\rangle
    + \langle\Bq_{\eq},\Bxi_{\eq}\rangle 
    \tag{\textsf{LSIP}-$\CV$}
    \label{eqn:dro-dual-LSIP}\\
    \mathrm{subject~to}\quad & 
    \psi_i  
    \ge {-\!\hspace{-0.6pt}\int_{\Omega_i}} \sup_{z\in\R}\!\big\{\hspace{-0.8pt}z\omega-f_i(z)\hspace{-0.8pt}\big\}\hspace{-0.6pt} \DIFFM{\mu_i}{\DIFF\omega}
    \hspace{-0.6pt}-\hspace{-0.6pt} \sum_{\Blambda\in\CV}\hspace{-0.6pt}f_i\big(\hspace{-0.8pt}\langle\Bp_i,\Blambda\rangle\hspace{-0.8pt}\big)\nu_{\Blambda} \;\;\; \forall f_i\hspace{-0.6pt}\in\hspace{-0.6pt}\CC_{\mathrm{cvx}}(\R),\;\hspace{-1pt} \forall 1\hspace{-0.6pt}\le\hspace{-0.6pt} i\hspace{-0.6pt}\le\hspace{-0.6pt} d,\nonumber\\
    & \Bc_1+\sum_{\Blambda\in\CV}\BQ^\TRANSP\Blambda\nu_{\Blambda}
    \ge \BL_{\ineq}^\TRANSP\Bxi_{\ineq} 
    + \BL_{\eq}^\TRANSP\Bxi_{\eq}, \qquad 
    \sum_{\Blambda\in\CV} \nu_{\Blambda}=1,\nonumber \\
    & \nu_{\Blambda}\in\R_+ \quad \forall \Blambda\in\CV, \qquad
    \Bxi_{\ineq}\in\R_{-}^{n_{\ineq}}, \qquad
    \Bxi_{\eq}\in \R^{n_{\eq}}.\nonumber
\end{align}
In particular, as a consequence of Theorem~\ref{thm:primaldual}\ref{thms:primaldual-dual-optimizer-existence},
the optimal values of \eqref{eqn:dro-dual} and \eqref{eqn:dro-dual-LSIP} are equal whenever $\CV\supseteq\CV_{\extremept}$. 
Observe that \eqref{eqn:dro-dual-LSIP} has infinitely many constraints of the form 
$\psi_i \ge {-\int_{\Omega_i}} \sup_{z\in\R}\big\{z\omega-f_i(z)\big\} \DIFFM{\mu_i}{\DIFF\omega} - \sum_{\Blambda\in\CV}f_i\big(\langle\Bp_i,\Blambda\rangle\big)\nu_{\Blambda}$ for $f_i\in\CC_{\mathrm{cvx}}(\R)$.
Despite this,
for any $\psi_1,\ldots,\psi_d\in\R$,
for any $\tau:=\sum_{\Blambda\in\CV}\nu_{\Blambda}\delta_{\Blambda}$ with $\nu_{\Blambda}\in\R_+$ $\forall\Blambda\in\CV$, 
$\sum_{\Blambda\in\CV}\nu_{\Blambda}=1$,
and for $i=1,\ldots,d$, 
we are able to efficiently determine a continuous piece-wise affine function $f^\star_{i,\tau}\in\CC_{\mathrm{cvx}}(\R)$ which leads to the largest violation of the constraint, 
i.e., $f^\star_{i,\tau}\in\argmin_{f\in\CC_{\mathrm{cvx}}(\R)}\big\{\psi_i  + \int_{\Omega_i} \sup_{z\in\R}\big\{z\omega-f(z)\big\} \DIFFM{\mu_i}{\DIFF\omega} + \sum_{\Blambda\in\CV}f\big(\langle\Bp_i,\Blambda\rangle\big)\nu_{\Blambda}\big\}$.
Specifically, for any $\tau\in\CP_1(\CS^*_2)$ with $\support(\tau)\subseteq\CV$ and for $i=1,\ldots,d$,
let $m_i:=\big|\big\{\langle\Bp_i,\Blambda\rangle:\Blambda\in\CV\big\}\big|$,
let $\{\kappa_{i,1},\ldots,\kappa_{i,m_i}\}$ be an enumeration of $\big\{\langle\Bp_i,\Blambda\rangle:\Blambda\in\CV\big\}$
with $\kappa_{i,1}<\kappa_{i,2}<\cdots<\nobreak\kappa_{i,m_i}$,
let $G^{-}_{\mu_i}:[0,1]\to\Omega_i$ be defined by (\ref{eqn:invcdf-def}),
and let $f^\star_{i,\tau}:\R\to\R$ be defined as follows:
\begin{align}
    f^\star_{i,\tau}(z):=\begin{cases}
        \underline{\omega}_i(z-\kappa_{i,1}) & 
        \hspace{82.6pt} \forall z\in(-\infty,\kappa_{i,1}], \\
        f^\star_{i,\tau}(\kappa_{i,j}) + G^{-}_{\mu_i}\big(\langle\Bp_i,\cdot\,\rangle\sharp\tau\big((-\infty,\kappa_{i,j}]\big)\big)(z - \kappa_{i,j}) \hspace{-60pt}  \\
        & \hspace{-10pt} \forall z\in (\kappa_{i,j}, \kappa_{i,j+1}],\; \forall 1\le j\le m_i-1,\\
        f^\star_{i,\tau}(\kappa_{i,m_i}) + \overline{\omega}_i(z-\kappa_{i,m_i}) & 
        \hspace{83.5pt} \forall z\in(\kappa_{i,m_i},\infty).
    \end{cases}
    \label{eqn:algo-potential-CPWA}
\end{align}
One can verify that $f^\star_{i,\tau}\in\argmax_{f\in\CC_{\mathrm{cvx}}(\R)}\big\{\int_{\Omega_i} \sup_{z\in\R}\big\{z\omega-f(z)\big\} \DIFFM{\mu_i}{\DIFF\omega} + \int_{\R}f\DIFFX{\langle\Bp_i,\cdot\,\rangle\sharp\tau}\big\}$; see
Lemma~\ref{lem:OT-1d}\ref{lems:OT-1d-potential-existence}.
Moreover, it holds that $\partial f^\star_{i,\tau}(z)\subseteq\Omega_i$ for all $z\in\R$.

Subsequently,
we develop a cutting-plane algorithm presented in Algorithm~\ref{alg:cuttingplane-LSIP} for solving \eqref{eqn:dro-dual-LSIP}
that is inspired by the Conceptual Algorithm~11.4.1 of \citet{goberna1998linear}.
Each iteration of Algorithm~\ref{alg:cuttingplane-LSIP} solves a linear programming (LP) relaxation \eqref{eqn:dro-dual-LSIP} and its dual LP problem defined as follows.
For any $\widehat{\CF}=(\widehat{\CF}_{1},\ldots,\widehat{\CF}_{d})$ where 
$\widehat{\CF}_{i}\subset \CC_{\mathrm{cvx}}(\R)$ is non-empty and finite for $i=1,\ldots,d$,
we consider the LP relaxation of \eqref{eqn:dro-dual-LSIP} formed via replacing 
each semi-infinite constraint in \eqref{eqn:dro-dual-LSIP} with respect to all $f_i\in\CC_{\mathrm{cvx}}(\R)$ with 
finitely many constraints with respect to all $f_i\in\widehat{\CF}_{i}$,
and denote it
by \eqref{eqn:dro-dual-LP-relax}:
\begin{align}
    \maximize_{\substack{(\nu_{\Blambda}),\,(\psi_i),\\ \Bxi_{\ineq},\,\Bxi_{\eq}}} \quad & 
    \Bigg(\sum_{\Blambda\in\CV}\langle\Bb,\Blambda\rangle\nu_{\Blambda}\Bigg) 
    - \Bigg(\sum_{i=1}^d \psi_i\Bigg)
    + \langle\Bq_{\ineq},\Bxi_{\ineq}\rangle
    + \langle\Bq_{\eq},\Bxi_{\eq}\rangle
    \tag{\textsf{LP}-$\CV$-$\widehat{\CF}$}
    \label{eqn:dro-dual-LP-relax}\\
    \mathrm{subject~to}\quad 
    & \psi_i 
    \ge {-\int_{\Omega_i}} \sup_{z\in\R}\big\{z\omega-f_i(z)\big\} \DIFFM{\mu_i}{\DIFF\omega}
    - \sum_{\Blambda\in\CV}f_i\big(\langle\Bp_i,\Blambda\rangle\big)\nu_{\Blambda}\quad \forall f_i\in\widehat{\CF}_{i},\; \forall 1\le i\le d,\nonumber\\
    & \Bc_1+\sum_{\Blambda\in\CV}\BQ^\TRANSP\Blambda\nu_{\Blambda}
    \ge \BL_{\ineq}^\TRANSP\Bxi_{\ineq} 
    + \BL_{\eq}^\TRANSP\Bxi_{\eq}, \qquad 
    \sum_{\Blambda\in\CV} \nu_{\Blambda}=1,\nonumber \\
    & \nu_{\Blambda}\in\R_+ \quad \forall \Blambda\in\CV, \qquad
    \Bxi_{\ineq}\in\R_{-}^{n_{\ineq}}, \qquad
    \Bxi_{\eq}\in \R^{n_{\eq}},\nonumber
\end{align}
and we 
denote the corresponding dual LP by \eqref{eqn:dro-dual-LP-relax-dual}:
\begin{align}
    \minimize_{\!\BIx,\,h_0,\,(y_{i,f_i})} \quad & \langle\Bc_1,\BIx\rangle-h_0+\sum_{i=1}^{d}\sum_{f_i\in\widehat{\CF}_{i}} y_{i,f_i}\int_{\Omega_i} \sup_{z\in\R}\big\{z\omega - f_i(z)\big\}\DIFFM{\mu_i}{\DIFF\omega}
    \tag{$\textsf{LP}^*$-$\CV$-$\widehat{\CF}$}
    \label{eqn:dro-dual-LP-relax-dual}\\
    \mathrm{subject~to} \quad & \langle\BQ\BIx+\Bb,\Blambda\rangle + h_0 + \sum_{i=1}^{d} \sum_{f_i\in\widehat{\CF}_{i}}y_{i,f_i}f_i\big(\langle\Bp_i,\Blambda\rangle\big) \le 0 \quad \forall \Blambda\in\CV,\nonumber\\
    & \BL_{\ineq}\BIx \le \Bq_{\ineq}, \qquad 
    \BL_{\eq}\BIx = \Bq_{\eq}, \qquad 
    \sum_{f_i\in\widehat{\CF}_{i}}y_{i,f_i}=1 \quad \forall 1\le i\le d,\nonumber \\
    & \BIx\in \R_{+}^{k_1}, \qquad 
    h_0\in\R, \qquad
    y_{i,f_i}\in\R_{+} \quad\hspace{2.4pt} \forall f_i\in\widehat{\CF}_{i},\; \forall 1\le i\le d.\nonumber
\end{align}
Setting~\ref{sett:cuttingplane-LSIP} below introduces the assumptions for the inputs of Algorithm~\ref{alg:cuttingplane-LSIP}
and discusses concrete implementation details about 
some parts of Algorithm~\ref{alg:cuttingplane-LSIP}.
Proposition~\ref{prop:cuttingplane-LSIP} establishes crucial theoretical properties of Algorithm~\ref{alg:cuttingplane-LSIP}.

\setlength{\textfloatsep}{6pt}
\begin{algorithm}[t]
	\caption{\textbf{Cutting-plane algorithm for solving \eqref{eqn:dro-dual-LSIP}.}}
	\label{alg:cuttingplane-LSIP}
	\KwIn{
		$\BOmega$, 
        $(\mu_i)_{i=1:d}$, 
        $\Bc_1$, 
        $\BL_{\ineq}$,
        $\Bq_{\ineq}$,
        $\BL_{\eq}$,
        $\Bq_{\eq}$,
        $\BQ$,
        $\BP$,
        $\Bb$,
        $\CV$, 
        $\varepsilon_{\OT}>0$.}
	\KwOut{$\hat{\tau}$,
    $\hat{\Bxi}_{\ineq}$, $\hat{\Bxi}_{\eq}$,
    $\hat{\BIx}$,
    $\hat{f}_1$, $\ldots$~, $\hat{f}_d$,
    $\widehat{C}_{\DRO}^{\mathrm{LB}}$.
    }
	\nl\label{alglin:cuttingplane-LSIP-initialmeas}%
    $\hat{\tau}^{(0)}\leftarrow \frac{1}{|\CV|}\sum_{\Blambda\in\CV}\delta_{\Blambda}$.

    \nl\label{alglin:cuttingplane-LSIP-initialfuncs-forloop}%
    \For{$i=1,\ldots,d$}{
        \nl\label{alglin:cuttingplane-LSIP-initialfuncs}%
        $\widehat{\CF}_{i}^{(1)}\leftarrow\big\{f^{\star}_{i,\hat{\tau}^{(0)}}\big\}$ where $f^{\star}_{i,\hat{\tau}^{(0)}}\in\CC_{\mathrm{cvx}}(\R)$ is defined by (\ref{eqn:algo-potential-CPWA}).
    }

    \nl\label{alglin:cuttingplane-LSIP-forloop}%
    \For{$t=1,2,\ldots$}{

        \nl\label{alglin:cuttingplane-LSIP-LP}%
        Solve \LPprimal{\widehat{\CF}^{(t)}} and 
        \LPdual{\widehat{\CF}^{(t)}} 
        with respect to $\widehat{\CF}^{(t)}\leftarrow \big(\widehat{\CF}_{1}^{(t)},\ldots,\widehat{\CF}_{d}^{(t)}\big)$
        to compute a pair of primal and dual optimizers
        $\big((\hat{\nu}_{\Blambda}^{(t)})_{\Blambda\in\CV},
        (\hat{\psi}_i^{(t)})_{i=1:d},
        \hat{\Bxi}_{\ineq}^{(t)}, \hat{\Bxi}_{\eq}^{(t)}\big)$ 
        and 
        $\big(\hat{\BIx}^{(t)},\hat{h}_0^{(t)},(\hat{y}_{i,f_i}^{(t)})_{f_i\in\widehat{\CF}_{i}^{(t)},\,i=1:d}\big)$, and
        denote the computed optimal value by $\widehat{C}_{\DRO}^{(t)}$.

        \nl\label{alglin:cuttingplane-LSIP-constructmeas}%
        $\hat{\tau}^{(t)}\leftarrow \sum_{\Blambda\in\CV}\hat{\nu}_{\Blambda}^{(t)}\delta_{\Blambda}$.

        \nl\label{alglin:cuttingplane-LSIP-newfuncs-forloop}%
        \For{$i=1,\ldots,d$}{

            \nl\label{alglin:cuttingplane-LSIP-newfuncs}%
            $\widehat{\CF}_{i}^{(t+1)}\leftarrow \widehat{\CF}_{i}^{(t)}\union \big\{f^{\star}_{i,\hat{\tau}^{(t)}}\big\}$ where $f^{\star}_{i,\hat{\tau}^{(t)}}\in\CC_{\mathrm{cvx}}(\R)$ is defined by (\ref{eqn:algo-potential-CPWA}).

            \nl\label{alglin:cuttingplane-LSIP-error}%
            $\hat{\varepsilon}_{\OT,i}^{(t)}\leftarrow {-\int_{\Omega_i}}\sup_{z\in\R}\big\{z\omega-f_{i,\hat{\tau}^{(t)}}^{\star}(z)\big\}\DIFFM{\mu_i}{\DIFF\omega} - \int_{\R}f_{i,\hat{\tau}^{(t)}}^{\star}\DIFFX{\langle\Bp_i,\cdot\,\rangle\sharp\hat{\tau}^{(t)}} - \hat{\psi}_i^{(t)}$.
        }

        \nl\label{alglin:cuttingplane-LSIP-termination}%
        \If{$\max_{1\le i\le d}\big\{\hat{\varepsilon}_{\OT,i}^{(t)}\}\le \varepsilon_{\OT}$}{
            \nl \textbf{break}.
        }
    }

    \nl\label{alglin:cuttingplane-LSIP-returnval}%
    $\hat{\tau}\leftarrow \hat{\tau}^{(t)}$,
    $\hat{\Bxi}_{\ineq}\leftarrow \hat{\Bxi}_{\ineq}^{(t)}$,
    $\hat{\Bxi}_{\eq}\leftarrow \hat{\Bxi}_{\eq}^{(t)}$,
    $\hat{\BIx}\leftarrow\hat{\BIx}^{(t)}$.

    \nl\label{alglin:cuttingplane-LSIP-returnfuncs-forloop}%
    \For{$i=1,\ldots,d$}{

        \nl\label{alglin:cuttingplane-LSIP-returnfuncs}%
        $\hat{f}_i\leftarrow \INDI_{\{i=1\}}\hat{h}_0^{(t)} + \sum_{f_i\in\widehat{\CF}_{i}^{(t)}}\hat{y}_{i,f_i}^{(t)}f_i$.
    }

    \nl\label{alglin:cuttingplane-LSIP-lowerbound}%
    $\widehat{C}_{\DRO}^{\mathrm{LB}}\leftarrow\widehat{C}_{\DRO}^{(t)}-\big(\sum_{i=1}^d\hat{\varepsilon}_{\OT,i}^{(t)}\big)$.

	\nl\Return $\hat{\tau}$,
    $\hat{\Bxi}_{\ineq}$, $\hat{\Bxi}_{\eq}$,
    $\hat{\BIx}$,
    $\hat{f}_1$, $\ldots$ , $\hat{f}_d$,
    $\widehat{C}_{\DRO}^{\mathrm{LB}}$.
\end{algorithm}

\begin{setting}[Algorithm~\ref{alg:cuttingplane-LSIP}]\label{sett:cuttingplane-LSIP}
    We make the following assumptions about the inputs of Algorithm~\ref{alg:cuttingplane-LSIP}.

    \begin{itemize}[beginpenalty=10000]
        \item 
        $\BOmega$, $(\mu_i)_{i=1:d}$, $\Bc_1$, $\BL_{\ineq}$, $\Bq_{\ineq}$, $\BL_{\eq}$, $\Bq_{\eq}$, $\BQ$, $\BP$, $\Bb$ are taken from the specification of the two-stage DRO problem \eqref{eqn:dro} in Setting~\ref{sett:dro-alt}.
        We assume in addition that $C_{\DRO}>-\infty$.
        
        \item
        $\CV$ is a non-empty finite subset of the set $\CS^*_2$ in Setting~\ref{sett:dro-alt}
        such that there exist $\Blambda\in\CV$, $\Bxi_{\ineq}\in\R_{-}^{n_{\ineq}}$, $\Bxi_{\eq}\in\R^{n_{\eq}}$ satisfying 
        $\Bc_1+\BQ^\TRANSP\Blambda\ge \BL_{\ineq}^\TRANSP\Bxi_{\ineq}+\BL_{\eq}^\TRANSP\Bxi_{\eq}$.
        Such a choice of $\CV$ guarantees that \eqref{eqn:dro-dual-LSIP} 
        is feasible.
        Under the assumption that $C_{\DRO}>-\infty$,
        Theorem~\ref{thm:primaldual}\ref{thms:primaldual-boundedness}
        shows that there exist $\Blambda_0\in\CS^*_2$, $\Bxi_{\ineq,0}\in\nobreak\R_{-}^{n_{\ineq}}$, $\Bxi_{\eq,0}\in\R^{n_{\eq}}$ 
        which satisfy $\Bc_1+\BQ^\TRANSP\Blambda_0\ge \BL_{\ineq}^\TRANSP\Bxi_{\ineq,0}+\BL_{\eq}^\TRANSP\Bxi_{\eq,0}$.
        Subsequently, one may set $\CV:=\{\Blambda_0\}$ as a particular choice of $\CV$ to guarantee the feasibility of \eqref{eqn:dro-dual-LSIP}.
        We remark that such $\Blambda_0$, $\Bxi_{\ineq,0}$, $\Bxi_{\eq,0}$  can be obtained 
        via computing a feasible solution of an LP problem.

        \item
        $\varepsilon_{\OT}>0$ is a user-specified tolerance parameter
        which controls the termination of the algorithm.
        As shown later in 
        Proposition~\ref{prop:cuttingplane-LSIP}\ref{props:cuttingplane-LSIP-termination}--\ref{props:cuttingplane-LSIP-primal}, 
        Algorithm~\ref{alg:cuttingplane-LSIP} is guaranteed to terminate 
        after finitely many iterations for any $\varepsilon_{\OT}>0$
        and the returned
        $\hat{\BIx},\hat{f}_1,\ldots,\hat{f}_d$ 
        satisfy the constraints
        (\ref{eqn:dual-suboptimality-primal-feasibility-potential-eq})--(\ref{eqn:dual-suboptimality-potential-subgradient}) with respect to $\varepsilon_{\OT}$.
    \end{itemize}

    In the following, we discuss the details of some lines in
    Algorithm~\ref{alg:cuttingplane-LSIP}.
    Line~\ref{alglin:cuttingplane-LSIP-initialfuncs} and Line~\ref{alglin:cuttingplane-LSIP-newfuncs}
    require one to compute a continuous piece-wise affine function 
    $f^\star_{i,\hat{\tau}^{(t)}}\in\CC_{\mathrm{cvx}}(\R)$ defined in (\ref{eqn:algo-potential-CPWA}),
    where $\hat{\tau}^{(t)}\in\CP_1(\CS^*_2)$ satisfies $\support(\hat{\tau}^{(t)})\subseteq\nobreak\CV$ due to 
    Line~\ref{alglin:cuttingplane-LSIP-initialmeas} and Line~\ref{alglin:cuttingplane-LSIP-constructmeas};
    see Remark~\ref{rmk:cuttingplane-LSIP-tractability} for our discussion about the tractability of computing $f^\star_{i,\hat{\tau}^{(t)}}$ and the integral $\int_{\Omega_i}\sup_{z\in\R}\big\{z\omega-f^\star_{i,\hat{\tau}^{(t)}}(z)\big\}\DIFFM{\mu_i}{\DIFF\omega}$.
    Moreover, since $\widehat{\CF}^{(t)}_1,\ldots,\widehat{\CF}^{(t)}_{d}$ contain only continuous piece-wise affine functions,
    we are able to tractably formulate a mixed-integer programming problem when solving the global maximization problem in Line~\ref{alglin:iterative-globalmax} of Algorithm~\ref{alg:iterative}; 
    see Proposition~\ref{prop:cuttingplane-LSIP}\ref{props:cuttingplane-LSIP-primal} and
    Proposition~\ref{prop:globalmax-MILP}.
    
    Line~\ref{alglin:cuttingplane-LSIP-LP} solves an LP relaxation of \eqref{eqn:dro-dual-LSIP} and its corresponding dual LP.
    This can be numerically carried out using, for example, 
    the simplex algorithm
    (see, e.g., \citep[Chapter~6.4]{vanderbei2020linear})
    or the interior point algorithm with a crossover phase
    (see, e.g., \citep{megiddo1991on}).
    We would like to remark that both \eqref{eqn:dro-dual-LP-relax} and \eqref{eqn:dro-dual-LP-relax-dual} are feasible; see
    Proposition~\ref{prop:cuttingplane-LSIP}\ref{props:cuttingplane-LSIP-LPfeasibility}.
    Thus, the strong duality between \eqref{eqn:dro-dual-LP-relax} and \eqref{eqn:dro-dual-LP-relax-dual} holds.
    Moreover, if the simplex algorithm is used, 
    the primal basis and dual basis computed by the simplex algorithm when solving 
    \LPprimal{\widehat{\CF}^{(t-1)}} and 
    \LPdual{\widehat{\CF}^{(t-1)}} 
    can be stored and used to ``warm start''
    the simplex algorithm when solving 
    \LPprimal{\widehat{\CF}^{(t)}} and 
    \LPdual{\widehat{\CF}^{(t)}}.
    This can speed up Algorithm~\ref{alg:cuttingplane-LSIP} considerably.
\end{setting}

\begin{proposition}[Properties of Algorithm~\ref{alg:cuttingplane-LSIP}]\label{prop:cuttingplane-LSIP}
    Let the inputs of Algorithm~\ref{alg:cuttingplane-LSIP} be specified according to
    Setting~\ref{sett:cuttingplane-LSIP}.
    Then, the following statements hold.
    \begin{enumerate}[label=(\roman*),beginpenalty=10000]
        \item\label{props:cuttingplane-LSIP-LPfeasibility}
        For any $\widehat{\CF}=(\widehat{\CF}_{1},\ldots,\widehat{\CF}_{d})$, where $\widehat{\CF}_{i}\subset\CC_{\mathrm{cvx}}(\R)$ is non-empty and finite for $i=1,\ldots,d$,
        both \eqref{eqn:dro-dual-LP-relax} and \eqref{eqn:dro-dual-LP-relax-dual} are feasible.
        In particular, both \LPprimal{\widehat{\CF}^{(t)}} and \LPdual{\widehat{\CF}^{(t)}} in Line~\ref{alglin:cuttingplane-LSIP-LP}
        of Algorithm~\ref{alg:cuttingplane-LSIP} are feasible. 

        \item\label{props:cuttingplane-LSIP-termination}
        Algorithm~\ref{alg:cuttingplane-LSIP} terminates after finitely many iterations.

        \item\label{props:cuttingplane-LSIP-dual}
        $(\hat{\tau},\hat{\Bxi}_{\ineq},\hat{\Bxi}_{\eq})$
        is feasible for \eqref{eqn:dro-dual} with objective $\widehat{C}_{\DRO}^{\mathrm{LB}}$ and
        $\support(\hat{\tau})\subseteq\CV$.
        
        \item\label{props:cuttingplane-LSIP-primal}
        $\hat{\BIx},\hat{f}_1,\ldots,\hat{f}_d$ 
        satisfy (\ref{eqn:dual-suboptimality-primal-feasibility-potential-eq})--(\ref{eqn:dual-suboptimality-potential-subgradient})
        with respect to 
        $\tau\leftarrow\hat{\tau}$, $\Bxi_{\ineq}\leftarrow\hat{\Bxi}_{\ineq}$, $\Bxi_{\eq}\leftarrow\hat{\Bxi}_{\eq}$,
        as well as the user-specified tolerance~$\varepsilon_{\OT}>0$,
        and $\hat{f}_1,\ldots,\hat{f}_d$ are continuous piece-wise affine.

    \end{enumerate}
\end{proposition}

\begin{remark}[Computational tractability of Algorithm~\ref{alg:cuttingplane-LSIP}]\label{rmk:cuttingplane-LSIP-tractability}
    In order to solve the LP problems
    \eqref{eqn:dro-dual-LP-relax} and \eqref{eqn:dro-dual-LP-relax-dual},
    one needs to evaluate the one-dimensional integral:
    $\int_{\Omega_i}\sup_{z\in\R}\big\{z\omega-f_i(z)\big\}\DIFFM{\mu_i}{\DIFF \omega}$ 
    for each $f_i\in\widehat{\CF}_{i}$,
    for $i=1,\ldots,d$.
    Since Line~\ref{alglin:cuttingplane-LSIP-initialfuncs}, Line~\ref{alglin:cuttingplane-LSIP-newfuncs}, and (\ref{eqn:algo-potential-CPWA}) guarantee that every $f_i\in\widehat{\CF}^{(t)}_{i}$ is continuous piece-wise affine,
    one checks that the integrand
    $\Omega_i\ni\omega\mapsto \sup_{z\in\R}\big\{z\omega-f_i(z)\big\}=\max_{1\le j\le m_i}\big\{\kappa_{i,j}\omega-f_i(\kappa_{i,j})\big\}\in\nobreak\R$
    is also continuous piece-wise affine, 
    and the integral can be tractably computed 
    as long as the integral of any continuous piece-wise affine function on $\Omega_i$ with respect to $\mu_i$ can be tractably computed.
    In particular, 
    when $\mu_i$ is the mixture of a discrete measure on $\Omega_i$ 
    and an absolutely continuous measure on $\Omega_i$,
    where the absolutely continuous measure has a density with respect to the Lebesgue measure that is continuous piece-wise affine,
    or given by a truncated mixture of Gaussian measures,
    one may check that the integral of any continuous piece-wise affine function on $\Omega_i$ with respect to $\mu_i$ can be tractably computed.
\end{remark}

\begin{remark}[Justification of the LSIP formulation]
    \label{rmk:justification-LSIP}
    We choose to tackle the nonlinear optimization problem \eqref{eqn:dro-dual} via the LSIP formulation \eqref{eqn:dro-dual-LSIP} since it allows us to
    develop a cutting-plane algorithm (i.e., Algorithm~\ref{alg:cuttingplane-LSIP}) 
    to utilize state-of-the-art LP solvers to efficiently compute approximate optimizers of \eqref{eqn:dro-dual}.
    More importantly, 
    each iteration of Algorithm~\ref{alg:cuttingplane-LSIP}
    computes a pair of primal and dual optimizers of \eqref{eqn:dro-dual-LP-relax} and \eqref{eqn:dro-dual-LP-relax-dual} 
    which then 
    yield $\hat{\BIx},\hat{f}_1,\ldots,\hat{f}_d$ 
    satisfying (\ref{eqn:dual-suboptimality-primal-feasibility-potential-eq})--(\ref{eqn:dual-suboptimality-potential-subgradient}).
    These conditions are crucial when we construct feasible and approximately optimal solutions of \eqref{eqn:dro-primal}, \eqref{eqn:dro-dual}, and \eqref{eqn:dro} in Algorithm~\ref{alg:iterative};
    see Theorem~\ref{thm:iterative}.
\end{remark}

Next, using Algorithm~\ref{alg:cuttingplane-LSIP} as a subroutine,
we develop Algorithm~\ref{alg:iterative} for computing approximate optimizers of \eqref{eqn:dro-primal}, \eqref{eqn:dro-dual}, \eqref{eqn:dro}, as well as an approximate worst-case probability measure with respect to the computed first-stage decision.
Setting~\ref{sett:iterative} introduces the assumptions for the inputs of Algorithm~\ref{alg:iterative}
and provides detailed explanations about some lines in Algorithm~\ref{alg:iterative}.
The theoretical properties of Algorithm~\ref{alg:iterative} are presented in Theorem~\ref{thm:iterative}.

\setlength{\textfloatsep}{6pt}
\begin{algorithm}[t]
    \caption{\textbf{Computation of approximate optimizers of \eqref{eqn:dro-primal}, \eqref{eqn:dro-dual}, and \eqref{eqn:dro}.}}
	\label{alg:iterative}
	\KwIn{
		$\BOmega$, 
        $(\mu_i)_{i=1:d}$, 
        $\Bc_1$, 
        $\BL_{\ineq}$,
        $\Bq_{\ineq}$,
        $\BL_{\eq}$,
        $\Bq_{\eq}$,
        $\BQ$,
        $\BP$,
        $\Bb$,
        $\BA_{\ineq}$,
        $\BA_{\eq}$,
        $\Bc_2$,
        $\big(\varepsilon_{\OT}^{(t)}\big)_{t\in\N}\subset(0,\infty)$,
        $\big(\varsigma^{(t)}\big)_{t\in\N}\subset \R_+$, 
        $\big(\vartheta^{(t)}\big)_{t\in\N}\subset \R_+$,
        $\varepsilon_{\mathrm{tol}}>\liminf_{t\to\infty} d\big(\varepsilon_{\OT}^{(t)}+\big(1+\vartheta^{(t)}\big)\varsigma^{(t)}\big)$.}
	\KwOut{$\hat{\tau}$,
    $\hat{\Bxi}_{\ineq}$, $\hat{\Bxi}_{\eq}$,
    $\hat{\BIx}$,
    $\tilde{f}_1$, $\ldots$~, $\tilde{f}_d$,
    $\hat{\mu}$,
    $\widehat{C}_{\DRO}^{\mathrm{LB}}$,
    $\widehat{C}_{\DRO}^{\mathrm{UB}}$,
    $\hat{\varepsilon}_{\mathrm{sub}}$,
    $\hat{\varepsilon}_{\mathrm{prob}}$.}
	
    \nl\label{alglin:iterative-initialization}%
    Find $\Blambda_0\in\CS^*_2$, $\Bxi_{\ineq,0}\in\R_{-}^{n_{\ineq}}$, $\Bxi_{\eq,0}\in\R^{n_{\eq}}$ that satisfy
    $\Bc_1+\BQ^\TRANSP\Blambda_0\ge \BL_{\ineq}^\TRANSP\Bxi_{\ineq,0}+\BL_{\eq}^\TRANSP\Bxi_{\eq,0}$.

    \nl\label{alglin:iterative-initial-support}%
    $\CV^{(1)}\leftarrow\{\Blambda_0\}$.

    \nl\label{alglin:iterative-forloop}%
    \For{$t=1,2,\ldots$}{

        \nl\label{alglin:iterative-LSIP}%
        Call Algorithm~\ref{alg:cuttingplane-LSIP} with inputs
        $\BOmega$, 
        $(\mu_i)_{i=1:d}$, 
        $\Bc_1$, 
        $\BL_{\ineq}$,
        $\Bq_{\ineq}$,
        $\BL_{\eq}$,
        $\Bq_{\eq}$,
        $\BQ$,
        $\BP$,
        $\Bb$,
        $\CV^{(t)}$, 
        $\varepsilon_{\OT}^{(t)}$,
        and denote the outputs by
        $\hat{\tau}^{(t)}$,
        $\hat{\Bxi}_{\ineq}^{(t)}$, $\hat{\Bxi}_{\eq}^{(t)}$,
        $\hat{\BIx}^{(t)}$,
        $\hat{f}_1^{(t)}$, $\ldots$~, $\hat{f}_d^{(t)}$,
        $\widehat{C}_{\DRO}^{\mathrm{LB}(t)}$.

        \nl\label{alglin:iterative-funcrelax-forloop}%
        \For{$i=1,\ldots,d$}{
            \nl\label{alglin:iterative-funcrelax}%
            Let $\overline{f}_i^{(t)}\in\CC_{\mathrm{cvx}}(\R)$ be a continuous piece-wise affine function that satisfies
            $\partial \overline{f}_i^{(t)}(z)\subseteq\Omega_i$ and
            $0\le \overline{f}_i^{(t)}(z)-\hat{f}_i^{(t)}(z)\le \varsigma^{(t)}$ for all $z\in\R$.

            \nl\label{alglin:iterative-primalobj}%
            $\hat{\beta}_i^{(t)}\leftarrow \int_{\Omega_i}\sup_{z\in\R}\big\{z\omega-\overline{f}_i^{(t)}(z)\big\}\DIFFM{\mu_i}{\DIFF\omega}$.
        }

        \nl\label{alglin:iterative-globalmax}%
        Solve the global maximization problem $\mathrm{vio}\big(\hat{\BIx}^{(t)},\overline{f}_{1}^{(t)},\ldots,\overline{f}_{d}^{(t)}\big)$ in (\ref{eqn:dual-suboptimality-vio-def})
        to compute $\hat{\Blambda}^{(t)}\in\CS^*_2$ and $\overline{\varepsilon}_{\mathrm{vio}}^{(t)}\in\R$ satisfying
        $\frac{\overline{\varepsilon}_{\mathrm{vio}}^{(t)}}{1+\vartheta^{(t)}}\le \!\big\langle\BQ\hat{\BIx}^{(t)}+\Bb,\hat{\Blambda}^{(t)}\big\rangle\!+\!\sum_{i=1}^d\overline{f}_{i}^{(t)}\!\big(\langle\Bp_i,\hat{\Blambda}^{(t)}\rangle\big)\!\le$\linebreak$\mathrm{vio}\big(\hat{\BIx}^{(t)},\overline{f}_{1}^{(t)},\ldots,\overline{f}_{d}^{(t)}\big) \!\le \overline{\varepsilon}_{\mathrm{vio}}^{(t)}$. 

        \nl\label{alglin:iterative-newsupport}%
        $\CV^{(t+1)}\leftarrow \CV^{(t)}\union \widehat{\CV}^{(t)}$ where 
        $\widehat{\CV}^{(t)}\subset\CS^*_2$ is a finite set of sub-optimizers of the global maximization problem $\mathrm{vio}\big(\hat{\BIx}^{(t)},\overline{f}_{1}^{(t)},\ldots,\overline{f}_{d}^{(t)}\big)$ with $\hat{\Blambda}^{(t)}\in\widehat{\CV}^{(t)}$.

        \nl\label{alglin:iterative-suboptimality}%
        $\hat{\varepsilon}_{\mathrm{sub}}^{(t)}\leftarrow \big\langle\Bc_1,\hat{\BIx}^{(t)}\big\rangle + \big(\sum_{i=1}^d\hat{\beta}_i^{(t)}\big) - \widehat{C}_{\DRO}^{\mathrm{LB}(t)} + \overline{\varepsilon}_{\mathrm{vio}}^{(t)}$.

        \nl\label{alglin:iterative-termination}%
        \If{$\hat{\varepsilon}_{\mathrm{sub}}^{(t)}\le \varepsilon_{\mathrm{tol}}$}{
            \nl\textbf{break}.
        }
    }

    \nl\label{alglin:iterative-OTerror-forloop}%
    \For{$i=1,\ldots,d$}{
        \nl\label{alglin:iterative-OTerror}%
        $\hat{\varepsilon}_{\OT,i}\leftarrow \OT_{\mu_i}\big(\langle\Bp_i,\cdot\,\rangle\sharp\hat{\tau}^{(t)}\big) + \int_{\Omega_i}\sup_{z\in\R}\big\{z\omega - \hat{f}_i^{(t)}(z)\big\}\DIFFM{\mu_i}{\DIFF\omega} + \int_{\R}\hat{f}_i^{(t)}\DIFFX{\langle\Bp_i,\cdot\,\rangle\sharp\hat{\tau}^{(t)}}$.
    }

    \nl\label{alglin:iterative-returnval}%
    $\hat{\tau}\leftarrow\hat{\tau}^{(t)}$,
    $\hat{\Bxi}_{\ineq}\leftarrow \hat{\Bxi}_{\ineq}^{(t)}$,
    $\hat{\Bxi}_{\eq}\leftarrow \hat{\Bxi}_{\eq}^{(t)}$,
    $\hat{\BIx}\leftarrow \hat{\BIx}^{(t)}$,
    $\tilde{f}_1\leftarrow \overline{f}_1^{(t)}-\overline{\varepsilon}_{\mathrm{vio}}^{(t)}$,
    $\tilde{f}_2\leftarrow \overline{f}_2^{(t)}$,
    $\ldots$ ,
    $\tilde{f}_d\leftarrow \overline{f}_d^{(t)}$.

    \nl\label{alglin:iterative-bounds}%
    $\widehat{C}_{\DRO}^{\mathrm{LB}}\leftarrow \widehat{C}_{\DRO}^{\mathrm{LB}(t)}$,
    $\widehat{C}_{\DRO}^{\mathrm{UB}}\leftarrow \widehat{C}_{\DRO}^{\mathrm{LB}(t)}\! + \hat{\varepsilon}_{\mathrm{sub}}^{(t)}$,
    $\hat{\varepsilon}_{\mathrm{sub}}\leftarrow \hat{\varepsilon}_{\mathrm{sub}}^{(t)}$,
    $\hat{\varepsilon}_{\mathrm{prob}}\leftarrow \overline{\varepsilon}_{\mathrm{vio}}^{(t)}\! + \sum_{i=1}^{d}\hat{\varepsilon}_{\OT,i}$,
    $\hat{\varsigma}\leftarrow\varsigma^{(t)}$.\!\!

    \nl\label{alglin:iterative-support}%
    $K\leftarrow\big|\CV^{(t)}\big|$.
    Enumerate $\CV^{(t)}=\{\Blambda_1,\ldots,\Blambda_K\}$.

    \nl\label{alglin:iterative-measprepare-forloop}%
    \For{$i=1,\ldots,d$}{
        \nl\label{alglin:iterative-measprepare}%
        Compute a permutation 
        $\sigma_i:\{1,\ldots,K\}\to\{1,\ldots,K\}$ 
        such that $\langle\Bp_i,\Blambda_{\sigma_i(1)}\rangle\le \langle\Bp_i,\Blambda_{\sigma_i(2)}\rangle\le\nobreak \cdots \le \langle\Bp_i,\Blambda_{\sigma_i(K)}\rangle$.
    }

    \nl\label{alglin:iterative-probspace}%
    Let $(\mathscr{X},\mathscr{F},\PROB)$ be a probability space
    and let $J:\mathscr{X}\to\{1,\ldots,K\}$, $U_1,\ldots,U_d:\mathscr{X}\to[0,1]$\!\!\!\!\!\!\!\!\!\!\!\!\!\! \linebreak
    be jointly independent random variables 
    where $\PROB[J=j]=\hat{\tau}\big(\{\Blambda_j\}\big)$ for $j=1,\ldots,K$
    and $U_i$ is uniformly distributed on $[0,1]$ for $i=1,\ldots,d$.

    \nl\label{alglin:iterative-randomvars-forloop}%
    \For{$i=1,\ldots,d$}{
        
        \nl\label{alglin:iterative-randomvars}%
        Define $Y_i:\mathscr{X}\to\Omega_i$ by $Y_i:=G^{-}_{\mu_i}\Big(U_i\hat{\tau}\big(\{\Blambda_{J}\}\big) + \sum_{j=1}^{\sigma_i^{-1}(J)-1}\hat{\tau}\big(\{\Blambda_{\sigma_i(j)}\}\big)\Big)$, where $G^{-}_{\mu_i}(\cdot)$ is defined in (\ref{eqn:invcdf-def}).
    }

    \nl\label{alglin:iterative-meas}%
    $\hat{\mu}\leftarrow$ the law of $(Y_1,\ldots,Y_d)^\TRANSP:\mathscr{X}\to\BOmega$.

	\nl\Return $\hat{\tau}$,
    $\hat{\Bxi}_{\ineq}$, $\hat{\Bxi}_{\eq}$,
    $\hat{\BIx}$,
    $\tilde{f}_1$, $\ldots$ , $\tilde{f}_d$,
    $\hat{\mu}$,
    $\widehat{C}_{\DRO}^{\mathrm{LB}}$,
    $\widehat{C}_{\DRO}^{\mathrm{UB}}$,
    $\hat{\varepsilon}_{\mathrm{sub}}$,
    $\hat{\varepsilon}_{\mathrm{prob}}$.%
\end{algorithm}%

\begin{setting}[Algorithm~\ref{alg:iterative}]\label{sett:iterative}
    We make the following assumptions about the inputs of Algorithm~\ref{alg:iterative}.
    \begin{itemize}[beginpenalty=10000]
        \item 
        $\BOmega$, $(\mu_i)_{i=1:d}$, $\Bc_1$, $\BL_{\ineq}$, $\Bq_{\ineq}$, $\BL_{\eq}$, $\Bq_{\eq}$, $\BQ$, $\BP$, $\Bb$, $\BA_{\ineq}$, $\BA_{\eq}$, $\Bc_2$
        specify the two-stage DRO problem \eqref{eqn:dro} in Setting~\ref{sett:dro-alt}.
        We assume in addition that $C_{\DRO}>-\infty$.
        
        \item 
        For each $t\in\N$,
        $\varepsilon_{\OT}^{(t)}>0$ is a user-specified tolerance parameter
        which controls the termination of Algorithm~\ref{alg:cuttingplane-LSIP}; see Setting~\ref{sett:cuttingplane-LSIP} and Proposition~\ref{prop:cuttingplane-LSIP}.
        A small value of~$\varepsilon_{\OT}^{(t)}$ increases the number of iterations needed for Algorithm~\ref{alg:cuttingplane-LSIP} to terminate and slows down Line~\ref{alglin:iterative-LSIP} of Algorithm~\ref{alg:iterative}.
        In practice, however, we have observed that Line~\ref{alglin:iterative-LSIP} remains efficient even when $\varepsilon_{\OT}^{(t)}$ is set to be close to~0;
        see our discussions about the numerical results in Section~\ref{ssec:experiments-discussions}.

        \item
        For each $t\in\N$,
        $\varsigma^{(t)}\ge 0$ is a user-specified tolerance parameter 
        controlling the supremum difference between 
        $\hat{f}_i^{(t)}$ and its upper approximation $\overline{f}_i^{(t)}$ in Line~\ref{alglin:iterative-funcrelax}. 
        The purpose of this upper approximation is to find
        $\overline{f}_i^{(t)}\in\CC_{\mathrm{cvx}}(\R)$ which admits a ``simpler'' representation than $\hat{f}_i^{(t)}$ in order to reduce the difficulty of the global maximization problem in Line~\ref{alglin:iterative-globalmax}; see our detailed discussion about Line~\ref{alglin:iterative-funcrelax} below.

        \item
        For each $t\in\N$,
        $\vartheta^{(t)}\ge 0$ is a user-specified tolerance parameter
        controlling the sub-optimality of the approximate maximizer $\hat{\Blambda}^{(t)}$ of the global maximization problem $\mathrm{vio}\big(\hat{\BIx}^{(t)},\overline{f}_{1}^{(t)},\ldots,\overline{f}_{d}^{(t)}\big)$ in Line~\ref{alglin:iterative-globalmax}; see our detailed discussion about Line~\ref{alglin:iterative-globalmax} below.
        A large value of $\vartheta^{(t)}$ reduces the cost of solving the global maximization problem in Line~\ref{alglin:iterative-globalmax}.

        \item
        $\varepsilon_{\mathrm{tol}}>\liminf_{t\to\infty} d\big(\varepsilon_{\OT}^{(t)}+\big(1+\vartheta^{(t)}\big)\varsigma^{(t)}\big)$ is a user-specified tolerance parameter
        which controls the sub-optimality of the computed approximate optimizers of \eqref{eqn:dro-primal}, \eqref{eqn:dro-dual}, and \eqref{eqn:dro}; see \mbox{Theorem~\ref{thm:iterative}\ref{thms:iterative-bounds}--\ref{thms:iterative-dro-worstmeas}}.
        In particular, 
        by setting 
        $\varepsilon_{\OT}^{(t)}=\varepsilon_{\OT}^{(1)}$,
        $\varsigma^{(t)}=\varsigma^{(1)}$,
        $\vartheta^{(t)}=\vartheta^{(1)}$ for all $t\in\N$ and 
        by requiring that $\hat{\Blambda}^{(t)}\in\CV_{\extremept}$ in Line~\ref{alglin:iterative-globalmax} in each iteration of Algorithm~\ref{alg:iterative},
        one can guarantee the termination of Algorithm~\ref{alg:iterative} in at most ${|\CV_{\extremept}|+1}$ iterations;
        see
        Theorem~\ref{thm:iterative}\ref{thms:iterative-termination}.
    
    \end{itemize}
    
    In the following, we discuss the details of
    Algorithm~\ref{alg:iterative}.
    As discussed in Setting~\ref{sett:cuttingplane-LSIP},
    $\Blambda_0$, $\Bxi_{\ineq,0}$, $\Bxi_{\eq,0}$ in Line~\ref{alglin:iterative-initialization} can be found by computing a feasible solution of an LP problem.
    In particular, this LP problem involves the constraints describing the polyhedral convex set~$\CS^*_2$ specified by (\ref{eqn:stage2-dual-feasible-set}) with respect to $\BA_{\ineq}$, $\BA_{\eq}$, and $\Bc_2$ in the inputs of Algorithm~\ref{alg:iterative}.
    Due to Line~\ref{alglin:iterative-initial-support} and Line~\ref{alglin:iterative-newsupport}, it holds that $\Blambda_0\in\CV^{(t)}$ for every iteration~$t\in\N$.
    Subsequently, 
    Theorem~\ref{thm:primaldual}\ref{thms:primaldual-boundedness}
    guarantees that the inputs of Algorithm~\ref{alg:cuttingplane-LSIP} in Line~\ref{alglin:iterative-LSIP} of Algorithm~\ref{alg:iterative} satisfy the assumptions in Setting~\ref{sett:cuttingplane-LSIP}.

    The motivation of Line~\ref{alglin:iterative-funcrelax} is explained as follows.
    Recall that Proposition~\ref{prop:cuttingplane-LSIP}\ref{props:cuttingplane-LSIP-primal} guarantees that $\hat{f}_1^{(t)},\ldots,\hat{f}_d^{(t)}$ are continuous piece-wise affine functions.
    For $i=1,\ldots,d$, let $m_i\in\N \intersection [3,\infty)$ and $-\infty<\kappa_{i,1}<\kappa_{i,2}<\cdots<\kappa_{i,m_i}<\infty$ be such that 
    $\hat{f}_i^{(t)}$ is piece-wise affine on 
    $(-\infty,\kappa_{i,1}],\allowbreak[\kappa_{i,1},\kappa_{i,2}],\ldots,\allowbreak[\kappa_{i,m_i-1},\kappa_{i,m_i}],\allowbreak[\kappa_{i,m_i},\infty)$.
    When the piece-wise representations of $\hat{f}_1^{(t)},\ldots,\hat{f}_d^{(t)}$ consist of a large number of intervals, 
    it becomes challenging to solve the global maximization problem $\mathrm{vio}\big(\hat{\BIx}^{(t)},\hat{f}_{1}^{(t)},\ldots,\hat{f}_{d}^{(t)}\big)$; 
    see the mixed-integer programming (MIP) formulation of $\mathrm{vio}\big(\hat{\BIx}^{(t)},\hat{f}_{1}^{(t)},\ldots,\hat{f}_{d}^{(t)}\big)$
    in Proposition~\ref{prop:globalmax-MILP}, which contains $\sum_{i=1}^d(m_i-2)$ binary-valued auxiliary variables and incurs up to $\prod_{i=1}^{d}(m_i-1)$ distinct feasible combinations of values of these auxiliary variables.
    Due to this reason, 
    Line~\ref{alglin:iterative-funcrelax} computes an upper approximation $\overline{f}_i^{(t)}\in\CC_{\mathrm{cvx}}(\R)$ of $\hat{f}_i^{(t)}$ such that $\sup_{z\in\R}\big\{\overline{f}_i^{(t)}(z) - \hat{f}_{i}^{(t)}(z)\big\}\le \varsigma^{(t)}$, 
    where $\overline{f}_i^{(t)}$ is continuous piece-wise affine with fewer intervals.
    A particular choice is to 
    let $\tilde{m}_i\in\N\intersection[3,\infty)$, 
    $\tilde{\kappa}_{i,1}<\tilde{\kappa}_{i,2}<\cdots<\tilde{\kappa}_{i,\tilde{m}_i}$ be such that 
    $\{\tilde{\kappa}_{i,1},\ldots,\tilde{\kappa}_{i,\tilde{m}_i}\}\subseteq \{\kappa_{i,1},\ldots,\kappa_{i,m_i}\}$,
    and such that the unique continuous function $\overline{f}_i^{(t)}$ that is piece-wise affine on $(-\infty,\tilde{\kappa}_{i,1}],\allowbreak[\tilde{\kappa}_{i,1},\tilde{\kappa}_{i,2}],\ldots,[\tilde{\kappa}_{i,\tilde{m}_i-1},\tilde{\kappa}_{i,\tilde{m}_i}],\allowbreak[\tilde{\kappa}_{i,\tilde{m}_i},\infty)$ with $\overline{f}_i^{(t)}(\tilde{\kappa}_{i,j})=\hat{f}_i^{(t)}(\tilde{\kappa}_{i,j})$ for $j=1,\ldots,\tilde{m}_i$,
    $\partial_{-}\overline{f}_i^{(t)}(\tilde{\kappa}_{i,1})=\partial_{-}\hat{f}_i^{(t)}(\kappa_{i,1})$,
    $\partial_{+}\overline{f}_i^{(t)}(\tilde{\kappa}_{i,\tilde{m}_i})=\partial_{+}\hat{f}_i^{(t)}(\kappa_{i,m_i})$ 
    satisfies $\sup_{z\in\R}\big\{\overline{f}_i^{(t)}(z) - \hat{f}_{i}^{(t)}(z)\big\}\le \varsigma^{(t)}$.
    If $\,\prod_{i=1}^{d}(\tilde{m}_i-1)$ is significantly smaller than $\prod_{i=1}^{d}(m_i-1)$, 
    then it is more computationally efficient to solve $\mathrm{vio}\big(\hat{\BIx}^{(t)},\overline{f}_{1}^{(t)},\ldots,\overline{f}_{d}^{(t)}\big)$ 
    compared to $\mathrm{vio}\big(\hat{\BIx}^{(t)},\hat{f}_{1}^{(t)},\ldots,\hat{f}_{d}^{(t)}\big)$.

    Line~\ref{alglin:iterative-globalmax} approximately solves 
    $\mathrm{vio}\big(\hat{\BIx}^{(t)},\overline{f}_{1}^{(t)},\ldots,\overline{f}_{d}^{(t)}\big)$,
    which corresponds to the maximization of a convex function over a polyhedral convex set, and can be done via, for example, a suitably chosen \emph{branch-and-bound} algorithm \citep{benson1990separable,benson1991branch}.
    In a typical branch-and-bound algorithm for globally maximizing a convex function over a convex set,
    one constructs a search tree to divide the underlying decision space such that each node in the search tree represents a polyhedral convex subset of the decision space. 
    At each node, the algorithm finds an affine function which is a local upper bound for the convex objective function, and subsequently solves the resulting local LP relaxation.
    The maximum of the optimal values of the local LP relaxations at the leaf nodes of the search tree then provides an upper bound for the global maximization problem,
    whereas the local LP maximizers encountered in the algorithm with the largest objective (that is called the \emph{incumbent solution})
    provides a lower bound for the global maximization problem.
    The algorithm iteratively subdivides the leaf nodes into further sub-nodes until the relative difference between the upper and lower bounds falls below a user-specified tolerance threshold. 
    In Line~\ref{alglin:iterative-globalmax}, 
    $\vartheta^{(t)}\ge0$ corresponds to the tolerance threshold for the branch-and-bound algorithm, 
    while
    $\hat{\Blambda}^{(t)}\in\CS^*_2$ and $\overline{\varepsilon}_{\mathrm{vio}}^{(t)}\in\R$ are the incumbent solution and the best upper bound when the algorithm terminates.
    The termination condition requires that
    $\frac{\overline{\varepsilon}_{\mathrm{vio}}^{(t)}}{1+\vartheta^{(t)}}\le \big\langle\BQ\hat{\BIx}^{(t)}+\Bb,\hat{\Blambda}^{(t)}\big\rangle+\sum_{i=1}^d\overline{f}_{i}^{(t)}\big(\langle\Bp_i,\hat{\Blambda}^{(t)}\rangle\big)\le\mathrm{vio}\big(\hat{\BIx}^{(t)},\overline{f}_{1}^{(t)},\ldots,\overline{f}_{d}^{(t)}\big) \le \overline{\varepsilon}_{\mathrm{vio}}^{(t)}$,
    i.e., the upper bound differs from the lower bound at most by a factor of $1+\vartheta^{(t)}$.
    Proposition~\ref{prop:globalmax-MILP} shows that the global maximization problem $\mathrm{vio}\big(\hat{\BIx}^{(t)},\overline{f}_{1}^{(t)},\ldots,\overline{f}_{d}^{(t)}\big)$ can be formulated into an MIP problem and subsequently solved by state-of-the-art MIP solvers utilizing specialized branch-and-bound algorithms such as Gurobi \citep{gurobi}.\footnote{See an overview of MIP and the branch-and-bound algorithm provided by \citet{gurobi}: \url{https://www.gurobi.com/resources/mixed-integer-programming-mip-a-primer-on-the-basics/}; accessed: July~4, 2025.}
    In particular,
    the MIP formulation in Proposition~\ref{prop:globalmax-MILP} possesses the property that any search tree formed by branching on the binary-valued auxiliary variables results in local LP relaxations 
    where the incumbent solution can be attained in $\CV_{\extremept}$.
    Thus, the assumption $\hat{\Blambda}^{(t)}\in\CV_{\extremept}$ in Theorem~\ref{thm:iterative} can be satisfied in practice when
    using the MIP formulation in Proposition~\ref{prop:globalmax-MILP}.
    We would like to remark that there are branch-and-bound algorithms specialized for maximizing non-concave separable piece-wise affine functions over convex sets such as
    \citep{hubner2025spatial}.
    However, the specialized algorithm in \citep{hubner2025spatial} only demonstrated superior performance over the MIP-based approaches when the number of breakpoints in each univariate function exceeds 100, i.e., when $\tilde{m}_i\ge 100$ for $i=1,\ldots,d$.
    Since the number of breakpoints in $\overline{f}_{1}^{(t)},\ldots,\overline{f}_{d}^{(t)}$ that we have encountered in our numerical experiments are well below~100 (see Section~\ref{sec:experiments}), the MIP formulation in Proposition~\ref{prop:globalmax-MILP} suffices for solving $\mathrm{vio}\big(\hat{\BIx}^{(t)},\overline{f}_{1}^{(t)},\ldots,\overline{f}_{d}^{(t)}\big)$ in our algorithm.

    Observe that the tolerance parameters $\varepsilon_{\OT}^{(t)},\varsigma^{(t)},\vartheta^{(t)}$ in Lines~\ref{alglin:iterative-LSIP}, \ref{alglin:iterative-funcrelax}, \ref{alglin:iterative-globalmax}
    are allowed to vary for each iteration~$t$.
    The purpose is to speed up early iterations of the for-loop in Line~\ref{alglin:iterative-forloop}
    in order to improve the overall efficiency of Algorithm~\ref{alg:iterative}.
    In our numerical experiments in Section~\ref{sec:experiments}, we fix the values of $\varepsilon_{\OT}^{(t)}$ to a constant close to~0 for all $t$,
    and we decrease the values of $\varsigma^{(t)}$ and $\vartheta^{(t)}$ gradually throughout the iterations.
    In particular, we have observed that the computational time of Algorithm~\ref{alg:iterative} is dominated by the global maximization problem in Line~\ref{alglin:iterative-globalmax}, which is sensitive to the parameter $\vartheta^{(t)}$.
    Because of this, 
    we set the value of $\vartheta^{(t)}$ to be large for early iterations and only decrease it in the final few iterations.
    
    In Lines~\ref{alglin:iterative-probspace}--\ref{alglin:iterative-meas},
    a probability measure $\hat{\mu}\in\CP(\BOmega)$ is constructed as the law of a random variable $(Y_1,\ldots,Y_d)^\TRANSP$ on a probability space $(\mathscr{X},\mathscr{F},\PROB)$ via the procedure in Corollary~\ref{cor:dual-suboptimality-primal-coupling}.
    Such a characterization allows one to efficiently generate samples from $\hat{\mu}$ as follows. 
    Given a sample size $M\in\N$, 
    one generates independent and identically distributed (i.i.d.\@) samples $\big\{J^{[m]}\big\}_{m=1:M}$ 
    from the categorical distribution with probabilities 
    $\PROB[J=j]=\hat{\tau}\big(\{\Blambda_j\}\big)$ for $j=1,\ldots,K$
    as well as i.i.d.\@ samples 
    $\big\{U_{i}^{[m]}\big\}_{i=1:d,\,m=1:M}$
    from the uniform distribution on $[0,1]$
    that are independent of $\big\{J^{[m]}\big\}_{m=1:M}$,
    and then computes $Y_i^{[m]}:=G^{-}_{\mu_i}\Big(U_i^{[m]}\hat{\tau}\big(\{\Blambda_{J^{[m]}}\}\big) + \sum_{j=1}^{\sigma_i^{-1}(J^{[m]})-1}\hat{\tau}\big(\{\Blambda_{\sigma_i(j)}\}\big)\Big)$ for $i=1,\ldots,d$, $m=1,\ldots,M$.
    For sufficiently large sample size~$M$, one can approximate $\hat{\mu}$ by the empirical measure
    $\frac{1}{M}\sum_{m=1}^M \delta_{(Y_1^{[m]},\ldots,Y_{d}^{[m]})^\TRANSP}\approx \hat{\mu}$.

    Line~\ref{alglin:iterative-bounds} computes two sub-optimality estimates $\hat{\varepsilon}_{\mathrm{sub}}$ and $\hat{\varepsilon}_{\mathrm{prob}}$ satisfying 
    $\hat{\varepsilon}_{\mathrm{sub}}\le \varepsilon_{\mathrm{tol}}$
    and 
    $\hat{\varepsilon}_{\mathrm{sub}}\le \hat{\varepsilon}_{\mathrm{prob}}\le \hat{\varepsilon}_{\mathrm{sub}} + d\hat{\varsigma}$; see Theorem~\ref{thm:iterative}\ref{thms:iterative-bounds}.
    Even though we are able to show in Theorem~\ref{thm:iterative}\ref{thms:iterative-dual}--\ref{thms:iterative-dro-worstmeas} that
    $(\hat{\tau},\hat{\Bxi}_{\ineq},\hat{\Bxi}_{\eq})$ is an $\hat{\varepsilon}_{\mathrm{sub}}$-optimizer of \eqref{eqn:dro-dual},
    $(\hat{\BIx},\tilde{f}_1,\ldots,\tilde{f}_d)$ is an $\hat{\varepsilon}_{\mathrm{sub}}$-optimizer of \eqref{eqn:dro-primal},
    and $\hat{\BIx}$ is an $\hat{\varepsilon}_{\mathrm{sub}}$-optimizer of \eqref{eqn:dro},
    we are only able to show that 
    $\hat{\mu}$ is an $\hat{\varepsilon}_{\mathrm{prob}}$-optimizer of the inner maximization problem in \eqref{eqn:dro} with respect to $\hat{\BIx}$.
    Nevertheless, both sub-optimality estimates can be controlled to be arbitrarily close to~0 with appropriate choices of the inputs $\varepsilon_{\mathrm{tol}}$ and $\big(\varsigma^{(t)}\big)_{t\in\N}$ of Algorithm~\ref{alg:iterative}; see Theorem~\ref{thm:iterative}\ref{thms:iterative-bounds}.
\end{setting}

\begin{theorem}[Properties of Algorithm~\ref{alg:iterative}]\label{thm:iterative}
    Let the inputs of Algorithm~\ref{alg:iterative} be specified according to
    Setting~\ref{sett:iterative}.
    \begin{enumerate}[label=(\roman*),beginpenalty=10000]
        \item\label{thms:iterative-globalmax}
        Then, in each iteration of Algorithm~\ref{alg:iterative}, 
        the global maximization problem in Line~\ref{alglin:iterative-globalmax} is bounded from above
        and a maximizer $\hat{\Blambda}^{(t)}$ can be attained in $\CV_{\extremept}$.
        
    \end{enumerate}
    In the following, let us assume in addition that
    $\hat{\Blambda}^{(t)}$ computed in Line~\ref{alglin:iterative-globalmax} satisfies $\hat{\Blambda}^{(t)}\in\CV_{\extremept}$ in each iteration of Algorithm~\ref{alg:iterative}.
    Then, the following statements hold.
    \begin{enumerate}[label=(\roman*),beginpenalty=10000]
        \setcounter{enumi}{1}

        \item\label{thms:iterative-termination}
        Algorithm~\ref{alg:iterative} terminates after finitely many iterations.
        Moreover, if 
        $\varepsilon_{\OT}^{(t)}=\nobreak\varepsilon_{\OT}^{(1)}$,
        $\varsigma^{(t)}=\varsigma^{(1)}$,
        $\vartheta^{(t)}=\vartheta^{(1)}$ for all $t\in\N$,
        then Algorithm~\ref{alg:iterative} terminates after
        at most
        $|\CV_{\extremept}|+\nobreak1$ iterations.
        
        \item\label{thms:iterative-bounds}
        It holds that
        $\widehat{C}_{\DRO}^{\mathrm{LB}}\le \eqref{eqn:dro} \le \widehat{C}_{\DRO}^{\mathrm{UB}}$, 
        $\hat{\varepsilon}_{\mathrm{sub}}=\widehat{C}_{\DRO}^{\mathrm{UB}}-\widehat{C}_{\DRO}^{\mathrm{LB}}\le \varepsilon_{\mathrm{tol}}$,
        and
        $\hat{\varepsilon}_{\mathrm{sub}}\le \hat{\varepsilon}_{\mathrm{prob}}\le\hat{\varepsilon}_{\mathrm{sub}}+d\hat{\varsigma}\le\varepsilon_{\mathrm{tol}}+d\sup_{t\in\N}\big\{\varsigma^{(t)}\big\}$.
        In particular, both $\hat{\varepsilon}_{\mathrm{sub}}$ and $\hat{\varepsilon}_{\mathrm{prob}}$ can be controlled to be arbitrarily close to~0 by setting the inputs $\varepsilon_{\mathrm{tol}}$ and $\big(\varsigma^{(t)}\big)_{t\in\N}$ sufficiently small.
        
        \item\label{thms:iterative-dual}
        $(\hat{\tau},\hat{\Bxi}_{\ineq},\hat{\Bxi}_{\eq})$ is an $\hat{\varepsilon}_{\mathrm{sub}}$-optimizer of \eqref{eqn:dro-dual} with objective $\widehat{C}_{\DRO}^{\mathrm{LB}}$.
        
        \item\label{thms:iterative-primal}
        $(\hat{\BIx},\tilde{f}_1,\ldots,\tilde{f}_d)$ is an $\hat{\varepsilon}_{\mathrm{sub}}$-optimizer of \eqref{eqn:dro-primal} with objective $\widehat{C}_{\DRO}^{\mathrm{UB}}$.
        
        \item\label{thms:iterative-dro-worstmeas}
        $\hat{\BIx}$ is an $\hat{\varepsilon}_{\mathrm{sub}}$-optimizer of \eqref{eqn:dro} and $\hat{\mu}$ is an $\hat{\varepsilon}_{\mathrm{prob}}$-optimizer of the inner maximization problem in \eqref{eqn:dro} with respect to $\hat{\BIx}$,
        i.e., $\hat{\mu}$ is an $\hat{\varepsilon}_{\mathrm{prob}}$-approximate worst-case probability measure with respect to the first-stage decision~$\hat{\BIx}$.
    \end{enumerate}
\end{theorem}

Since $\overline{f}_1^{(t)},\ldots,\overline{f}_d^{(t)}$ computed in Line~\ref{alglin:iterative-funcrelax} of Algorithm~\ref{alg:iterative} are continuous piece-wise affine functions,
the global maximization problem $\mathrm{vio}\big(\hat{\BIx}^{(t)},\overline{f}_1^{(t)},\ldots,\overline{f}_{d}^{(t)}\big)$ solved in Line~\ref{alglin:iterative-globalmax} of Algorithm~\ref{alg:iterative} can be formulated into a linear mixed-integer programming (MIP) problem via the so-called incremental model; see, e.g.,
\citep[Section~3.4]{vielma2010mixed}.
We provide the details in the following proposition.
\begin{proposition}[MIP formulation of $\mathrm{vio}(\BIx,f_1,\ldots,f_d)$]\label{prop:globalmax-MILP}
    Under Setting~\ref{sett:dro-alt}, let $\BIx\in\nobreak\CS_1$, and 
    for $i=1,\ldots,d$, let $m_i\in\N\intersection [3,\infty)$, $-\infty<\kappa_{i,1}<\kappa_{i,2}<\cdots<\kappa_{i,m_i}<\infty$,
    let $f_i\in\CC_{\mathrm{cvx}}(\R)$ be piece-wise affine 
    on $(-\infty,\kappa_{i,1}],\allowbreak[\kappa_{i,1},\kappa_{i,2}],\ldots,[\kappa_{i,m_i-1},\kappa_{i,m_i}],\allowbreak[\kappa_{i,m_i},\infty)$ 
    satisfying $\partial_{-}f_i(\kappa_{i,1})=\nobreak\underline{\omega}_i$, 
    $\partial_{+}f_i(\kappa_{i,m_i})=\overline{\omega}_i$.
    Let us consider the following linear MIP problem:
    \begin{align}
        \begin{split}
            \maximize_{\!\!\!\!\substack{\Blambda,\,(\zeta_{i,j}),\,(\iota_{i,j})}} \quad & \langle\BQ\BIx+\Bb,\Blambda\rangle + \sum_{i=1}^d \Bigg(f_i(\kappa_{i,1}) + \sum_{j=1}^{m_i-1}\big(f_i(\kappa_{i,j+1}) - f_i(\kappa_{i,j})\big)\zeta_{i,j}\Bigg) \\
            \mathrm{subject~to} \quad & \Blambda\in\CS^*_2, \mathrm{~and~for~}i=1,\ldots,d:\\
            & \begin{cases}
                \zeta_{i,j} \in \R_+ & \forall 1\le j\le m_i-1,\\
                \iota_{i,j} \in \{0, 1\} & \forall 1\le j\le m_i-2, \\
                \zeta_{i,1} \le 1, \\
                \zeta_{i,j+1} \le \iota_{i,j} \le \zeta_{i,j} & \forall 1\le j\le m_i-2, \\
                \kappa_{i,1} + \sum_{j=1}^{m_i-1} (\kappa_{i,j+1}-\kappa_{i,j})\zeta_{i,j} = \langle\Bp_i,\Blambda\rangle.
            \end{cases}
        \end{split}
        \label{eqn:globalmax-MILP}
    \end{align}
    If we assume that $\kappa_{i,1}\le \min_{\Blambda\in\CV_{\extremept}}\big\{\langle\Bp_i,\Blambda\rangle\big\}$ and $\kappa_{i,m_i}\ge \max_{\Blambda\in\CV_{\extremept}}\big\{\langle\Bp_i,\Blambda\rangle\big\}$ for ${i=1,\ldots,d}$,
    then (\ref{eqn:globalmax-MILP}) is an equivalent formulation of the global maximization problem $\mathrm{vio}(\BIx,f_1,\ldots,f_d)$ in (\ref{eqn:dual-suboptimality-vio-def}) and the following statements hold.
    \begin{enumerate}[label=(\roman*), beginpenalty=10000]
        \item\label{props:globalmax-MILP-optimalval}
        The optimal values of (\ref{eqn:globalmax-MILP}) and (\ref{eqn:dual-suboptimality-vio-def}) are equal.
        
        \item\label{props:globalmax-MILP-suboptimizer}
        For any $\varepsilon>0$ and any $\varepsilon$-optimizer $\big(\Blambda_{\varepsilon},(\zeta_{i,j,\varepsilon})_{j=1:m_i-1,\,i=1:d}, (\iota_{i,j,\varepsilon})_{j=1:m_i-2,\,i=1:d}\big)$ of (\ref{eqn:globalmax-MILP}),
        $\Blambda_{\varepsilon}$ is an $\varepsilon$-optimizer of (\ref{eqn:dual-suboptimality-vio-def}).

        \item\label{props:globalmax-MILP-extremept}
        There exists an optimizer 
        $\big(\Blambda^\star,(\zeta^\star_{i,j})_{j=1:m_i-1,\,i=1:d}, (\iota^\star_{i,j})_{j=1:m_i-2,\,i=1:d}\big)$
        of (\ref{eqn:globalmax-MILP}) where $\Blambda^\star\in\nobreak\CV_{\extremept}$.
        
        \item\label{props:globalmax-MILP-optimizer}
        For every optimizer 
        $\big(\Blambda^\star,(\zeta^\star_{i,j})_{j=1:m_i-1,\,i=1:d},\allowbreak (\iota^\star_{i,j})_{j=1:m_i-2,\,i=1:d}\big)$ of (\ref{eqn:globalmax-MILP}),
        $\Blambda^\star$ is an optimizer of (\ref{eqn:dual-suboptimality-vio-def}).

    \end{enumerate}
\end{proposition}

In our implementation of Algorithm~\ref{alg:iterative},
we adopt the MIP formulation of the global maximization problem $\mathrm{vio}\big(\hat{\BIx}^{(t)},\overline{f}_1^{(t)},\ldots,\overline{f}_{d}^{(t)}\big)$ in Proposition~\ref{prop:globalmax-MILP}
and numerically solve it by the MIP solver of Gurobi \citep{gurobi}.
Notice that (\ref{eqn:globalmax-MILP}) contains $\sum_{i=1}^{d}(m_i-2)$ binary-valued auxiliary variables, which can take up to $\prod_{i=1}^{d}(m_i-1)$ distinct feasible combinations of values. 
We remark that there exist various alternative MIP formulations of $\mathrm{vio}(\BIx,f_1,\ldots,f_d)$, some involving logarithmically many binary-valued auxiliary variables; see, e.g., \citep[Sections~3.1--3.3]{vielma2010mixed}.
We have chosen the formulation in Proposition~\ref{prop:globalmax-MILP} due to its simplicity and its good empirical performance.
In fact, in our numerical experiments in Section~\ref{sec:experiments},
the function $\overline{f}_i^{(t)}$ computed by Line~\ref{alglin:iterative-funcrelax} of Algorithm~\ref{alg:iterative} rarely has more than 10 breakpoints.
This means that MIP formulations with logarithmically many binary-valued auxiliary variables or specialized branch-and-bound algorithms such as \citep{hubner2025spatial} are unlikely to provide better computational efficiency than the incremental model in Proposition~\ref{prop:globalmax-MILP} for the problem instances under our settings; 
see, e.g., \citep[Section~6]{hubner2025spatial}.

\section{Numerical experiments}
\label{sec:experiments}

In this section, we perform three numerical experiments to empirically test the performance of Algorithm~\ref{alg:iterative}.
The code used in the experiments is available at: \url{https://github.com/qikunxiang/TwoStageDROwMarginals}; please refer to the code for further details.

\subsection{Experiment~1: task scheduling}
\label{ssec:experiments-scheduling}
\begin{figure}[t]
\centering
\includegraphics[width=0.48\linewidth]{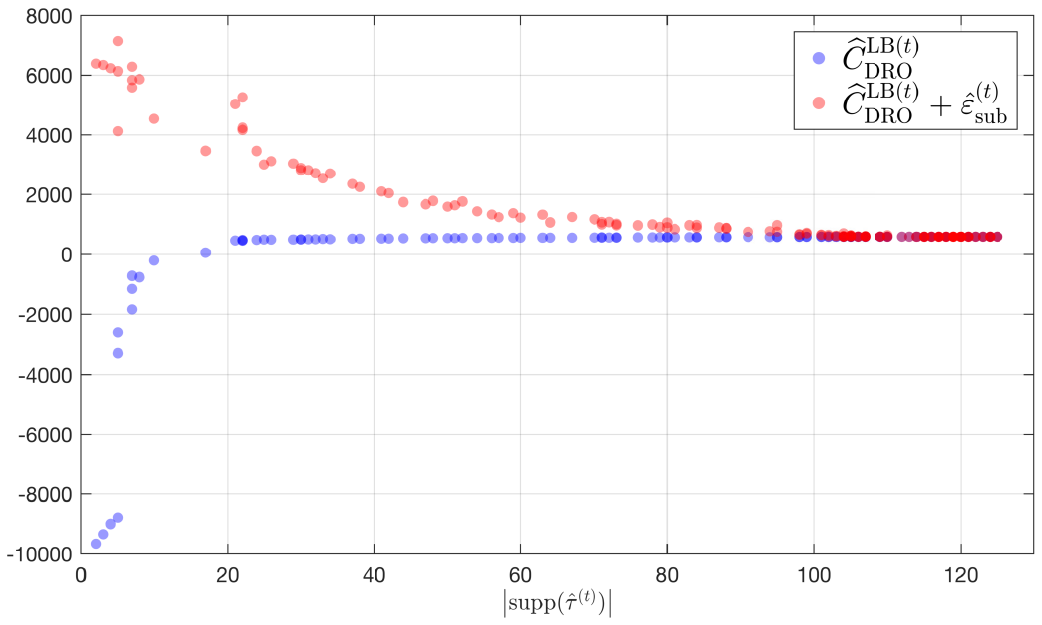}
~
\includegraphics[width=0.48\linewidth]{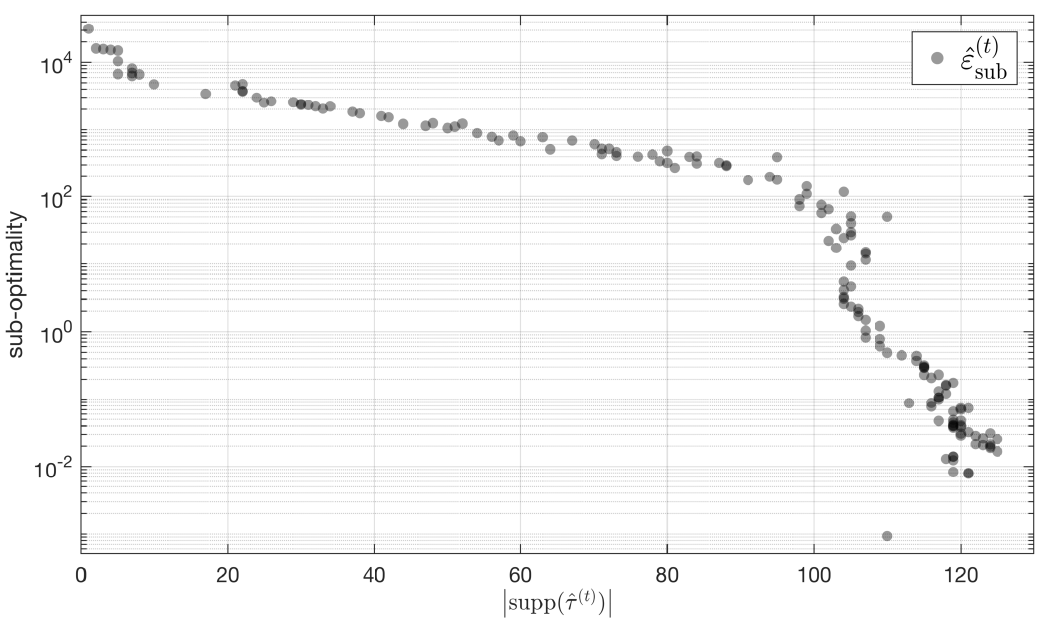}

\vspace{2pt}

\includegraphics[width=0.48\linewidth]{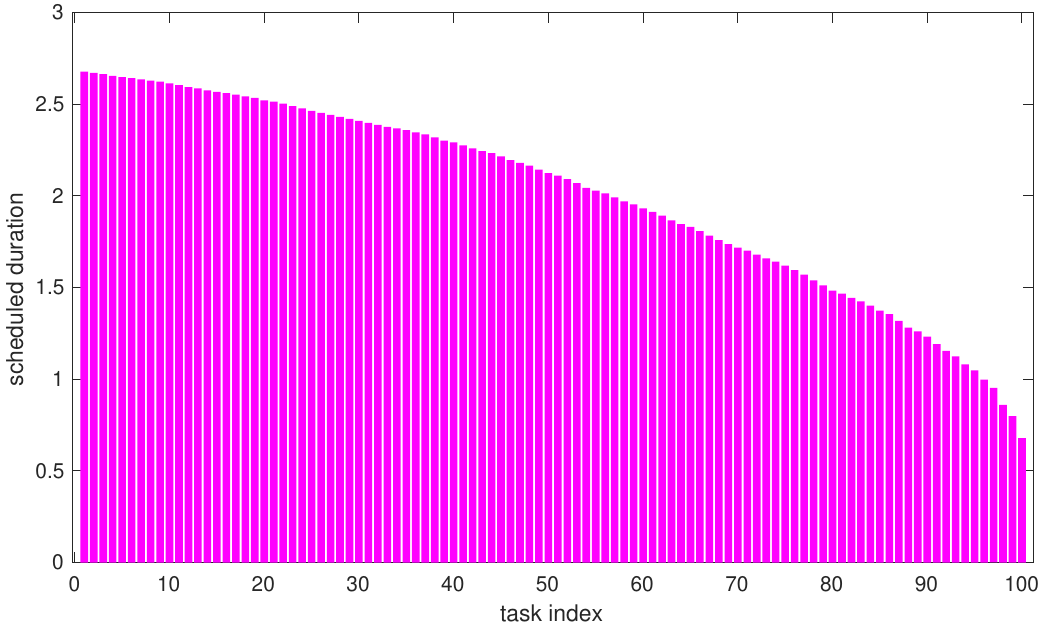}
\caption{\textbf{Experiment~1.} 
\textit{Top-left:} 
the computed lower bounds $\widehat{C}_{\DRO}^{\mathrm{LB}(t)}$ and upper bounds $\widehat{C}_{\DRO}^{\mathrm{LB}(t)}+\hat{\varepsilon}_{\mathrm{sub}}^{(t)}$ plotted against $\big|\support(\hat{\tau}^{(t)})\big|$. 
\textit{Top-right:} 
the computed sub-optimality estimates $\hat{\varepsilon}_{\mathrm{sub}}^{(t)}$ plotted against $\big|\support(\hat{\tau}^{(t)})\big|$. 
\textit{Bottom:} the computed approximately optimal scheduled task durations $\hat{\BIx}$.}%
\label{fig:scheduling-errors}%
\end{figure}

\begin{figure}[t]
\centering%
\includegraphics[width=0.80\linewidth]{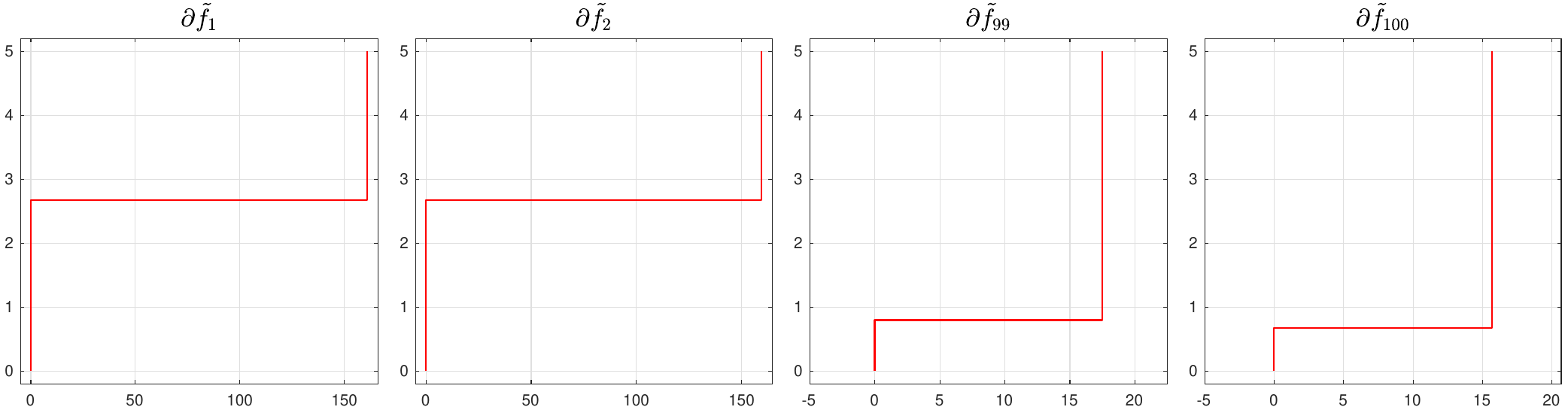}
\caption{\textbf{Experiment~1.} 
The graphs of the computed subdifferential mappings $\partial\tilde{f}_1,\partial\tilde{f}_2,\allowbreak\partial\tilde{f}_{99},\partial\tilde{f}_{100}$.}
\label{fig:scheduling-potentials}
\end{figure}

\begin{figure}[t]
\centering
\includegraphics[width=0.58\linewidth]{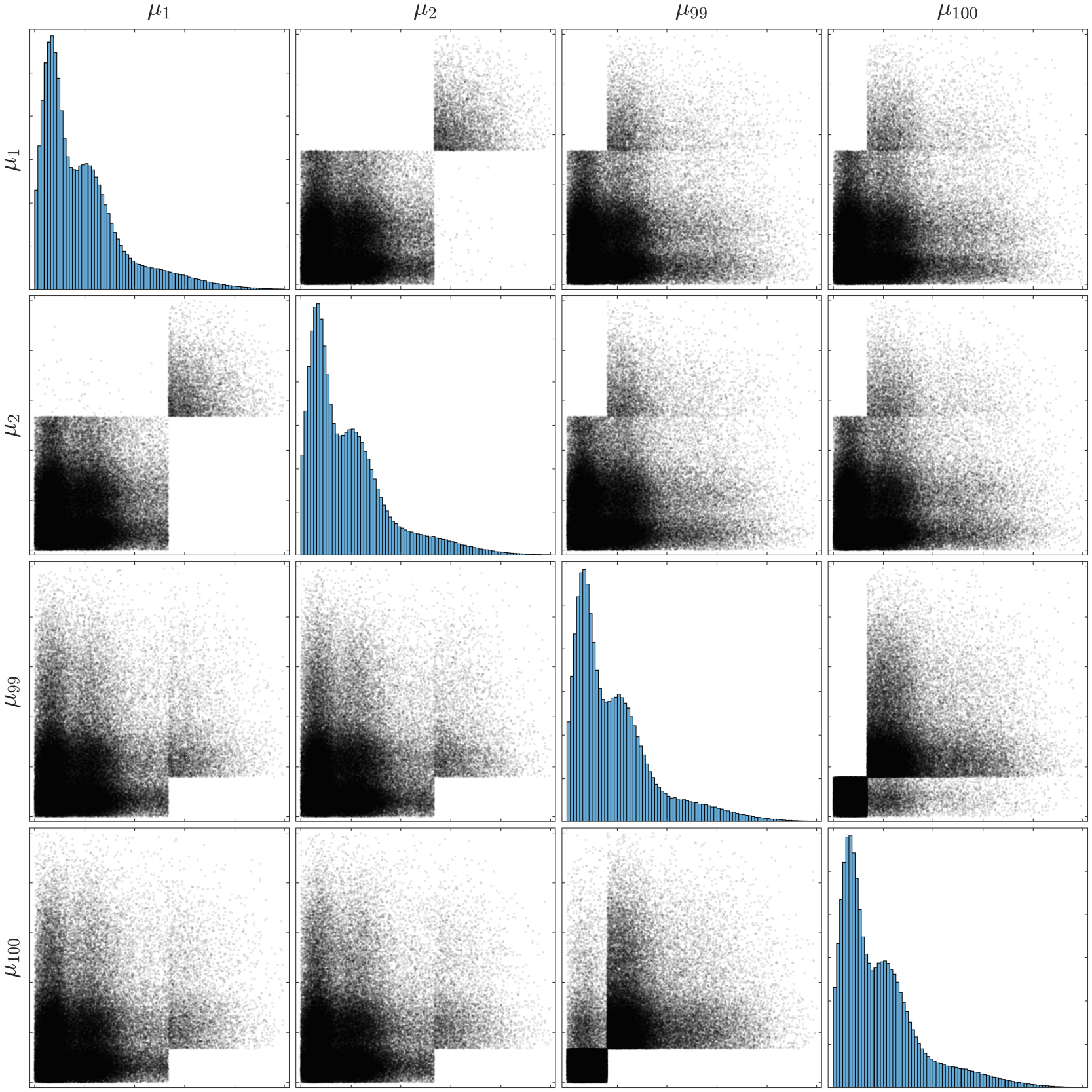}
\caption{\textbf{Experiment~1.} 
Pair-wise scatter plots of 4 marginals of the computed $\hat{\varepsilon}_{\mathrm{prob}}$-approximate worst-case probability measure $\hat{\mu}$ with respect to $\hat{\BIx}$.}
\label{fig:scheduling-scatter}
\end{figure}

In this experiment, we solve the distributionally robust task scheduling problem in Example~\ref{exp:scheduling}, where we consider the scheduling of $d=100$ tasks in a time window $[0,T]$ where $T=200$. 
We let $\mu_1=\cdots=\mu_d\in\CP([0,5])$, where $\mu_1$ is a mixture of three Gaussian measures truncated to $[0,5]$. 
The mean of $\mu_1$ is equal to~$0.9554$ and histograms of samples from $\mu_1$ can be seen in the diagonal panels of Figure~\ref{fig:scheduling-scatter}.
The weights $c_{2,1},\ldots,c_{2,d}\in\R_+$ are randomly generated.
Under this setting, we apply Algorithm~\ref{alg:iterative} to compute approximate optimizers of \eqref{eqn:dro-primal}, \eqref{eqn:dro-dual}, and \eqref{eqn:dro},
where we set 
$\varepsilon_{\mathrm{tol}}=10^{-3}$,
$\varepsilon_{\OT}^{(t)}=10^{-7}$ for all $t\in\N$,
and vary the parameters $\varsigma^{(t)}$ and $\vartheta^{(t)}$ as discussed in Setting~\ref{sett:iterative}.
Algorithm~\ref{alg:iterative} terminated after 161~iterations
and computed 
the lower bound
$\widehat{C}_{\DRO}^{\mathrm{LB}}=594.3361$,
the upper bound
$\widehat{C}_{\DRO}^{\mathrm{UB}}=594.3371$,
and the sub-optimality estimates
$\hat{\varepsilon}_{\mathrm{sub}}=9.3388\times 10^{-4}$,
$\hat{\varepsilon}_{\mathrm{prob}}=9.3389\times 10^{-4}$.
The relative errors are 
$\frac{\hat{\varepsilon}_{\mathrm{sub}}}{|\widehat{C}_{\DRO}^{\mathrm{UB}}|}=1.5713\times 10^{-6}$,
$\frac{\hat{\varepsilon}_{\mathrm{prob}}}{|\widehat{C}_{\DRO}^{\mathrm{UB}}|}=1.5713\times 10^{-6}$.
These results show that the approximate optimizers of \eqref{eqn:dro-primal}, \eqref{eqn:dro-dual}, and \eqref{eqn:dro}
computed by Algorithm~\ref{alg:iterative}
are very close to optimality.
In the following, let us discuss some further results produced by Algorithm~\ref{alg:iterative} shown in 
Figures~\ref{fig:scheduling-errors}--\ref{fig:scheduling-scatter}.

\setlength{\floatsep}{4pt plus 2.0pt minus 3.0pt}

The top panels of Figure~\ref{fig:scheduling-errors} contain scatter plots showing 
the values of the lower bounds 
$\widehat{C}_{\DRO}^{\mathrm{LB}(t)}$,
the upper bounds 
$\widehat{C}_{\DRO}^{\mathrm{LB}(t)}+\hat{\varepsilon}_{\mathrm{sub}}^{(t)}$,
and the sub-optimality estimates 
$\hat{\varepsilon}_{\mathrm{sub}}^{(t)}$
computed by Algorithm~\ref{alg:iterative}
throughout the iterations,
where the horizontal axes show $\big|\support(\hat{\tau}^{(t)})\big|$.
Even though the lower bound $\widehat{C}_{\DRO}^{\mathrm{LB}(t)}$ is non-decreasing in~$t$,
the upper bound $\widehat{C}_{\DRO}^{\mathrm{LB}(t)}+\hat{\varepsilon}_{\mathrm{sub}}^{(t)}$ is not necessarily decreasing in~$t$,
and we find the results more readable when plotting 
$\widehat{C}_{\DRO}^{\mathrm{LB}(t)}$,
$\widehat{C}_{\DRO}^{\mathrm{LB}(t)}+\hat{\varepsilon}_{\mathrm{sub}}^{(t)}$, $\hat{\varepsilon}_{\mathrm{sub}}^{(t)}$
against $\big|\support(\hat{\tau}^{(t)})\big|$ rather than $t$.
The top-left panel of Figure~\ref{fig:scheduling-errors}
shows an increasing trend of $\widehat{C}_{\DRO}^{\mathrm{LB}(t)}$ against $\big|\support(\hat{\tau}^{(t)})\big|$ and shows a decreasing trend of
$\widehat{C}_{\DRO}^{\mathrm{LB}(t)}+\hat{\varepsilon}_{\mathrm{sub}}^{(t)}$ against $\big|\support(\hat{\tau}^{(t)})\big|$.
Moreover, it shows that the decrease of the upper bound is more gradual than the increase of the lower bound.
The top-right panel of Figure~\ref{fig:scheduling-errors}
shows that $\hat{\varepsilon}_{\mathrm{sub}}^{(t)}$ exhibits a steep downward trend when $\big|\support(\hat{\tau}^{(t)})\big|\ge100$.
Notably, the sub-optimality $\hat{\varepsilon}_{\mathrm{sub}}=9.3388\times 10^{-4}$ was achieved with respect to $\hat{\tau}$ where $|\support(\hat{\tau})|=110$.
This shows that in the approximate optimizer $(\hat{\tau},\hat{\Bxi}_{\ineq},\hat{\Bxi}_{\eq})$ of \eqref{eqn:dro-dual} computed by Algorithm~\ref{alg:iterative},
$\hat{\tau}$ has high sparsity relative to the problem dimension $d=k^*_2=100$.
In fact, the maximum value of $\big|\support(\hat{\tau}^{(t)})\big|$ encountered during the execution of Algorithm~\ref{alg:iterative} was~125.
The bottom panel of Figure~\ref{fig:scheduling-errors} shows the approximate optimizer $\hat{\BIx}$ of \eqref{eqn:dro} computed by Algorithm~\ref{alg:iterative}, which corresponds to the scheduled duration of the 100 tasks.
The result shows a decreasing pattern where earlier tasks are allocated longer durations compared to later tasks. 
This is due to the possibility that the delay of an early task could lead to a chain reaction causing later tasks to be further delayed.

Figure~\ref{fig:scheduling-potentials} shows the graphs of the subdifferential mappings of four computed functions $\tilde{f}_1,\allowbreak\tilde{f}_2,\allowbreak\tilde{f}_{99},\allowbreak\tilde{f}_{100}$,
which are related to the first two and the last two tasks out of the 100~tasks being scheduled.
Recall that for any
$\bar{f}_i\in\argmax_{f\in\CC_{\mathrm{cvx}}(\R)}\big\{\int_{\Omega_i} \sup_{z\in\R}\big\{z\omega-f(z)\big\} \DIFFM{\mu_i}{\DIFF\omega} + \int_{\R}f\DIFFX{\langle\Bp_i,\cdot\,\rangle\sharp\hat{\tau}}\big\}$,
the optimal coupling of $\mu_i$ and $\langle\Bp_i,\cdot\,\rangle\sharp\hat{\tau}$ with respect to the cost function $\Omega_i\times\R\ni(\omega,z)\mapsto{-z\omega}\in\R$
is necessarily concentrated on the graph of 
$\partial\bar{f}_i$; see, e.g., Lemma~\ref{lem:OT-1d}\ref{lems:OT-1d-optimality}. 
Thus, since $\tilde{f}_i$ is an approximate maximizer of the above problem (see Lemma~\ref{lem:weak-duality}), 
the graph of $\partial \tilde{f}_i$ can be seen as approximately containing the support of the optimal coupling of $\mu_i$ and $\langle\Bp_i,\cdot\,\rangle\sharp\hat{\tau}$.
It can be observed from Figure~\ref{fig:scheduling-potentials} that 
$\tilde{f}_1,\allowbreak\tilde{f}_2,\allowbreak\tilde{f}_{99},\allowbreak\tilde{f}_{100}$ are continuous piece-wise affine functions each containing 2~breakpoints (i.e., 3~pieces).
In fact, we have found that each of
$\tilde{f}_1,\ldots,\tilde{f}_{d}$ contains at most 3~breakpoints, 
and the average number of breakpoints is~2.14.

Moreover, we generate $10^5$ independent samples from the first two and the last two marginals of the $\hat{\varepsilon}_{\mathrm{prob}}$-approximate worst-case probability measure $\hat{\mu}$ with respect to $\hat{\BIx}$ computed by Algorithm~\ref{alg:iterative} and show the pair-wise scatter plots in Figure~\ref{fig:scheduling-scatter},
where the diagonal panels contain the marginal histograms. 
These scatter plots exhibit block-diagonal patterns.
We have also created pair-wise scatter plots of various other choices of marginals and observed similar block-diagonal patterns.
These patterns indicate that under $\hat{\mu}$, the actual durations of the tasks are positively correlated.
Based on these observations,
we theorize that, 
if $\bar{\BIx}=(\bar{x}_1,\ldots,\bar{x}_d)^\TRANSP$ is an optimizer of \eqref{eqn:dro} and $\bar{\mu}$ is a worst-case probability measure with respect to $\bar{\BIx}$,
then it holds that
\begin{align*}
    \bar{\mu}\Big(\big\{(\omega_1,\ldots,\omega_d)^\TRANSP\in\BOmega:\exists 1\le i<j\le d,\;\omega_i > \bar{x}_i,\; \omega_j \le \bar{x}_j\big\}\Big)=0.
\end{align*}
This means that, under a worst-case probability measure,
whenever the actual duration of a task exceeds its scheduled duration, then the actual durations of all subsequent tasks will also exceed their respective scheduled durations,
which is a sensible adversarial choice of $\bar{\mu}$.
We empirically confirm this by generating $10^6$ independent samples from the computed approximate worst-case probability measure $\hat{\mu}$,
denoting the collection of samples by $\widetilde{\BOmega}$, 
denoting $\hat{\BIx}=(\hat{x}_1,\ldots,\hat{x}_d)^\TRANSP$,
and defining
$\widetilde{\BOmega}_{\zeta}:=\big\{(\omega_1,\ldots,\omega_d)^\TRANSP\in\widetilde{\BOmega}:\exists 1\le i<j\le d,\allowbreak\omega_i > \hat{x}_i+\zeta,\;\allowbreak \omega_j \le \hat{x}_j-\zeta\big\}$ for all $\zeta\ge0$.
We have found that, out of the $10^6$ generated samples, we have
$\big|\widetilde{\BOmega}_{0}\big|=265$ and
$\big|\widetilde{\BOmega}_{10^{-4}}\big|=0$.
Therefore, $\hat{\BIx}$ and $\hat{\mu}$ computed by Algorithm~\ref{alg:iterative} approximately satisfy the theorized property.

Lastly, we would like to remark that \citep[Section~5.1]{chen2018distributionally} has shown that the distributionally robust task scheduling problem is polynomial time solvable under the additional assumptions that the weights $c_{2,1},\ldots,c_{2,d}$ are positive integers and/or the probability measures $\mu_1,\ldots,\mu_d$ have finite supports.
However, despite that we have set $\mu_1,\ldots,\mu_d$ to be non-discrete and chosen non-integer weights $c_{2,1},\ldots,c_{2,d}$ in this experiment, 
our algorithm remains computationally efficient and capable of computing an approximate optimizer of \eqref{eqn:dro} whose sub-optimality is extremely close to~0.
Moreover, we have experimented with problem instances with integer and non-integer valued weights $c_{2,1},\ldots,c_{2,d}$, and we have not observed any discernible difference in their computational time.
This suggests that it might be possible to strengthen the computational complexity results of \citet{chen2018distributionally} to more general assumptions.

\subsection{Experiment~2: multi-product assembly}
\label{ssec:experiments-assembly}

In this experiment, we solve the distributionally robust multi-product assembly problem in Example~\ref{exp:multiproduct},
where we consider a manufacturer that produces $d=50$ products with $k_1=100$ parts. 
For $i=1,\ldots,d$, 
the law of the random demand $\mu_i\in\CP\big([2,10]\big)$ is a mixture of three equally weighted Gaussian measures with randomly generated parameters truncated to $[2,10]$.
The mean values of $\mu_1,\ldots,\mu_d$ range between 
2.2955 and 7.1648, with an average of 3.8005.
The 
per-unit revenues $r_1,\ldots,r_d$ of the products, 
per-unit prices $c_{1,1},\ldots,c_{1,k_1}$ of the parts,
and the per-unit salvage values $s_1,\ldots,s_{k_1}$ of the parts 
are all randomly generated.
The matrix $\BU\in\R^{d\times k_1}$ that represents the types and amounts of parts needed for producing each unit of product is a randomly generated sparse matrix with random entries. 
The 50 products each uses an average of 19.68 parts,
and the 100 parts are each used by an average of 9.84 products.

Under this setting, we apply Algorithm~\ref{alg:iterative} to compute approximate optimizers of \eqref{eqn:dro-primal}, \eqref{eqn:dro-dual}, and \eqref{eqn:dro},
where we set 
$\varepsilon_{\mathrm{tol}}=10$,
$\varepsilon_{\OT}^{(t)}=2\times10^{-4}$ for all $t\in\N$,
and vary the parameters $\varsigma^{(t)}$ and $\vartheta^{(t)}$ as discussed in Setting~\ref{sett:iterative}.
Algorithm~\ref{alg:iterative} terminated after 10307~iterations
and computed 
the lower bound
$\widehat{C}_{\DRO}^{\mathrm{LB}}=-3719.1723$,
the upper bound
$\widehat{C}_{\DRO}^{\mathrm{UB}}=-3715.6591$,
and the sub-optimality estimates
$\hat{\varepsilon}_{\mathrm{sub}}=3.5132$,
$\hat{\varepsilon}_{\mathrm{prob}}=4.6196$.
The relative errors are 
$\frac{\hat{\varepsilon}_{\mathrm{sub}}}{|\widehat{C}_{\DRO}^{\mathrm{UB}}|}=0.0946\%$,
$\frac{\hat{\varepsilon}_{\mathrm{prob}}}{|\widehat{C}_{\DRO}^{\mathrm{UB}}|}=0.1243\%$.
These results show that Algorithm~\ref{alg:iterative} has produced high-quality approximate optimizers of \eqref{eqn:dro-primal}, \eqref{eqn:dro-dual}, and \eqref{eqn:dro} with small sub-optimality.

\begin{figure}[t]
\centering
\includegraphics[width=0.49\linewidth]{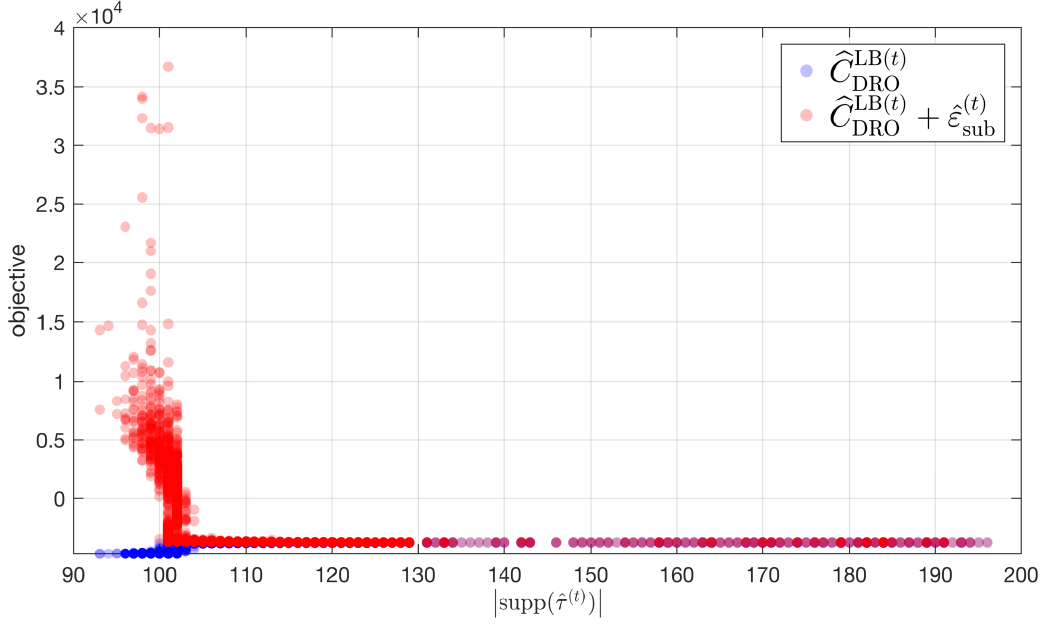}~
\includegraphics[width=0.49\linewidth]{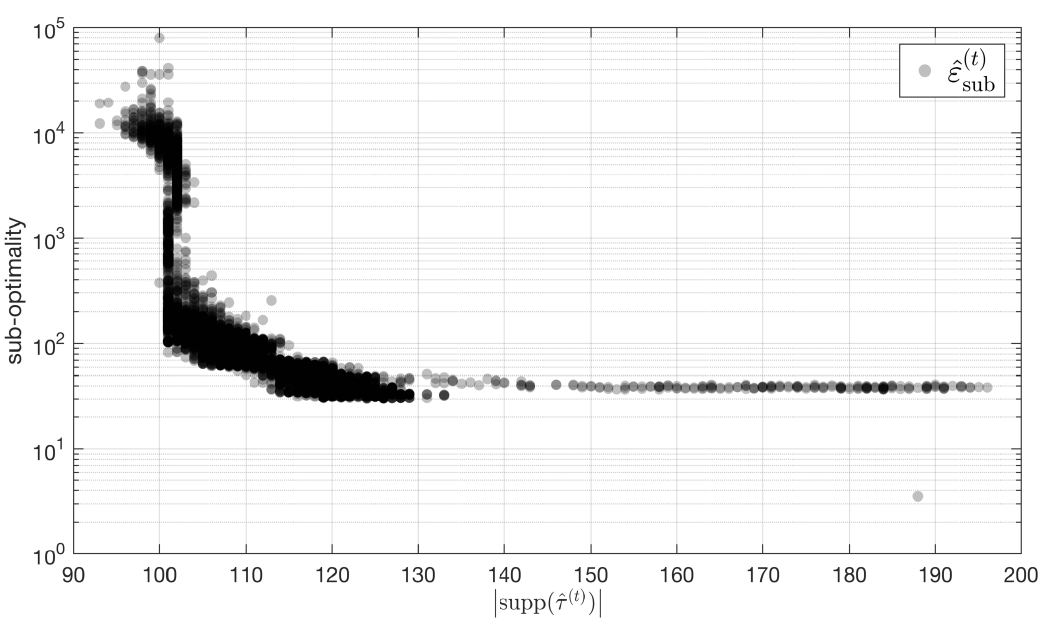}
\caption{\textbf{Experiment~2.} 
\textit{Left:} 
the computed lower bounds $\widehat{C}_{\DRO}^{\mathrm{LB}(t)}$ and upper bounds $\widehat{C}_{\DRO}^{\mathrm{LB}(t)}+\hat{\varepsilon}_{\mathrm{sub}}^{(t)}$ plotted against $\big|\support(\hat{\tau}^{(t)})\big|$. 
\textit{Right:} 
the computed sub-optimality estimates $\hat{\varepsilon}_{\mathrm{sub}}^{(t)}$ plotted against $\big|\support(\hat{\tau}^{(t)})\big|$.}
\label{fig:assembly-errors}
\end{figure}

In Figure~\ref{fig:assembly-errors}, the left panel
shows an increasing trend of $\widehat{C}_{\DRO}^{\mathrm{LB}(t)}$ against $\big|\support(\hat{\tau}^{(t)})\big|$ and a decreasing trend of
$\widehat{C}_{\DRO}^{\mathrm{LB}(t)}+\hat{\varepsilon}_{\mathrm{sub}}^{(t)}$ against $\big|\support(\hat{\tau}^{(t)})\big|$, and
the right panel
shows that $\hat{\varepsilon}_{\mathrm{sub}}^{(t)}$ exhibits a slow downward trend when $\big|\support(\hat{\tau}^{(t)})\big|\ge110$.
The sub-optimality $\hat{\varepsilon}_{\mathrm{sub}}=3.5132$ was achieved with respect to $\hat{\tau}$ where $|\support(\hat{\tau})|=188$,
as indicated by the lowest point in the right panel of Figure~\ref{fig:assembly-errors}.
Here, the sub-optimality is significantly lower in the last iteration of Algorithm~\ref{alg:iterative} compared to previous iterations since it uses a much smaller value of the parameter $\vartheta^{(t)}$; see our discussion in Setting~\ref{sett:iterative} about the strategy for setting the parameters $\varsigma^{(t)}$ and $\vartheta^{(t)}$.
Similar to the results of Experiment~1, 
$\hat{\tau}$ here has high sparsity relative to the problem dimensions $d=50$, $k^*_2=150$.
In fact, the maximum value of $\big|\support(\hat{\tau}^{(t)})\big|$ encountered during the execution of Algorithm~\ref{alg:iterative} was~196.

\begin{figure}[t]
\centering
\includegraphics[width=0.80\linewidth]{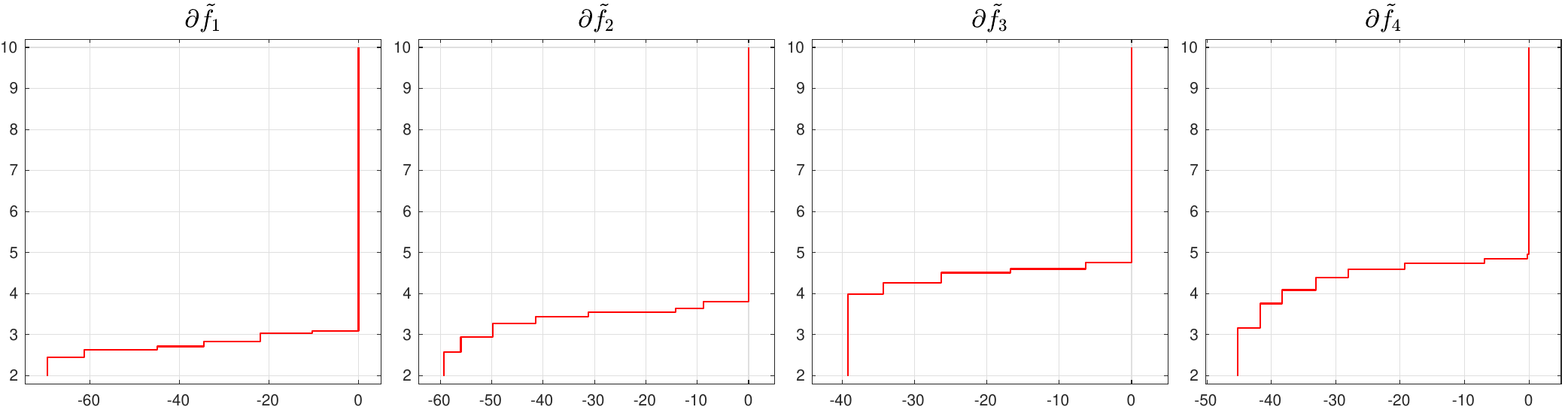}
\caption{\textbf{Experiment~2.} 
The graphs of the computed subdifferential mappings $\partial\tilde{f}_1,\partial\tilde{f}_2,\partial\tilde{f}_3,\partial\tilde{f}_{4}$.}
\label{fig:assembly-potentials}
\end{figure}

Figure~\ref{fig:assembly-potentials} shows the graphs of the subdifferential mappings of four computed functions $\tilde{f}_1,\allowbreak\tilde{f}_2,\allowbreak\tilde{f}_{3},\allowbreak\tilde{f}_{4}$.
It can be observed that the graphs in Figure~\ref{fig:assembly-potentials} exhibit more sophisticated structures compared to Figure~\ref{fig:scheduling-potentials} in Experiment~1.
Specifically, each of
$\tilde{f}_1,\ldots,\tilde{f}_{d}$ contains at least 3~breakpoints, at most 9~breakpoints, 
and the average number of breakpoints is~6.40.

\begin{figure}[t]
\centering
\includegraphics[width=0.58\linewidth]{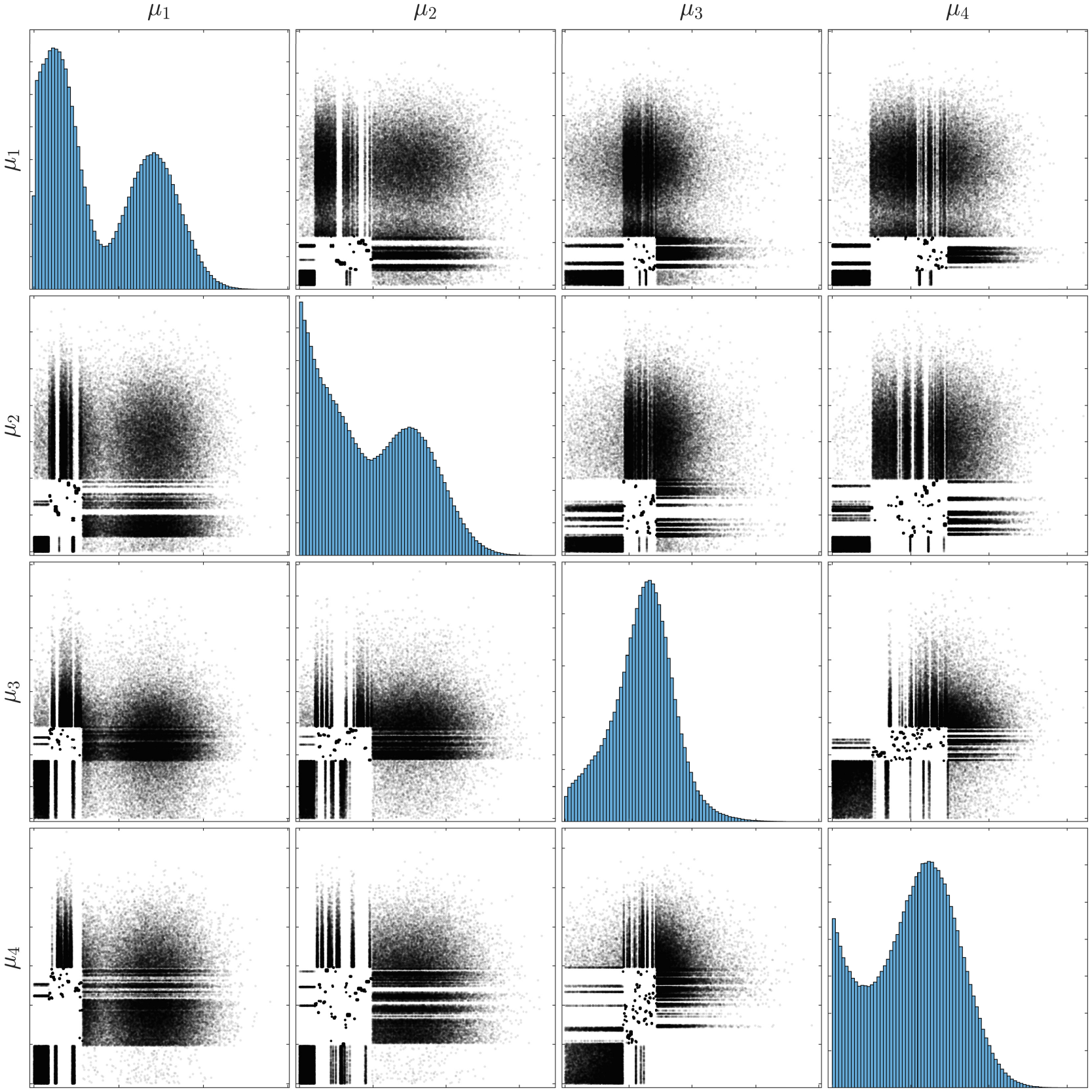}
\caption{\textbf{Experiment~2.} 
Pair-wise scatter plots of 4 marginals of the computed $\hat{\varepsilon}_{\mathrm{prob}}$-approximate worst-case probability measure $\hat{\mu}$ with respect to $\hat{\BIx}$.}
\label{fig:assembly-scatter}
\end{figure}

We also generate $10^5$ independent samples from the first four marginals of the $\hat{\varepsilon}_{\mathrm{prob}}$-approximate worst-case probability measure $\hat{\mu}$ with respect to $\hat{\BIx}$ computed by Algorithm~\ref{alg:iterative} and show the pair-wise scatter plots in Figure~\ref{fig:assembly-scatter}.
It seems that
the samples from $\hat{\mu}$ are organized in a large number of distinct clusters.
Moreover, the scatter plots lack simple and interpretable patterns like the block-diagonal patterns observed in Figure~\ref{fig:scheduling-scatter}.

\subsection{Experiment~3: supply chain network design with edge failure}
\label{ssec:experiments-network}

\begin{figure}[t]
\centering
\includegraphics[width=\linewidth]{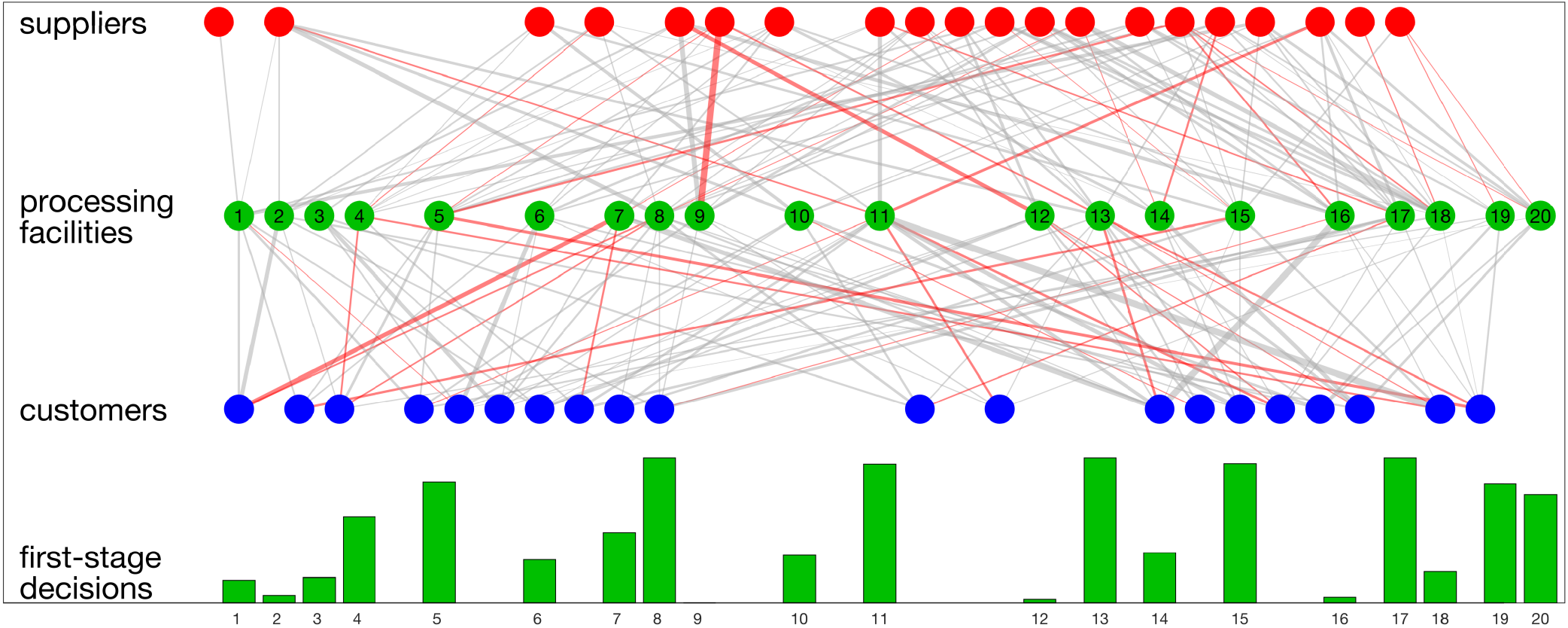}
\caption{\textbf{Experiment~3.}
The supply chain network configuration along with the computed approximately optimal processing capabilities $\hat{\BIx}$ in a bar plot.}
\label{fig:network}
\end{figure}

In this experiment, we solve the distributionally robust supply chain network design problem with edge failure introduced in Example~\ref{exp:supplychain}. 
We consider a supply chain network $(\TTV,\TTE)$ of 
$|\TTS|=20$ suppliers, 
$|\TTP|=20$ processing facilities, 
and $|\TTC|=20$ customers. 
We randomly generate $|\TTE_{\TTS,\TTP}|=100$ edges between the suppliers and the processing facilities
as well as $|\TTE_{\TTP,\TTC}|=100$ edges between the processing facilities and the customers.
For each customer $\TTc\in\TTC$, the law of the random demand 
$\omega_{\TTc}\in\CP\big([1, 2]\big)$ is a mixture of three equally weighted Gaussian measures with randomly generated parameters truncated to $[1,2]$.
The mean values of $(\mu_{\TTc})_{\TTc\in\TTC}$ range between 
1.2698 and 1.6165, with an average of 1.4539.
Moreover, we fix the maximum processing capability $\overline{x}_{\TTp}=2$ for all $\TTp\in\TTP$,
and we randomly generate 
the supplies $(b_{\TTs})_{\TTs\in\TTS}$,
the investment costs $(c_{1,\TTp})_{\TTp\in\TTP}$,
the transportation/processing costs $\big(c_{2,\TTs,\TTp}\big)_{(\TTs,\TTp)\in\TTE_{\TTS,\TTP}}$, $\big(c_{2,\TTp,\TTc}\big)_{(\TTp,\TTc)\in\TTE_{\TTP,\TTC}}$, 
and the maximum transportation capacities $(b_{\TTs,\TTp})_{(\TTs,\TTp)\in\TTE_{\TTS,\TTP}}$, $(b_{\TTp,\TTc})_{(\TTp,\TTc)\in\TTE_{\TTP,\TTC}}$.
Subsequently, we take $|\widetilde{\TTE}_{\TTS,\TTP}|=20$ edges from $\TTE_{\TTS,\TTP}$ with the least costs and $|\widetilde{\TTE}_{\TTP,\TTC}|=20$ edges from $\TTE_{\TTP,\TTC}$ with the least costs and set them to be susceptible to failure,
where the failure probabilities are randomly generated. 
The generated failure probabilities range between
0.0033 and 0.4448, with an average of 0.1112.
The configuration of this supply chain network is illustrated in Figure~\ref{fig:network}, in which 
the suppliers are represented by red circles, 
the processing facilities are represented by green circles, 
the customers are represented by blue circles,
the maximum transportation capacities of the edges are represented by their widths,
and the edges that are susceptible to failure are colored red.
In summary, this two-stage DRO model in Setting~\ref{sett:dro-alt} has $k_1=220$, 
$d=60$, and $k^*_2=280$.

\begin{figure}[t]
\centering
\includegraphics[width=0.49\linewidth]{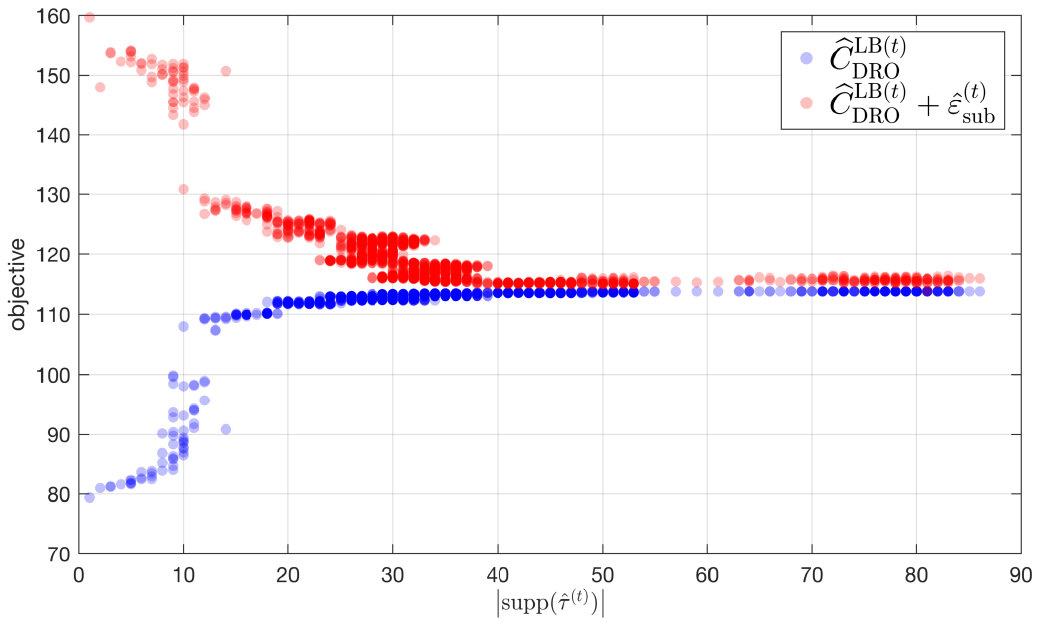}~
\includegraphics[width=0.49\linewidth]{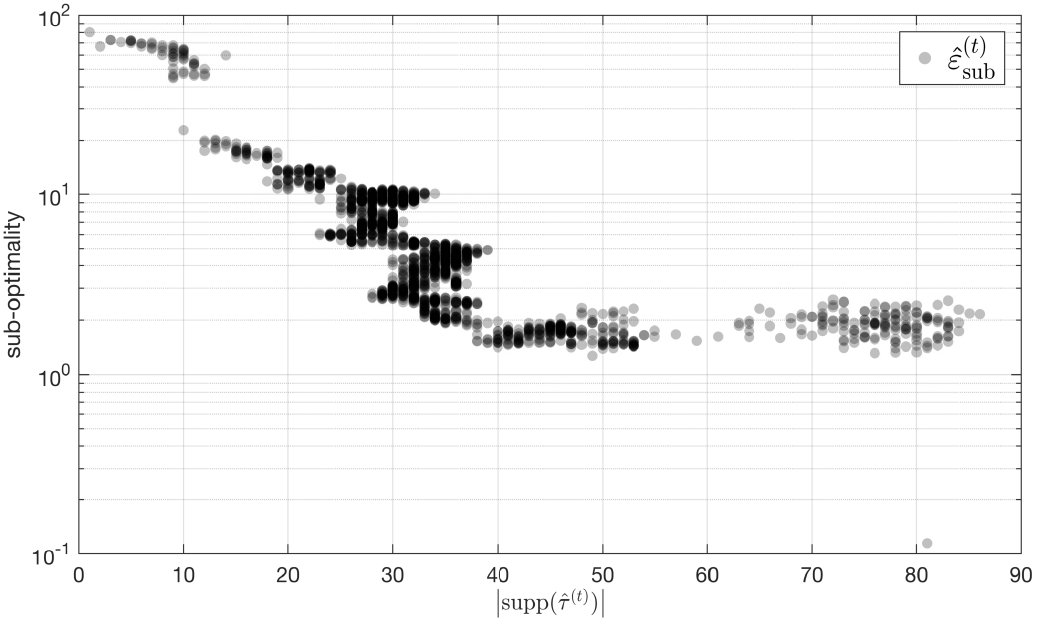}
\caption{\textbf{Experiment~3.} 
\textit{Left:} 
the computed lower bounds $\widehat{C}_{\DRO}^{\mathrm{LB}(t)}$ and upper bounds $\widehat{C}_{\DRO}^{\mathrm{LB}(t)}+\hat{\varepsilon}_{\mathrm{sub}}^{(t)}$ plotted against $\big|\support(\hat{\tau}^{(t)})\big|$. 
\textit{Right:} 
the computed sub-optimality estimates $\hat{\varepsilon}_{\mathrm{sub}}^{(t)}$ plotted against $\big|\support(\hat{\tau}^{(t)})\big|$.}
\label{fig:network-errors}
\end{figure}

We apply Algorithm~\ref{alg:iterative} to compute approximate optimizers of \eqref{eqn:dro-primal}, \eqref{eqn:dro-dual}, and \eqref{eqn:dro},
where we set 
$\varepsilon_{\mathrm{tol}}=0.5$,
$\varepsilon_{\OT}^{(t)}=\frac{1}{60000}$ for all $t\in\N$,
and vary the parameters $\varsigma^{(t)}$ and $\vartheta^{(t)}$ as discussed in Setting~\ref{sett:iterative}.
Algorithm~\ref{alg:iterative} terminated after 2203~iterations
and computed 
the lower bound
$\widehat{C}_{\DRO}^{\mathrm{LB}}=113.7691$,
the upper bound
$\widehat{C}_{\DRO}^{\mathrm{UB}}=113.8840$,
and the sub-optimality estimates
$\hat{\varepsilon}_{\mathrm{sub}}=0.1149$,
$\hat{\varepsilon}_{\mathrm{prob}}=0.1717$.
The relative errors are 
$\frac{\hat{\varepsilon}_{\mathrm{sub}}}{|\widehat{C}_{\DRO}^{\mathrm{UB}}|}=0.1009\%$,
$\frac{\hat{\varepsilon}_{\mathrm{prob}}}{|\widehat{C}_{\DRO}^{\mathrm{UB}}|}=0.1508\%$.
The computed $\hat{\varepsilon}_{\mathrm{sub}}$-optimizer $\hat{\BIx}$ of \eqref{eqn:dro} is shown in the bottom row of Figure~\ref{fig:network} as a bar plot.
These results show that Algorithm~\ref{alg:iterative} has produced high-quality approximate optimizers of \eqref{eqn:dro-primal}, \eqref{eqn:dro-dual}, and \eqref{eqn:dro} with small sub-optimality.

Figure~\ref{fig:network-errors} shows trends that are similar to Figure~\ref{fig:assembly-errors}.
The sub-optimality $\hat{\varepsilon}_{\mathrm{sub}}=0.1149$ was achieved with respect to $\hat{\tau}$ where $|\support(\hat{\tau})|=81$,
as indicated by the lowest point in the right panel of Figure~\ref{fig:network-errors}.
Here, the significantly lower sub-optimality in the last iteration of Algorithm~\ref{alg:iterative} is due to the same reason as in Experiment~2.
$\hat{\tau}$ has high sparsity relative to the problem dimensions $d=60$, $k^*_2=280$.
In fact, the maximum value of $\big|\support(\hat{\tau}^{(t)})\big|$ encountered during the execution of Algorithm~\ref{alg:iterative} was~$86$.
We have also examined the graphs of $\tilde{f}_1,\ldots,\tilde{f}_d$ 
and the pair-wise scatter plots of the marginals of $\hat{\tau}$ as in Experiments~1 and 2, but we omit them here since we have found that they look similar to those in Figure~\ref{fig:assembly-potentials} and Figure~\ref{fig:assembly-scatter} and do not provide additional insights.
We would like to remark that each of $\tilde{f}_1,\ldots,\tilde{f}_d$ contains at least 2 breakpoints, at most 13 breakpoints, and the average number of breakpoints is 5.53.

\subsection{Discussions}
\label{ssec:experiments-discussions}

In conclusion, the three numerical experiments in this section have revealed the following empirical findings about the primal and dual formulations \eqref{eqn:dro-primal} and \eqref{eqn:dro-dual} as well as our numerical algorithm (i.e., Algorithm~\ref{alg:iterative}).
\begin{itemize}
    \item The motivation of Algorithm~\ref{alg:iterative} is the presumption that \eqref{eqn:dro-dual} admits an approximate optimizer $(\tau,\Bxi_{\ineq},\Bxi_{\eq})$ where $\support(\tau)$ is finite and of moderate size, even though $|\CV_{\extremept}|$ is prohibitively large.
    We have observed that this is indeed the case for the three problem instances solved in our numerical experiments. 
    Despite the vastly different structures of the three problems, Algorithm~\ref{alg:iterative} was able to compute high-quality approximate optimizers $(\hat{\tau},\hat{\Bxi}_{\ineq},\hat{\Bxi}_{\eq})$ of \eqref{eqn:dro-dual} where $|\support(\hat{\tau})|$ is small relative to the problem dimensions in $d$ and $k^*_2$.
    
    \item 
    Theorem~\ref{thm:iterative}\ref{thms:iterative-primal} states that $(\hat{\BIx},\tilde{f}_1,\ldots,\tilde{f}_d)$ computed by Algorithm~\ref{alg:iterative} is an $\hat{\varepsilon}_{\mathrm{sub}}$-optimizer of \eqref{eqn:dro-primal}.
    In our numerical results, we have observed that $\tilde{f}_1,\ldots,\tilde{f}_d$ are continuous piece-wise affine functions each with a small number of breakpoints (typically no more than~10).
    This explains why we have chosen the incremental model in Proposition~\ref{prop:globalmax-MILP} to formulate the global maximization problem 
    $\mathrm{vio}\big(\hat{\BIx}^{(t)},\overline{f}_{1}^{(t)},\ldots,\overline{f}_{d}^{(t)}\big)$
    in Line~\ref{alglin:iterative-globalmax} of Algorithm~\ref{alg:iterative} into an MIP problem, as more sophisticated formulations and approaches are unlikely to be more computationally efficient when each $\overline{f}_{i}^{(t)}$ contains only a handful of breakpoints.

    \item Algorithm~\ref{alg:iterative} computes not only a high-quality approximate optimizer $\hat{\BIx}$ of \eqref{eqn:dro} but also an approximate worst-case probability measures $\hat{\mu}\in\CP(\BOmega)$ for the inner maximization problem in \eqref{eqn:dro} with respect to $\hat{\BIx}$, from which independent samples can be efficiently generated.
    These samples can be used to approximately compute quantities of interest related to $\hat{\mu}$ which can provide valuable insights about the worst-case probability measure to the decision maker.
    
    \item 
    In each of the three experiments, we fixed the tolerance parameter $\varepsilon_{\OT}^{(t)}$ to be a constant that is close to~0 for all $t\in\N$, and we observed that the computation of Line~\ref{alglin:iterative-LSIP} in Algorithm~\ref{alg:iterative} was highly efficient. 
    This was thanks to the warm starting technique when using the simplex algorithm to solve \LPprimal{\widehat{\CF}^{(t)}} and 
    \LPdual{\widehat{\CF}^{(t)}}
    in Line~\ref{alglin:cuttingplane-LSIP-LP} of Algorithm~\ref{alg:cuttingplane-LSIP}, 
    as discussed in Setting~\ref{sett:cuttingplane-LSIP}.
    We have observed in the three experiments that
    Algorithm~\ref{alg:cuttingplane-LSIP} terminated after at most 53 iterations.

    \item In the distributionally robust task scheduling problem (i.e., Example~\ref{exp:scheduling}) studied in Experiment~1,
    the so-called separation problem in \citep[Eq.~(10)]{chen2018distributionally} can be seen as a special case of the global maximization problem (\ref{eqn:dual-suboptimality-vio-def}) that we solve in Algorithm~\ref{alg:iterative}, with respect to $\BP\leftarrow\BI_d$.
    \citet[Section~5.1.1 \& Corollary~6]{chen2018distributionally} have established the polynomial solvability of this separation problem under the additional assumption that the weights $c_{2,1},\ldots,c_{2,d}$ are all positive integers, by explicitly identifying an LP reformulation.
    In Experiment~1, we have observed that the global maximization problem in Line~\ref{alglin:iterative-globalmax} of Algorithm~\ref{alg:iterative} was highly efficient even when the weights $c_{2,1},\ldots,c_{2,d}$ are non-integers, 
    and that the computational time of Algorithm~\ref{alg:iterative} was dominated by Line~\ref{alglin:iterative-LSIP}, which solves \LSIPprimal{\CV^{(t)}} by Algorithm~\ref{alg:cuttingplane-LSIP}.
    This demonstrates that the MIP formulation of (\ref{eqn:dual-suboptimality-vio-def}) in Proposition~\ref{prop:globalmax-MILP} along with Gurobi's \citep{gurobi} implementation of the branch-and-bound algorithm can solve (\ref{eqn:dual-suboptimality-vio-def}) efficiently
    without the need of explicit reformulation.

    \item 
    In contrast to Experiment~1, the computational time of Algorithm~\ref{alg:iterative} in Experiments~2 and 3 
    was dominated by the global maximization problem in Line~\ref{alglin:iterative-globalmax}.
    This is not surprising since global maximization problems of the form (\ref{eqn:dual-suboptimality-vio-def}) have been known to be NP-hard in general (see, e.g., \citep[Theorem~2]{chen2018distributionally}).
    Despite this,
    by adopting the approximation technique in Line~\ref{alglin:iterative-funcrelax} of Algorithm~\ref{alg:iterative},
    the MIP formulation of (\ref{eqn:dual-suboptimality-vio-def}) in Proposition~\ref{prop:globalmax-MILP},
    and by utilizing the state-of-the-art MIP solver of Gurobi \citep{gurobi}, we managed to solve moderately sized problem instances ($d=50$ in Experiment~2 and $d=60$ in Experiment~3, respectively) 
    and obtain high-quality approximate optimizers of \eqref{eqn:dro}.
\end{itemize}

\newpage

\appendix

\section{Proof of theoretical results}
\label{apx:proof}

\subsection{Results about one-dimensional semi-discrete optimal transport}
\label{sapx:OT-1d}

The following lemma recalls some key properties of one-dimensional semi-discrete optimal transport which are used throughout the paper.
\begin{lemma}[One-dimensional semi-discrete optimal transport]\label{lem:OT-1d}
    Let $-\infty<\underline{\omega}\le \overline{\omega}<\nobreak\infty$, let $\mu\in\CP([\underline{\omega},\overline{\omega}])$,
    and define 
    \begin{align}
        \OT_{\mu}(\nu):=\inf_{\gamma\in\Gamma(\mu,\nu)}\bigg\{\int_{[\underline{\omega},\overline{\omega}]\times\R}-z\omega\DIFFM{\gamma}{\DIFF\omega,\DIFF z}\bigg\} \qquad \forall \nu\in\CP_1(\R).
        \label{eqn:OT-1d-def}
    \end{align}
    Then, the following statements hold.\widowpenalties-1000
    \begin{enumerate}[label=(\roman*),beginpenalty=10000]
        \item\label{lems:OT-1d-duality}
        The function $\CP_1(\R)\ni\nu\mapsto \OT_{\mu}(\nu)\in\R$ is 
        weakly continuous, and 
        it holds that%
        \begin{align}%
            \hspace{20pt}\OT_{\mu}(\nu)&= \max_{f\in\CC_{\mathrm{cvx}}(\R)}\bigg\{{-\int_{[\underline{\omega},\overline{\omega}]}}\sup_{z\in\R}\big\{z\omega-f(z)\big\}\DIFFM{\mu}{\DIFF\omega}- \int_{\R}f\DIFFX{\nu}\bigg\} \quad \forall \nu\in\CP_1(\R).\label{eqn:OT-1d-dual-discside}%
        \end{align}%

        \item\label{lems:OT-1d-optimality}
        For any $\nu\in\CP_1(\R)$ and any $f\in\CC_{\mathrm{cvx}}(\R)$, 
        $f$ maximizes (\ref{eqn:OT-1d-dual-discside}) 
        if and only if there exists $\gamma\in\Gamma(\mu,\nu)$ which satisfies $\omega \in \partial f(z)$ for $\gamma$-almost every $(\omega,z)\in [\underline{\omega},\overline{\omega}]\times\R$.
        Moreover, such $\gamma\in\Gamma(\mu,\nu)$ also minimizes 
        (\ref{eqn:OT-1d-def}).

    \end{enumerate}
    Next, for $m\in\N$ and $-\infty<\kappa_1<\cdots<\kappa_{m}<\infty$,
    consider arbitrary $\nu_1,\ldots,\nu_m\in\R_{+}$ satisfying $\sum_{j=1}^m\nu_j=1$,
    and define $\nu:=\sum_{j=1}^m\nu_j\delta_{\kappa_j}\in\CP_1(\R)$.
    Let $G^{-}_{\mu}:[0,1]\to[\underline{\omega},\overline{\omega}]$ and $G^{-}_{\nu}:[0,1]\to\R$ be defined as follows:
    \begin{align*}
        G^{-}_{\mu}(u)&:=\inf\big\{t\in[\underline{\omega},\overline{\omega}]: \mu\big([\underline{\omega},t]\big)\ge u\big\} \hspace{87.5pt} \forall u\in[0,1], \\
        G^{-}_{\nu}(u)&:=\begin{cases}
            \kappa_1 & \hspace{137.6pt} \forall u\in \big[0, \nu\big((-\infty,\kappa_1]\big)\big], \\
            \kappa_j & \hspace{20pt} \forall u\in \big(\nu\big((-\infty,\kappa_{j})\big), \nu\big((-\infty,\kappa_{j}]\big)\big],\; \forall 2\le j\le m.
        \end{cases}
    \end{align*}
    Then, the following statement holds.
    \begin{enumerate}[label=(\roman*),beginpenalty=10000]
        \setcounter{enumi}{2}
        \item\label{lems:OT-1d-coupling}
        Let $\mathscr{L}_{[0,1]}$ denote the Lebesgue measure on $[0,1]$ and define $\gamma:=\big[G^{-}_{\mu},G^{-}_{\nu}\big]\sharp\mathscr{L}_{[0,1]}\in\CP\big([\underline{\omega},\overline{\omega}]\times\nobreak \R\big)$.
        Then, $\gamma\in\Gamma(\mu,\nu)$ and $\gamma$ minimizes (\ref{eqn:OT-1d-def}).

        \item\label{lems:OT-1d-potential-existence}
        Let $f_1\in\R$ be arbitrary and
        let $f:\R\to\R$ be defined as follows:
        \begin{align}
            \qquad\;\; f(z):=\begin{cases}
                f_1+\underline{\omega}(z-\kappa_1) & 
                \hspace{84pt} \forall z\in(-\infty,\kappa_1], \\
                f(\kappa_{j}) + G^{-}_{\mu}\big(\nu\big((-\infty,\kappa_{j}]\big)\big)(z - \kappa_{j})   & \forall z\in (\kappa_{j}, \kappa_{j+1}],\; \forall 1\le j\le m-1,\\
                f(\kappa_{m}) + \overline{\omega}(z-\kappa_{m}) & 
                \hspace{88pt} \forall z\in(\kappa_{m},\infty).
            \end{cases}\!\!\!
            \label{eqn:OT-1d-potential-CPWA}
        \end{align}%
        Then, $f\in\CC_{\mathrm{cvx}}(\R)$ is a maximizer of (\ref{eqn:OT-1d-dual-discside}).

        \item\label{lems:OT-1d-potential-compact}
        Let $B:=\Big\{(\tilde{f}_{1},\ldots,\tilde{f}_{m})^\TRANSP\in\R^{m}:\tilde{f}_{1}=0,\; \underline{\omega} \le \frac{\tilde{f}_{2}-\tilde{f}_{1}}{\kappa_{2}-\kappa_{1}}\le \cdots \le \frac{\tilde{f}_{m}-\tilde{f}_{m-1}}{\kappa_{m}-\kappa_{m-1}}\le \overline{\omega}\Big\}$, and let us consider the following maximization problem over $B$:
        \begin{align}
            \max_{(\tilde{f}_{1},\ldots,\tilde{f}_{m})^\TRANSP\in B}\Bigg\{{-\int_{[\underline{\omega},\overline{\omega}]}}\max_{1\le j\le m}\big\{\kappa_{j}\omega- \tilde{f}_{j}\big\}\DIFFM{\mu}{\DIFF\omega}-\Bigg(\sum_{j=1}^{m}\tilde{f}_{j}\nu_{j}\Bigg)\Bigg\}.
            \label{eqn:OT-1d-potential-compact-representation}
        \end{align}
        Then, it holds that the optimal value of (\ref{eqn:OT-1d-potential-compact-representation}) is equal to $\OT_{\mu}(\nu)$,
        and that $(\tilde{f}_{1},\ldots,\tilde{f}_{m})^\TRANSP\in B$ maximizes (\ref{eqn:OT-1d-potential-compact-representation}) only if
        there exists $f\in\CC_{\mathrm{cvx}}(\R)$ that maximizes (\ref{eqn:OT-1d-dual-discside}) where $f(\kappa_j)=\tilde{f}_j$ for $j=1,\ldots,m$ and $\partial f(z)\subseteq[\underline{\omega},\overline{\omega}]$ for all $z\in\R$.
    \end{enumerate}
\end{lemma}

\begin{proof}[Proof of Lemma~\ref{lem:OT-1d}]
    Statement~\ref{lems:OT-1d-duality} follows directly from classical results about the stability of optimal transport and the Kantorovich duality; see, e.g., \citep[Theorem~5.20 \& Theorem~5.10(iii)]{villani2008optimal}.
    Statement~\ref{lems:OT-1d-optimality} follows from 
    the equivalence between (a) and (d) in \citep[Theorem~5.10(ii)]{villani2008optimal} as well as \citep[Theorem~5.10(iii)]{villani2008optimal}.
    
    Throughout the rest of this proof,
    we will fix arbitrary $\nu_1,\ldots,\nu_m\in\R_{+}$ with $\sum_{j=1}^m\nu_j=1$, and 
    let $\nu:=\sum_{j=1}^m\nu_j\delta_{\kappa_j}$.
    Observe from the definitions of $G_{\mu}^{-}$ and $G_{\nu}^{-}$ that 
    $G_{\mu}^{-}$ is the left-continuous generalized inverse of the distribution function $[\underline{\omega},\overline{\omega}]\ni \omega\mapsto \mu\big([\underline{\omega},\omega]\big)\in[0,1]$ of~$\mu$, and
    $G^{-}_{\nu}(u)=\inf\big\{t\in[\kappa_1,\infty):\nu\big((-\infty,t]\big)\ge u\big\}$ $\forall u\in[0,1]$,
    which coincides with the left-continuous generalized inverse of the distribution function
    $\R\ni z\mapsto \nu\big([\kappa_1,z]\big)\in[0,1]$ of~$\nu$.
    Subsequently, the properties of $G_{\mu}^{-}$, $G_{\nu}^{-}$ in \citep[Proposition~A.3(iv)]{mcneil2005quantitative} imply that $\gamma\in\Gamma(\mu, \nu)$.
    Moreover, since $\int_{[\underline{\omega},\overline{\omega}]\times\R} {-z\omega}\DIFFM{\gamma}{\DIFF \omega,\DIFF z}=\int_{(0,1)}{-G^{-}_{\mu}}(u)G^{-}_{\nu}(u)\DIFFM{\mathscr{L}_{[0,1]}}{\DIFF u}$,
    it follows from \citep[Eq.~(3.1.4)]{rachev1998mass} that $\gamma$ minimizes (\ref{eqn:OT-1d-def}).
    This completes the proof of statement~\ref{lems:OT-1d-coupling}.
    
    Next,
    we define $f:\R\to\R$ by (\ref{eqn:OT-1d-potential-CPWA})
    and prove statement~\ref{lems:OT-1d-potential-existence}.
    Since $\underline{\omega}\le G^{-}_{\mu}\big(\nu\big((-\infty,\kappa_{j}]\big)\big)\le G^{-}_{\mu}\big(\nu\big((-\infty,\kappa_{j+1}]\big)\big)\le \overline{\omega}$ for $j=1,\ldots,m-1$, 
    it holds that $f\in\CC_{\mathrm{cvx}}(\R)$. 
    To show that $f$ is a maximizer of (\ref{eqn:OT-1d-dual-discside}), let $\gamma:=\big[G^{-}_{\mu},G^{-}_{\nu}\big]\sharp\mathscr{L}_{[0,1]}\in\Gamma(\mu,\nu)$ be defined as in statement~\ref{lems:OT-1d-coupling}.
    For any $u\in[0,1]$, it follows from the definition of $G^{-}_{\nu}$ and the non-decreasing property of $G^{-}_{\mu}$ that
    \begin{align*}
        G^{-}_{\nu}(u)=\kappa_j \quad &\Rightarrow \quad u\le \nu\big((-\infty,\kappa_j]\big)\\
        &\Rightarrow \quad G^{-}_{\mu}(u)\le G^{-}_{\mu}\big(\nu\big((-\infty,\kappa_j]\big)\big)=\partial_{+}f(\kappa_j) \qquad \forall 1\le j\le m-1,\\
        G^{-}_{\nu}(u)=\kappa_j \quad &\Rightarrow \quad u> \nu\big((-\infty,\kappa_{j-1}]\big) \\
        &\Rightarrow \quad G^{-}_{\mu}(u) \ge G^{-}_{\mu}\big(\nu\big((-\infty,\kappa_{j-1}]\big)\big)=\partial_{-}f(\kappa_j) \qquad \hspace{8.2pt} \forall 2\le j\le m.
    \end{align*}
    The above properties show that
    \begin{align*}
        \big\{u\in[0,1]:G_{\nu}^{-}(u)=\kappa_j\big\} &= \Big\{u\in[0,1]:G_{\nu}^{-}(u)=\kappa_j,\; G_{\mu}^{-}(u)\in\big[\partial_{-}f(\kappa_j),\partial_{+}f(\kappa_j)\big]\Big\}\\
        & \hspace{220pt} \forall 1\le j\le m.
    \end{align*}
    We hence get
    \begin{align*}
        &\hspace{-20pt}\gamma\Big(\big\{(\omega,z)\in[\underline{\omega},\overline{\omega}]\times\R: \omega \in \partial f(z)\big\}\Big) \\
        &= \mathscr{L}_{[0,1]}\Bigg(\bigcup_{j=1}^{m}\Big\{u\in[0,1]:G_{\nu}^{-}(u)=\kappa_j,\; G_{\mu}^{-}(u)\in\big[\partial_{-}f(\kappa_j),\partial_{+}f(\kappa_j)\big]\Big\}\Bigg)\allowdisplaybreaks\\
        &= \mathscr{L}_{[0,1]}\Bigg(\bigcup_{j=1}^{m}\big\{u\in[0,1]:G_{\nu}^{-}(u)=\kappa_j\big\}\Bigg)=G_{\nu}^{-}\sharp\mathscr{L}_{[0,1]}\big(\{\kappa_1,\ldots,\kappa_m\}\big)= 1,
    \end{align*}
    and it thus holds for $\gamma$-almost every $(\omega,z)\in[\underline{\omega},\overline{\omega}]\times\R$ that $\omega\in \partial f(z)$.
    We conclude by statement~\ref{lems:OT-1d-optimality} that $f$ is a maximizer of (\ref{eqn:OT-1d-dual-discside}).
    This completes the proof of statement~\ref{lems:OT-1d-potential-existence}.

    Lastly, to prove statement~\ref{lems:OT-1d-potential-compact}, let us fix an arbitrary maximizer $(\tilde{f}_{1},\ldots,\tilde{f}_{m})^\TRANSP\in B$ of (\ref{eqn:OT-1d-potential-compact-representation}), and let 
    $f$ be the unique continuous function on $\R$ which is piece-wise affine on $(-\infty,\kappa_1],\allowbreak[\kappa_1,\kappa_2],\ldots,[\kappa_{m-1},\kappa_{m}],\allowbreak[\kappa_{m},\infty)$ and satisfies $f(\kappa_j)=\tilde{f}_j$ for $j=1,\ldots,m$, as well as $\partial_{-}f(\kappa_1)=\nobreak\underline{\omega}$, $\partial_{+}f(\kappa_m)=\overline{\omega}$.
    One verifies that $f\in\CC_{\mathrm{cvx}}(\R)$, 
    $\partial f(z)\subseteq[\underline{\omega},\overline{\omega}]$ for all $z\in\R$, 
    and that ${-\int_{[\underline{\omega},\overline{\omega}]}}\sup_{z\in\R}\big\{z\omega-f(z)\big\}\DIFFM{\mu}{\DIFF\omega}-\int_{\R}f\DIFFX{\nu}={-\int_{[\underline{\omega},\overline{\omega}]}}\max_{1\le j\le m}\big\{\kappa_j\omega-\tilde{f}_j\big\}\DIFFM{\mu}{\DIFF\omega}-\big(\sum_{j=1}^{m}\tilde{f}_j\nu_j\big)$.
    Statement~\ref{lems:OT-1d-duality} then shows that $\eqref{eqn:OT-1d-potential-compact-representation}\le \OT_{\mu}(\nu)$.
    To show that $\eqref{eqn:OT-1d-potential-compact-representation}\ge \OT_{\mu}(\nu)$, let us fix an arbitrary maximizer $f\in\CC_{\mathrm{cvx}}(\R)$ of (\ref{eqn:OT-1d-dual-discside}). 
    Let us denote 
    $\underline{j}:=\min\big\{j\in\{1,\ldots,m\}:\nu_{j}>\nobreak0\big\}$,
    $\overline{j}:=\max\big\{j\in\{1,\ldots,m\}:\nu_{j}>\nobreak0\big\}$.
    By statement~\ref{lems:OT-1d-optimality}, 
    there exists $\gamma\in\Gamma(\mu,\nu)$ such that $\omega\in\partial f(z)$ for $\gamma$-almost every $(\omega,z)\in [\underline{\omega},\overline{\omega}]\times\R$,
    and it hence holds that 
    $\partial_{+}f(\kappa_{\underline{j}})\ge \underline{\omega}$, 
    $\partial_{-}f(\kappa_{\overline{j}})\le\nobreak \overline{\omega}$.
    Let $\tilde{f}$ be the unique continuous function on $\R$ which is piece-wise affine on $(-\infty,\kappa_{\underline{j}}],\allowbreak[\kappa_{\underline{j}},\kappa_{\underline{j}+1}],\ldots,\allowbreak[\kappa_{\overline{j}-1},\kappa_{\overline{j}}],[\kappa_{\overline{j}},\infty)$ and satisfies
    $\tilde{f}(\kappa_{j})=f(\kappa_{j})$ for $j=\underline{j},\ldots,\overline{j}$,
    as well as
    $\partial_{-}\tilde{f}(\kappa_{\underline{j}})=\underline{\omega}$,
    $\partial_{+}\tilde{f}(\kappa_{\overline{j}})=\overline{\omega}$.
    One verifies that $\tilde{f}\in\CC_{\mathrm{cvx}}(\R)$ 
    and that $\omega\in\partial \tilde{f}(z)$ holds for $\gamma$-almost every $(\omega,z)\in [\underline{\omega},\overline{\omega}]\times\R$.
    Thus, $\tilde{f}$ is also a maximizer of (\ref{eqn:OT-1d-dual-discside}).
    Moreover, let $\tilde{f}_j:=\tilde{f}(\kappa_j)-\tilde{f}(\kappa_1)$ for $j=1,\ldots,m$.
    Observe that 
    $(\tilde{f}_1,\ldots,\tilde{f}_m)^\TRANSP\in B$
    and that ${-\int_{[\underline{\omega},\overline{\omega}]}}\max_{1\le j\le m}\big\{\kappa_j\omega-\tilde{f}_j\big\}\DIFFM{\mu}{\DIFF\omega}-\big(\sum_{j=1}^{m}\tilde{f}_j\nu_j\big)={-\int_{[\underline{\omega},\overline{\omega}]}}\sup_{z\in\R}\big\{z\omega-\tilde{f}(z)\big\}\DIFFM{\mu}{\DIFF\omega}-\int_{\R}\tilde{f}\DIFFX{\nu}=\OT_{\mu}(\nu)$.
    It thus holds that $\eqref{eqn:OT-1d-potential-compact-representation}\ge \OT_{\mu}(\nu)$.
    The above arguments also complete the proof of the rest of statement~\ref{lems:OT-1d-potential-compact}.
    The proof is now complete.
\end{proof}

\clearpage
\subsection{Proof of results in Section~\ref{sec:dro}}\label{sapx:proof-dro}
\nopagebreak
\begin{proof}[Proof of Lemma~\ref{lem:stage2-dual}]
    It follows from the assumption \ref{setts:dro-stage2} that $C_2(\BIx,\Bomega)$ corresponds to a feasible and bounded LP problem for all $\BIx\in \CS_1$ and all $\Bomega\in\BOmega$. 
    Statement~\ref{lems:stage2-dual-duality} then follows from the strong duality of LP problems. 
    
    In the following, let 
    \begin{align*}
    \begin{split}
    \BQ:=\left[\begin{tabular}{l}
    $\BQ_{\ineq}$ \\
    $\BQ_{\eq}$
    \end{tabular}
    \right]\in\R^{(m_{\ineq}+m_{\eq})\times k_1}, \;\;
    \BP:=\left[\begin{tabular}{l}
    $\BP_{\ineq}$ \\
    $\BP_{\eq}$
    \end{tabular}
    \right]\in\R^{(m_{\ineq}+m_{\eq})\times d},\;\;
    \Bb:=\left[\begin{tabular}{l}
    $\Bb_{\ineq}$ \\
    $\Bb_{\eq}$
    \end{tabular}
    \right]\in\R^{m_{\ineq}+m_{\eq}},
    \end{split}
    \end{align*}
    for notational simplicity. 
    It then follows from statement~\ref{lems:stage2-dual-duality} that
    \begin{align}
    C_2(\BIx,\Bomega)=\max_{\Blambda\in \CS_2^*}\big\{\langle\BQ\BIx+\BP\Bomega+\Bb,\Blambda\rangle\big\} \qquad \forall \BIx\in\CS_1,\; \forall \Bomega\in\BOmega.
    \label{eqn:stage2-dual-simplified}
    \end{align}
    Subsequently, it follows the assumption \ref{setts:dro-stage2} that the problem (\ref{eqn:stage2-dual-simplified}) is bounded above by $\alpha\in\R$. 
    Therefore, for every $\Brho\in\recession(\CS^*_2)$, every $\BIx\in \CS_1$, and every $\Bomega\in\BOmega$, it holds that $\langle\BQ\BIx+\BP\Bomega+\Bb,\Brho\rangle\le 0$. 
    Let $\CV_{\extremept}:=\extremept\big(\CS^*_2\intersection \lineality(\CS^*_2)^{\perp}\big)$.
    By \citep[p.170 \& Theorem~18.5 \& Theorem~19.1]{rockafellar1970convex}, the polyhedral convex set $\CS^*_2$ can be expressed as $\CS^*_2=\conv(\CV_{\extremept})+\recession(\CS^*_2)$. 
    Consequently, one observes that $\CS^*_2$ in (\ref{eqn:stage2-dual-simplified}) can be replaced by any finite set $\CV$ that satisfies $\CV_{\extremept}\subseteq\CV\subset\CS^*_2$ without changing the value of the supremum. 
    The proof is now complete.
\end{proof}

\subsection{Proof of results in Section~\ref{sec:primaldual}}
\label{sapx:proof-primaldual}


Before proving Theorem~\ref{thm:primaldual}, let us first establish Lemmas~\ref{lem:primal-objective}--\ref{lem:weak-duality} below.

\begin{lemma}\label{lem:primal-objective}
    Under Setting~\ref{sett:dro-alt}, let $(\BIx,f_1,\ldots,f_d)$ be feasible for \eqref{eqn:dro-primal}.
    Then, it holds that $\sum_{i=1}^d \sup_{z\in\R}\big\{z\omega_i - f_i(z)\big\}\ge C_2(\BIx,\Bomega)$ for all $\Bomega=(\omega_1,\ldots,\omega_d)^\TRANSP\in\BOmega$ and
    \begin{align*}
        \sum_{i=1}^d\int_{\Omega_i}\sup_{z\in\R}\big\{z\omega - f_i(z)\big\}\DIFFM{\mu_i}{\DIFF\omega} &\ge \sup_{\mu\in\Gamma(\mu_1,\ldots,\mu_d)}\bigg\{\int_{\BOmega}C_2(\BIx,\Bomega)\DIFFM{\mu}{\DIFF\Bomega}\bigg\}.
    \end{align*}
\end{lemma}

\begin{proof}[Proof of Lemma~\ref{lem:primal-objective}]
    Since $\langle\BQ\BIx+\Bb,\Blambda\rangle+\sum_{i=1}^{d}f_i\big(\langle\Bp_i,\Blambda\rangle\big)\le 0$ for all $\Blambda\in\CS^*_2$,
    it holds for all $\Bomega=(\omega_1,\ldots,\omega_d)^\TRANSP\in\BOmega$ that 
    \begin{align*}
        \begin{split}
            C_2(\BIx,\Bomega)&=\max_{\Blambda\in\CS_2^*}\big\{\langle\BQ\BIx+\BP\Bomega+\Bb,\Blambda\rangle\big\}= \max_{\Blambda\in\CS^*_2}\Bigg\{\langle\BQ\BIx+\Bb,\Blambda\rangle+\sum_{i=1}^d \langle\Bp_i,\Blambda\rangle\omega_i\Bigg\} \\
            &\le \sup_{\Blambda\in\CS_2^*}\left\{\sum_{i=1}^d\langle\Bp_i,\Blambda\rangle \omega_i - f_i\big(\langle\Bp_i,\Blambda\rangle\big)\right\}\le \sum_{i=1}^d \sup_{z\in\R}\big\{z\omega_i - f_i(z)\big\}.
        \end{split}
    \end{align*}
    Integrating the above inequalities with respect to any $\mu\in\Gamma(\mu_1,\ldots,\mu_d)$ and then taking the supremum over all $\mu\in\Gamma(\mu_1,\ldots,\mu_d)$ completes the proof.
\end{proof}

\begin{lemma}\label{lem:dual-gluing}
    Under Setting~\ref{sett:dro-alt}, the following statements hold.
    \begin{enumerate}[label=(\roman*),beginpenalty=10000]
        \item\label{lems:dual-gluing-construction}
        For any $\tau\in\CP_1(\CS^*_2)$,
        there exists $\pi\in\CP(\BOmega\times\CS^*_2)$ such that 
        the marginal of $\pi$ on $\CS^*_2$ is $\tau$,
        and for $i=1,\ldots,d$, the marginal of $\pi$ on $\Omega_i$ is $\mu_i$, the marginal $\pi_i$ of $\pi$ on $\Omega_i\times\CS^*_2$ satisfies 
        $\int_{\Omega_i\times\CS^*_2}{-\langle\Bp_i,\Blambda\rangle\omega}\DIFFM{\pi_i}{\DIFF\omega,\DIFF\Blambda}=\OT_{\mu_i}\big(\langle\Bp_i,\cdot\,\rangle\sharp\tau\big)$.
        
        \item\label{lems:dual-gluing-equalities}
        The following equalities hold for all $\BIx\in\CS_1$:
        \begin{subequations}%
        \label{eqn:dual-gluing-equalities}%
        \begin{align}
            &\sup_{\mu\in\Gamma(\mu_1,\ldots,\mu_d)}\bigg\{\int_{\BOmega}C_2(\BIx,\Bomega)\DIFFM{\mu}{\DIFF\Bomega}\bigg\}\label{eqn:dual-gluing-equalities-MMOT} \\
            &\qquad= \sup_{\tau\in\CP(\CV_{\extremept})}\Bigg\{\int_{\CS^*_2}\langle\BQ\BIx+\Bb,\Blambda\rangle\DIFFM{\tau}{\DIFF\Blambda} -\Bigg(\sum_{i=1}^d\OT_{\mu_i}\big(\langle\Bp_i,\cdot\,\rangle\sharp\tau\big)\Bigg)\Bigg\}\label{eqn:dual-gluing-equalities-finitemeas}\\
            &\qquad= \sup_{\tau\in\CP_1(\CS^*_2)}\Bigg\{\int_{\CS^*_2}\langle\BQ\BIx+\Bb,\Blambda\rangle\DIFFM{\tau}{\DIFF\Blambda} -\Bigg(\sum_{i=1}^d\OT_{\mu_i}\big(\langle\Bp_i,\cdot\,\rangle\sharp\tau\big)\Bigg)\Bigg\}.
            \label{eqn:dual-gluing-equalities-meas}
        \end{align}%
        \end{subequations}%

        \item\label{lems:dual-gluing-optimality}
        Let $\BIx\in\CS_1$ be arbitrary,
        let $\bar{\tau}\in\CP_1(\CS^*_2)$ be a maximizer of (\ref{eqn:dual-gluing-equalities-meas}), 
        let $\bar{\pi}\in\CP(\BOmega\times\nobreak\CS^*_2)$ be the probability measure in statement~\ref{lems:dual-gluing-construction} with $\tau\leftarrow\bar{\tau}$,
        and let $\bar{\mu}$ denote the marginal of $\bar{\pi}$ on~$\BOmega$.
        Then, it holds that 
        $\langle\BQ\BIx+\BP\Bomega+\Bb,\Blambda\rangle=C_2(\BIx,\Bomega)$ for $\bar{\pi}$-almost every $(\Bomega,\Blambda)\in\BOmega\times\CS^*_2$,
        and that 
        $\bar{\mu}\in\argmax_{\mu\in\Gamma(\mu_1,\ldots,\mu_d)}\big\{\int_{\BOmega}C_2(\BIx,\Bomega)\DIFFM{\mu}{\DIFF\BIx}\big\}$.
        
    \end{enumerate}
\end{lemma}

\begin{proof}[Proof of Lemma~\ref{lem:dual-gluing}]
    Let us fix an arbitrary $\tau\in\CP_1(\CS_2^*)$ and prove statement~\ref{lems:dual-gluing-construction}.
    Observe that $\langle\Bp_i,\cdot\,\rangle\sharp\tau\in\CP_1(\R)$ for $i=1,\ldots,d$.
    We let $\CZ_i:=\R$ denote the space on which $\langle\Bp_i,\cdot\,\rangle\sharp\tau$ lives in order to differentiate multiple copies of the real line. 
    Moreover, let $\gamma_i\in\Gamma\big(\mu_i,\langle\Bp_i,\cdot\,\rangle\sharp\tau\big)\subset\CP(\Omega_i\times\CZ_i)$ be an optimal coupling of $\mu_i$ and $\langle\Bp_i,\cdot\,\rangle\sharp\tau$ with respect to the cost function 
    $\Omega_i\times\CZ_i\ni(\omega,z)\mapsto -z\omega\in\nobreak\R$, i.e., $\int_{\Omega_i\times\CZ_i}-z\omega\DIFFM{\gamma_i}{\DIFF \omega,\DIFF z}=\OT_{\mu_i}\big(\langle\Bp_i,\cdot\,\rangle\sharp\tau\big)$.
    We subsequently apply the ``gluing'' lemma (see, e.g., \citep[Lemma~7.6]{villani2003topics}) inductively to ``glue together'' 
    $\big(\mathrm{Id},\langle\Bp_1,\cdot\,\rangle,\ldots,\langle\Bp_d,\cdot\,\rangle\big)\sharp\tau\in\CP(\CS_2^*\times\CZ_1\times\cdots\times\CZ_d)$ and $\gamma_1\in\CP(\Omega_1\times\CZ_1),$ $\ldots,$ $\gamma_d\in\CP(\Omega_d\times\CZ_d)$, where $\mathrm{Id}:\CS_2^*\to\CS_2^*$ denotes the identity mapping.
    Specifically, let $\gamma^{(0)}:=\big(\mathrm{Id},\langle\Bp_1,\cdot\,\rangle,\ldots,\langle\Bp_d,\cdot\,\rangle\big)\sharp\tau\in\Gamma\big(\tau,\langle\Bp_1,\cdot\,\rangle\sharp\tau,\ldots,\langle\Bp_d,\cdot\,\rangle\sharp\tau\big)$.
    For $i=1,\ldots,d$, 
    let $\gamma^{(i)}\in\Gamma\big(\mu_1,\ldots,\mu_{i},\tau,\langle\Bp_1,\cdot\,\rangle\sharp\tau,\ldots,\langle\Bp_d,\cdot\,\rangle\sharp\tau\big)$ be constructed 
    by ``gluing together'' $\gamma^{(i-1)}\in\Gamma\big(\mu_1,\ldots,\mu_{i-1},\tau,\langle\Bp_1,\cdot\,\rangle\sharp\tau,\ldots,\langle\Bp_d,\cdot\,\rangle\sharp\tau\big)$ and $\gamma_i\in\Gamma\big(\mu_i,\langle\Bp_i,\cdot\,\rangle\sharp\tau\big)$,
    that is,
    we apply \citep[Lemma~7.6]{villani2003topics} with
    $X_1\leftarrow\Omega_1\times\cdots\times\Omega_{i-1}\times\CS^*_2\times\CZ_1\times\cdots\times\CZ_{i-1}\times\CZ_{i+1}\times\cdots\times\CZ_{d}$,
    $X_2\leftarrow \CZ_i$,
    $X_3\leftarrow \Omega_i$,
    $\mu_1\leftarrow$ the marginal of $\gamma^{(i-1)}$ on $\Omega_1\times\cdots\times\Omega_{i-1}\times\CS^*_2\times\CZ_1\times\cdots\times\CZ_{i-1}\times\CZ_{i+1}\times\cdots\times\CZ_{d}$,
    $\mu_2\leftarrow \langle\Bp_i,\cdot\,\rangle\sharp\tau$,
    $\mu_3\leftarrow \mu_i$,
    $\pi_{12}\leftarrow \gamma^{(i-1)}$, 
    $\pi_{23}\leftarrow \gamma_i$,
    $\pi\leftarrow \gamma^{(i)}$ 
    in the notation of \citep[Lemma~7.6]{villani2003topics}.
    After $\gamma^{(i)}$ has been constructed for $i=1,\ldots,d$, let 
    $\pi$ be the marginal of $\gamma^{(d)}$ on $\BOmega\times\CS^*_2$.
    One checks that $\pi\in\Gamma(\mu_1,\ldots,\mu_d,\tau)$ satisfies the properties in statement~\ref{lems:dual-gluing-construction}.

    To prove statement~\ref{lems:dual-gluing-equalities},
    let us fix an arbitrary $\BIx\in\CS_1$ 
    and show that $\eqref{eqn:dual-gluing-equalities-MMOT}\le \eqref{eqn:dual-gluing-equalities-finitemeas}\le \eqref{eqn:dual-gluing-equalities-meas}\le \eqref{eqn:dual-gluing-equalities-MMOT}$.
    We first fix an arbitrary $\mu\in\Gamma(\mu_1,\ldots,\mu_d)$ and show that $\eqref{eqn:dual-gluing-equalities-MMOT}\le \eqref{eqn:dual-gluing-equalities-finitemeas}$. 
    Let $\Blambda_{\max}:\BOmega\to\CV_{\extremept}$ be a Borel measurable function that satisfies $C_2(\BIx,\Bomega)=\langle\BQ\BIx+\BP\Bomega+\nobreak\Bb,\Blambda_{\max}(\Bomega)\rangle$ for all $\Bomega\in\nobreak\BOmega$,
    which exists due to Lemma~\ref{lem:stage2-dual}\ref{lems:stage2-dual-finiteset}. 
    Moreover, let $\tau:=\Blambda_{\max}\sharp\mu\in\CP(\CV_{\extremept})$.
    Hence, $\langle\Bp_i,\cdot\,\rangle\sharp\tau\in\CP_1(\R)$ for $i=1,\ldots,d$, and 
    we have $\gamma_i:=\big(\mathrm{proj}_i,\langle\Bp_i,\cdot\Blambda_{\max}(\cdot)\rangle\big)\sharp\mu\in\Gamma\big(\mu_i,\langle\Bp_i,\cdot\,\rangle\sharp\tau\big)$, where $\mathrm{proj}_i:\BOmega\to\Omega_i$ denotes the function that projects its input onto the $i$-th coordinate. 
    Consequently, we get
    \begin{align*}
        \int_{\BOmega}C_2(\BIx,\Bomega)\DIFFM{\mu}{\DIFF\Bomega} &= \int_{\BOmega} \langle\BQ\BIx+\BP\Bomega+\Bb,\Blambda_{\max}(\Bomega)\rangle\DIFFM{\mu}{\DIFF\Bomega}\\ 
        &=\int_{\CS^*_2}\langle\BQ\BIx+\Bb,\Blambda\rangle\DIFFM{\tau}{\DIFF\Blambda} + \Bigg(\sum_{i=1}^d \int_{\Omega_i\times\R}z\omega\DIFFM{\gamma_i}{\DIFF \omega,\DIFF z}\Bigg)\\
        &\le\int_{\CS^*_2}\langle\BQ\BIx+\Bb,\Blambda\rangle\DIFFM{\tau}{\DIFF\Blambda} - \Bigg(\sum_{i=1}^d\OT_{\mu_i}\big(\langle\Bp_i,\cdot\,\rangle\sharp\tau\big)\Bigg)\le \eqref{eqn:dual-gluing-equalities-finitemeas}.
    \end{align*}
    Taking the supremum over all $\mu\in\Gamma(\mu_1,\ldots,\mu_d)$ in the inequalities above proves $\eqref{eqn:dual-gluing-equalities-MMOT}\le \eqref{eqn:dual-gluing-equalities-finitemeas}$. 
    Next, $\eqref{eqn:dual-gluing-equalities-finitemeas}\le \eqref{eqn:dual-gluing-equalities-meas}$ follows due to $\CV_{\extremept}\subset\CS^*_2$.
    To show that $\eqref{eqn:dual-gluing-equalities-meas}\le \eqref{eqn:dual-gluing-equalities-MMOT}$, 
    let us fix an arbitrary $\tau\in\CP_1(\CS^*_2)$, construct $\pi\in\CP(\BOmega\times\CS^*_2)$ by statement~\ref{lems:dual-gluing-construction},
    and let $\mu$ denote the marginal of $\pi$ on $\BOmega$.
    We subsequently get $\mu\in\Gamma(\mu_1,\ldots,\mu_d)$ as well as
    \begin{align}
        \begin{split}
            &\hspace{-20pt}\int_{\CS^*_2}\langle\BQ\BIx+\Bb,\Blambda\rangle\DIFFM{\tau}{\DIFF\Blambda} - \Bigg(\sum_{i=1}^d\OT_{\mu_i}\big(\langle\Bp_i,\cdot\,\rangle\sharp\tau\big)\Bigg) \\
            & = \int_{\BOmega\times\CS_2^*} \langle\BQ\BIx+\Bb,\Blambda\rangle + {\sum_{i=1}^d}\langle\Bp_i,\Blambda\rangle\omega_i \DIFFM{\pi}{\DIFF \omega_1,\ldots,\DIFF \omega_d,\DIFF\Blambda}\\
            & \le \int_{\BOmega} \max_{\Blambda\in\CS_2^*} \Big\{\langle\BQ\BIx+\Bb,\Blambda\rangle + {\textstyle\sum_{i=1}^d}\langle\Bp_i,\Blambda\rangle\omega_i \Big\} \DIFFM{\mu}{\DIFF\omega_1,\ldots,\DIFF\omega_d}\\
            & =\int_{\BOmega} C_2(\BIx,\Bomega)\DIFFM{\mu}{\DIFF\Bomega} \le \eqref{eqn:dual-gluing-equalities-MMOT}.
        \end{split}
        \label{eqn:dual-gluing-proof-ineq}
    \end{align}
    Taking the supremum over all $\tau\in\CP_1(\CS^*_2)$ in the inequalities above proves $\eqref{eqn:dual-gluing-equalities-meas}\le \eqref{eqn:dual-gluing-equalities-MMOT}$ and completes the proof of statement~\ref{lems:dual-gluing-equalities}.
    
    Lastly, to prove statement~\ref{lems:dual-gluing-optimality},
    observe that (\ref{eqn:dual-gluing-proof-ineq}) and statement~\ref{lems:dual-gluing-equalities} yield
    \begin{align*}
        \eqref{eqn:dual-gluing-equalities-meas}&=\int_{\CS^*_2}\langle\BQ\BIx+\Bb,\Blambda\rangle\DIFFM{\bar{\tau}}{\DIFF\Blambda} - \Bigg(\sum_{i=1}^d\OT_{\mu_i}\big(\langle\Bp_i,\cdot\,\rangle\sharp\bar{\tau}\big)\Bigg)\\
        &=\int_{\BOmega\times\CS^*_2}\langle\BQ\BIx+\BP\Bomega+\Bb,\Blambda\rangle\DIFFM{\bar{\pi}}{\DIFF\Bomega,\DIFF\Blambda}\le \int_{\BOmega} C_2(\BIx,\Bomega)\DIFFM{\bar{\mu}}{\DIFF\Bomega} \le \eqref{eqn:dual-gluing-equalities-MMOT}=\eqref{eqn:dual-gluing-equalities-meas},
    \end{align*}
    which shows that $\langle\BQ\BIx+\BP\Bomega+\Bb,\Blambda\rangle=C_2(\BIx,\Bomega)$ for $\bar{\pi}$-almost every $(\Bomega,\Blambda)\in\BOmega\times\CS^*_2$,
    and 
    $\bar{\mu}\in\argmax_{\mu\in\Gamma(\mu_1,\ldots,\mu_d)}\big\{\int_{\BOmega}C_2(\BIx,\Bomega)\DIFFM{\mu}{\DIFF\BIx}\big\}$.
    The proof is now complete.
\end{proof}

\begin{lemma}\label{lem:OT-function-convex}
    Under Setting~\ref{sett:dro-alt},
    let $K:=\nobreak|\CV_{\extremept}|<\infty$, let us enumerate $\CV_{\extremept}=\{\Blambda_1,\ldots,\Blambda_K\}$,
    and define $\Delta_K:=\big\{(\nu_1,\ldots,\nu_K)^\TRANSP\in\nobreak\R^K_{+}:\sum_{j=1}^K\nu_j=\nobreak1\big\}$.
    Moreover, let $\Phi:\R^K\to\R\union\{\infty\}$ be defined as follows:
    \begin{align*}
        \Phi(\nu_1,\ldots,\nu_{K})&:=\begin{cases}
            \sum_{i=1}^d\OT_{\mu_i}\big({\textstyle\sum_{j=1}^{K}\nu_j\delta_{\langle\Bp_i,\Blambda_j\rangle}}\big) & \mathrm{if~}(\nu_1,\ldots,\nu_{K})^\TRANSP\in\Delta_{K} \\
            \infty & \mathrm{if~}(\nu_1,\ldots,\nu_{K})^\TRANSP\notin\Delta_{K}
        \end{cases} \\
        & \hspace{180pt} \forall (\nu_1,\ldots,\nu_{K})^\TRANSP\in\R^{K}.
    \end{align*}
    Then, 
    the following statements hold.
    \begin{enumerate}[label=(\roman*),beginpenalty=10000]
        \item\label{lems:OT-function-convex-lscconvex}
        $\Phi$ is lower semi-continuous and convex.
        Moreover, for any $(\nu_1,\ldots,\nu_K)^\TRANSP\in\Delta_K$ and
        for $\tau:=\sum_{j=1}^{K}\nu_j\delta_{\Blambda_j}$, it holds that
        $\Phi(\nu_1,\ldots,\nu_{K})=\sum_{i=1}^d\OT_{\mu_i}\big(\langle\Bp_i,\cdot\,\rangle\sharp\tau\big)$.
        
        \item\label{lems:OT-function-convex-subgradient}
        Let $(\nu_1,\ldots,\nu_{K})^\TRANSP\in\Delta_{K}$, let $\tau:=\sum_{j=1}^{K}\nu_j\delta_{\Blambda_j}$, and let 
        $(g_1,\ldots,g_{K})^\TRANSP \in\partial\Phi(\nu_1,\ldots,\nu_{K})$.
        Then, there exist $f_1,\ldots,f_d\in\CC_{\mathrm{cvx}}(\R)$ satisfying ${-\int_{\Omega_i}}\sup_{z\in\R}\big\{z\omega-f_i(z)\big\}\DIFFM{\mu_i}{\DIFF\omega}-\int_{\R}f_i\DIFFX{\langle\Bp_i,\cdot\,\rangle\sharp\tau}=\OT_{\mu_i}\big(\langle\Bp_i,\cdot\,\rangle\sharp\tau\big)$, $\partial f_i(z)\subseteq\Omega_i$ $\forall z\in\R$ for $i=1,\ldots,d$, where 
        $g_j = {-\sum_{i=1}^{d}}f_i\big(\langle\Bp_i,\Blambda_j\rangle\big)$ for every $j$ with $\nu_j>0$,
        $g_j \le {-\sum_{i=1}^{d}}f_i\big(\langle\Bp_i,\Blambda_j\rangle\big)$ for every $j$ with $\nu_j=0$.
    \end{enumerate}
\end{lemma}

\begin{proof}[Proof of Lemma~\ref{lem:OT-function-convex}]
    Firstly, $\Phi$ is lower semi-continuous and convex due to the convexity of $\Delta_{K}$ and the representation of each summand in $\sum_{i=1}^d\OT_{\mu_i}\big({\textstyle\sum_{j=1}^{K}\nu_j\delta_{\langle\Bp_i,\Blambda_j\rangle}}\big)$ by Lemma~\ref{lem:OT-1d}\ref{lems:OT-1d-duality}.
    Moreover, for any $(\nu_1,\ldots,\nu_K)^\TRANSP\in\Delta_K$ and
    for $\tau:=\sum_{j=1}^{K}\nu_j\delta_{\Blambda_j}$,
    it holds for $i=1,\ldots,d$ that $\sum_{j=1}^{K}\nu_j\delta_{\langle\Bp_i,\Blambda_j\rangle}=\langle\Bp_i,\cdot\,\rangle\sharp\big(\sum_{j=1}^K\nu_j\delta_{\Blambda_j}\big)=\langle\Bp_i,\cdot\,\rangle\sharp\tau$,
    which shows that $\Phi(\nu_1,\ldots,\nu_{K})=\sum_{i=1}^d\OT_{\mu_i}\big(\langle\Bp_i,\cdot\,\rangle\sharp\tau\big)$.
    Thus, statement~\ref{lems:OT-function-convex-lscconvex} holds.

    Let us now prove statement~\ref{lems:OT-function-convex-subgradient}.
    For $i=1,\ldots,d$, let 
    $m_i:=\big|\big\{\langle\Bp_i,\Blambda\rangle:\Blambda\in\CV_{\extremept}\big\}\big|$, 
    enumerate 
    $\big\{\langle\Bp_i,\Blambda\rangle:\Blambda\in\CV_{\extremept}\big\}=\{\kappa_{i,1},\ldots,\kappa_{i,m_i}\}$
    where $\kappa_{i,1}<\kappa_{i,2}<\cdots<\kappa_{i,m_i}$,
    define $\Delta_{m_i}:=\big\{(\nu_1,\ldots,\nu_{m_i})^\TRANSP\in\R^{m_i}_{+}:\sum_{j=1}^{m_i}\nu_j=\nobreak1\big\}$,
    and define the compact set $B_i\subset\R^{m_i}$ as follows:
    \begin{align*}
        B_{i}:=\bigg\{(\tilde{f}_{i,1},\ldots,\tilde{f}_{i,m_i})^\TRANSP\in\R^{m_i}:\tilde{f}_{i,1}=0,\; \underline{\omega}_i \le \frac{\tilde{f}_{i,2}-\tilde{f}_{i,1}}{\kappa_{i,2}-\kappa_{i,1}}\le \cdots \le \frac{\tilde{f}_{i,m_i}-\tilde{f}_{i,m_i-1}}{\kappa_{i,m_i}-\kappa_{i,m_i-1}}\le \overline{\omega}_i\bigg\}.
    \end{align*}
    Moreover, for any $(\nu_{i,1},\ldots,\nu_{i,m_i})^\TRANSP\in\R^{m_i}$,
    let $\widetilde{\OT}_{\mu_i}(\nu_{i,1},\ldots,\nu_{i,m_i})\in\R$ denote the optimal value of the following concave maximization problem: 
    \begin{align}
        \max_{(\tilde{f}_{i,1},\ldots,\tilde{f}_{i,m_i})^\TRANSP\in B_i}\Bigg\{{-\int_{\Omega_i}}\max_{1\le j\le m_i}\big\{\kappa_{i,j}\omega- \tilde{f}_{i,j}\big\}\DIFFM{\mu_i}{\DIFF\omega}-\Bigg(\sum_{j=1}^{m_i}\tilde{f}_{i,j}\nu_{i,j}\Bigg)\Bigg\},
        \label{eqn:OT-function-convex-proof-maxrepresentation}
    \end{align}
    and let $B^\star_i(\nu_{i,1},\ldots,\nu_{i,m_i})\subset\R^{m_i}$ denote the set of maximizers,
    which is non-empty, convex, and compact.
    Thus, $\widetilde{\OT}_{\mu_i}:\R^{m_i}\to\R$ is a continuous and convex function.
    Furthermore, let $\BPi_i\in\R^{m_i\times K}$ be an $m_i$-by-$K$ matrix with entries $[\BPi_i]_{l,j}:=\INDI_{\big\{\kappa_{i,l}=\langle\Bp_i,\Blambda_j\rangle\big\}}$ for $l=1,\ldots,m_i$, $j=1,\ldots,K$. 
    Notice that for any $\Bnu=(\nu_1,\ldots,\nu_K)^\TRANSP\in\Delta_K$,
    it holds that $(\nu_{i,1},\ldots,\nu_{i,m_i})^\TRANSP:=\BPi_i\Bnu$ satisfies 
    $(\nu_{i,1},\ldots,\nu_{i,m_i})^\TRANSP\in\Delta_{m_i}$
    and $\sum_{l=1}^{m_i}\nu_{i,l}\delta_{\kappa_{i,l}}=\sum_{j=1}^K\nu_j\delta_{\langle\Bp_i,\Blambda_j\rangle}=\langle\Bp_i,\cdot\,\rangle\sharp\tau$.
    Consequently, Lemma~\ref{lem:OT-1d}\ref{lems:OT-1d-potential-compact} implies that we can re-express 
    $\Phi:\R^{K}\to\R\union\{\infty\}$ by 
    \begin{align*}
        \Phi(\Bnu)=\Bigg(\sum_{i=1}^{d}\widetilde{\OT}_{\mu_i}(\BPi_i\Bnu)\Bigg)+\chi_{\Delta_K}(\Bnu) \qquad \forall \Bnu\in\R^K,
    \end{align*}
    where 
    \begin{align*}
        \chi_{\Delta_{K}}(\Bnu)
        &:=\begin{cases}
            0 & \mathrm{if~}\Bnu\in\Delta_{K}\\
            \infty & \mathrm{if~}\Bnu \notin \Delta_{K}
        \end{cases} \hspace{41pt} \qquad \forall \Bnu\in\R^{K}.
    \end{align*}
    Since the term being maximized in (\ref{eqn:OT-function-convex-proof-maxrepresentation}) is jointly continuous in $\tilde{f}_{i,1},\ldots,\tilde{f}_{i,m_i},\nu_{i,1},\ldots,\nu_{i,m_i}$, we apply a classical result about the subdifferential of the point-wise supremum of convex functions (see, e.g., \citep{hantoute2008subdifferential}) to get 
    $\partial\widetilde{\OT}_{\mu_i}(\Bnu_i) 
    = {\clos\big(\conv\big({-B_i^\star}(\Bnu_i)\big)\big)}
    = {-B_i^\star}(\Bnu_i)$ for all $\Bnu_i\in\Delta_{m_i}$.
    Furthermore, since $\R^{K}\ni\Bnu\mapsto\widetilde{\OT}_{\mu_i}(\BPi_i\Bnu)\in\R$ is $\R$-valued on the entire $\R^K$ for $i=1,\ldots,d$, 
    \citep[Theorem~23.8 \& Theorem~23.9]{rockafellar1970convex}
    imply that 
    \begin{align}
        \partial \Phi(\Bnu) &= \Bigg(\sum_{i=1}^d\Big\{{\BPi_i^\TRANSP\BIg_i}:\BIg_i\in \partial\widetilde{\OT}_{\mu_i}(\BPi_i\Bnu)\Big\} \Bigg)
        + \partial \chi_{\Delta_{K}}(\Bnu) \label{eqn:OT-function-convex-proof-subdifferential-addition}\\
        &=\Bigg(\sum_{i=1}^d\Big\{{-\BPi_i^\TRANSP\tilde{\BIf}_i}:\tilde{\BIf}_i\in B^\star_{i}(\BPi_i\Bnu)\Big\} \Bigg)
        + \big\{\beta\vecone_{K}:\beta\in\R\big\} + \sum_{j\in\{j':\nu_{j'}=0\}}\big\{\eta_j\BIe_{j}:\eta_j\in\R_{-}\big\} \nonumber \\
        &\hspace{280pt} \forall \Bnu=(\nu_1,\ldots,\nu_K)^\TRANSP\in\Delta_{K},\nonumber%
    \end{align}%
    where $\BIe_{j}$ denotes the $j$-th standard basis vector in $\R^{K}$, that is, the $j$-th component of $\BIe_j$ is equal to~1 and all other components of $\BIe_j$ are equal to~0.

    In the following, let us fix an arbitrary $\Bnu=(\nu_1,\ldots,\nu_K)^\TRANSP\in\Delta_K$, 
    define $\tau:=\sum_{j=1}^K\nu_j\delta_{\Blambda_j}$,
    denote $\BPi_i\Bnu=(\nu_{i,1},\ldots,\nu_{i,m_i})^\TRANSP\in\Delta_{m_i}$ for $i=1,\ldots,d$,
    and fix an arbitrary $(g_1,\ldots,g_K)^\TRANSP\in \partial\Phi(\nu_1,\ldots,\nu_K)$.
    By (\ref{eqn:OT-function-convex-proof-subdifferential-addition}), there exist 
    $\tilde{\BIf}_i=(\tilde{f}_{i,1},\ldots,\tilde{f}_{i,m_i})^\TRANSP\in B_i^\star(\BPi_i\Bnu)$ for $i=1,\ldots,d$, 
    $\beta\in\R$, 
    and $\eta_j\in\R_{-}$ for every $j$ with $\nu_j=0$, 
    such that for $j=1,\ldots,K$, 
    $g_j=\Big(\sum_{i=1}^d -\big[\BPi_i^\TRANSP\tilde{\BIf}_i\big]_j\Big)+\beta$ if $\nu_j>0$, 
    $g_j=\Big(\sum_{i=1}^d -\big[\BPi_i^\TRANSP\tilde{\BIf}_i\big]_j\Big)+\beta+\eta_j$ if $\nu_j=0$.
    For $i=1,\ldots,d$,
    Lemma~\ref{lem:OT-1d}\ref{lems:OT-1d-potential-compact} implies the existence of $f_i\in\CC_{\mathrm{cvx}}(\R)$ such that $f_i(\kappa_{i,l})=\tilde{f}_{i,l}-\frac{\beta}{d}$ for $l=1,\ldots,m_i$, 
    $\partial f_i(z)\subseteq\Omega_i$ for all $z\in\R$, 
    and 
    ${-\int_{\Omega_i}}\sup_{z\in\R}\big\{z\omega-f_i(z)\big\}\DIFFM{\mu_i}{\DIFF\omega}-\int_{\R}f_i\DIFFX{\langle\Bp_i,\cdot\,\rangle\sharp\tau}=\OT_{\mu_i}\big(\langle\Bp_i,\cdot\,\rangle\sharp\tau\big)$.
    Moreover, it holds for $i=1,\ldots,d$, $j=1,\ldots,K$ that
    $\big[\BPi_i^\TRANSP\tilde{\BIf}_i\big]_j=\tilde{f}_{i,l}$ where $l\in\{1,\ldots,m_i\}$ satisfies $\langle\Bp_i,\Blambda_j\rangle=\kappa_{i,l}$.
    It follows that 
    $\Big(\sum_{i=1}^d -\big[\BPi_i^\TRANSP\tilde{\BIf}_i\big]_j\Big)+\beta={-\sum_{i=1}^{d}}f_i\big(\langle\Bp_i,\Blambda_j\rangle\big)$ 
    for $j=1,\ldots,K$.
    We have thus shown that
    $g_j={-\sum_{i=1}^{d}}f_i\big(\langle\Bp_i,\Blambda_j\rangle\big)$ for every $j$ with $\nu_j>0$, and
    $g_j\le {-\sum_{i=1}^{d}}f_i\big(\langle\Bp_i,\Blambda_j\rangle\big)$ for every $j$ with $\nu_j=0$.
    The proof is now complete.
\end{proof}

\begin{lemma}\label{lem:global-simplification}
    Under Setting~\ref{sett:dro-alt}, for $i=1,\ldots,d$, let $f_i\in\CC_{\mathrm{cvx}}(\R)$ satisfy $\partial f_i(z)\subseteq\Omega_i$ for all $z\in\R$.
    Then, it holds that
    \begin{align*}
        \sup_{\Bzeta\in\CS^*_2}\Bigg\{\langle\BQ\BIx+\Bb,\Bzeta\rangle + \sum_{i=1}^{d}f_i\big(\langle\Bp_i,\Bzeta\rangle\big)\Bigg\}=\max_{\Blambda\in\CV_{\extremept}}\Bigg\{\langle\BQ\BIx+\Bb,\Blambda\rangle + \sum_{i=1}^{d}f_i\big(\langle\Bp_i,\Blambda\rangle\big)\Bigg\} \quad\;\; \forall \BIx\in\CS_1.
    \end{align*}
\end{lemma}

\begin{proof}[Proof of Lemma~\ref{lem:global-simplification}]
    Let us fix an arbitrary $\BIx\in\CS_1$. 
    It follows from \citep[p.170 \& Theorem~18.5 \& Theorem~19.1]{rockafellar1970convex} that the polyhedral convex set $\CS^*_2$ can be expressed as $\CS^*_2=\conv(\CV_{\extremept})+\recession(\CS^*_2)$.
    Consequently, let us express an arbitrary $\Bzeta\in\CS^*_2$ as 
    $\Bzeta=\Blambda+\Brho$ where $\Blambda\in\conv(\CV_{\extremept})$ and $\Brho\in\recession(\CS^*_2)$. 
    Taking $\omega_i\in\partial f_i\big(\langle\Bp_i,\Bzeta\rangle\big)\subseteq\Omega_i$ for $i=1,\ldots,d$, 
    we get $f_i\big(\langle\Bp_i,\Bzeta\rangle\big)\le f_i\big(\langle\Bp_i,\Blambda\rangle\big) + \langle\Bp_i,\Brho\rangle\omega_i$ for $i=1,\ldots,d$.
    Denoting $\Bomega:=(\omega_1,\ldots,\omega_d)^\TRANSP\in\BOmega$, 
    and recalling that the proof of Lemma~\ref{lem:stage2-dual}\ref{lems:stage2-dual-finiteset} has shown that $\langle\BQ\BIx+\BP\Bomega+\Bb,\Brho\rangle\le 0$, we subsequently get
    $\langle\BQ\BIx+\Bb,\Bzeta\rangle + \sum_{i=1}^{d} f_i\big(\langle\Bp_i,\Bzeta\rangle\big) \le \langle\BQ\BIx+\Bb,\Blambda\rangle + \langle\BQ\BIx+\Bb,\Brho\rangle + \sum_{i=1}^{d} f_i\big(\langle\Bp_i,\Blambda\rangle\big)+\langle\Bp_i,\Brho\rangle\omega_i 
    \le \langle\BQ\BIx+\Bb,\Blambda\rangle + \sum_{i=1}^{d} f_i\big(\langle\Bp_i,\Blambda\rangle\big)$.
    This shows that 
    $\sup_{\Bzeta\in\CS^*_2}\big\{\langle\BQ\BIx+\Bb,\Bzeta\rangle + \sum_{i=1}^{d}f_i\big(\langle\Bp_i,\Bzeta\rangle\big)\big\}=\max_{\Blambda\in\conv(\CV_{\extremept})}\big\{\langle\BQ\BIx+\Bb,\Blambda\rangle + \sum_{i=1}^{d}f_i\big(\langle\Bp_i,\Blambda\rangle\big)\big\}$, where the latter maximum can be attained in $\CV_{\extremept}$ due to the convexity of $\conv(\CV_{\extremept})\ni\Blambda\mapsto \langle\BQ\BIx+\Bb,\Blambda\rangle + \sum_{i=1}^{d}f_i\big(\langle\Bp_i,\Blambda\rangle\big)\in\R$.
    The proof is now complete.
\end{proof}

\begin{lemma}\label{lem:weak-duality}
    Under Setting~\ref{sett:dro-alt}, it holds for any $(\BIx,f_1,\ldots,f_d)$ feasible for \eqref{eqn:dro-primal} and any $(\tau,\Bxi_{\ineq},\Bxi_{\eq})$ feasible for \eqref{eqn:dro-dual} that
    \begin{align*}
        &\phantom{=}\,\;\Bigg(\langle\Bc_1,\BIx\rangle+\sum_{i=1}^d\int_{\Omega_i}\sup_{z\in\R}\big\{z\omega-f_{i}(z)\big\}\DIFFM{\mu_i}{\DIFF\omega}\Bigg) \\
        & \qquad - \Bigg(\int_{\CS^*_2}\langle\Bb,\Blambda\rangle\DIFFM{\tau}{\DIFF\Blambda} 
        - \Bigg(\sum_{i=1}^d\OT_{\mu_i}\big(\langle\Bp_i,\cdot\,\rangle\sharp\tau\big)\Bigg)
        + \langle\Bq_{\ineq},\Bxi_{\ineq}\rangle+\langle\Bq_{\eq},\Bxi_{\eq}\rangle\Bigg)\\
        &= \langle\BL_{\ineq}\BIx-\Bq_{\ineq}, \Bxi_{\ineq}\rangle + \big\langle\Bc_1+{\textstyle\int_{\CS^*_2}\BQ^\TRANSP\Blambda}\DIFFM{\tau}{\DIFF\Blambda}-\BL_{\ineq}^\TRANSP\Bxi_{\ineq}-\BL_{\eq}^{\TRANSP}\Bxi_{\eq},\BIx\big\rangle\\
        & \qquad - \int_{\CS^*_2}\langle\BQ\BIx+\Bb,\Blambda\rangle+{\textstyle\sum_{i=1}^{d}}f_i\big(\langle\Bp_i,\Blambda\rangle\big)\DIFFM{\tau}{\DIFF\Blambda}\\
        & \qquad + \sum_{i=1}^{d}\bigg(\OT_{\mu_i}\big(\langle\Bp_i,\cdot\,\rangle\sharp\tau\big)+\int_{\Omega_i}\sup_{z\in\R}\big\{z\omega-f_i(z)\big\}\DIFFM{\mu_i}{\DIFF\omega}+\int_{\R}f_i\DIFFX{\langle\Bp_i,\cdot\,\rangle\sharp\tau}\bigg),
    \end{align*}
    where 
    $\langle\BL_{\ineq}\BIx-\Bq_{\ineq}, \Bxi_{\ineq}\rangle\ge 0$, 
    $\big\langle\Bc_1+{\textstyle\int_{\CS^*_2}\BQ^\TRANSP\Blambda}\DIFFM{\tau}{\DIFF\Blambda}-\BL_{\ineq}^\TRANSP\Bxi_{\ineq}-\BL_{\eq}^{\TRANSP}\Bxi_{\eq},\BIx\big\rangle\ge 0$,
    $\int_{\CS^*_2}\langle\BQ\BIx+\nobreak\Bb,\Blambda\rangle+{\textstyle\sum_{i=1}^{d}}f_i\big(\langle\Bp_i,\Blambda\rangle\big)\DIFFM{\tau}{\DIFF\Blambda}\le 0$, 
    and, for $i=1,\ldots,d$,
    $\OT_{\mu_i}\big(\langle\Bp_i,\cdot\,\rangle\sharp\tau\big)+\int_{\Omega_i}\sup_{z\in\R}\big\{z\omega-f_i(z)\big\}\DIFFM{\mu_i}{\DIFF\omega}+\int_{\R}f_i\DIFFX{\langle\Bp_i,\cdot\,\rangle\sharp\tau}\ge 0$.
\end{lemma}

\begin{proof}[Proof of Lemma~\ref{lem:weak-duality}]
    The identity can be directly verified by reorganizing the right-hand side. 
    Moreover, 
    $\langle\BL_{\ineq}\BIx-\Bq_{\ineq}, \Bxi_{\ineq}\rangle\ge 0$, 
    $\big\langle\Bc_1+{\textstyle\int_{\CS^*_2}\BQ^\TRANSP\Blambda}\DIFFM{\tau}{\DIFF\Blambda}-\BL_{\ineq}^\TRANSP\Bxi_{\ineq}-\BL_{\eq}^{\TRANSP}\Bxi_{\eq},\BIx\big\rangle\ge 0$, and
    $\int_{\CS^*_2}\langle\BQ\BIx+\nobreak\Bb,\Blambda\rangle+{\textstyle\sum_{i=1}^{d}}f_i\big(\langle\Bp_i,\Blambda\rangle\big)\DIFFM{\tau}{\DIFF\Blambda}\le 0$
    hold due to the constraints in \eqref{eqn:dro-primal} and \eqref{eqn:dro-dual}.
    Lastly, for $i=1,\ldots,d$, 
    $\OT_{\mu_i}\big(\langle\Bp_i,\cdot\,\rangle\sharp\tau\big)+\int_{\Omega_i}\sup_{z\in\R}\big\{z\omega-f_i(z)\big\}\DIFFM{\mu_i}{\DIFF\omega}+\int_{\R}f_i\DIFFX{\langle\Bp_i,\cdot\,\rangle\sharp\tau}\ge \nobreak0$ follows from Lemma~\ref{lem:OT-1d}\ref{lems:OT-1d-duality} with respect to $\mu\leftarrow\mu_i$, $\nu\leftarrow\langle\Bp_i,\cdot\,\rangle\sharp\tau$.
    The proof is complete.
\end{proof}

\begin{proof}[Proof of Theorem~\ref{thm:primaldual}]
    Let us begin the proof of statement~\ref{thms:primaldual-strongduality} by showing that $\eqref{eqn:dro-primal}\le C_{\DRO}$.
    Since Lemma~\ref{lem:stage2-dual}\ref{lems:stage2-dual-finiteset}
    has shown that $C_2(\BIx,\Bomega)=\max_{\Blambda\in \CV_{\extremept}}\big\{\langle\BQ\BIx+\BP\Bomega+\nobreak\Bb,\Blambda\rangle\big\}$ 
    for all \mbox{$\BIx\in\CS_1$} and all $\Bomega\in\BOmega$,
    it follows that the function $\BOmega\ni \Bomega\mapsto C_2(\BIx,\cdot\,)\in\R$ is continuous for every $\BIx\in\CS_1$. 
    Subsequently, it follows from an extension of the classical Kantorovich duality (see, e.g., \citep[Theorem~1.3]{villani2003topics}) to the multi-marginal case that
    \begin{align}
        \sup_{\mu\in\Gamma(\mu_1,\ldots,\mu_d)}\bigg\{\int_{\BOmega} C_2(\BIx,\Bomega)\DIFFM{\mu}{\DIFF\Bomega}\bigg\}
        &=\inf_{\substack{h_i\in\CL^1(\mu_i)\,\forall 1\le i\le d, \\ \bigoplus_{i=1}^d h_i \ge C_2(\BIx,\,\cdot\,)}}\Bigg\{\sum_{i=1}^d \int_{\Omega_i}h_i\DIFFX{\mu_i}\Bigg\} \qquad \forall \BIx\in\CS_1.
        \label{eqn:primaldual-proof-kantorovichduality}
    \end{align}
    where
    \begin{align*}
        \CL^1(\mu_i):=\big\{h:\Omega_i\to\R\union\{\infty\}: h \text{ is Borel measurable and } \textstyle\int_{\Omega_i}|h|\DIFFX{\mu_i}<\infty\big\} \qquad \forall 1\le i\le d,
    \end{align*}
    and $\bigoplus_{i=1}^d h_i(\Bomega):=\sum_{i=1}^d h_i(\omega_i)$ for all $\Bomega=(\omega_1,\ldots,\omega_d)^\TRANSP\in\BOmega$.
    Now, let us fix an arbitrary $\BIx\in\nobreak\CS_1$ and fix arbitrary $h_1\in\CL^1(\mu_1),\ldots,h_d\in\CL^1(\mu_d)$ such that $\bigoplus_{i=1}^d h_i\ge C_2(\BIx,\cdot\,)$.
    For $i=1,\ldots,d$, let $h_i^*(z):=\sup_{\omega_i\in\Omega_i}\big\{z\omega_i-h_i(\omega_i)\big\}$ $\forall z\in\R$ denote the convex conjugate of $h_i$.
    It thus holds that $h_i^*\in\CC_{\mathrm{cvx}}(\R)$.
    Observe that 
    \begin{align*}
        \sum_{i=1}^d h_i(\omega_i)\ge \langle\BQ\BIx+\Bb,\Blambda\rangle + \sum_{i=1}^d \langle\Bp_i,\Blambda\rangle\omega_i \qquad \forall \Bomega=(\omega_1,\ldots,\omega_d)^\TRANSP\in\BOmega,\; \forall \Blambda\in\CS^*_2.
    \end{align*}
    Hence, it holds that
    \begin{align*}
        0\ge \langle\BQ\BIx+\Bb,\Blambda\rangle+\sum_{i=1}^d\sup_{\omega_i\in\Omega_i}\big\{\langle\Bp_i,\Blambda\rangle\omega_i-h_i(\omega_i)\big\} = \langle\BQ\BIx+\Bb,\Blambda\rangle+\sum_{i=1}^d h_i^*\big(\langle\Bp_i,\Blambda\rangle\big) \quad \forall \Blambda\in\CS^*_2.
    \end{align*}
    Since $\BIx$ satisfies $\BL_{\ineq}\BIx\le \Bq_{\ineq}$ and $\BL_{\eq}\BIx=\Bq_{\eq}$, this shows that $(\BIx,h_1^*,\ldots,h_d^*)$ is feasible for \eqref{eqn:dro-primal}. 
    Moreover, it holds that
    \begin{align*}
        \sup_{z\in\R}\big\{z\omega_i-h_i^*(z)\big\} &= \sup_{z\in\R}\bigg\{z\omega_i-\sup_{\omega_i'\in\Omega_i}\big\{z\omega_i'-h_i(\omega_i')\big\}\bigg\}\le h_i(\omega_i) \qquad \forall \omega_i\in\Omega_i,\;\forall 1\le i\le d,
    \end{align*}
    which yields
    \begin{align*}
        \langle\Bc_1,\BIx\rangle+\sum_{i=1}^d\int_{\Omega_i}\sup_{z\in\R}\big\{z\omega_i-h_i^*(z)\big\} \DIFFM{\mu_i}{\DIFF\omega_i}\le \langle\Bc_1,\BIx\rangle+ \sum_{i=1}^d \int_{\Omega_i}h_i\DIFFX{\mu_i}.
    \end{align*}
    Taking the infimum on both sides of the above inequality first over all $h_1\in\CL^1(\mu_1),\ldots,h_d\in\CL^1(\mu_d)$ satisfying $\bigoplus_{i=1}^d h_i\ge C_2(\BIx,\cdot\,)$ and then over all $\BIx\in\CS_1$ leads to
    \begin{align*}
        \eqref{eqn:dro-primal} &\le \inf_{\substack{\BIx\in\CS_1,\,h_i\in\CL^1(\mu_i)\,\forall 1\le i\le d, \\ \bigoplus_{i=1}^d h_i \ge C_2(\BIx,\,\cdot\,)}} \Bigg\{\langle\Bc_1,\BIx\rangle+\sum_{i=1}^d\int_{\Omega_i}\sup_{z\in\R}\big\{z\omega_i-h_i^*(z)\big\} \DIFFM{\mu_i}{\DIFF\omega_i}\Bigg\}\\
        &\le \inf_{\BIx\in\CS_1}\left\{\langle\Bc_1,\BIx\rangle+ \inf_{\substack{h_i\in\CL^1(\mu_i)\,\forall 1\le i\le d, \\ \bigoplus_{i=1}^d h_i \ge C_2(\BIx,\,\cdot\,)}}\Bigg\{\sum_{i=1}^d \int_{\Omega_i}h_i\DIFFX{\mu_i}\Bigg\}\right\}.
    \end{align*}
    Applying (\ref{eqn:primaldual-proof-kantorovichduality}) then proves $\eqref{eqn:dro-primal}\le C_{\DRO}$.

    Conversely, $\eqref{eqn:dro-primal}\ge C_{\DRO}$ follows directly from Lemma~\ref{lem:primal-objective}.
    We hence conclude that the optimal value of \eqref{eqn:dro-primal} is equal to $C_{\DRO}$.
    Furthermore, 
    Let $\alpha:=\sup_{\BIx\in\CS_1,\,\Bomega\in\BOmega}\big\{C_2(\BIx,\Bomega)\big\}<\nobreak\infty$ be given by the assumption~\ref{setts:dro-stage2-alt} in Setting~\ref{sett:dro-alt}.
    Subsequently, let us take an arbitrary $\BIx\in\CS_1$, take an arbitrary $\Bomega=(\omega_1,\ldots,\omega_d)^\TRANSP\in\BOmega$, 
    and define $f_i(z):=\omega_iz$ $\forall z\in\R$ for $i=1,\ldots,d-1$, $f_d(z):=\omega_dz-\alpha$ $\forall z\in\R$.
    Then, $f_1,\ldots,f_d\in\CC_{\mathrm{cvx}}(\R)$, and it holds for all $\Blambda\in\CS^*_2$ that
    $\langle\BQ\BIx+\nobreak\Bb,\Blambda\rangle+\sum_{i=1}^df_i\big(\langle\Bp_i,\Blambda\rangle\big)=\langle\BQ\BIx+\BP\Bomega+\nobreak\Bb,\Blambda\rangle-\alpha\le C_2(\BIx,\Bomega)-\alpha\le 0$.
    Hence, $(\BIx,f_1,\ldots,f_d)$ is feasible for \eqref{eqn:dro-primal}, showing that $C_{\DRO}<\infty$.

    Next, to show that $\eqref{eqn:dro-dual}\le C_{\DRO}$, let us fix 
    an arbitrary $\BIx\in\CS_1$, as well as arbitrary $\tau\in\CP_1(\CS^*_2)$, $\Bxi_{\ineq}\in\R^{n_{\ineq}}_{-}$, $\Bxi_{\eq}\in\R^{n_{\eq}}$ that satisfy 
    $\Bc_1+\int_{\CS^*_2}\BQ^\TRANSP\Blambda\DIFFM{\tau}{\DIFF\Blambda}\ge \BL_{\ineq}^\TRANSP\Bxi_{\ineq}+\BL_{\eq}^\TRANSP\Bxi_{\eq}$.
    It follows from the equality $\eqref{eqn:dual-gluing-equalities-MMOT}=\eqref{eqn:dual-gluing-equalities-meas}$ in Lemma~\ref{lem:dual-gluing}\ref{lems:dual-gluing-equalities} that
    \begin{align}
        \begin{split}
        &\hspace{-20pt}\langle\Bc_1,\BIx\rangle+\sup_{\mu\in\Gamma(\mu_1,\ldots,\mu_d)}\bigg\{\int_{\BOmega}C_2(\BIx,\Bomega)\DIFFM{\mu}{\DIFF\Bomega}\bigg\}\\
        &\ge \langle\Bc_1,\BIx\rangle+\int_{\CS^*_2}\langle\BQ\BIx+\Bb,\Blambda\rangle\DIFFM{\tau}{\DIFF\Blambda} -\Bigg(\sum_{i=1}^d\OT_{\mu_i}\big(\langle\Bp_i,\cdot\,\rangle\sharp\tau\big)\Bigg)\\
        &=\big\langle\Bc_1+{\textstyle\int_{\CS^*_2}\BQ^\TRANSP\Blambda\DIFFM{\tau}{\DIFF\Blambda}},\BIx\big\rangle + \int_{\CS^*_2}\langle\Bb,\Blambda\rangle\DIFFM{\tau}{\DIFF\Blambda} - \Bigg(\sum_{i=1}^d\OT_{\mu_i}\big(\langle\Bp_i,\cdot\,\rangle\sharp\tau\big)\Bigg).
        \end{split}
        \label{eqn:primaldual-proof-weakduality-step1}
    \end{align}
    Since $\BIx$ satisfies $\BIx\in\R_{+}^{k_1}$, $\BL_{\ineq}\BIx\le \Bq_{\ineq}$, $\BL_{\eq}\BIx=\Bq_{\eq}$, we have $\big\langle\Bc_1+{\textstyle\int_{\CS^*_2}\BQ^\TRANSP\Blambda\DIFFM{\tau}{\DIFF\Blambda}},\BIx\big\rangle \ge \langle\BL_{\ineq}^\TRANSP\Bxi_{\ineq}+\BL_{\eq}^\TRANSP\Bxi_{\eq},\BIx\rangle\ge \langle\Bq_{\ineq},\Bxi_{\ineq}\rangle + \langle\Bq_{\eq},\Bxi_{\eq}\rangle$.
    Combining this with (\ref{eqn:primaldual-proof-weakduality-step1}), and then taking the infimum over all $\BIx\in\CS_1$ and the supremum over all $\tau\in\CP_1(\CS^*_2)$, $\Bxi_{\ineq}\in\R^{n_{\ineq}}_{-}$, $\Bxi_{\eq}\in\R^{n_{\eq}}$ satisfying
    $\Bc_1+\int_{\CS^*_2}\BQ^\TRANSP\Blambda\DIFFM{\tau}{\DIFF\Blambda}\ge \BL_{\ineq}^\TRANSP\Bxi_{\ineq}+\BL_{\eq}^\TRANSP\Bxi_{\eq}$ 
    proves $\eqref{eqn:dro-dual}\le C_{\DRO}$.

    It remains to show that $\eqref{eqn:dro-dual}\ge C_{\DRO}$.
    Observe that (\ref{eqn:dual-gluing-equalities-finitemeas}) in Lemma~\ref{lem:dual-gluing}\ref{lems:dual-gluing-equalities} is a maximization over probability measures on the finite set $\CV_{\extremept}$.
    Let us denote $K:=|\CV_{\extremept}|$, enumerate $\CV_{\extremept}=\{\Blambda_1,\ldots,\Blambda_K\}$,
    and define $\Delta_K\subset\R^K_{+}$, $\Phi:\R^K\to\R\union\{\infty\}$ by Lemma~\ref{lem:OT-function-convex}.
    It holds by Lemma~\ref{lem:OT-function-convex}\ref{lems:OT-function-convex-lscconvex} that $\Phi$ is a lower semi-continuous and convex function, and it holds that $\Phi(\nu_1,\ldots,\nu_K)=\sum_{i=1}^{d}\OT_{\mu_i}\big(\langle\Bp_i,\cdot\,\rangle\sharp\tau\big)$ for every $(\nu_1,\ldots,\nu_K)^\TRANSP\in\nobreak\Delta_K$ and for $\tau:=\sum_{j=1}^K \nu_j\delta_{\Blambda_j}$.
    Hence, the mapping 
    $\CS_1\times \Delta_K \ni (\BIx,\nu_1,\ldots,\nu_K) \mapsto \langle\Bc_1,\BIx\rangle+\big(\sum_{j=1}^K\langle\BQ\BIx+\nobreak\Bb,\Blambda_j\rangle\nu_{j}\big) - \Phi(\nu_1,\ldots,\nu_K)\in\R$
    is affine in $\BIx$ for every $(\nu_1,\ldots,\nu_K)^\TRANSP\in\Delta_K$ and is upper semi-continuous and concave in $(\nu_1,\ldots,\nu_K)^\TRANSP$ for every $\BIx\in\CS_1$. 
    Therefore, using Lemma~\ref{lem:dual-gluing}\ref{lems:dual-gluing-equalities}, the compactness of $\Delta_K$, and Sion's minimax theorem \citep{sion1958on} leads to
    \begin{align}
        C_{\DRO}&=\!\inf_{\BIx\in\CS_1}\! \Bigg\{\!\langle\Bc_1,\BIx\rangle + \sup_{\tau\in\CP(\CV_{\extremept})}\!\Bigg\{\int_{\CS^*_2}\langle\BQ\BIx+\Bb,\Blambda\rangle\DIFFM{\tau}{\DIFF\Blambda} -\Bigg(\sum_{i=1}^d\OT_{\mu_i}\big(\langle\Bp_i,\cdot\,\rangle\sharp\tau\big)\Bigg)\Bigg\}\!\Bigg\} \nonumber\\
        \begin{split}
        &=\!\inf_{\BIx\in\CS_1}\!\Bigg\{\!\sup_{(\nu_1,\ldots,\nu_K)^\TRANSP\in\Delta_K}\!\Bigg\{\langle\Bc_1,\BIx\rangle+\Bigg(\sum_{j=1}^{K}\langle\BQ\BIx+\Bb,\Blambda\rangle\nu_{j}\Bigg)-\Phi(\nu_1,\ldots,\nu_K)\Bigg\}\!\Bigg\}\\
        &=\!\sup_{(\nu_1,\ldots,\nu_K)^\TRANSP\in\Delta_K}\!\Bigg\{\!\inf_{\BIx\in\CS_1}\!\Bigg\{\langle\Bc_1,\BIx\rangle+\Bigg(\sum_{j=1}^{K}\langle\BQ\BIx+\Bb,\Blambda\rangle\nu_{j}\Bigg)-\Phi(\nu_1,\ldots,\nu_K)\Bigg\}\!\Bigg\}
        \end{split}\label{eqn:primaldual-proof-dualminimax}\\
        &=\!\sup_{\tau\in\CP(\CV_{\extremept})}\!\!\Bigg\{\!\!\int_{\CS^*_2}\!\langle\Bb,\Blambda\rangle\DIFFM{\tau}{\DIFF\Blambda}-\!\Bigg(\!\sum_{i=1}^{d}\!\OT_{\mu_i}\!\big(\langle\Bp_i,\cdot\,\rangle\sharp\tau\big)\!\!\Bigg)\!+\!\inf_{\BIx\in\CS_1}\!\!\Big\{\!\big\langle\Bc_1+{\textstyle\int_{\CS^*_2}\BQ^\TRANSP\Blambda\DIFFM{\tau}{\DIFF\Blambda}},\BIx\big\rangle\!\Big\}\!\Bigg\}.\nonumber
    \end{align}
    Since $\CS_1$ is non-empty by the assumption~\ref{setts:dro-stage1} in Setting~\ref{sett:dro}, 
    for every $\tau\in\CP(\CV_{\extremept})$, it holds that $\inf_{\BIx\in\CS_1}\Big\{\big\langle\Bc_1+{\textstyle\int_{\CS^*_2}\BQ^\TRANSP\Blambda\DIFFM{\tau}{\DIFF\Blambda}},\BIx\big\rangle\Big\}$ corresponds to a feasible LP problem:
    \begin{align}
        \label{eqn:primaldual-proof-innerLPmin}
        \begin{split}
            \minimize_{\BIx}\quad & \big\langle\Bc_1+{\textstyle\int_{\CS^*_2}\BQ^\TRANSP\Blambda\DIFFM{\tau}{\DIFF\Blambda}},\BIx\big\rangle \\
            \mathrm{subject~to} \quad & \BL_{\ineq}\BIx\le \Bq_{\ineq}, \quad \BL_{\eq}\BIx=\Bq_{\eq}, \quad \BIx\in\R_{+}^{k_1},
        \end{split}
    \end{align}
    whose dual LP problem is:
    \begin{align}
        \label{eqn:primaldual-proof-innerLPmax}
        \begin{split}
            \maximize_{\Bxi_{\ineq},\,\Bxi_{\eq}}\quad & \langle\Bq_{\ineq},\Bxi_{\ineq}\rangle + \langle\Bq_{\eq},\Bxi_{\eq}\rangle \\
            \mathrm{subject~to} \quad & \BL_{\ineq}^\TRANSP\Bxi_{\ineq} + \BL_{\eq}^\TRANSP\Bxi_{\eq} \le \Bc_1+{\textstyle\int_{\CS^*_2}\BQ^\TRANSP\Blambda\DIFFM{\tau}{\DIFF\Blambda}}, \\
            & \Bxi_{\ineq} \in\R_{-}^{n_{\ineq}}, \quad \Bxi_{\eq} \in \R^{n_{\eq}}.
        \end{split}
    \end{align}
    It then follows from the feasibility of (\ref{eqn:primaldual-proof-innerLPmin}) and the strong duality of LP problems that the optimal values of (\ref{eqn:primaldual-proof-innerLPmin}) and (\ref{eqn:primaldual-proof-innerLPmax}) are identical; note that $\eqref{eqn:primaldual-proof-innerLPmin}=\eqref{eqn:primaldual-proof-innerLPmax}=-\infty$ is possible.
    Now, consider 
    the following restriction of \eqref{eqn:dro} where the maximization over $\tau$ is restricted to $\CP(\CV_{\extremept})$:
    \begin{align}
        \begin{split}
            \maximize_{\tau,\,\Bxi_{\ineq},\,\Bxi_{\eq}} \quad & \int_{\CS^*_2}\langle\Bb,\Blambda\rangle\DIFFM{\tau}{\DIFF\Blambda} 
            - \Bigg(\sum_{i=1}^d\OT_{\mu_i}\big(\langle\Bp_i,\cdot\,\rangle\sharp\tau\big)\Bigg)
            +\langle\Bq_{\ineq},\Bxi_{\ineq}\rangle
            +\langle\Bq_{\eq},\Bxi_{\eq}\rangle  \\
            \mathrm{subject~to} \quad & \Bc_1+\int_{\CS_2^*}\BQ^\TRANSP\Blambda \DIFFM{\tau}{\DIFF\Blambda} \ge \BL_{\ineq}^\TRANSP\Bxi_{\ineq} + \BL_{\eq}^\TRANSP\Bxi_{\eq},\\
            & \tau\in\CP(\CV_{\extremept}), \qquad \Bxi_{\ineq} \in \R^{n_{\ineq}}_{-},\qquad \Bxi_{\eq}\in\R^{n_{\eq}}.
        \end{split}
        \label{eqn:primaldual-proof-dual-finitesupport} 
    \end{align}
    Substituting the equality $\eqref{eqn:primaldual-proof-innerLPmin}=\eqref{eqn:primaldual-proof-innerLPmax}$ into (\ref{eqn:primaldual-proof-dualminimax}) proves $C_{\DRO}=\eqref{eqn:primaldual-proof-dual-finitesupport}\le \eqref{eqn:dro-dual}$.
    The proof of statement~\ref{thms:primaldual-strongduality} is complete.
    In particular, notice that this also proves $\eqref{eqn:primaldual-proof-dual-finitesupport}=\eqref{eqn:dro-dual}$ 
    as well as
    \begin{subequations}%
        \label{eqn:primaldual-proof-minimax}
        \begin{align}
            C_{\DRO}&=\inf_{\BIx\in\CS_1}\sup_{\tau\in\CP_1(\CS^*_2)}\! \Bigg\{\!\langle\Bc_1,\BIx\rangle + \int_{\CS^*_2}\langle\BQ\BIx+\Bb,\Blambda\rangle\DIFFM{\tau}{\DIFF\Blambda} -\!\Bigg(\!\sum_{i=1}^d\OT_{\mu_i}\big(\langle\Bp_i,\cdot\,\rangle\sharp\tau\big)\!\Bigg)\!\Bigg\}
            \label{eqn:primaldual-proof-minimax-infsup}\\
            &=\sup_{\tau\in\CP_1(\CS^*_2)}\inf_{\BIx\in\CS_1}\! \Bigg\{\!\langle\Bc_1,\BIx\rangle + \int_{\CS^*_2}\langle\BQ\BIx+\Bb,\Blambda\rangle\DIFFM{\tau}{\DIFF\Blambda} -\!\Bigg(\!\sum_{i=1}^d\OT_{\mu_i}\big(\langle\Bp_i,\cdot\,\rangle\sharp\tau\big)\!\Bigg)\!\Bigg\}.
            \label{eqn:primaldual-proof-minimax-supinf}
        \end{align}
    \end{subequations}

    To prove statement~\ref{thms:primaldual-boundedness},
    notice that since $\eqref{eqn:dro-dual}=C_{\DRO}$, we have $C_{\DRO}>-\infty$ if and only if \eqref{eqn:dro-dual} is feasible. 
    Let us first assume that \eqref{eqn:dro-dual} admits a feasible solution $(\tau,\Bxi_{\ineq},\Bxi_{\eq})$, and 
    let $\hat{\Blambda}:=\int_{\CS^*_2}\Blambda\DIFFM{\tau}{\DIFF\Blambda}$.
    It holds that $\hat{\Blambda}\in\CS^*_2$ due to the convexity of $\CS^*_2$. 
    Subsequently, we get 
    $\Bc_1+\BQ^\TRANSP\hat{\Blambda}
    =\Bc_1+\int_{\CS^*_2}\BQ^\TRANSP\Blambda\DIFFM{\delta_{\hat{\Blambda}}}{\DIFF\Blambda}
    =\Bc_1+\int_{\CS^*_2}\BQ^\TRANSP\Blambda\DIFFM{\tau}{\DIFF\Blambda}
    \ge \BL_{\ineq}^\TRANSP\Bxi_{\ineq}+\BL_{\eq}^\TRANSP\Bxi_{\eq}$.
    In particular, this shows that $(\delta_{\hat{\Blambda}},\Bxi_{\ineq},\Bxi_{\eq})$ is also feasible for \eqref{eqn:dro-dual}.
    Conversely, let us assume that $\hat{\Blambda}\in\CS^*_2$, $\Bxi_{\ineq}\in\R_{-}^{n_{\ineq}}$, $\Bxi_{\eq}\in\R^{n_{\eq}}$ satisfy 
    $\Bc_1+\BQ^\TRANSP\hat{\Blambda}\ge \BL_{\ineq}^\TRANSP\Bxi_{\ineq}+\BL_{\eq}^\TRANSP\Bxi_{\eq}$.
    Then, $(\delta_{\hat{\Blambda}},\Bxi_{\ineq},\Bxi_{\eq})$ is feasible for \eqref{eqn:dro-dual}.
    This proves statement~\ref{thms:primaldual-boundedness}.

    In the following, let us assume in addition that $C_{\DRO}>-\infty$ and prove the remaining statements.
    Let us defer the proof of statement~\ref{thms:primaldual-primal-optimizer-existence} until we have proved 
    statement~\ref{thms:primaldual-dual-optimizer-existence}.
    To prove statement~\ref{thms:primaldual-primal-optimizer-properties}, 
    let $(\bar{\BIx},\bar{f}_1,\ldots,\bar{f}_d)$ be an arbitrary optimizer of \eqref{eqn:dro-primal}.
    It follows from Lemma~\ref{lem:primal-objective} and statement~\ref{thms:primaldual-strongduality} that 
    $\sum_{i=1}^d \bar{f}^*_i(\omega_i)\ge C_2(\bar{\BIx},\Bomega)$ for all $\Bomega=(\omega_1,\ldots,\omega_d)^\TRANSP\in\BOmega$, and
    \begin{align*}
        \langle\Bc_1,\bar{\BIx}\rangle+ \sup_{\mu\in\Gamma(\mu_1,\ldots,\mu_d)}\bigg\{\int_{\BOmega}C_2(\bar{\BIx},\Bomega)\DIFFM{\mu}{\DIFF\Bomega}\bigg\}
        &\le \langle\Bc_1,\bar{\BIx}\rangle + \sum_{i=1}^d \int_{\Omega_i}\bar{f}^*_i\DIFFM{\mu_i}{\DIFF\omega}\\
        &= \eqref{eqn:dro-primal} = C_{\DRO},
    \end{align*}
    which shows that $\bar{\BIx}$ is an optimizer of \eqref{eqn:dro}.
    This also shows that $(\bar{f}^*_1,\ldots,\bar{f}^*_d)$ minimizes (\ref{eqn:primaldual-proof-kantorovichduality}) with $\BIx\leftarrow\bar{\BIx}$, and thus
    $\mu\in\argmax_{\mu'\in\Gamma(\mu_1,\ldots,\mu_d)}\big\{\!\int_{\BOmega}C_2(\bar{\BIx},\Bomega)\DIFFM{\mu'}{\DIFF\Bomega}\big\}$
    if and only if 
    $\sum_{i=1}^{d}\bar{f}_i^*(\omega_i)=C_2(\bar{\BIx},\Bomega)$ for 
    $\mu$-almost every $\Bomega=(\omega_1,\ldots,\omega_d)^\TRANSP\in\BOmega$.
    The proof of statement~\ref{thms:primaldual-primal-optimizer-properties} is complete. 

    To prove statement~\ref{thms:primaldual-primal-approxoptimizer}, we fix an arbitrary $\varepsilon>0$ as well as an arbitrary $(\BIx_{\varepsilon},f_{1,\varepsilon},\ldots,f_{d,\varepsilon})$ that is $\varepsilon$-optimal for \eqref{eqn:dro-primal}. 
    Thus, it follows from Lemma~\ref{lem:primal-objective} and statement~\ref{thms:primaldual-strongduality} that
    \begin{align*}
        \langle\Bc_1,\BIx_{\varepsilon}\rangle+ \sup_{\mu\in\Gamma(\mu_1,\ldots,\mu_d)}\bigg\{\int_{\BOmega}C_2(\BIx_{\varepsilon},\Bomega)\DIFFM{\mu}{\DIFF\Bomega}\bigg\}
        &\le \langle\Bc_1,\BIx_{\varepsilon}\rangle + \sum_{i=1}^d \int_{\Omega_i}\sup_{z\in\R}\big\{z\omega-f_{i,\varepsilon}(z)\big\}\DIFFM{\mu_i}{\DIFF\omega}\\
        &\le \eqref{eqn:dro-primal} + \varepsilon = C_{\DRO} + \varepsilon,
    \end{align*}
    which shows that $\BIx_{\varepsilon}$ is an $\varepsilon$-optimizer of \eqref{eqn:dro} and proves statement~\ref{thms:primaldual-primal-approxoptimizer}.

    Let us prove statement~\ref{thms:primaldual-dual-optimizer-existence}.
    Recall that we have shown in the proof of statement~\ref{thms:primaldual-strongduality} that the optimal values of \eqref{eqn:primaldual-proof-dual-finitesupport} and \eqref{eqn:dro-dual} are equal.
    Hence, it suffices to show that (\ref{eqn:primaldual-proof-dual-finitesupport}) admits an optimizer.
    To that end, let us denote $K:=|\CV_{\extremept}|$, enumerate $\CV_{\extremept}=\{\lambda_1,\ldots,\lambda_K\}$, and define $\Delta_K\subset\R^K_{+}$, $\Phi:\R^K\to\R$ as in the proof of statement~\ref{thms:primaldual-strongduality}.
    Moreover, let us define $h:\R^K\times\R^{n_{\ineq}}\times\R^{n_{\eq}}\to\R\union\{-\infty\}$ as follows:
    \begin{align}
        \begin{split}
            h(\nu_1,\ldots,\nu_K,\Bxi_{\ineq},\Bxi_{\eq})&:=\Bigg(\sum_{j=1}^{K}\langle\Bb,\Blambda_j\rangle\nu_j\Bigg) 
            - \Phi(\nu_1,\ldots,\nu_K)
            + \langle\Bq_{\ineq},\Bxi_{\ineq}\rangle + \langle\Bq_{\eq},\Bxi_{\eq}\rangle \\
            & \hspace{70pt} \forall (\nu_1,\ldots,\nu_K)^\TRANSP\in\R^K,\; \forall \Bxi_{\ineq}\in\R^{n_{\ineq}}, \; \forall \Bxi_{\eq}\in\R^{n_{\eq}}.
        \end{split}
        \label{eqn:primaldual-proof-dual-objective}
    \end{align}
    One may check that (\ref{eqn:primaldual-proof-dual-finitesupport}) is equivalent to the maximization of the upper semi-continuous and concave function $h$ over the polyhedral convex set: 
    $D:=\big\{(\nu_1,\ldots,\nu_K,\Bxi_{\ineq}^\TRANSP,\Bxi_{\eq}^\TRANSP)^\TRANSP\in\R^K\times\R_{-}^{n_{\ineq}}\times\R^{n_{\eq}}:
    \Bc_1+\sum_{j=1}^K\nu_j\BQ^\TRANSP\Blambda_j\ge \BL_{\ineq}^\TRANSP\Bxi_{\ineq}+\BL_{\eq}^\TRANSP\Bxi_{\eq}\big\}$.
    Since $\eqref{eqn:primaldual-proof-dual-finitesupport}<\infty$ and 
    $h$ is affine in every direction of recession of $h$ (see \citep[p.69]{rockafellar1970convex} for the definition of the directions of recession of convex functions), 
    it follows from 
    \citep[Corollary~27.3.1]{rockafellar1970convex} that 
    $h$ admits a maximizer over $D$.
    Therefore, (\ref{eqn:primaldual-proof-dual-finitesupport}) admits an optimizer $(\bar{\tau},\bar{\Bxi}_{\ineq},\bar{\Bxi}_{\eq})$ which is also an optimizer of \eqref{eqn:dro-dual}.

    Next, let us prove statement~\ref{thms:primaldual-primal-optimizer-existence}.
    To that end, let us fix an optimizer $(\bar{\tau},\bar{\Bxi}_{\ineq},\bar{\Bxi}_{\eq})$ of (\ref{eqn:primaldual-proof-dual-finitesupport})
    and denote $\bar{\nu}_j:=\bar{\tau}\big(\{\Blambda_j\}\big)$ for $j=1,\ldots,K$.
    Hence, $(\bar{\nu}_1,\ldots,\bar{\nu}_{K},\bar{\Bxi}_{\ineq},\bar{\Bxi}_{\eq})$ maximizes the function $h$ defined in (\ref{eqn:primaldual-proof-dual-objective}) subject to the affine inequality constraint $\Bc_1+\sum_{j=1}^K\bar{\nu}_j\BQ^\TRANSP\Blambda_j\ge \BL_{\ineq}^\TRANSP\bar{\Bxi}_{\ineq}+\BL_{\eq}^\TRANSP\bar{\Bxi}_{\eq}$.
    Subsequently, the Karush--Kuhn--Tucker optimality conditions imply that there exist $\bar{\BIx}\in\R^{k_1}_{+}$, $\bar{\Btheta}\in\R_{+}^{n_{\ineq}}$, and $(\bar{g}_1,\ldots,\bar{g}_K)^\TRANSP\in\R^K$ such that:
    \begin{subequations}\label{eqn:primaldual-proof-KKT}%
        \begin{align}%
            (\bar{g}_1,\ldots,\bar{g}_K)^\TRANSP &\in \partial \Phi(\bar{\nu}_1,\ldots,\bar{\nu}_{K})\label{eqn:primaldual-proof-KKT-subgradient1}\allowdisplaybreaks\\
            \langle\BQ\bar{\BIx}+\Bb,\Blambda_j\rangle &=\bar{g}_j \hspace{50pt}\qquad\forall 1\le j\le K,\label{eqn:primaldual-proof-KKT-subgradient2}\allowdisplaybreaks\\
            \BL_{\ineq}\bar{\BIx}+\bar{\Btheta}&=\Bq_{\ineq},\label{eqn:primaldual-proof-KKT-in} \allowdisplaybreaks\\
            \BL_{\eq}\bar{\BIx}&=\Bq_{\eq},\label{eqn:primaldual-proof-KKT-eq} \allowdisplaybreaks\\
            \!\!\big\langle\Bc_1+\big({\textstyle\sum_{j=1}^{K}}\bar{\nu}_j\BQ^\TRANSP\Blambda_j\big) - \BL_{\ineq}^\TRANSP\bar{\Bxi}_{\ineq} - \BL_{\eq}^\TRANSP\bar{\Bxi}_{\eq}, \bar{\BIx}\big\rangle &= 0,\label{eqn:primaldual-proof-KKT-complementarity-x}\allowdisplaybreaks\\
            \langle\bar{\Btheta},\bar{\Bxi}_{\ineq}\rangle&=0. \label{eqn:primaldual-proof-KKT-complementarity-in}%
        \end{align}%
    \end{subequations} 
    Applying Lemma~\ref{lem:OT-function-convex}\ref{lems:OT-function-convex-subgradient} to (\ref{eqn:primaldual-proof-KKT-subgradient1}) and (\ref{eqn:primaldual-proof-KKT-subgradient2}) then shows the existence of $\bar{f}_1,\ldots,\bar{f}_d\in\CC_{\mathrm{cvx}}(\R)$ that satisfy
    ${-\int_{\Omega_i}}\sup_{z\in\R}\big\{z\omega-\bar{f}_i(z)\big\}\DIFFM{\mu_i}{\DIFF\omega}-\int_{\R}\bar{f}_i\DIFFX{\langle\Bp_i,\cdot\,\rangle\sharp\bar{\tau}}=\OT_{\mu_i}\big(\langle\Bp_i,\cdot\,\rangle\sharp\bar{\tau}\big)$, 
    $\partial \bar{f}_i(z)\subseteq\Omega_i$ $\forall z\in\R$
    for $i=1,\ldots,d$ 
    as well as 
    $\langle\BQ\bar{\BIx}+\Bb,\Blambda_j\rangle+\sum_{i=1}^{d}\bar{f}_i\big(\langle\Bp_i,\Blambda_j\rangle\big)=0$ for every $j$ with $\bar{\nu}_j>0$,
    $\langle\BQ\bar{\BIx}+\nobreak\Bb,\Blambda_j\rangle+\sum_{i=1}^{d}\bar{f}_i\big(\langle\Bp_i,\Blambda_j\rangle\big)\le0$ for every $j$ with $\bar{\nu}_j=0$.
    Lemma~\ref{lem:global-simplification} guarantees that $\langle\BQ\bar{\BIx}+\Bb,\Blambda\rangle+\sum_{i=1}^{d}\bar{f}_i\big(\langle\Bp_i,\Blambda\rangle\big)\le 0$ for all $\Blambda\in\CS^*_2$.
    It subsequently follows from (\ref{eqn:primaldual-proof-KKT-in}) and (\ref{eqn:primaldual-proof-KKT-eq}) that 
    $(\bar{\BIx},\bar{f}_1,\ldots,\bar{f}_d)$ is feasible for \eqref{eqn:dro-primal}.
    Moreover, 
    since $\sum_{j=1}^{K}\bar{\nu}_j\BQ^\TRANSP\Blambda_j=\int_{\CS^*_2}\BQ^\TRANSP\Blambda\DIFFM{\bar{\tau}}{\DIFF\Blambda}$,
    (\ref{eqn:primaldual-proof-KKT-complementarity-x})
    and (\ref{eqn:primaldual-proof-KKT-complementarity-in}) guarantee that
    $\langle\BL_{\ineq}\bar{\BIx}-\Bq_{\ineq},\bar{\Bxi}_{\ineq}\rangle=0$,
    $\big\langle\Bc_1+\int_{\CS^*_2}\BQ^\TRANSP\Blambda\DIFFM{\bar{\tau}}{\DIFF\Blambda} - \BL_{\ineq}^\TRANSP\bar{\Bxi}_{\ineq} - \BL_{\eq}^\TRANSP\bar{\Bxi}_{\eq}, \bar{\BIx}\big\rangle = 0$.
    Combining these properties with the established properties
    $\int_{\CS^*_2}\langle\BQ\bar{\BIx}+\nobreak\Bb,\Blambda\rangle+{\textstyle\sum_{i=1}^{d}}\bar{f}_i\big(\langle\Bp_i,\Blambda\rangle\big)\DIFFM{\bar{\tau}}{\DIFF\Blambda}= 0$ 
    and 
    $\OT_{\mu_i}\big(\langle\Bp_i,\cdot\,\rangle\sharp\bar{\tau}\big) +\int_{\Omega_i}\sup_{z\in\R}\big\{z\omega-\bar{f}_i(z)\big\}\DIFFM{\mu_i}{\DIFF\omega}+\int_{\R}\bar{f}_i\DIFFX{\langle\Bp_i,\cdot\,\rangle\sharp\bar{\tau}}=0$ for $i=1,\ldots,d$,
    we conclude by Lemma~\ref{lem:weak-duality} and statement~\ref{thms:primaldual-strongduality} that $(\bar{\BIx},\bar{f}_1,\ldots,\bar{f}_d)$ is an optimizer of \eqref{eqn:dro-primal}.
    This completes the proof of statement~\ref{thms:primaldual-primal-optimizer-existence}.

    To prove statement~\ref{thms:primaldual-dual-optimizer-properties}, 
    let us fix an arbitrary optimizer $(\bar{\tau},\bar{\Bxi}_{\ineq},\bar{\Bxi}_{\eq})$ of \eqref{eqn:dro-dual}, 
    let $\bar{\pi}\in\CP(\BOmega\times\CS^*_2)$ satisfy the properties in Lemma~\ref{lem:dual-gluing}\ref{lems:dual-gluing-construction} with respect to $\tau\leftarrow\bar{\tau}$,
    and let $\bar{\mu}$ denote the marginal of $\bar{\pi}$ on~$\BOmega$.
    Recall the equalities (\ref{eqn:primaldual-proof-minimax}) established in the proof of statement~\ref{thms:primaldual-strongduality}.
    For any optimizer $\bar{\BIx}\in\CS_1$ of \eqref{eqn:dro},
    it holds that 
    $(\bar{\BIx},\bar{\tau})$ is a saddle point of the minimax problem (\ref{eqn:primaldual-proof-minimax})
    where $\bar{\BIx}$ is a minimizer of the inf-sup problem (\ref{eqn:primaldual-proof-minimax-infsup})
    and $\bar{\tau}$ is a maximizer of the sup-inf problem (\ref{eqn:primaldual-proof-minimax-supinf}).
    Consequently, $\bar{\tau}$ also maximizes (\ref{eqn:dual-gluing-equalities-meas}) in Lemma~\ref{lem:dual-gluing}\ref{lems:dual-gluing-equalities} with respect to $\BIx\leftarrow\bar{\BIx}$,
    and thus statement~\ref{thms:primaldual-dual-optimizer-properties} follows from Lemma~\ref{lem:dual-gluing}\ref{lems:dual-gluing-optimality}.

    Lastly, statement~\ref{thms:primaldual-saddle} is a direct consequence of Lemma~\ref{lem:weak-duality} and statement~\ref{thms:primaldual-strongduality}.
    The proof is now complete.
\end{proof}


\begin{proof}[Proof of Corollary~\ref{cor:dual-optimizers}]
    For $i=1,\ldots,d$,
    it follows from Theorem~\ref{thm:primaldual}\ref{thms:primaldual-saddle} that 
    ${-\int_{\Omega_i}}\bar{f}_i^*\DIFFX{\mu_i}-\int_{\R}\bar{f}_i\DIFFX{\langle\Bp_i,\cdot\,\rangle\sharp\bar{\tau}}= \OT_{\mu_i}\big(\langle\Bp_i,\cdot\,\rangle\sharp\bar{\tau}\big)$.
    Lemma~\ref{lem:OT-1d}\ref{lems:OT-1d-optimality} with $\mu\leftarrow\mu_i$, $\nu\leftarrow\langle\Bp_i,\cdot\,\rangle\sharp\bar{\tau}$ then implies that there exists $\bar{\gamma}_i\in\Gamma\big(\mu_i,\langle\Bp_i,\cdot\,\rangle\sharp\bar{\tau}\big)$ such that
    $\omega\in\partial \bar{f}_i(z)$ and thus $z\in\partial \bar{f}_i^*(\omega)$ for $\bar{\gamma}_i$-almost every $(\omega,z)\in\Omega_i\times\R$.
    Consequently, 
    $\bar{\gamma}_i$ is concentrated on the graph $\big\{(\omega,z):\omega\in\Omega_i,\; z\in \partial \bar{f}_i^*(\omega)\big\}$ of $\partial \bar{f}_i^*$.
    Since $\mu_i$ is absolutely continuous with respect to the Lebesgue measure, it holds that
    $\partial \bar{f}_i^*(\omega)$ is a singleton for $\mu_i$-almost every $\omega\in\Omega_i$; see, e.g., \citep[Theorem~25.5]{rockafellar1970convex}.
    Thus,
    $\bar{\gamma}_i$ is concentrated on the graph $\big\{(\omega,\nabla\bar{f}_i^*(\omega)):\omega\in\Omega_i\big\}$ of $\nabla \bar{f}_i^*:\Omega_i\to\R$,
    and hence $\langle\Bp_i,\cdot\,\rangle\sharp\bar{\tau}=\nabla \bar{f}_i^*\sharp\mu_i$.
    In particular, this implies that 
    $\langle\Bp_i,\cdot\,\rangle\sharp\bar{\tau} = \langle\Bp_i,\cdot\,\rangle\sharp\bar{\tau}'$ holds for $i=1,\ldots,d$ and for any two optimizers $(\bar{\tau},\bar{\Bxi}_{\ineq},\bar{\Bxi}_{\eq}),(\bar{\tau}',\bar{\Bxi}_{\ineq}',\bar{\Bxi}_{\eq}')$ of \eqref{eqn:dro-dual}.
    The proof is now complete.
\end{proof}


\begin{proof}[Proof of Theorem~\ref{thm:dual-suboptimality}]
    To begin, let us fix $\BIx\in\CS_1$, $f_1,\ldots,f_d\in\CC_{\mathrm{cvx}}(\R)$ that satisfy (\ref{eqn:dual-suboptimality-primal-feasibility-potential-eq})--(\ref{eqn:dual-suboptimality-potential-subgradient}), assuming that they exist.
    It follows from (\ref{eqn:dual-suboptimality-primal-feasibility-potential-eq}) and $\support(\tau)\subset\CS^*_2$ that
    $\mathrm{vio}(\BIx,f_1,\ldots,f_d)\ge 0$.
    The rest of statement~\ref{thms:dual-suboptimality-vio-positivity} follows from
    (\ref{eqn:dual-suboptimality-potential-subgradient}) and Lemma~\ref{lem:global-simplification}.
    
    To prove statement~\ref{thms:dual-suboptimality-suboptimality}, 
    let us fix $\BIx\in\CS_1$, $f_1,\ldots,f_d\in\CC_{\mathrm{cvx}}(\R)$ that satisfy (\ref{eqn:dual-suboptimality-primal-feasibility-potential-eq})--(\ref{eqn:dual-suboptimality-potential-subgradient}), 
    define
    $\varepsilon_{\mathrm{vio}}:=\mathrm{vio}(\BIx,f_1,\ldots,f_d)$,
    $\tilde{f}_1:=f_1-\varepsilon_{\mathrm{vio}}$,
    $\tilde{f}_i:=f_i$ for $i=2,\ldots,d$,
    and examine the properties of $\BIx,\tilde{f}_1,\ldots,\tilde{f}_d$.
    Firstly, 
    we have $\sup_{\Blambda\in\CS^*_2}\big\{\langle\BQ\BIx+\Bb,\Blambda\rangle+\sum_{i=1}^d \tilde{f}_i\big(\langle\Bp_i,\Blambda\rangle\big)\big\}= \sup_{\Blambda\in\CS^*_2}\big\{\langle\BQ\BIx+\nobreak\Bb,\Blambda\rangle+\sum_{i=1}^d f_i\big(\langle\Bp_i,\Blambda\rangle\big)\big\} - \varepsilon_{\mathrm{vio}}= 0$,
    which shows that $(\BIx,\tilde{f}_1,\ldots,\tilde{f}_d)$ is feasible for \eqref{eqn:dro-primal}. 
    Secondly,
    (\ref{eqn:dual-suboptimality-primal-feasibility-potential-eq}) implies that
    $\int_{\CS^*_2}\langle\BQ\BIx+\Bb,\Blambda\rangle+\sum_{i=1}^{d}\tilde{f}_i\big(\langle\Bp_i,\Blambda\rangle\big)\DIFFM{\tau}{\DIFF\Blambda}={-\varepsilon_{\mathrm{vio}}}$.
    Thirdly, 
    one observes from (\ref{eqn:dual-suboptimality-potential-optimizer}) that 
    $\OT_{\mu_i}\big(\langle\Bp_i,\cdot\,\rangle\sharp\tau\big)+\int_{\Omega_i}\sup_{z\in\R}\big\{z\omega-\tilde{f}_i(z)\big\}\DIFFM{\mu_i}{\DIFF\omega}+\int_{\R}\tilde{f}_i\DIFFX{\langle\Bp_i,\cdot\,\rangle\sharp\tau}\le \varepsilon_{\OT}$ for $i=1,\ldots,d$.
    It subsequently follows from (\ref{eqn:dual-suboptimality-primal-complementarity}),
    (\ref{eqn:dual-suboptimality-dual-complementarity}), 
    and Lemma~\ref{lem:weak-duality} that
    the objective of $(\BIx,\tilde{f}_1,\ldots,\tilde{f}_d)$ with respect to \eqref{eqn:dro-primal} is at most $\varepsilon_{\mathrm{vio}}+d\varepsilon_{\OT}$ larger than 
    the objective of $(\tau,\Bxi_{\ineq},\Bxi_{\eq})$ with respect to \eqref{eqn:dro-dual}.
    It then follows from Theorem~\ref{thm:primaldual}\ref{thms:primaldual-strongduality} that 
    $(\tau,\Bxi_{\ineq},\Bxi_{\eq})$ is an $\varepsilon_{\mathrm{sub}}$-optimizer of \eqref{eqn:dro-dual}.
    Lastly, in the case where $\varepsilon_{\mathrm{vio}}>0$, 
    (\ref{eqn:dual-suboptimality-primal-feasibility-potential-eq})
    and (\ref{eqn:dual-suboptimality-primal-feasibility-potential-le}) guarantee that
    $\big\{\Blambda\in\CS^*_2:\langle\BQ\BIx+\Bb, \Blambda\rangle+\sum_{i=1}^d f_i\big(\langle\Bp_i,\Blambda\rangle\big)>\nobreak0\big\}$ 
    is disjoint from~$\CV$.
    The proof is now complete.
\end{proof}


\begin{proof}[Proof of Corollary~\ref{cor:dual-suboptimality-primal-coupling}]
    Let us define $\tilde{f}_1:=f_1-\varepsilon_{\mathrm{vio}}$ and define $\tilde{f}_i:=f_i$ for $i=2,\ldots,d$.
    Recall that the proof of Theorem~\ref{thm:dual-suboptimality}\ref{thms:dual-suboptimality-suboptimality} has already shown that 
    $(\BIx,\tilde{f}_1,\ldots,\tilde{f}_d)$ is an 
    $\varepsilon_{\mathrm{sub}}$-optimizer of \eqref{eqn:dro-primal}.
    It subsequently follows from Theorem~\ref{thm:primaldual}\ref{thms:primaldual-primal-approxoptimizer} that $\BIx$ is an $\varepsilon_{\mathrm{sub}}$-optimizer of \eqref{eqn:dro}.
    This proves statement~\ref{cors:dual-suboptimality-primal-coupling-primal}.

    Now, let us prove statement~\ref{cors:dual-suboptimality-primal-coupling-meas}.
    For $i=1,\ldots,d$, let us denote 
    $m_i:=\big|\big\{\langle\Bp_i,\Blambda\rangle:\Blambda\in\CV\big\}\big|$,
    enumerate 
    $\big\{\langle\Bp_i,\Blambda\rangle:\Blambda\in\CV\big\}=\{\kappa_{i,1},\ldots,\kappa_{i,m_i}\}$
    where $\kappa_{i,1}<\kappa_{i,2}<\cdots<\kappa_{i,m_i}$,
    and
    define $G^{-}_{\langle\Bp_i,\cdot\,\rangle\sharp\tau}:[0,1]\to \R$ as follows:
    \begin{align*}
        G^{-}_{\langle\Bp_i,\cdot\,\rangle\sharp\tau}(u):=\begin{cases}
            \kappa_{i,1} & \hspace{174.9pt} \forall u\in\big[0, \langle\Bp_i,\cdot\,\rangle\sharp\tau\big((-\infty,\kappa_{i,1}]\big)\big], \\
            \kappa_{i,l} & \hspace{20pt} \forall u\in \big(\langle\Bp_i,\cdot\,\rangle\sharp\tau\big((-\infty,\kappa_{i,l})\big),\langle\Bp_i,\cdot\,\rangle\sharp\tau\big((-\infty,\kappa_{i,l}]\big)\big],\; \forall 2\le l\le m_i.
        \end{cases}
    \end{align*}
    Moreover, for $i=1,\ldots,d$,
    let us denote
    $\eta_j:=\sum_{l=1}^{j}\tau\big(\{\Blambda_{\sigma_i(l)}\}\big)$ 
    for $j=0,1,\ldots,K$,
    and define the random variable 
    $W_i:\mathscr{X}\to[0,1]$ by
    $W_i:=U_i\tau\big(\{\Blambda_{J}\}\big) + \sum_{j=1}^{\sigma_i^{-1}(J)-1}\tau\big(\{\Blambda_{\sigma_i(j)}\}\big)=U_i\tau\big(\{\Blambda_{J}\}\big) + \eta_{\sigma_i^{-1}(J)-1}$.
    Furthermore, let $\mathscr{L}_{[0,1]}$ denote the Lebesgue measure on $[0,1]$,
    let $\CJ:=\big\{j\in\{1,\ldots,K\}:\tau\big(\{\Blambda_j\}\big)>\nobreak0\big\}$,
    and define the random variable $\BILambda:\mathscr{X}\to\CS^*_2$ as $\BILambda:=\Blambda_{J}$.
    Observe that the law of $\BILambda$ is~$\tau$.
    In the following, let us fix an arbitrary $i\in\{1,\ldots,d\}$ and show that $W_i$ is uniformly distributed on $[0,1]$, i.e., the law of $W_i$ is $\mathscr{L}_{[0,1]}$,
    and that $G^{-}_{\langle\Bp_i,\cdot\,\rangle\sharp\tau}(W_i)=\langle\Bp_i,\BILambda\rangle$ $\PROB$-almost surely. 

    Let us fix an arbitrary $\hat{j}\in\{1,\ldots,K\}$ such that $\sigma_i(\hat{j})\in\CJ$, and observe that 
    $\tau\big(\{\Blambda_{\sigma_i(\hat{j})}\}\big)+\eta_{\hat{j}-1}
    =\nobreak\eta_{\hat{j}}$.
    On the event $\big\{J=\sigma_i(\hat{j})\big\}$, 
    it holds that $W_i=U_i\tau\big(\{\Blambda_{\sigma_i(\hat{j})}\}\big)+\eta_{\hat{j}-1}$.
    This shows that $W_i$ is uniformly distributed on $[\eta_{\hat{j}-1},\eta_{\hat{j}}]$ conditional on $J=\sigma_i(\hat{j})$,
    where $\eta_{\hat{j}}-\eta_{\hat{j}-1}=\tau\big(\{\Blambda_{\sigma_i(\hat{j})}\}\big)=\PROB\big[J=\sigma_i(\hat{j})\big]$.
    Since $0=\eta_{0}\le \eta_{1} \le \cdots \le \eta_{K-1}\le\eta_{K}=1$,
    marginalizing the joint distribution of $(W_i,J)$ over $J$ yields that $W_i$ is uniformly distributed on $[0,1]$.
    Moreover, there exists $l\in\{1,\ldots,m_i\}$ with $\kappa_{i,l}=\langle\Bp_i,\Blambda_{\sigma_i(\hat{j})}\rangle$.
    Let $\underline{j}:=\min\big\{j\in\{1,\ldots,K\}:\langle\Bp_i,\Blambda_{\sigma_{i}(j)}\rangle=\kappa_{i,l}\big\}\le \hat{j}$ and 
    $\overline{j}:=\max\big\{j\in\{1,\ldots,K\}:\langle\Bp_i,\Blambda_{\sigma_{i}(j)}\rangle=\nobreak\kappa_{i,l}\big\}\ge \hat{j}$.
    Subsequently,
    we get
    \begin{align}
        \langle\Bp_i,\cdot\,\rangle\sharp\tau\hspace{-1pt}\big((-\infty,\kappa_{i,l}]\big)\!&= \!\sum_{j=1}^{K}\INDI_{\{\langle\Bp_i,\Blambda_{\sigma_i(j)}\rangle\le \kappa_{i,l}\}}\tau\hspace{-1pt}\big(\hspace{-1pt}\{\Blambda_{\sigma_i(j)}\}\hspace{-1pt}\big)\!=\! \sum_{j=1}^{\overline{j}}\tau\hspace{-1pt}\big(\hspace{-1pt}\{\Blambda_{\sigma_i(j)}\}\hspace{-1pt}\big)\!=\eta_{\overline{j}}\ge \eta_{\hat{j}},
        \label{eqn:dual-suboptimality-primal-coupling-proof-ineq1}
    \end{align}
    as well as
    \begin{align}
        \!\langle\Bp_i,\cdot\,\rangle\sharp\tau\hspace{-1pt}\big((-\infty,\kappa_{i,l})\big)\!&= \!\sum_{j=1}^{K}\INDI_{\{\langle\Bp_i,\Blambda_{\sigma_i(j)}\rangle< \kappa_{i,l}\}}\tau\hspace{-1pt}\big(\hspace{-1pt}\{\Blambda_{\sigma_i(j)}\}\hspace{-1pt}\big)\!= \!\sum_{j=1}^{\underline{j}-1}\tau\big(\hspace{-1pt}\{\Blambda_{\sigma_i(j)}\}\hspace{-1pt}\big)\!=\eta_{\underline{j}-1}\le \eta_{\hat{j}-1}.\!
        \label{eqn:dual-suboptimality-primal-coupling-proof-ineq2}
    \end{align}
    Since $\PROB\big[\eta_{\hat{j}-1}<W_i\le \eta_{\hat{j}}|J=\sigma_i(\hat{j})\big]=\nobreak1$,
    applying (\ref{eqn:dual-suboptimality-primal-coupling-proof-ineq1}), (\ref{eqn:dual-suboptimality-primal-coupling-proof-ineq2}),
    and the definition of $G^{-}_{\langle\Bp_i,\cdot\,\rangle\sharp\tau}$ 
    shows that 
    $\PROB\big[G^{-}_{\langle\Bp_i,\cdot\,\rangle\sharp\tau}(W_i)=\langle\Bp_i,\Blambda_{\sigma_i(\hat{j})}\rangle|J=\sigma_i(\hat{j})\big]=\nobreak1$.
    Therefore, we conclude that 
    $G^{-}_{\langle\Bp_i,\cdot\,\rangle\sharp\tau}(W_i)=\langle\Bp_i,\Blambda_{J}\rangle=\langle\Bp_i,\BILambda\rangle$ $\PROB$-almost surely. 

    In summary,
    we have shown that for $i=1,\ldots,d$, $W_i$ is uniformly distributed on $[0,1]$ and $G^{-}_{\langle\Bp_i,\cdot\,\rangle\sharp\tau}(W_i)=\langle\Bp_i,\BILambda\rangle$ $\PROB$-almost surely. 
    For $i=1,\ldots,d$, let $\gamma_i\in\CP(\Omega_i\times\R)$ denote the law of the random variable $\big(Y_i,\langle\Bp_i,\BILambda\rangle\big)^\TRANSP:\mathscr{X}\to\Omega_i\times\R$.
    Since the law of $W_i$ is $\mathscr{L}_{[0,1]}$ and $Y_i=G^{-}_{\mu_i}(W_i)$, $\langle\Bp_i,\BILambda\rangle=G^{-}_{\langle\Bp_i,\cdot\,\rangle\sharp\tau}(W_i)$ $\PROB$-almost surely, 
    we get $\gamma_i=\big[G^{-}_{\mu_i},G^{-}_{\langle\Bp_i,\cdot\,\rangle\sharp\tau}\big]\sharp\mathscr{L}_{[0,1]}$.
    Lemma~\ref{lem:OT-1d}\ref{lems:OT-1d-coupling} with $\mu\leftarrow\mu_i$, $\nu\leftarrow \langle\Bp_i,\cdot\,\rangle\sharp\tau$ then implies that $\gamma_i\in\Gamma(\mu_i,\langle\Bp_i,\cdot\,\rangle\sharp\tau)$ and
    \begin{align}
        \OT_{\mu_i}\big(\langle\Bp_i,\cdot\,\rangle\sharp\tau\big)=\int_{\Omega_i\times\R}{-z\omega}\DIFFM{\gamma_i}{\DIFF\omega,\DIFF z}={-\EXP}\big[\langle\Bp_i,\BILambda\rangle Y_i\big] \qquad \forall 1\le i\le d.
        \label{eqn:dual-suboptimality-primal-coupling-proof-optimal-coupling}
    \end{align}
    Let us denote $\BIY:=(Y_1,\ldots,Y_d)^\TRANSP:$ $\mathscr{X}\to \BOmega$.
    Since the law of $\BIY$ is $\mu_{\tau}$, we have $\mu_{\tau}\in\Gamma(\mu_1,\ldots,\mu_d)$,
    and
    (\ref{eqn:dual-suboptimality-primal-coupling-proof-optimal-coupling}) implies that
    \begin{align}
        \int_{\BOmega}C_2(\BIx,\Bomega)\DIFFM{\mu_{\tau}}{\DIFF\Bomega}
        &=\EXP\bigg[\max_{\Blambda\in\CS^*_2}\big\{\langle\BQ\BIx+\BP\BIY+\Bb,\Blambda\rangle\big\}\bigg]\ge\EXP\big[\langle\BQ\BIx+\Bb,\BILambda\rangle\big] + \sum_{i=1}^d \EXP\big[\langle\Bp_i,\BILambda\rangle Y_i\big]\nonumber\\
        &=\int_{\CS^*_2}\langle\BQ\BIx+\Bb,\Blambda\rangle\DIFFM{\tau}{\DIFF\Blambda}-\Bigg(\sum_{i=1}^d \OT_{\mu_i}\big(\langle\Bp_i,\cdot\,\rangle\sharp\tau\big)\Bigg).
        \label{eqn:dual-suboptimality-primal-coupling-proof-lb}
    \end{align}
    Moreover, 
    recall from the proof of Theorem~\ref{thm:dual-suboptimality}\ref{thms:dual-suboptimality-suboptimality} that
    $\int_{\CS^*_2}\langle\BQ\BIx+\Bb,\Blambda\rangle+\sum_{i=1}^{d}\tilde{f}_i\big(\langle\Bp_i,\Blambda\rangle\big)\DIFFM{\tau}{\DIFF\Blambda}={-\varepsilon_{\mathrm{vio}}}$.
    Substituting this and (\ref{eqn:dual-suboptimality-potential-optimizer}) into (\ref{eqn:dual-suboptimality-primal-coupling-proof-lb}) and then applying Lemma~\ref{lem:primal-objective} yields
    \begin{align*}
        \int_{\BOmega}C_2(\BIx,\Bomega)\DIFFM{\mu_{\tau}}{\DIFF\Bomega}
        &\ge \Bigg(\sum_{i=1}^{d}\int_{\Omega_i}\sup_{z\in\R}\big\{z\omega-\tilde{f}_i(z)\big\}\DIFFM{\mu_i}{\DIFF\omega}\Bigg)-\varepsilon_{\mathrm{sub}}\\
        &\ge \sup_{\mu\in\Gamma(\mu_1,\ldots,\mu_d)}\bigg\{\int_{\BOmega}C_2(\BIx,\Bomega) \DIFFM{\mu}{\DIFF\Bomega}\bigg\} - \varepsilon_{\mathrm{sub}}.
    \end{align*}
    The proof is now complete.
\end{proof}

\subsection{Proof of results in Section~\ref{sec:algorithm}}
\label{sapx:proof-algorithm}


\begin{proof}[Proof of Proposition~\ref{prop:cuttingplane-LSIP}]
    To prove statement~\ref{props:cuttingplane-LSIP-LPfeasibility},
    let us fix an arbitrary $\widehat{\CF}=(\widehat{\CF}_{1},\ldots,\widehat{\CF}_{d})$, where $\widehat{\CF}_{i}\subset\CC_{\mathrm{cvx}}(\R)$ is non-empty and finite for $i=1,\ldots,d$.
    On the one hand, \eqref{eqn:dro-dual-LP-relax} is feasible since it is a relaxation of \eqref{eqn:dro-dual-LSIP}, 
    which is assumed to be feasible under Setting~\ref{sett:cuttingplane-LSIP}.
    On the other hand,
    let us take an arbitrary $\hat{\BIx}\in\CS_1$,
    and, for $i=1,\ldots,d$, take an
    arbitrary $\hat{f}_i\in\widehat{\CF}_{i}$,
    define $\hat{y}_{i,\hat{f}_i}:=1$,
    and define $\hat{y}_{i,f_i}:=0$ for every $f_i\in\widehat{\CF}_{i}\setminus\{\hat{f}_i\}$.
    Subsequently, let us define 
    $\hat{h}_0:={-\max_{\Blambda\in\CV}}\big\{\langle\BQ\hat{\BIx}+\nobreak\Bb,\Blambda\rangle+\sum_{i=1}^{d}\hat{f}_i\big(\langle\Bp_i,\Blambda\rangle\big)\big\}$.
    One checks that 
    $\big(\hat{\BIx},\hat{h}_0, (\hat{y}_{i,f_i})_{f_i\in\widehat{\CF}_{i},\,i=1:d}\big)$
    is feasible for \eqref{eqn:dro-dual-LP-relax-dual}.
    This proves statement~\ref{props:cuttingplane-LSIP-LPfeasibility}.

    To prove statement~\ref{props:cuttingplane-LSIP-termination},
    let us suppose for the sake of contradiction that
    Algorithm~\ref{alg:cuttingplane-LSIP} loops forever
    and generates an infinite sequence $(\hat{\tau}^{(t)})_{t\in\N}\subset\CP(\CV)$ as well as infinite sequences
    $(\hat{\psi}_i^{(t)})_{t\in\N}\subset\nobreak\R$ for $i=1,\ldots,d$.
    Let us fix an arbitrary $i\in\{1,\ldots,d\}$ and show that 
    $\liminf_{t\to\infty}\hat{\varepsilon}_{\OT,i}^{(t)}\le0$ to get a contradiction.
    For each $t\in\N$,
    Line~\ref{alglin:cuttingplane-LSIP-LP}, Line~\ref{alglin:cuttingplane-LSIP-constructmeas},
    and the objective and constraints of \LPprimal{\widehat{\CF}^{(t)}} imply that 
    $\hat{\psi}_i^{(t)}=\max_{f_i\in\widehat{\CF}_{i}^{(t)}}\Big\{ {-\int_{\Omega_i}} \sup_{z\in\R}\big\{z\omega-f_i(z)\big\} \DIFFM{\mu_i}{\DIFF\omega}-\int_{\R}f_i\DIFFX{\langle\Bp_i,\cdot\,\rangle\sharp\hat{\tau}^{(t)}}\Big\}$.
    Observe that for each $t\in\N$, $\hat{\tau}^{(t)}$ can be regarded as an element of $\CP(\CV)$,
    and that $\CP(\CV)$ is compact with respect to weak convergence.
    Consequently, if follows from Lemma~\ref{lem:OT-1d}\ref{lems:OT-1d-duality} with $\mu\leftarrow\mu_i$ that
    \begin{align}
        \begin{split}
            \hat{\psi}_i^{(t)} &\le \OT_{\mu_i}\big(\langle\Bp_i,\cdot\,\rangle\sharp\hat{\tau}^{(t)}\big)
            \le \sup_{\tau\in\CP(\CV)}\big\{\OT_{\mu_i}\big(\langle\Bp_i,\cdot\,\rangle\sharp\tau\big)\big\}<\infty \qquad 
            \forall t\in\N.
        \end{split}
        \label{eqn:cuttingplane-LSIP-proof-psi-UB}
    \end{align}
    Moreover, since Line~\ref{alglin:cuttingplane-LSIP-initialfuncs} and Line~\ref{alglin:cuttingplane-LSIP-newfuncs} guarantee that $f^{\star}_{i,\hat{\tau}^{(0)}}$ belongs to $\widehat{\CF}_{i}^{(t)}$ for every $t\in\N$,
    we get
    \begin{align}
        \begin{split}
            \hat{\psi}_i^{(t)} &\ge {-\int_{\Omega_i}} \sup_{z\in\R}\big\{z\omega-f^{\star}_{i,\hat{\tau}^{(0)}}(z)\big\} \DIFFM{\mu_i}{\DIFF\omega}-\int_{\R}f^{\star}_{i,\hat{\tau}^{(0)}}\DIFFX{\langle\Bp_i,\cdot\,\rangle\sharp\hat{\tau}^{(t)}}\\
            &\ge {-\int_{\Omega_i}} \sup_{z\in\R}\big\{z\omega-f^{\star}_{i,\hat{\tau}^{(0)}}(z)\big\} \DIFFM{\mu_i}{\DIFF\omega}
            -\max_{\Blambda\in\CV}\big\{f^{\star}_{i,\hat{\tau}^{(0)}}\big(\langle\Bp_i,\Blambda\rangle\big)\big\}>{-\infty} \qquad \forall t\in\N.
        \end{split}
        \label{eqn:cuttingplane-LSIP-proof-psi-LB}
    \end{align}
    Due to the compactness of $\CP(\CV)$ with respect to weak convergence, 
    and due to (\ref{eqn:cuttingplane-LSIP-proof-psi-UB}) and (\ref{eqn:cuttingplane-LSIP-proof-psi-LB}),
    there exists a subsequence $(t_r)_{r\in\N}\subseteq\N$ such that
    as $r\to\infty$,
    $\big(\hat{\tau}^{(t_r)}\big)_{r\in\N}$ converges weakly to some $\hat{\tau}^{(\infty)}\in\CP(\CV)$
    and 
    $\big(\hat{\psi}_i^{(t_r)}\big)_{r\in\N}$ converges to some $\hat{\psi}_i^{(\infty)}\in\R$.
    Next, for any $r\in\N$, observe from Line~\ref{alglin:cuttingplane-LSIP-newfuncs} that $f^{\star}_{i,\hat{\tau}^{(t_r)}}\in\widehat{\CF}_{i}^{(t_s)}$ for all $s\in\N\intersection[r+1,\infty)$.
    Thus, 
    Line~\ref{alglin:cuttingplane-LSIP-constructmeas} and
    the constraints of \LPprimal{\widehat{\CF}^{(t_s)}} guarantee that
    $\hat{\psi}_i^{(t_s)} + \int_{\Omega_i} \sup_{z\in\R}\big\{z\omega-f^{\star}_{i,\hat{\tau}^{(t_r)}}(z)\big\} \DIFFM{\mu_i}{\DIFF\omega} + \int_{\R}f^{\star}_{i,\hat{\tau}^{(t_r)}}\DIFFX{\langle\Bp_i,\cdot\,\rangle\sharp\hat{\tau}^{(t_s)}} \ge 0$ for all $s\in\N\intersection[r+1,\infty)$,
    and taking $s\to\infty$ yields
    \begin{align}
        \hat{\psi}_i^{(\infty)} + \int_{\Omega_i} \sup_{z\in\R}\big\{z\omega-f^{\star}_{i,\hat{\tau}^{(t_r)}}(z)\big\} \DIFFM{\mu_i}{\DIFF\omega} + \int_{\R}f^{\star}_{i,\hat{\tau}^{(t_r)}}\DIFFX{\langle\Bp_i,\cdot\,\rangle\sharp\hat{\tau}^{(\infty)}}\ge 0 \qquad \forall r\in\N.
        \label{eqn:cuttingplane-LSIP-proof-limitviolation}
    \end{align}
    Combining (\ref{eqn:cuttingplane-LSIP-proof-limitviolation}) with Line~\ref{alglin:cuttingplane-LSIP-error} then leads to
    \begin{align}
        \hat{\varepsilon}_{\OT,i}^{(t_r)}\!&\le {-\int_{\Omega_i}} \sup_{z\in\R}\big\{z\omega-f^{\star}_{i,\hat{\tau}^{(t_r)}}(z)\big\} \DIFFM{\mu_i}{\DIFF\omega}
        - \int_{\R}f^{\star}_{i,\hat{\tau}^{(t_r)}}\DIFFX{\langle\Bp_i,\cdot\,\rangle\sharp\hat{\tau}^{(t_r)}} - \hat{\psi}_i^{(t_r)} \nonumber \\
        &\qquad + \bigg(\!\hat{\psi}_i^{(\infty)}  + \int_{\Omega_i} \sup_{z\in\R}\big\{z\omega-f^{\star}_{i,\hat{\tau}^{(t_r)}}(z)\big\} \DIFFM{\mu_i}{\DIFF\omega} + \int_{\R}f^{\star}_{i,\hat{\tau}^{(t_r)}}\DIFFX{\langle\Bp_i,\cdot\,\rangle\sharp\hat{\tau}^{(\infty)}}\!\bigg) 
        \label{eqn:cuttingplane-LSIP-proof-errorlimit}\\
        &\le \!\bigg|\!\int_{\R}f^{\star}_{i,\hat{\tau}^{(t_r)}}\DIFFX{\langle\Bp_i,\cdot\,\rangle\sharp\hat{\tau}^{(\infty)}} - \int_{\R}f^{\star}_{i,\hat{\tau}^{(t_r)}}\DIFFX{\langle\Bp_i,\cdot\,\rangle\sharp\hat{\tau}^{(t_r)}}\bigg| + \big|\hat{\psi}_i^{(\infty)}-\hat{\psi}_i^{(t_r)}\big| \nonumber\\
        &\le \max_{\Blambda\in\CV}\!\big\{\big|f^{\star}_{i,\hat{\tau}^{(t_r)}}\big(\langle\Bp_i,\Blambda\rangle\big)\big|\big\}\!\Bigg(\!\sum_{\Blambda\in\CV}\!\left|\hat{\tau}^{(\infty)}\big(\{\Blambda\}\big) - \hat{\tau}^{(t_r)}\big(\{\Blambda\}\big)\right|\!\Bigg)\! + \big|\hat{\psi}_i^{(\infty)}-\hat{\psi}_i^{(t_r)}\big| \quad\; \forall r\in\N. \nonumber
    \end{align}
    Since (\ref{eqn:algo-potential-CPWA}) guarantees that $f^{\star}_{i,\hat{\tau}^{(t_r)}}$ satisfies $f^{\star}_{i,\hat{\tau}^{(t_r)}}(\kappa_{i,1})=0$ and $\partial f^{\star}_{i,\hat{\tau}^{(t_r)}}(z)\subseteq[\underline{\omega}_i,\overline{\omega}_i]$ for all $z\in\R$,
    it holds that $\max_{r\in\N}\big\{\max_{\Blambda\in\CV}\big\{\big|f^{\star}_{i,\hat{\tau}^{(t_r)}}\big(\langle\Bp_i,\Blambda\rangle\big)\big|\big\}\big\}<\infty$.
    Thus, taking $r\to\infty$ in (\ref{eqn:cuttingplane-LSIP-proof-errorlimit}) leads to
    $\liminf_{r\to\infty}\hat{\varepsilon}_{\OT,i}^{(t_r)}\le 0$,
    which shows that the termination condition in Line~\ref{alglin:cuttingplane-LSIP-termination} will be satisfied after finitely many iterations of Algorithm~\ref{alg:cuttingplane-LSIP}.
    This proves statement~\ref{props:cuttingplane-LSIP-termination}.

    Now, let us assume that the termination condition in Line~\ref{alglin:cuttingplane-LSIP-termination} is satisfied in iteration~$t\in\nobreak\N$ and prove statements~\ref{props:cuttingplane-LSIP-dual} and \ref{props:cuttingplane-LSIP-primal}.
    It holds by Lines~\ref{alglin:cuttingplane-LSIP-returnval}--\ref{alglin:cuttingplane-LSIP-lowerbound} that
    $\hat{\tau}=\nobreak \hat{\tau}^{(t)}$,
    $\hat{\Bxi}_{\ineq}=\nobreak \hat{\Bxi}_{\ineq}^{(t)}$,
    $\hat{\Bxi}_{\eq}=\nobreak \hat{\Bxi}_{\eq}^{(t)}$,
    $\hat{\BIx}=\nobreak\hat{\BIx}^{(t)}$,
    $\hat{f}_i=\nobreak \INDI_{\{i=1\}}\hat{h}_0^{(t)} + \sum_{f_i\in\widehat{\CF}_{i}^{(t)}}\hat{y}_{i,f_i}^{(t)}f_i$ for $i=1,\ldots,d$,
    and
    $\widehat{C}_{\DRO}^{\mathrm{LB}}=\widehat{C}_{\DRO}^{(t)}-\big(\sum_{i=1}^d\hat{\varepsilon}_{\OT,i}^{(t)}\big)$.
    Observe that 
    $\support\big(\hat{\tau}^{(t)}\big)\subseteq\CV$,
    $\int_{\CS^*_2}\BQ^\TRANSP\Blambda\DIFFM{\hat{\tau}^{(t)}}{\DIFF\Blambda}=\sum_{\Blambda\in\CV}\BQ^\TRANSP\Blambda\hat{\nu}_{\Blambda}^{(t)}$ by
    Line~\ref{alglin:cuttingplane-LSIP-constructmeas},
    and hence Line~\ref{alglin:cuttingplane-LSIP-LP} guarantees that
    $(\hat{\tau}^{(t)},\Bxi_{\ineq}^{(t)},\Bxi_{\ineq}^{(t)})$ is feasible for \eqref{eqn:dro-dual}.
    Moreover, 
    Lines~\ref{alglin:cuttingplane-LSIP-newfuncs}--\ref{alglin:cuttingplane-LSIP-termination}, (\ref{eqn:algo-potential-CPWA}),
    and Lemma~\ref{lem:OT-1d}\ref{lems:OT-1d-potential-existence} with $\mu\leftarrow\mu_i$, $\nu\leftarrow\langle\Bp_i,\cdot\,\rangle\sharp\hat{\tau}^{(t)}$ imply that
    \begin{align}
        \begin{split}
            \varepsilon_{\OT}\ge \hat{\varepsilon}_{\OT,i}^{(t)}&={-\int_{\Omega_i}}\sup_{z\in\R}\big\{z\omega-f^{\star}_{i,\hat{\tau}^{(t)}}(z)\big\}\DIFFM{\mu_i}{\DIFF\omega} - \int_{\R}f^{\star}_{i,\hat{\tau}^{(t)}}\DIFFX{\langle\Bp_i,\cdot\,\rangle\sharp\hat{\tau}^{(t)}} - \hat{\psi}_i^{(t)}\\
            &=\OT_{\mu_i}\big(\langle\Bp_i,\cdot\,\rangle\sharp\hat{\tau}^{(t)}\big)-\hat{\psi}_i^{(t)}.
        \end{split}
        \label{eqn:cuttingplane-LSIP-proof-potential-ineq1}
    \end{align}
    Thus, 
    (\ref{eqn:cuttingplane-LSIP-proof-potential-ineq1}), 
    Line~\ref{alglin:cuttingplane-LSIP-LP}, and Line~\ref{alglin:cuttingplane-LSIP-lowerbound}
    imply that
    \begin{align*}
        &\hspace{-20pt}\int_{\CS^*_2}\langle\Bb,\Blambda\rangle \DIFFM{\hat{\tau}^{(t)}}{\DIFF\Blambda} 
        - \Bigg(\sum_{i=1}^d\OT_{\mu_i}\big(\langle\Bp_i,\cdot\,\rangle\sharp\hat{\tau}^{(t)}\big)\Bigg) 
        + \big\langle\Bq_{\ineq},\hat{\Bxi}_{\ineq}^{(t)}\big\rangle + \big\langle\Bq_{\eq},\hat{\Bxi}_{\eq}^{(t)}\big\rangle \\
        &=\Bigg(\sum_{\Blambda\in\CV}\langle\Bb,\Blambda\rangle\hat{\nu}_{\Blambda}^{(t)}\Bigg) 
        - \Bigg(\sum_{i=1}^{d} \hat{\psi}_{i}^{(t)} + \hat{\varepsilon}_{\OT,i}^{(t)}\Bigg)
        + \big\langle\Bq_{\ineq},\hat{\Bxi}_{\ineq}^{(t)}\big\rangle + \big\langle\Bq_{\eq},\hat{\Bxi}_{\eq}^{(t)}\big\rangle\\
        &= \widehat{C}_{\DRO}^{(t)} - \Bigg(\sum_{i=1}^{d} \hat{\varepsilon}_{\OT,i}^{(t)}\Bigg) = \widehat{C}_{\DRO}^{\mathrm{LB}}.
    \end{align*}
    We have completed the proof of statement~\ref{props:cuttingplane-LSIP-dual}.
    Furthermore, Line~\ref{alglin:cuttingplane-LSIP-LP} also guarantees that
    $\hat{\BIx}^{(t)}\in\nobreak\CS_1$, and that
    $\big((\hat{\nu}_{\Blambda}^{(t)})_{\Blambda\in\CV},\allowbreak (\hat{\psi}_i^{(t)})_{i=1:d}, \allowbreak\hat{\Bxi}_{\ineq}^{(t)},\allowbreak\hat{\Bxi}_{\eq}^{(t)}\big)$ and
    $\big(\hat{\BIx}^{(t)},\allowbreak\hat{h}_0^{(t)},\allowbreak(\hat{y}_{i,f_i}^{(t)})_{f_i\in\widehat{\CF}_{i}^{(t)},\,i=1:d}\big)$
    satisfy the following complementary slackness conditions:
    \begin{subequations}\label{eqn:cuttingplane-LSIP-proof-complementary-slackness}
        \begin{align}
            \sum_{\Blambda\in\CV}\hat{\nu}_{\Blambda}^{(t)}\Bigg(\langle\BQ\hat{\BIx}^{(t)}+\Bb,\Blambda\rangle + \hat{h}_0^{(t)} + \sum_{i=1}^{d} \sum_{f_i\in\widehat{\CF}_{i}^{(t)}}\hat{y}_{i,f_i}^{(t)} f_i\big(\langle\Bp_i,\Blambda\rangle\big)\Bigg) &= 0,
            \label{eqn:cuttingplane-LSIP-proof-complementary-slackness-funcdominance} \allowdisplaybreaks\\
            \langle\BL_{\ineq}\hat{\BIx}^{(t)}-\Bq_{\ineq}, \hat{\Bxi}_{\ineq}^{(t)}\rangle&=0,
            \label{eqn:cuttingplane-LSIP-proof-complementary-slackness-primal-in}\allowdisplaybreaks\\
            \sum_{i=1}^{d}\sum_{f_i\in\widehat{\CF}_i^{(t)}}\hat{y}_{i,f_i}^{(t)}\Bigg(\hat{\psi}_i^{(t)} + \int_{\Omega_i} \sup_{z\in\R}\big\{z\omega-f_i(z)\big\} \DIFFM{\mu_i}{\DIFF\omega} + \sum_{\Blambda\in\CV}f_i\big(\langle\Bp_i,\Blambda\rangle\big)\hat{\nu}_{\Blambda}^{(t)}\Bigg)&=0,
            \label{eqn:cuttingplane-LSIP-proof-complementary-slackness-funcweight}\allowdisplaybreaks\\
            \big\langle\Bc_1+\big({\textstyle\sum_{\Blambda\in\CV}}\BQ^\TRANSP\Blambda\hat{\nu}_{\Blambda}^{(t)}\big)-\BL_{\ineq}^\TRANSP\hat{\Bxi}_{\ineq}^{(t)}-\BL_{\eq}^\TRANSP\hat{\Bxi}_{\eq}^{(t)},\hat{\BIx}^{(t)}\big\rangle&=0.
            \label{eqn:cuttingplane-LSIP-proof-complementary-slackness-dual-in}
        \end{align}
    \end{subequations}
    Specifically, combining (\ref{eqn:cuttingplane-LSIP-proof-complementary-slackness-funcdominance})
    with the constraints of \LPdual{\widehat{\CF}^{(t)}}, 
    Line~\ref{alglin:cuttingplane-LSIP-constructmeas},
    and Line~\ref{alglin:cuttingplane-LSIP-returnfuncs}
    leads to
    $\big\langle\BQ\hat{\BIx}+\Bb,\Blambda\big\rangle+\sum_{i=1}^{d}\hat{f}_i\big(\langle\Bp_i,\Blambda\rangle\big)\le 0$
    $\forall \Blambda\in\CV$,
    as well as 
    $\big\langle\BQ\hat{\BIx}+\Bb,\Blambda\big\rangle+\sum_{i=1}^{d}\hat{f}_i\big(\langle\Bp_i,\Blambda\rangle\big)=\nobreak0$
    $\forall \Blambda\in\support(\hat{\tau})$.
    Since
    (\ref{eqn:cuttingplane-LSIP-proof-complementary-slackness-primal-in})
    and 
    (\ref{eqn:cuttingplane-LSIP-proof-complementary-slackness-dual-in})
    both hold,
    one may check that 
    $\hat{\BIx},\hat{f}_1,\ldots,\hat{f}_d$ 
    satisfy (\ref{eqn:dual-suboptimality-primal-feasibility-potential-eq})--(\ref{eqn:dual-suboptimality-dual-complementarity})
    with respect to 
    $\tau\leftarrow\hat{\tau}$, $\Bxi_{\ineq}\leftarrow\hat{\Bxi}_{\ineq}$, $\Bxi_{\eq}\leftarrow\hat{\Bxi}_{\eq}$.
    Moreover, for $i=1,\ldots,d$,
    since Line~\ref{alglin:cuttingplane-LSIP-initialfuncs}, Line~\ref{alglin:cuttingplane-LSIP-newfuncs}, and (\ref{eqn:algo-potential-CPWA}) guarantee that every $f_i\in\widehat{\CF}_i^{(t)}$ satisfies (\ref{eqn:dual-suboptimality-potential-subgradient}),
    and since the constraints of \LPdual{\widehat{\CF}^{(t)}} guarantee
    that $\hat{y}_{i,f_i}^{(t)}\ge 0$ $\forall f_i\in\widehat{\CF}_{i}^{(t)}$ and $\sum_{f_i\in\widehat{\CF}_{i}^{(t)}}\hat{y}_{i,f_i}^{(t)}=1$,
    Line~\ref{alglin:cuttingplane-LSIP-returnfuncs} implies that $\hat{f}_i$ satisfies 
    (\ref{eqn:dual-suboptimality-potential-subgradient}).
    On the other hand,
    (\ref{eqn:cuttingplane-LSIP-proof-complementary-slackness-funcweight}) and Line~\ref{alglin:cuttingplane-LSIP-constructmeas} imply that
    \begin{align*}
        {-\int_{\Omega_i}}\sup_{z\in\R}\big\{z\omega-\hat{f}_i(z)\big\}\DIFFM{\mu_i}{\DIFF\omega}
        &=\INDI_{\{i=0\}}\hat{h}_0^{(t)} - \int_{\Omega_i}\sup_{z\in\R}\Big\{z\omega - {\textstyle\sum_{f_i\in\widehat{\CF}_{\CV,i}^{(t)}}}\hat{y}_{i,f_i}^{(t)}f_i(z)\Big\}\DIFFM{\mu_i}{\DIFF\omega} \\
        &\ge \INDI_{\{i=0\}}\hat{h}_0^{(t)} -\sum_{f_i\in\widehat{\CF}_{i}^{(t)}}\hat{y}_{i,f_i}^{(t)} \int_{\Omega_i}\sup_{z\in\R}\big\{z\omega-f_i(z)\big\}\DIFFM{\mu_i}{\DIFF\omega} \\
        &= \INDI_{\{i=0\}}\hat{h}_0^{(t)} + \sum_{f_i\in\widehat{\CF}_{i}^{(t)}}\hat{y}_{i,f_i}^{(t)} \Bigg(\hat{\psi}_i^{(t)} + \sum_{\Blambda\in\CV}f_i\big(\langle\Bp_i,\Blambda\rangle\big)\hat{\nu}_{\Blambda}^{(t)}\Bigg)\\
        &= \hat{\psi}_i^{(t)} + \int_{\R} \hat{f}_i\DIFFX{\langle\Bp_i,\cdot\,\rangle\sharp\hat{\tau}}.
    \end{align*}
    Combining this with (\ref{eqn:cuttingplane-LSIP-proof-potential-ineq1}) and 
    Line~\ref{alglin:cuttingplane-LSIP-termination} 
    shows that 
    $\hat{f}_i$ satisfies (\ref{eqn:dual-suboptimality-potential-optimizer})
    with respect to $\hat{\tau}$.
    Lastly, 
    since each $f^{\star}_{i,\hat{\tau}^{(t)}}$ in Line~\ref{alglin:cuttingplane-LSIP-initialfuncs} and Line~\ref{alglin:cuttingplane-LSIP-newfuncs} is continuous piece-wise affine, 
    it holds that $\widehat{\CF}_{i}$ contains only continuous piece-wise affine functions for $i=1,\ldots,d$,
    and Line~\ref{alglin:cuttingplane-LSIP-returnfuncs} and Line~\ref{alglin:cuttingplane-LSIP-returnval} subsequently guarantee that
    $\hat{f}_1,\ldots,\hat{f}_d$ are also continuous piece-wise affine.
    The proof is now complete.
\end{proof}


\begin{proof}[Proof of Theorem~\ref{thm:iterative}]
    Let us begin by examining the properties of 
    $\hat{\tau}^{(t)}$,
    $\hat{\Bxi}_{\ineq}^{(t)}$,
    $\hat{\Bxi}_{\eq}^{(t)}$,
    $\hat{\BIx}^{(t)}$,
    $\hat{f}_1^{(t)}$,~$\ldots$~,~$\hat{f}_d^{(t)}$,
    $\overline{f}_1^{(t)}$,~$\ldots$~,~$\overline{f}_d^{(t)}$,
    $\widehat{C}_{\DRO}^{\mathrm{LB}(t)}$,
    $\overline{\varepsilon}_{\mathrm{vio}}^{(t)}$,
    $\hat{\varepsilon}_{\mathrm{sub}}^{(t)}$ 
    for each iteration~$t\in\N$.
    Line~\ref{alglin:iterative-LSIP} and
    Proposition~\ref{prop:cuttingplane-LSIP}\ref{props:cuttingplane-LSIP-dual}--\ref{props:cuttingplane-LSIP-primal} imply that
    $\big(\hat{\tau}^{(t)},\hat{\Bxi}_{\ineq}^{(t)},\hat{\Bxi}_{\eq}^{(t)}\big)$ is feasible for \eqref{eqn:dro-dual} 
    with objective $\widehat{C}_{\DRO}^{\mathrm{LB}(t)}$,
    $\support\big(\hat{\tau}^{(t)}\big)\subseteq\CV^{(t)}$,
    and
    $\hat{\BIx}^{(t)},\hat{f}_1^{(t)},\ldots,\hat{f}_d^{(t)}$ satisfy (\ref{eqn:dual-suboptimality-primal-feasibility-potential-eq})--(\ref{eqn:dual-suboptimality-potential-subgradient})
    with respect to $\tau\leftarrow\hat{\tau}^{(t)}$, $\Bxi_{\ineq}\leftarrow\hat{\Bxi}_{\ineq}^{(t)}$, $\Bxi_{\eq}\leftarrow\hat{\Bxi}_{\eq}^{(t)}$, $\CV\leftarrow\CV^{(t)}$, $\varepsilon_{\OT}\leftarrow\varepsilon_{\OT}^{(t)}$.
    Firstly, it follows from
    (\ref{eqn:dual-suboptimality-primal-feasibility-potential-eq}),
    (\ref{eqn:dual-suboptimality-primal-feasibility-potential-le}),
    Line~\ref{alglin:iterative-funcrelax},
    and Line~\ref{alglin:iterative-globalmax}
    that
    \begin{align}
        \big\langle\BQ\hat{\BIx}^{(t)}+\Bb,\Blambda\big\rangle + \sum_{i=1}^{d}\hat{f}_i^{(t)}\big(\langle\Bp_i,\Blambda\rangle\big) &\le 0 \hspace{31.0pt} \qquad \forall \Blambda\in\CV^{(t)},
        \label{eqn:iterative-proof-iter-potential-le}\allowdisplaybreaks\\
        \big\langle\BQ\hat{\BIx}^{(t)}+\Bb,\Blambda\big\rangle + \sum_{i=1}^{d}\hat{f}_i^{(t)}\big(\langle\Bp_i,\Blambda\rangle\big) &= 0 \qquad \forall \Blambda\in\support\big(\hat{\tau}^{(t)}\big),
        \label{eqn:iterative-proof-iter-potential-eq}\allowdisplaybreaks\\
        \big\langle\BQ\hat{\BIx}^{(t)}+\Bb,\Blambda\big\rangle + \sum_{i=1}^{d}\overline{f}_i^{(t)}\big(\langle\Bp_i,\Blambda\rangle\big) &\le \overline{\varepsilon}_{\mathrm{vio}}^{(t)} \hspace{25.6pt} \qquad \forall \Blambda\in\CS^*_2.
        \label{eqn:iterative-proof-iter-primal-feasibility}
    \end{align}
    Secondly, 
    observe that (\ref{eqn:iterative-proof-iter-primal-feasibility})
    shows that
    $\big(\hat{\BIx}^{(t)},\overline{f}_1^{(t)}-\overline{\varepsilon}_{\mathrm{vio}}^{(t)},\overline{f}_2^{(t)},\ldots,\overline{f}_{d}^{(t)}\big)$ is feasible for \eqref{eqn:dro-primal},
    whereas Line~\ref{alglin:iterative-suboptimality} and Line~\ref{alglin:iterative-primalobj} imply that 
    the objective of
    $\big(\hat{\BIx}^{(t)},\overline{f}_1^{(t)}-\overline{\varepsilon}_{\mathrm{vio}}^{(t)},\overline{f}_2^{(t)},\ldots,\overline{f}_{d}^{(t)}\big)$
    with respect to \eqref{eqn:dro-primal} is given by
    \begin{align}
        \big\langle\Bc_1,\hat{\BIx}^{(t)}\big\rangle 
        + \Bigg(\sum_{i=1}^{d} \int_{\Omega_i}\sup_{z\in\R}\big\{z\omega-\overline{f}_i^{(t)}(z)\big\}\DIFFM{\mu_i}{\DIFF\omega}\Bigg)
        + \overline{\varepsilon}_{\mathrm{vio}}^{(t)} 
        &= \widehat{C}_{\DRO}^{\mathrm{LB}(t)} + \hat{\varepsilon}_{\mathrm{sub}}^{(t)}.
        \label{eqn:iterative-proof-iter-primal-objective}
    \end{align}
    Thirdly, 
    integrating both sides of (\ref{eqn:iterative-proof-iter-potential-eq}) with respect to $\hat{\tau}^{(t)}$ and then combining the resulting equality with 
    (\ref{eqn:dual-suboptimality-primal-complementarity}),
    (\ref{eqn:dual-suboptimality-dual-complementarity}),
    and Line~\ref{alglin:iterative-OTerror} yields
    \begin{align}
        \hspace{20pt}&\hspace{-20pt}\big\langle\Bc_1,\hat{\BIx}^{(t)}\big\rangle
        + \sum_{i=1}^{d} \int_{\Omega_i}\sup_{z\in\R}\big\{z\omega-\hat{f}_i^{(t)}(z)\big\}\DIFFM{\mu_i}{\DIFF\omega} \nonumber\\
        \begin{split}
            &= \int_{\CS^*_2}\langle\Bb,\Blambda\rangle\DIFFM{\hat{\tau}^{(t)}}{\DIFF\Blambda} 
            + \big\langle\Bq_{\ineq},\hat{\Bxi}_{\ineq}^{(t)}\big\rangle
            + \big\langle\Bq_{\eq},\hat{\Bxi}_{\eq}^{(t)}\big\rangle \\
            &\qquad + \sum_{i=1}^{d}\!\int_{\Omega_i}\sup_{z\in\R}\big\{z\omega-\hat{f}_i^{(t)}(z)\big\}\DIFFM{\mu_i}{\DIFF\omega}
            + \int_{\R}\hat{f}_{i}^{(t)}\DIFFX{\langle\Bp_i,\cdot\,\rangle\sharp\hat{\tau}^{(t)}}
        \end{split}\label{eqn:iterative-proof-iter-LSIP-objective}\\
        &= \widehat{C}_{\DRO}^{\mathrm{LB}(t)} + \sum_{i=1}^{d}\OT_{\mu_i}\big(\langle\Bp_i,\cdot\,\rangle\sharp\hat{\tau}^{(t)}\big)
        + \int_{\Omega_i}\sup_{z\in\R}\big\{z\omega-\hat{f}_i^{(t)}(z)\big\}\DIFFM{\mu_i}{\DIFF\omega}
        + \int_{\R}\hat{f}_{i}^{(t)}\DIFFX{\langle\Bp_i,\cdot\,\rangle\sharp\hat{\tau}^{(t)}}. \nonumber
    \end{align}
    Lastly, combining
    Line~\ref{alglin:iterative-funcrelax}, 
    (\ref{eqn:iterative-proof-iter-primal-objective}),
    (\ref{eqn:iterative-proof-iter-LSIP-objective}),
    and (\ref{eqn:dual-suboptimality-potential-optimizer})
    leads to
    \begin{align}
        \hat{\varepsilon}_{\mathrm{sub}}^{(t)}
        &\le \big\langle\Bc_1,\hat{\BIx}^{(t)}\big\rangle 
        + \!\Bigg(\!\sum_{i=1}^{d} \int_{\Omega_i}\sup_{z\in\R}\big\{\hspace{-0.1em}z\omega-\hat{f}_i^{(t)}(z)\hspace{-0.1em}\big\}\DIFFM{\mu_i}{\DIFF\omega}\!\Bigg)\!
        - \widehat{C}_{\DRO}^{\mathrm{LB}(t)} + \overline{\varepsilon}_{\mathrm{vio}}^{(t)}\label{eqn:iterative-proof-iter-suboptimality-bound1}\\
        &=\overline{\varepsilon}_{\mathrm{vio}}^{(t)} 
        + \!\sum_{i=1}^{d}\OT_{\mu_i}\big(\langle\Bp_i,\cdot\,\rangle\sharp\hat{\tau}^{(t)}\big)
        + \!\int_{\Omega_i}\sup_{z\in\R}\big\{\hspace{-0.1em}z\omega-\hat{f}_i^{(t)}(z)\hspace{-0.1em}\big\}\DIFFM{\mu_i}{\DIFF\omega}
        + \!\int_{\R}\hat{f}_{i}^{(t)}\DIFFX{\langle\Bp_i,\cdot\,\rangle\sharp\hat{\tau}^{(t)}}\nonumber\\
        &\le \overline{\varepsilon}_{\mathrm{vio}}^{(t)} + d\varepsilon_{\OT}^{(t)},\nonumber
    \end{align}
    as well as
    \begin{align}
        \hat{\varepsilon}_{\mathrm{sub}}^{(t)}
        &\ge \big\langle\Bc_1,\hat{\BIx}^{(t)}\big\rangle 
        + \!\Bigg(\!\sum_{i=1}^{d} \int_{\Omega_i}\sup_{z\in\R}\Big\{\hspace{-0.1em}z\omega-\big(\hat{f}_i^{(t)}(z)+\varsigma^{(t)}\big)\hspace{-0.1em}\Big\}\DIFFM{\mu_i}{\DIFF\omega}\!\Bigg)\!
        - \widehat{C}_{\DRO}^{\mathrm{LB}(t)} + \overline{\varepsilon}_{\mathrm{vio}}^{(t)}\label{eqn:iterative-proof-iter-suboptimality-bound2}\\
        &=\overline{\varepsilon}_{\mathrm{vio}}^{(t)} - d\varsigma^{(t)} + \!\sum_{i=1}^{d}\!\OT_{\mu_i}\!\big(\hspace{-1pt}\langle\Bp_i,\cdot\,\rangle\sharp\hat{\tau}^{(t)}\hspace{-1pt}\big)\!
        + \!\int_{\Omega_i}\sup_{z\in\R}\big\{\hspace{-0.1em}z\omega-\hat{f}_i^{(t)}(z)\hspace{-0.1em}\big\}\DIFFM{\mu_i}{\DIFF\omega}\!
        + \!\int_{\R}\hat{f}_{i}^{(t)}\DIFFX{\langle\Bp_i,\cdot\,\rangle\sharp\hat{\tau}^{(t)}}.\nonumber
    \end{align}

    Now, observe that in every iteration~$t\in\N$, Line~\ref{alglin:iterative-funcrelax} guarantees
    $\overline{f}_{1}^{(t)},\ldots,\overline{f}_{d}^{(t)}\in\CC_{\mathrm{cvx}}(\R)$ and $\partial \overline{f}_i^{(t)}(z)\subseteq\Omega_i$ $\forall z\in\R$ for $i=1,\ldots,d$.
    Subsequently, statement~\ref{thms:iterative-globalmax} follows from Lemma~\ref{lem:global-simplification}.
    To prove statement~\ref{thms:iterative-termination},
    let us suppose for the sake of contradiction that the termination condition in Line~\ref{alglin:iterative-termination} is never satisfied and consider the terms generated by the for-loop in Line~\ref{alglin:iterative-forloop} for all $t\in\N$.
    Let $\CT:=\big\{\tilde{t}\in\N:\hat{\Blambda}^{(\tilde{t})}\in\CV^{(\tilde{t})}\big\}$,
    and let $\tilde{t}_0:=\min \CT$.
    Observe that Line~\ref{alglin:iterative-newsupport} and the assumption that $\hat{\Blambda}^{(t)}\in\CV_{\extremept}$ in each iteration of Algorithm~\ref{alg:iterative} guarantee that 
    $|\N\setminus\CT|\le |\CV_{\extremept}|<\infty$ and thus $\tilde{t}_0\le |\CV_{\extremept}|+1$.
    Subsequently, 
    (\ref{eqn:iterative-proof-iter-potential-le}) and Line~\ref{alglin:iterative-newsupport} imply that
    $\big\langle\BQ\hat{\BIx}^{(\tilde{t})}+\Bb,\hat{\Blambda}^{(\tilde{t})}\big\rangle+\sum_{i=1}^d\hat{f}_{i}^{(\tilde{t})}\big(\langle\Bp_i,\hat{\Blambda}^{(\tilde{t})}\rangle\big)\le \nobreak0$ for all $\tilde{t}\in\CT$.
    Thus,
    we obtain from 
    (\ref{eqn:iterative-proof-iter-suboptimality-bound1}),
    Line~\ref{alglin:iterative-globalmax}, 
    and Line~\ref{alglin:iterative-funcrelax}, 
    that
    \begin{align}
        \begin{split}
        \hat{\varepsilon}_{\mathrm{sub}}^{(\tilde{t})}&\le \overline{\varepsilon}_{\mathrm{vio}}^{(\tilde{t})}+d\varepsilon_{\OT}^{(\tilde{t})} 
        \le d\varepsilon_{\OT}^{(\tilde{t})} + \big(1+\vartheta^{(\tilde{t})}\big)\Bigg(\big\langle\BQ\hat{\BIx}^{(\tilde{t})}+\Bb,\hat{\Blambda}^{(\tilde{t})}\big\rangle+\sum_{i=1}^d\overline{f}_{i}^{(\tilde{t})}\big(\langle\Bp_i,\hat{\Blambda}^{(\tilde{t})}\rangle\big)\Bigg)\\
        &\le d\big(\varepsilon_{\OT}^{(\tilde{t})}+\big(1+\vartheta^{(\tilde{t})}\big)\varsigma^{(\tilde{t})}\big) + \big(1+\vartheta^{(\tilde{t})}\big)\Bigg(\big\langle\BQ\hat{\BIx}^{(\tilde{t})}+\Bb,\hat{\Blambda}^{(\tilde{t})}\big\rangle+\sum_{i=1}^d\hat{f}_{i}^{(\tilde{t})}\big(\langle\Bp_i,\hat{\Blambda}^{(\tilde{t})}\rangle\big)\Bigg)\\
        &\le d\big(\varepsilon_{\OT}^{(\tilde{t})}+\big(1+\vartheta^{(\tilde{t})}\big)\varsigma^{(\tilde{t})}\big) \hspace{220pt} \forall \tilde{t}\in\CT.
        \end{split}
        \label{eqn:iterative-proof-termination-suboptimality}
    \end{align}
    Since Setting~\ref{sett:iterative} requires that $\varepsilon_{\mathrm{tol}}>\liminf_{t\to\infty} d\big(\varepsilon_{\OT}^{(t)}+\big(1+\vartheta^{(t)}\big)\varsigma^{(t)}\big)$,
    (\ref{eqn:iterative-proof-termination-suboptimality})
    shows that $\liminf_{\tilde{t}\to\infty}\hat{\varepsilon}_{\mathrm{sub}}^{(\tilde{t})} < \varepsilon_{\mathrm{tol}}$ and thus
    the termination condition in Line~\ref{alglin:iterative-termination} is satisfied after finitely many iterations.
    Moreover, under the additional assumption that $\varepsilon_{\OT}^{(t)}=\nobreak\varepsilon_{\OT}^{(1)}$,
    $\varsigma^{(t)}=\varsigma^{(1)}$,
    $\vartheta^{(t)}=\vartheta^{(1)}$ for all $t\in\N$,
    (\ref{eqn:iterative-proof-termination-suboptimality})
    shows that
    $\hat{\varepsilon}_{\mathrm{sub}}^{(\tilde{t}_0)} < \varepsilon_{\mathrm{tol}}$ and that 
    the termination condition in Line~\ref{alglin:iterative-termination} is satisfied at least once within the first $\tilde{t}_0\le |\CV_{\extremept}|+\nobreak1$ iterations.
    The proof of statement~\ref{thms:iterative-termination} is complete.
    
    Next, one observes from
    (\ref{eqn:iterative-proof-iter-primal-objective}),
    Line~\ref{alglin:iterative-returnval}, 
    Line~\ref{alglin:iterative-bounds},
    the termination condition in Line~\ref{alglin:iterative-termination},
    and Corollary~\ref{cor:dual-optimizers}\ref{cors:dual-suboptimality-primal-coupling-primal}
    that
    $(\hat{\tau},\hat{\Bxi}_{\ineq},\hat{\Bxi}_{\eq})$ is feasible for \eqref{eqn:dro-dual} with objective $\widehat{C}_{\DRO}^{\mathrm{LB}}$, 
    $(\hat{\BIx},\tilde{f}_1,\ldots,\tilde{f}_d)$ is feasible for \eqref{eqn:dro-primal} with objective $\widehat{C}_{\DRO}^{\mathrm{UB}}=\widehat{C}_{\DRO}^{\mathrm{LB}}+\hat{\varepsilon}_{\mathrm{sub}}$,
    and $\hat{\varepsilon}_{\mathrm{sub}}\le \varepsilon_{\mathrm{tol}}$.
    Moreover, observe from 
    (\ref{eqn:iterative-proof-iter-suboptimality-bound1}), 
    (\ref{eqn:iterative-proof-iter-suboptimality-bound2}),
    Line~\ref{alglin:iterative-OTerror},
    and Line~\ref{alglin:iterative-bounds}
    that
    $\hat{\varepsilon}_{\mathrm{sub}}\le \hat{\varepsilon}_{\mathrm{prob}}\le \hat{\varepsilon}_{\mathrm{sub}} + d\hat{\varsigma}\le \varepsilon_{\mathrm{tol}}+d\sup_{t\in\N}\big\{\varsigma^{(t)}\big\}$.
    Subsequently, applying Theorem~\ref{thm:primaldual}\ref{thms:primaldual-strongduality} proves 
    statements~\ref{thms:iterative-bounds}--\ref{thms:iterative-primal}.

    Lastly, to prove statement~\ref{thms:iterative-dro-worstmeas},
    observe that the construction of $\hat{\mu}$ in Lines~\ref{alglin:iterative-support}--\ref{alglin:iterative-meas} satisfies the assumptions of Corollary~\ref{cor:dual-suboptimality-primal-coupling} 
    with respect to 
    $\tau\leftarrow\hat{\tau}$,
    $\Bxi_{\ineq}\leftarrow\hat{\Bxi}_{\ineq}$,
    $\Bxi_{\eq}\leftarrow\hat{\Bxi}_{\eq}$,
    $\BIx\leftarrow\hat{\BIx}$,
    $f_1\leftarrow\nobreak\hat{f}_1^{(t)}$,~$\ldots$~,~$f_d\leftarrow\hat{f}_d^{(t)}$, 
    and $\mu_{\tau}\leftarrow\hat{\mu}$.
    Subsequently, it follows from
    (\ref{eqn:dual-suboptimality-primal-coupling-proof-lb}) in the proof of Corollary~\ref{cor:dual-suboptimality-primal-coupling} that
    \begin{align}
            \int_{\BOmega} C_2\big(\hat{\BIx},\Bomega\big)\DIFFM{\hat{\mu}}{\DIFF\Bomega}&\ge \int_{\CS^*_2}\langle\BQ\hat{\BIx}+\Bb,\Blambda\rangle\DIFFM{\hat{\tau}}{\DIFF\Blambda}-\Bigg(\sum_{i=1}^d \OT_{\mu_i}\big(\langle\Bp_i,\cdot\,\rangle\sharp\hat{\tau}\big)\Bigg).
        \label{eqn:iterative-proof-worstmeas-ineq1}
    \end{align}
    Integrating both sides of (\ref{eqn:iterative-proof-iter-potential-eq}) with respect to $\hat{\tau}=\hat{\tau}^{(t)}$, and then substituting the resulting equality together with
    Lines~\ref{alglin:iterative-OTerror-forloop}--\ref{alglin:iterative-bounds}
    into (\ref{eqn:iterative-proof-worstmeas-ineq1}) leads to
    \begin{align}
        \begin{split}
            \int_{\BOmega} C_2\big(\hat{\BIx},\Bomega\big)\DIFFM{\hat{\mu}}{\DIFF\Bomega}&\ge {-\Bigg(\sum_{i=1}^d \OT_{\mu_i}\big(\langle\Bp_i,\cdot\,\rangle\sharp\hat{\tau}\big)+ \int_{\R}\hat{f}_i^{(t)}\DIFFX{\langle\Bp_i,\cdot\,\rangle\sharp\hat{\tau}}\Bigg)}\\
            &= \sum_{i=1}^{d}\int_{\Omega_i}\sup_{z\in\R}\big\{z\omega-\hat{f}_i^{(t)}(z)\big\}\DIFFM{\mu_i}{\DIFF\omega} - \hat{\varepsilon}_{\OT,i} \\
            &\ge \sum_{i=1}^{d}\int_{\Omega_i}\sup_{z\in\R}\big\{z\omega-\overline{f}_i^{(t)}(z)\big\}\DIFFM{\mu_i}{\DIFF\omega} - \hat{\varepsilon}_{\OT,i} \\
            &=\Bigg(\sum_{i=1}^{d}\int_{\Omega_i}\sup_{z\in\R}\big\{z\omega-\tilde{f}_i(z)\big\}\DIFFM{\mu_i}{\DIFF\omega}\Bigg) - \hat{\varepsilon}_{\mathrm{prob}}.
        \end{split}
        \label{eqn:iterative-proof-worstmeas-ineq2}
    \end{align}
    Since $(\hat{\BIx},\tilde{f}_1,\ldots,\tilde{f}_d)$ is feasible for \eqref{eqn:dro-primal},
    applying Lemma~\ref{lem:primal-objective} to (\ref{eqn:iterative-proof-worstmeas-ineq2}) proves statement~\ref{thms:iterative-dro-worstmeas}.
    The proof is now complete.
\end{proof}


\begin{proof}[Proof of Proposition~\ref{prop:globalmax-MILP}]
    Let us remark that 
    since $\partial f_i(z)\subseteq\Omega_i$ for all $z\in\R$ and for $i=1,\ldots,d$,
    Lemma~\ref{lem:global-simplification} implies that $\mathrm{vio}(\BIx,f_1,\ldots,f_d)<\infty$ and 
    that (\ref{eqn:dual-suboptimality-vio-def}) attains a maximizer in $\CV_{\extremept}$.
    Let us first fix an arbitrary maximizer $\Blambda^\star\in\CV_{\extremept}$ of (\ref{eqn:dual-suboptimality-vio-def}).
    For $i=1,\ldots,d$,
    let $\iota_{i,j}^\star:=\INDI_{\{\langle\Bp_i,\Blambda^\star\rangle>\kappa_{i,j+1}\}}$ for $j=1,\ldots,m_i-2$,
    and let $\zeta_{i,j}^\star:= \frac{(\langle\Bp_i,\Blambda^\star\rangle-\kappa_{i,j})^+}{\kappa_{i,j+1}-\kappa_{i,j}}\wedge 1$ for $j=1,\ldots,m_i-1$. 
    Since $\kappa_{i,1}\le \min_{\Blambda\in\CV_{\extremept}}\big\{\langle\Bp_i,\Blambda\rangle\big\}$ and $\kappa_{i,m_i}\ge \max_{\Blambda\in\CV_{\extremept}}\big\{\langle\Bp_i,\Blambda\rangle\big\}$ for ${i=1,\ldots,d}$ by assumption,
    one verifies that 
    $\big(\Blambda^\star,(\zeta^\star_{i,j})_{j=1:m_i-1,\,i=1:d},\allowbreak (\iota^\star_{i,j})_{j=1:m_i-2,\,i=1:d}\big)$ is feasible for (\ref{eqn:globalmax-MILP}).
    Moreover, for $i=1,\ldots,d$, 
    since $f_i$ is piece-wise affine 
    on $(-\infty,\kappa_{i,1}],\allowbreak[\kappa_{i,1},\kappa_{i,2}],\ldots,[\kappa_{i,m_i-1},\kappa_{i,m_i}],\allowbreak[\kappa_{i,m_i},\infty)$
    and $\kappa_{i,1} + \sum_{j=1}^{m_i-1} (\kappa_{i,j+1}-\kappa_{i,j})\zeta^\star_{i,j} = \langle\Bp_i,\Blambda^\star\rangle$,
    it holds that 
    $f_i(\kappa_{i,1}) + \sum_{j=1}^{m_i-1}\big(f_i(\kappa_{i,j+1}) - f_i(\kappa_{i,j})\big)\zeta^\star_{i,j}=f_i\big(\langle\Bp_i,\Blambda^\star\rangle\big)$.
    Thus, the objective of $\big(\Blambda^\star,(\zeta^\star_{i,j})_{j=1:m_i-1,\,i=1:d},\allowbreak (\iota^\star_{i,j})_{j=1:m_i-2,\,i=1:d}\big)$ with respect to (\ref{eqn:globalmax-MILP}) 
    is equal to $\langle\BQ\BIx+\nobreak\Bb,\Blambda^\star\rangle+\sum_{i=1}^df_i\big(\langle\Bp_i,\Blambda^\star\rangle\big)=\mathrm{vio}(\BIx,f_1,\ldots,f_d)$,
    which shows that $(\ref{eqn:globalmax-MILP})\ge (\ref{eqn:dual-suboptimality-vio-def})$.
    
    Conversely, let us fix an arbitrary feasible solution 
    $\big(\Blambda,(\zeta_{i,j})_{j=1:m_i-1,\,i=1:d},\allowbreak (\iota_{i,j})_{j=1:m_i-2,\,i=1:d}\big)$ of (\ref{eqn:globalmax-MILP}).
    One checks that $f_i(\kappa_{i,1}) + \sum_{j=1}^{m_i-1}\big(f_i(\kappa_{i,j+1}) - f_i(\kappa_{i,j})\big)\zeta_{i,j}=f_i\big(\langle\Bp_i,\Blambda\rangle\big)$
    for $i=1,\ldots,d$,
    and thus the objective of $\big(\Blambda,(\zeta_{i,j})_{j=1:m_i-1,\,i=1:d},\allowbreak (\iota_{i,j})_{j=1:m_i-2,\,i=1:d}\big)$ with respect to (\ref{eqn:globalmax-MILP}) is equal to 
    $\langle\BQ\BIx+\nobreak\Bb,\Blambda\rangle+\sum_{i=1}^df_i\big(\langle\Bp_i,\Blambda\rangle\big)$.
    This shows that $(\ref{eqn:dual-suboptimality-vio-def})\ge (\ref{eqn:globalmax-MILP})$.
    Combining all results established above then completes the proof of statements~\ref{props:globalmax-MILP-optimalval}--\ref{props:globalmax-MILP-optimizer}.
\end{proof}

\newpage

\bibliographystyle{abbrvnat}
\bibliography{references}

\end{document}